%% file: MAIN.tex
\documentclass{article}

\usepackage{JOEconfig}

\usepackage[T1]{fontenc}

\usepackage{graphicx}
\usepackage{subcaption}
\usepackage[svgnames]{xcolor}

\usepackage{hyperref}
\hypersetup{
    linktocpage,  
    colorlinks,
    citecolor=DarkBlue,
    filecolor=blue,
    linkcolor=blue,
    urlcolor=blue
}

\usepackage{accents} 
\usepackage{bm} 

\usepackage{amsmath, amsthm} 
\usepackage{array}
\usepackage{longtable,tabularx} 
\usepackage{enumitem} 
\setlist[itemize,1]{left=3ex,labelsep=0.75ex,itemsep=0pt,topsep=1pt,label={\footnotesize{\textbullet}}}
\setlist[itemize,2]{itemsep=0pt,topsep=0pt}
\setlist[enumerate,1]{left=3ex,labelsep=0.75ex,itemsep=0pt,topsep=1pt}

\usepackage{stackengine}
\usepackage{xfrac}
\usepackage{tensor}
\usepackage{mathtools} 
\usepackage{bigints}
\usepackage{cancel} 

 \usepackage[slantedGreek,libertine,vvarbb]{newtx}
    
\DeclareRobustCommand{\textsb}[1]{{\textbf{#1}}} 

\usepackage[frak=euler,scr=cm]{mathalpha} 

\usepackage{upgreek} 

\usepackage[artemisia]{textgreek}

\newsavebox\fooboxx
\newcommand{\negslantbox}[2][-.25]{\mbox{%
        \sbox{\fooboxx}{#2}%
        \hskip\wd\fooboxx
        \pdfsave
        \pdfsetmatrix{1 0 #1 1}%
        \llap{\usebox{\fooboxx}}%
        \pdfrestore
}}

\newcommand{\upell}[0]{\negslantbox{$\mkern1mu\ell\mkern-1.5mu$}}

\DeclareMathSizes{10}{9.5}{7}{5.5} 
\DeclareMathSizes{9}{9.5}{7}{5}
 \DeclareMathSizes{8}{8}{6}{5}
 \DeclareMathSizes{7}{7}{5}{4}

\begin{document}

\title{\vspace{-1.0cm}Projective Transformations for Regularized Central-Force Dynamics: Hamiltonian Formulation} 

\author{Joseph T.A. Peterson\footnote{Department of Aerospace Engineering, Texas A\&M University, College Station, Texas 77843.}, 
Manoranjan Majji\footnote{Department of Aerospace Engineering, Texas A\&M University, College Station, Texas 77843.}, 
and John L.~Junkins\footnote{Department of Aerospace Engineering, Texas A\&M University, College Station, Texas 77843.}
}







\maketitle

\input{JOEcommands} 
\input{Mysecs_prj/00_abstract}

\renewcommand{\contentsname}{}
\tableofcontents  
\newpage

\input{Mysecs_prj/00_notation}
\input{Mysecs_prj/00_INTRO2} 

\input{Mysecs_prj/new_prjSUM_noElems} 

\input{Mysecs_prj/new_0_CTnonmin}

\input{Mysecs_prj/new_2_prjFrad}

\input{Mysecs_prj/new_3_2BP}

\input{Mysecs_prj/new_3.3_prjPerifocal}

\input{Mysecs_prj/new_3.4_prjSTM}

\input{Mysecs_prj/new_5_prjJ2}

\input{Mysecs_prj/99_conclusion}
\begin{appendices}
\renewcommand{\appendixname}{Appx.}
\titleformat{\section}{\Large\bfseries\sffamily}{\textmd{\large{Appendix \thesection}}}{0em}{\vskip .1\baselineskip}

\input{Mysecs_prj/new_1apx_prjCoord}

\input{Mysecs_prj/apx_angMoment}

\input{Mysecs_prj/apx_nonCon_coord}

\input{Mysecs_prj/apx_Pext_coord}

\input{Mysecs_prj/old_prjPlots}

\end{appendices}

\newpage
\clearpage
\addcontentsline{toc}{section}{References}
\begin{small}

\input{Ref_MAIN.bbl}
\end{small}

\end{document}

%% file: JOEcommands.tex


\newcommand{\pgrf}[1]{\noindent\textbf{#1}}
\newcommand\blfootnote[1]{%
  \begingroup
  \renewcommand\thefootnote{}\footnote{#1}%
  \addtocounter{footnote}{-1}%
  \endgroup
}

\newcommand{\cf}[0]{f}%
\newcommand{\cff}[0]{{\sl{f}}}

\newcommand{\tx}[1]{{\text{#1}}}
\newcommand{\txi}[1]{{\textnormal{\textit{#1}}}}
\newcommand{\tbf}[1]{{\textbf{\textup{\textrm{#1}}}}}
\newcommand{\trm}[1]{{\textup{\textrm{#1}}}}
\newcommand{\txnonb}[1]{{\textmd{#1}}} 
\newcommand{\txup}[1]{{\textup{#1}}} 
\newcommand{\bfi}[1]{{\textbf{\textit{#1}}}}
\newcommand{\Bfi}[1]{{ \scalebox{1.18}{\textbf{\textit{#1}}}  }}

\newcommand{\rmsb}[1]{{\textup{\textsb{#1}}}}%
\newcommand{\itsb}[1]{{\textit{\textsb{#1}}}}

\renewcommand{\sl}[1]{{\textnormal{\textrm{\textsl{#1}}}}}
\newcommand{\slb}[1]{{\textbf{\textsl{#1}}}}

\newcommand{\txsfb}[1]{\textbf{\textit{\textsf{#1}}}}
\newcommand{\sfi}[1]{{\textit{\textsf{#1}}}}%
\newcommand{\sfup}[1]{{\textup{\textsf{#1}}}}
\newcommand{\sfbup}[1]{ {\textbf{\textup{\textsf{#1}}}} }%

\renewcommand{\sc}[1]{{\textsc{#1}}}
\newcommand*{\rmsc}[1]{{\textrm{\textup{\textsc{#1}}}}}
\newcommand*{\bsc}[1]{{\textbf{\textsc{#1}}}}
\newcommand*{\rmbsc}[1]{{\textbf{\textrm{\textsc{#1}}}}}
\newcommand*{\itsc}[1]{{\textit{\textrm{\textsc{#1}}}}}
\newcommand*{\bisc}[1]{{\textit{\textbf{\textsc{#1}}}}}
\newcommand*{\sfsc}[1]{\textsf{\textit{\textsc{#1}}}}
\newcommand*{\sfbsc}[1]{\textbf{\textsf{\textit{\textsc{#1}}}}}
\newcommand{\scg}[0]{{\rmsc{g}}}
\newcommand{\sck}[0]{{\rmsc{k}}}
\newcommand{\scm}[0]{{\rmsc{m}}}
\newcommand{\scb}[0]{{\rmsc{b}}}
\newcommand{\qsc}[0]{{\itsc{q}}}

\newcommand{\kconst}[0]{k}

\newcommand{\isc}[1]{{\scriptstyle{\textit{\textrm{\textsc{#1}}}}}} 
\newcommand{\ssc}[1]{{\scalebox{0.6}{$#1$}}} 

\newcommand{\bemph}[1]{{\small{\itsb{#1}}}} 
 \newcommand{\sbemph}[1]{\rmsb{#1}} 


\newcommand{\smsize}[1]{ {\text{\small{#1}}} }
\newcommand*{\fnsize}[1]{ {\text{\footnotesize{#1}}} }%
\newcommand{\scrsize}[1]{ {\text{\scriptsize{#1}}} }
\newcommand{\nssize}[1]{ {\text{\normalsize{#1}}} }

\newcommand*{\nmsz}[1]{{\displaystyle{#1}}}
\newcommand*{\scrsz}[1]{{\scriptstyle{#1}}}
\newcommand*{\ssz}[1]{{\scriptstyle{#1}}}
\newcommand*{\ii}[1]{{\scriptscriptstyle{#1}}} 
\renewcommand*{\ss}[1]{ {\scalebox{0.7}{$#1$}} }

\newcommand*{\smsz}[1]{\text{\small{$#1$}} }
\newcommand*{\fnsz}[1]{\text{\footnotesize{$#1$}} }
\newcommand*{\ns}[1]{{\scalebox{1.15}{$#1$}} }
\renewcommand*{\lg}[1]{{\scalebox{1.25}{$#1$}} }

\newcommand{\mscale}[2][.7]{ {\text{\scalebox{#1}{$#2$}}} } 

\newcommand{\ttfrac}[2]{{\scriptstyle{\frac{\bs{#1}}{\bs{#2}}}}}


\newcommand{\til}[1]{\tilde{#1}}
\newcommand{\wt}[1]{ \widetilde{#1}  }
\newcommand{\wh}[1]{\widehat{#1}}
\newcommand{\ol}[1]{\overline{#1}}

\newcommand{\uln}[1]{{\underline{\smash{#1}\mkern-1mu}{\mkern1mu}}}
\newcommand{\ubar}[1]{ \underbar{$#1\hspace{-0.15mm}$}\hspace{0.15mm} }
\newcommand{\ub}[1]{{\underaccent{\text{—}}{\smash{#1}}}} 

\newcommand{\dt}[1]{\accentset{\mbox{\small\textbf{.}}}{#1}}
\newcommand*{\ddt}[1]{ \accentset{\mbox{\small\bfseries..}}{#1} }
\newcommand*{\Dt}[1]{\accentset{\mbox{\normalsize\bfseries .}}{#1}}
\newcommand*{\DDt}[1]{ \accentset{\mbox{\bfseries .\hspace{-0.03ex}.}}{#1} }

\newcommand{\rng}[1]{ \mathring{#1} }
\newcommand{\rrng}[1]{{\rlap{$\,\mathring{\vphantom{#1}}\mkern4.5mu\mathring{\vphantom{#1}}$}}#1}
\newcommand{\rring}[1]{ \accentset{\circ\circ}{#1} }
\newcommand{\ringp}[1]{\accentset{\circ'}{#1} }
\newcommand{\rringp}[1]{\accentset{\circ\circ'}{#1}}

\newcommand{\pdt}[1]{ \acute{#1} }
\newcommand{\pddt}[1]{ \rlap{$\;\,\acute{\vphantom{#1}}\!\acute{\vphantom{#1}}$}#1 }
\newcommand{\pdot}[1]{{\'{#1}}}
\newcommand{\pddot}[1]{{\H{#1}}}

\newcommand{\tridot}[1]{ \accentset{\triangle}{#1} }
\newcommand{\triddot}[1]{ \accentset{\triangle\triangle}{#1} }
\newcommand{\tridt}[1]{ \accentset{\blacktriangle}{#1} }
\newcommand{\triddt}[1]{ \accentset{\blacktriangle\!\blacktriangle}{#1} }
\newcommand{\deldot}[1]{ \accentset{\triangledown}{#1} }
\newcommand{\delddot}[1]{ \accentset{\triangledown\triangledown}{#1} }

\newcommand{\bxdot}[1]{ \accentset{\square}{#1} }
\newcommand{\bxddot}[1]{ \accentset{\square\square}{#1} }
\newcommand{\sqdt}[1]{ \accentset{\blacksquare}{#1} }
\newcommand{\sqddt}[1]{ \accentset{\blacksquare\blacksquare}{#1} }

\newcommand{\hdot}[1]{ \accentset{\mbox{\hspace{0.1mm}.}}{#1} }
\newcommand{\hddot}[2][.2ex]{\ddot{\raisebox{0pt}[\dimexpr\height+#1][\depth]{$#2$}}}
\newcommand{\ggvec}[2][-2pt]{\vec{\raisebox{0pt}[\dimexpr\height+#1][\depth]{$#2$}}}

\newcommand{\harp}[1]{ {\accentset{\rightharpoonup}{#1}} }
\newcommand{\varvec}[1]{ {\accentset{\to}{#1}} }


\let\oldcdot\cdot
\renewcommand*{\cdot}{{\mkern1.5mu\oldcdot\mkern1.5mu}}
\newcommand{\slot}[0]{{\oldcdot}} 
\newcommand{\cdt}[0]{ \pmb{\cdot} }
\newcommand{\bdt}[0]{ \bs{\cdot} }

\newcommand{\txblt}[0]{{\text{\small{\textbullet\hspace{.7ex}}}}}
\newcommand*{\sblt}[0]{ \raisebox{0.25ex}{$\scriptscriptstyle{\bullet}$}  }
\newcommand*{\bltt}[0]{ \text{\large{${\,\,\bullet\,\,}$}} }
\newcommand*{\nmblt}[0]{ \raisebox{0.15ex}{$\scriptstyle{\bullet}$}  }

\newcommand{\cddt}[0]{{\mkern2mu\pmb{\cdot}\mkern2mu}}
\newcommand{\cddot}[0]{ {\mkern2mu\uln{\pmb{\cdot}}\mkern2mu} }

\newcommand*{\otms}[0]{ {\raisebox{0.15ex}{$\scriptstyle{\,\otimes\,}$}}  }
\newcommand*{\opls}[0]{ \raisebox{0.15ex}{$\scriptstyle{\,\oplus\,}$}  }
\newcommand*{\wdg}[0]{ \raisebox{0.15ex}{$\scriptstyle{\,\wedge\,}$}  }
\newcommand*{\tms}[0]{ \raisebox{0.1ex}{$\scriptstyle{\,\times\,}$}  }
\newcommand{\botimes}[0]{ \scalebox{1.3}{$\otimes$} }
\newcommand{\boplus}[0]{ \scalebox{1.3}{$\oplus$} }
\newcommand{\bwedge}[1]{ {\scalebox{1}[1.3]{$\wedge$}}^{\!#1} }

\newcommand*{\ssqr}[0]{ \raisebox{0.15ex}{$\scriptstyle{\blacksquare}$}  }
\newcommand*{\iisqr}[0]{ \raisebox{0.25ex}{$\,\scriptscriptstyle{\blacksquare}\,$}  }
\newcommand*{\sbx}[0]{ \raisebox{0.15ex}{$\scriptstyle{\square}$}  }
\newcommand*{\iibx}[0]{ \raisebox{0.25ex}{$\,\scriptscriptstyle{\square}\,$}  }

\newcommand{\dimd}[0]{\diamond}%
\newcommand*{\tri}[0]{\triangle}%
\newcommand*{\trid}[0]{\triangledown} %
\newcommand*{\rtri}[0]{\triangleright} %
\newcommand*{\ltri}[0]{\triangleleft} %
\newcommand*{\stri}[0]{ \raisebox{0.25ex}{$\scriptstyle{\,\triangle\,}$} }
\newcommand*{\strid}[0]{\raisebox{0.25ex}{$\scriptstyle{\,\triangledown\,}$}}%
\newcommand*{\srtri}[0]{ \raisebox{0.25ex}{$\scriptstyle{\,\triangleright\,}$} }
\newcommand*{\slri}[0]{\raisebox{0.25ex}{$\scriptstyle{\,\triangleleft\,}$}}%

\newcommand*{\rnlarrow}[0]{ \overset{\nLeftarrow\;}{\Rightarrow}}
\newcommand*{\Rnlarrow}[0]{ \overset{\Leftarrow \!/ \!=\,}{\Longrightarrow}}
\newcommand*{\lnrarrow}[0]{ \overset{\,\nRightarrow}{\Leftarrow}}
\newcommand*{\Lnrarrow}[0]{ \overset{=\! /\!\Rightarrow\,}{\Longleftarrow}}
\newcommand{\nDownarrow}[0]{\rotatebox[origin=c]{-90}{$\nRightarrow$}}
\newcommand{\nUparrow}[0]{\rotatebox[origin=c]{90}{$\nRightarrow$}}


\newcommand*{\shrp}[0]{{\scalebox{0.55}{$\sharp$}}}
\newcommand*{\flt}[0]{{\scalebox{0.6}{$\flat$}}}

\newcommand{\iio}[0]{ {\scalebox{0.45}{$\bs{\circ}$}} }
\newcommand{\iix}[0]{{\scriptscriptstyle{\times}} }
\newcommand{\iistr}[0]{{\scriptscriptstyle{\star}} }
\newcommand{\str}[0]{{\scriptscriptstyle{\star}} }
\newcommand{\drg}[1]{ #1^\ss{\dagger} } 

\newcommand{\hdge}[1]{ #1^{\scriptscriptstyle{\star}}  }
\newcommand{\hhdge}[1]{ #1^{\scriptscriptstyle{\star\star}}  }
\newcommand{\hodge}[0]{\star}
\newcommand{\varhodge}[0]{ {\,^{\scriptscriptstyle{\bar{\star}}}} }
\newcommand{\ax}[1]{ #1^{\scriptscriptstyle{\times}} }
\newcommand{\axx}[1]{ #1^{\scriptstyle{\times}} }

\newcommand{\inv}[1]{ #1^{\scriptscriptstyle{-\!1}} }
\newcommand{\negg}[0]{{\textnormal{-}}}

 \newcommand{\trn}[1]{ #1^{\ii{\textsf{\textup{T}}}}  }
\newcommand{\invtrn}[1]{ #1^{\ii{\textsf{\textup{--T}}}} }
 
\newcommand*{\zer}[0]{ {\textrm{\tiny{\textit{0}}}} } %
\newcommand*{\zr}[0]{ {\ss{0}} }%
\newcommand{\nozer}[0]{ {\scriptscriptstyle{\emptyset}} }

\newcommand*{\sodot}[0]{ \raisebox{0.12ex}{$\scriptstyle{\odot}$}  }

\newcommand*{\nc}[0]{{\ii{\mrm{nc}}}}


\newcommand*{\mrm}[1]{ {\mathrm{#1}} } 
\newcommand*{\iimrm}[1]{ {\ii{\mathrm{#1}}} }
\newcommand*{\iirm}[1]{ {\ii{\textrm{#1}}} }

\newcommand*{\mbf}[1]{{\mathbf{#1}}} 

\newcommand{\rmb}[1]{\tbf{#1}}


\newcommand{\bs}[1]{{\boldsymbol{#1}}}
\newcommand{\bsi}[1]{\boldsymbol{\mit{#1}}} 
\newcommand{\accalign}[2]{ \mkern2mu #1{\mkern-2mu #2} }
\newcommand{\tbs}[1]{ {\accalign{\tilde}{\bs{#1}}} }%
\newcommand{\wtbs}[1]{ {\accalign{\widetilde}{\bs{#1}}} }%
\newcommand{\hbs}[1]{ {\accalign{\hat}{\bs{#1}}} }%
\newcommand{\whbs}[1]{ {\accalign{\widehat}{\bs{#1}}} }%
\newcommand{\barbs}[1]{ {\accalign{\bar}{\bs{#1}}} }%
\newcommand{\vecbs}[1]{ {\accalign{\vec}{\bs{#1}}} }%
\newcommand*{\dtbs}[1]{ {\accalign{\dt}{\bs{#1}}} }%
\newcommand*{\ddtbs}[1]{ {\accalign{\ddt}{\bs{#1}}} }%

\newcommand{\tbm}[1]{\tilde{\bm{#1}}}
\newcommand{\wtbm}[1]{\widetilde{\bm{#1}}}
\newcommand{\hbm}[1]{\hat{\bm{#1}}}
\newcommand{\whbm}[1]{\widehat{\bm{#1}}}
\newcommand{\barbm}[1]{\bar{\bm{#1}}}

\newcommand*{\nbs}[1]{ \scalebox{1.15}{$\boldsymbol{#1}$} }
\newcommand*{\lbs}[1]{ \scalebox{1.4}{$\boldsymbol{#1}$} }
\newcommand*{\nbm}[1]{ \scalebox{1.2}{$\bm{#1}$} }
\newcommand*{\lbm}[1]{ \scalebox{1.4}{$\bm{#1}$} }


\newcommand{\msfb}[1]{{\mathsfb{#1}}} 

\newcommand{\sfb}[1]{{\txsfb{#1}}} 

\renewcommand{\accalign}[2]{ \mkern1.5mu #1{\mkern-1.5mu #2} }
\newcommand{\tsfb}[1]{ {\accalign{\tilde}{\sfb{#1}}} }%
\newcommand{\wtsfb}[1]{ {\accalign{\widetilde}{\sfb{#1}}} }%
\newcommand{\hsfb}[1]{ {\accalign{\hat}{\sfb{#1}}} }%
\newcommand{\whsfb}[1]{ {\accalign{\widehat}{\sfb{#1}}} }%
\newcommand{\barsfb}[1]{ {\accalign{\bar}{\sfb{#1}}} }%
\newcommand{\bsfb}[1]{ {\accalign{\bar}{\sfb{#1}}} }%
\newcommand{\vecsfb}[1]{ {\accalign{\vec}{\sfb{#1}}} }%
\newcommand*{\dtsfb}[1]{ {\accalign{\dt}{\sfb{#1}}} }%
\newcommand*{\ddtsfb}[1]{ {\accalign{\ddt}{\sfb{#1}}} }%

\newcommand*{\sfg}[0]{{\bs{g}}}


\newcommand{\tn}[1]{{\sfb{#1}}} 

\newcommand*{\mbb}[1]{ {\mathbb{#1}} } 
\newcommand*{\tmbb}[1]{ \tilde{\mbb{#1}} }
\newcommand*{\wtmbb}[1]{ \widetilde{\mbb{#1}} }
\newcommand*{\hmbb}[1]{ \hat{\mbb{#1}} }
\newcommand*{\whmbb}[1]{ \widehat{\mbb{#1}} }
\newcommand*{\bbk}[0]{\Bbbk} 

\newcommand*{\mfrak}[1]{ {\mathfrak{#1}} }
\newcommand*{\bfrak}[1]{ {\bs{\mathfrak{#1}}} }
\newcommand*{\tmfrak}[1]{ \tilde{\mfrak{#1}} }
\newcommand*{\wtmfrak}[1]{ \widetilde{\mfrak{#1}} }
\newcommand*{\hmfrak}[1]{ \hat{\mfrak{#1}} }
\newcommand*{\whmfrak}[1]{ \widehat{\mfrak{#1}} }


\newcommand*{\mscr}[1]{{\mathscr{#1}}}%
\newcommand*{\bscr}[1]{ \bs{\mathscr{#1}} }
\newcommand*{\tscr}[1]{\tilde{\mathscr{#1}}}
\newcommand*{\wtscr}[1]{\widetilde{\mathscr{#1}}}
\newcommand*{\hscr}[1]{ \hat{\mathscr{#1}} }
\newcommand*{\whscr}[1]{ \widehat{\mathscr{#1}} }
\newcommand*{\barscr}[1]{ \bar{\mathscr{#1}} }
\newcommand*{\olscr}[1]{ \overline{\mathscr{#1}} }

\newcommand*{\mcal}[1]{ {\mathcal{#1}} }
\newcommand*{\bcal}[1]{ \bs{\mathcal{#1}} }
\newcommand*{\tcal}[1]{\tilde{\mathcal{#1}}}
\newcommand*{\wtcal}[1]{\widetilde{\mathcal{#1}}}
\newcommand*{\hcal}[1]{ \hat{\mathcal{#1}} }
\newcommand*{\whcal}[1]{ \widehat{\mathcal{#1}} }
\newcommand*{\barcal}[1]{ \bar{\mathcal{#1}} }
\newcommand*{\olcal}[1]{ \overline{\mathcal{#1}} }

\newcommand{\iical}[1]{{\scriptscriptstyle{\mcal{#1}}}}
\newcommand{\sscal}[1]{{{}^{{}_{\mathcal{#1}}}}} 
\newcommand{\iiscr}[1]{ {\scriptscriptstyle{\mscr{#1}}} }
\newcommand*{\sscr}[1]{ {\scalebox{0.7}{$\mscr{#1}$}} } 



\newcommand{\en}[0]{{\itsc{n}}}
\newcommand{\emm}[0]{{\itsc{m}}} 
\newcommand{\envec}[0]{ \hbe_\en }
\newcommand{\enform}[0]{\,\hat{\!\bs{\epsilon}}^\en}
\newcommand{\tilp}[0]{{\mkern2mu \til{\mkern-2mu\smash{p}}}}

\newcommand{\pu}[0]{\mu}  
\newcommand{\varpu}[0]{\pi}  

\newcommand*{\tp}[1]{{\vphantom{#1}\uln{#1}}} 
\newcommand*{\btp}[1]{ {\bar{\tp{#1}}} }
\newcommand*{\bartp}[1]{ {\bar{\tp{#1}}} }
\newcommand*{\tiltp}[1]{ {\til{\tp{#1}}} }
\newcommand*{\ttp}[1]{ {\til{\tp{#1}}} }
\newcommand*{\wttp}[1]{ {\wt{\tp{#1}}} }
\newcommand*{\htp}[1]{ {\hat{\tp{#1}}} }
\newcommand*{\whtp}[1]{ {\wh{\tp{#1}}} }
\newcommand{\dottp}[1]{{\dot{\tp{#1}}} }

\newcommand*{\tup}[1]{{\bs{\mathrm{#1}}}}
 \newcommand*{\bartup}[1]{ {\bar{\tup{#1}}} } 
\newcommand*{\tiltup}[1]{ {\til{\tup{#1}}} }
\newcommand*{\ttup}[1]{ {\til{\tup{#1}}} }
\newcommand*{\wttup}[1]{ {\wt{\tup{#1}}} }
\newcommand*{\htup}[1]{ {\hat{\tup{#1}}} }
\newcommand*{\whtup}[1]{ {\wh{\tup{#1}}} }
\newcommand{\dottup}[1]{{\dot{\tup{#1}}} }
\newcommand{\nrmtup}[1]{ {#1} } %


\newcommand*{\pt}[1]{{#1}}
\newcommand{\tilpt}[1]{ {\til{\pt{#1}}} }
\newcommand*{\barpt}[1]{ {\bar{\pt{#1}}} }
\newcommand{\hatpt}[1]{ {\hat{\pt{#1}}} }
\newcommand{\hpt}[1]{ {\hat{\pt{#1}}} }
\newcommand*{\ptvec}[1]{ {\vec{\pt{#1}}} }
\newcommand*{\dtpt}[1]{\dt{\pt{#1}}}


\renewcommand{\d}[0]{\mrm{d}}
\newcommand*{\rmd}[0]{{\text{\rmsb{d}}}}
\newcommand*{\dif}[0]{ {\textbf{\textsf{\textsl{d}}}} }
\newcommand*{\Dif}[0]{{ \scriptstyle{\textbf{\textsf{\textsl{D}}}} }}

\newcommand{\difbar}[0]{{ \txsfb{\kern0.6ex\ooalign{\hidewidth\raisebox{0.63ex}{--}\hidewidth\cr{\kern-0.6ex\textsf{d}}\cr}} \vphantom{d} }}
\newcommand{\itdifbar}{%
   {\mkern4.5mu \txsfb{\ooalign{\hidewidth\raisebox{0.65ex}{--}\hidewidth\cr$\mkern-4.5mu \bs{d}$\cr}} \vphantom{d} }%
}

\newcommand*{\hdif}[0]{  \hspace{1mm}\hat{\hspace{-1mm}\dif} }
\newcommand*{\tdif}[0]{  \hspace{1mm}\tilde{\hspace{-1mm}\dif} }
\newcommand*{\deldif}[0]{\bs{\delta}}

\newcommand{\tinywedge}[0]{ \scalebox{0.45}{$\pmb{\wedge}$}}
\newcommand{\exd}[0]{\dif_{\hspace{-0.25mm}{\tinywedge}}}
\newcommand{\exdbar}[0]{\difbar_{\hspace{-0.25mm}{\tinywedge}}}

\newcommand{\codif}[0]{ {\exd^{\scriptscriptstyle{\star}}} }

\newcommand*{\fibd}[0]{{\sfb{F}}}%

\newcommand{\pderiv}[2]{ \tfrac{\uppartial #1}{\uppartial #2} }
\newcommand{\ppderiv}[3]{\tfrac{\uppartial^2 #1}{\uppartial #2 \uppartial #3}}

\newcommand{\pder}[2]{ \tfrac{ \displaystyle{\uppartial #1}}{\displaystyle{\uppartial #2}} }  
\newcommand{\pderr}[2]{ \tfrac{ \textstyle{\uppartial #1}}{\textstyle{\uppartial #2}} }  
\newcommand{\ppder}[3]{ \tfrac{\displaystyle{\uppartial^2 #1}}{\displaystyle{\uppartial #2 \uppartial #3}} }

\newcommand{\Pderiv}[2]{\frac{\uppartial #1}{\uppartial #2}}
\newcommand{\PPderiv}[3]{\frac{\uppartial^2 #1}{\uppartial #2 \uppartial #3}}

 \newcommand*{\diff}[2]{ \tfrac{\mathrm{d} #1}{\mathrm{d} #2} }
\newcommand{\ddiff}[2]{ \tfrac{\mathrm{d}^2{#1}}{\mathrm{d} {#2}^2} }
\newcommand{\Diff}[2]{ \frac{\mathrm{d}#1}{\mathrm{d}\,#2} }
\newcommand{\DDiff}[2]{ \frac{\mathrm{d}^2{#1}}{\mathrm{d}\,{#2}^2} }

\newcommand{\fdiff}[1]{ \,\tfrac{ {\!\!\!\!}^{ \text{\footnotesize{$#1$}} }\d }{\;\d\, t} }
\newcommand{\fDiff}[3]{ \,\frac{ \!\!\!\!^{ \text{\footnotesize{$#3$}} }\d #1 }{\;\d\, #2} }
\newcommand{\fddiff}[1]{ \,\tfrac{ \!\!\!\!^{ \text{\scriptsize{$#1$}} }\d^2 }{\;\d\, t^2} }
\newcommand{\fDDiff}[3]{ \,\frac{ \!\!\!\!^{ \text{\footnotesize{$#3$}} }\d^2 #1 }{\;\d\, #2^2} }

\newcommand{\mdiff}[2]{ \tfrac{ \mathrm{D} #1}{\mathrm{d} #2}   }
\newcommand{\nabdiff}[2]{ \tfrac{\nabla #1}{\mathrm{d} #2}   }
\newcommand{\nabpderiv}[2]{ \tfrac{\nabla #1}{\uppartial #2}   }

\newcommand{\difscale}[1]{\scalebox{0.72}{$#1$}}
\newcommand{\lderiv}[1]{ {\textit{\textrm{\pounds}}}_{\difscale{#1}} } 
\newcommand{\nab}[1][\mkern1mu]{ \nabla_{\mkern-4mu \difscale{#1}\mkern1mu} }
\newcommand{\del}[1]{ \nabla_{\hspace{-0.6mm} {#1}}  }


\newcommand{\pmat}[1]{ \begin{pmatrix} #1 \end{pmatrix} }
\newcommand{\bmat}[1]{ \begin{bmatrix} #1 \end{bmatrix} }
\newcommand{\fnpmat}[1]{ \fnsz{\begin{pmatrix} #1 \end{pmatrix}} }
\newcommand{\smpmat}[1]{ \smsz{\begin{pmatrix} #1 \end{pmatrix}} }
\newcommand{\fnbmat}[1]{ \fnsz{\begin{bmatrix} #1 \end{bmatrix}} }

\newcommand*{\kd}[0]{1}
\newcommand*{\lc}[0]{{\itsc{e}}}
\newcommand{\ibase}[0]{\hat{\rmsb{\oldstylenums{1}}}}
\newcommand*{\imat}[0]{{\text{1}}\hspace{-.2mm}}
\newcommand*{\jmat}[0]{J}

\newcommand{\Id}[0]{{\mathrm{Id}}}%
\newcommand{\iden}[0]{\sfb{1}}
\newcommand*{\gvol}[0]{\epsilon}
\newcommand*{\spvol}[0]{\sigma} 
\newcommand*{\lagvol}[0]{\varsigma}
\newcommand*{\spform}[0]{\theta}  
\newcommand*{\spformup}[0]{\thetaup}
\newcommand*{\lagform}[0]{\vartheta} 
\newcommand*{\lagformup}[0]{\varthetaup}

\newcommand{\vf}[0]{v} 
\newcommand{\pf}[0]{p} 
\newcommand{\vlin}[0]{{\scalebox{0.7}{$\mathcal{V}$}}}
\newcommand{\plin}[0]{{\pi}} 
\newcommand{\uplin}[0]{{\piup}} 
\newcommand{\tilplin}[0]{{\mkern2mu \til{\mkern-2mu\smash{\pi}}}} 


\newcommand{\rmat}[1]{ \mathrm{M}^{#1}_\ii{\mbb{R}} }
\newcommand{\mats}[1]{ \rmat{#1} }
\newcommand{\mat}[2]{ \mathrm{M}^{#1}_\ii{#2} }
\newcommand{\lgrp}[3]{\mathrm{#1}^{#2}_\ii{\mbb{#3}}}
\newcommand{\lalg}[3]{\mathfrak{#1}^{#2}_\ii{\mbb{#3}}}
\newcommand{\Afmat}[1]{ \mathrm{Af}^{#1}_\ii{\mbb{R}} }
\newcommand{\afmat}[1]{ \mathfrak{af}^{#1}_\ii{\mbb{R}} }
\newcommand{\Glmat}[1]{ \mathrm{Gl}^{#1}_\ii{\mbb{R}} }
\newcommand{\Slmat}[1]{ \mathrm{Sl}^{#1}_\ii{\mbb{R}} }
\newcommand{\glmat}[1]{ \mathfrak{gl}^{#1}_\ii{\mbb{R}} }
\newcommand{\slmat}[1]{ \mathfrak{sl}^{#1}_\ii{\mbb{R}} }
\newcommand{\Spmat}[1]{ \mathrm{Sp}^{#1}_\ii{\mbb{R}} }
\newcommand{\spmat}[1]{ \mathfrak{sp}^{#1}_\ii{\mbb{R}} }
\newcommand{\Omat}[1]{ \mathrm{O}^{#1}_\ii{\mbb{R}} }
\newcommand{\Somat}[1]{ \mathrm{SO}^{#1}_\ii{\mbb{R}} }
\newcommand{\omat}[1]{ \mathfrak{o}^{#1}_\ii{\mbb{R}} }
\newcommand{\somat}[1]{ \mathfrak{so}^{#1}_\ii{\mbb{R}} }
\newcommand{\Umat}[1]{ \mathrm{U}^{#1}_\ii{\mbb{C}} }
\newcommand{\Sumat}[1]{ \mathrm{SU}^{#1}_\ii{\mbb{C}} }
\newcommand{\umat}[1]{ \mathfrak{u}^{#1}_\ii{\mbb{C}} }
\newcommand{\sumat}[1]{ \mathfrak{su}^{#1}_\ii{\mbb{C}} }
\newcommand{\Emat}[1]{ \mathrm{E}^{#1}_\ii{\mbb{R}} }
\newcommand{\Semat}[1]{ \mathrm{SE}^{#1}_\ii{\mbb{R}} }
\newcommand{\emat}[1]{ \mathfrak{e}^{#1}_\ii{\mbb{R}} }
\newcommand{\semat}[1]{ \mathfrak{se}^{#1}_\ii{\mbb{R}} }

\newcommand{\Aften}[0]{ \mathrm{Af} }
\newcommand{\aften}[0]{ \mathfrak{af} }
\newcommand{\Glten}[0]{ \mathrm{Gl} }
\newcommand{\Slten}[0]{ \mathrm{Sl} }
\newcommand{\glten}[0]{ \mathfrak{gl} }
\newcommand{\slten}[0]{ \mathfrak{sl} }
\newcommand{\Spten}[0]{ \mathrm{Sp}}
\newcommand{\spten}[0]{ \mathfrak{sp} }
\newcommand{\Oten}[0]{ \mathrm{O} }
\newcommand{\Soten}[0]{ \mathrm{SO} }
\newcommand{\oten}[0]{ \mathfrak{o} }
\newcommand{\soten}[0]{ \mathfrak{so} }
\newcommand{\Eten}[0]{ \mathrm{E} }
\newcommand{\Seten}[0]{ \mathrm{SE} }
\newcommand{\eten}[0]{ \mathfrak{e} }
\newcommand{\seten}[0]{ \mathfrak{se} }

\newcommand{\Spism}[0]{\mathscr{S}\mkern-1mu p}
\newcommand{\Dfism}[0]{\mathscr{D}\mkern-2mu f}
\newcommand{\Hmism}[0]{\mathscr{H}m }
\newcommand{\Ctism}[0]{\mathscr{C}t}
\newcommand{\Isom}[0]{\mathcal{I}\mkern-3mu s}


\newcommand{\man}[1]{ \mathcal{#1} }  
\newcommand{\tman}[1]{ \til{\man{#1}} } 
\newcommand{\barman}[1]{ \bar{\man{#1}} }
\newcommand{\bman}[1]{ \bar{\man{#1}} }
\newcommand{\hman}[1]{ \hat{\man{#1}} }
\newcommand{\whman}[1]{ \widehat{\man{#1}} }
\newcommand{\sman}[1]{\sl{#1}} 
\newcommand{\tint}[0]{\mathcal{I}} 

\newcommand{\chart}[2]{ \man{#1}_{\hspace{-0.2mm}\textup{\tiny(} {\scriptscriptstyle{#2}} \textup{\tiny)}} }
\newcommand{\rchart}[2]{ \mathbb{#1}_{\hspace{-0.2mm}\textup{\tiny(} {\scriptscriptstyle{#2}} \textup{\tiny)}} }
\newcommand{\chrt}[2]{ {#1}_{\hspace{-0.2mm}\textup{\tiny(} {\scriptscriptstyle{#2}} \textup{\tiny)}} }

\newcommand{\vsp}[1]{ \mathbb{#1} } 
\newcommand{\aff}[1]{\man{#1}} 
\newcommand{\tvsp}[1]{ \tilde{\vsp{#1}} } 
\newcommand{\bvsp}[1]{ \bar{\vsp{#1}} } 
\newcommand{\taff}[1]{ \tilde{\aff{#1}} } 
\newcommand{\baff}[1]{ \bar{\aff{#1}} } 

\newcommand{\affE}[0]{\man{E}}
\newcommand{\vecE}[0]{\vsp{E}}
\newcommand{\Eaf}[0]{\man{E}}
\newcommand{\Evec}[0]{\vsp{E}}
\newcommand{\emet}[0]{I} 
\newcommand{\bemet}[0]{\sfb{I}}
\newcommand{\rfun}[0]{{r}} 
\newcommand{\lang}[0]{{\ell}} 
\newcommand{\Lang}[0]{{L}}
\newcommand{\langb}[0]{\bm{\ell}} 
\newcommand{\Langb}[0]{\sfb{L}} 
\newcommand{\slang}[0]{{\ell}} 
\newcommand{\slangup}[0]{{\upell}}
\newcommand{\Slang}[0]{{L}} 
\newcommand{\slangb}[0]{\bm{\ell}}
\newcommand{\Slangb}[0]{\sfb{L}}

\newcommand{\blang}[0]{\langb} 
\newcommand{\bLang}[0]{\Langb} 
\newcommand{\bslang}[0]{\slangb} 
\newcommand{\bSlang}[0]{\Slangb}

\newcommand*{\fun}[0]{{\scalebox{0.85}[1]{$\mathcal{F}\mkern-2mu$}}}  
\newcommand*{\tens}[0]{ \mathscr{T} } 
\newcommand*{\forms}[0]{ \Lambdaup }
\newcommand{\formsex}[0]{ \Lambdaup_\ii{\mrm{ex}\!}}
\newcommand{\formscl}[0]{ \Lambdaup_\ii{\mrm{cl}\!}}

\newcommand*{\veckl}[0]{ \mathfrak{X}_{\scriptscriptstyle{\!\mathfrak{i\!s}}\!} }
\newcommand*{\vect}[0]{ \mathfrak{X} } 
\newcommand*{\vecsp}[0]{ \mathfrak{X}_{\scriptscriptstyle{\!\mathfrak{s\!p}}\!\!} } 
\newcommand*{\vechm}[0]{ \mathfrak{X}_{\scriptscriptstyle{\!\mathfrak{h\!\!m}}\!\!} } 
\newcommand{\ham}[1]{{\mathscr{#1}}}
\newcommand{\tham}[1]{\widetilde{\mathscr{#1}}}

 \newcommand{\Tspacefont}[1]{ {\mathrm{#1}} } 
\newcommand{\Tan}[0]{ \Tspacefont{T}   }
\newcommand{\Tanh}[0]{  \Tspacefont{H}  }
\newcommand{\Tanv}[0]{  \Tspacefont{V}  }
\newcommand{\tsp}[1][\mkern2mu]{ \Tspacefont{T}_{\mkern-2mu #1} }
\newcommand{\tspv}[1][\mkern2mu]{\Tspacefont{V}_{\!#1}}
\newcommand{\tsph}[1][]{\Tspacefont{H}_{#1}}
\newcommand{\tspn}[1][\hspace{.5mm}]{\Tspacefont{N}_{\hspace{-.2mm}#1}}
\newcommand{\cotsp}[1][\hspace{1mm}]{ {\Tspacefont{T}^* {\hspace{-1ex}}_{\hspace{-.4mm}#1\hspace{0.5mm}}} } 
\newcommand{\cotspv}[1][\hspace{1mm}]{ \Tspacefont{V}^* {\hspace{-1ex}}_{\!#1\hspace{0.5mm}} }
\newcommand{\cotsph}[1][\hspace{1mm}]{ \Tspacefont{H}^* {\hspace{-1ex}}_{#1\hspace{0.1mm}} }
\newcommand{\cotspn}[1][\hspace{1mm}]{ \Tspacefont{N}^* {\hspace{-1ex}}_{#1\hspace{0.1mm}} }

\newcommand{\prj}[0]{ \scalebox{0.7}{$\Pi$} }  
 \newcommand{\tpr}[0]{ \scalebox{0.5}[0.68]{$\bs{\mcal{T}}$} }
 \newcommand{\copr}[0]{{ \mkern1mu\hat{\mkern-1mu\smash{\tpr}} }}

\newcommand{\vfun}[1][]{{V^{#1}}}
\newcommand{\pfun}[1][]{{P^{#1}}}
\newcommand{\varpfun}[1][]{{l^{#1}}}

\newcommand*{\tlift}[0]{ T\hspace{-.3mm} }
\newcommand*{\colift}[0]{ \hat{T}\hspace{-.3mm} }
\newcommand*{\uptlift}[0]{ \mrm{T}\hspace{-.3mm} }
\newcommand*{\upcolift}[0]{ \hat{\mrm{T}}\hspace{-.3mm} }

\newcommand{\formsh}[0]{\Lambdaup_\ii{\mrm{h}\!}} 
\newcommand{\formsv}[0]{\Lambdaup_\ii{\mrm{v}\!}} 
\newcommand{\vectv}[0]{\mathfrak{X}_\ii{\mrm{v}\!}} 
\newcommand{\vecth}[0]{\mathfrak{X}_\ii{\mrm{h}\!}} 
\newcommand{\formsbh}[0]{\Lambdaup_\ii{\mrm{bh}\!}} 
\newcommand{\vectbv}[0]{\mathfrak{X}_\ii{\mrm{bv}\!}}

\newcommand{\lft}[1]{ #1^{\hspace{-0.2mm}\scriptscriptstyle\upharpoonright} } 
\newcommand{\lift}[1]{ #1^{\hspace{-0.2mm}\scriptscriptstyle\uparrow} } 
\newcommand{\cotlft}[1]{  #1^{\hspace{-0.2mm}\hat{\scriptscriptstyle{\uparrow}}} }
\newcommand{\invlift}[1]{ #1^{\hspace{-0.2mm}\scriptscriptstyle\downarrow} }

\newcommand{\vlift}[1]{ #1^{\hspace{-0.2mm}\scriptscriptstyle\uparrow} }
\newcommand{\Vlift}[1]{  #1^{\scriptscriptstyle\Uparrow} }
\newcommand{\coVlift}[1]{  #1^{\hat{\scriptscriptstyle{\Uparrow}}} }%

\newcommand{\covlift}[1]{ #1^{\hat{\textup{\tiny{\textsf{V}}}}}}
\newcommand{\verlift}[1]{ #1^{\textup{\tiny{\textsf{V}}}} }%
\newcommand{\coverlift}[1]{ #1^{\hat{\textup{\tiny{\textsf{V}}}}}}

\newcommand{\hlift}[1]{ #1^{\textup{\tiny{\textsf{H}}}} }%
\newcommand{\horlift}[1]{ #1^{\textup{\tiny{\textsf{H}}}} }%
\newcommand{\cohlift}[1]{ #1^{\hat{\textup{\tiny{\textsf{H}}}}} }%

 \newcommand{\dubuparrow}[0]{\rotatebox[origin=c]{90}{$\ii{\twoheadrightarrow}$}}
\newcommand{\tailuparrow}[0]{\rotatebox[origin=c]{90}{$\ii{\rightarrowtail}$}}
\newcommand{\lftt}[1]{#1^{\dubuparrow}}
\newcommand{\llft}[1]{#1^{\tailuparrow}}

\newcommand{\rhk}[0]{\hookrightarrow}
\newcommand{\lhk}[0]{\hooklefttarrow}
\newcommand{\subman}[0]{ \,\accentset{\scriptstyle{\subset}}{\scriptstyle{\hookrightarrow}}\, }

\newcommand{\pbrak}[2]{\{#1  ,  #2 \}}
\newcommand{\Pbrak}[2]{\big\{#1\;,\;#2\big\}}
\newcommand{\ppbrak}[2]{\{ \hspace{-0.8mm} \{#1 ,  #2 \} \hspace{-0.8mm} \}}

\newcommand{\lbrak}[2]{ [ #1 , #2 ] }
\newcommand{\Lbrak}[2]{ \big[ #1\;,\;#2 \big] }
\newcommand{\llbrak}[2]{ [\hspace{-0.8mm}[#1 ,  #2 ]\hspace{-0.8mm}] }

\newcommand{\altbrak}[2]{ {\lfloor #1 ,  #2 \rfloor} }
\newcommand{\aaltbrak}[2]{ {\lfloor \hspace{-0.8mm} \lfloor #1 ,  #2 \rfloor \hspace{-0.8mm} \rfloor} }
\newcommand{\lagbrak}[2]{ {\lfloor #1 ,  #2 \rfloor}  }
\newcommand{\llagbrak}[2]{ {\lfloor \hspace{-0.8mm} \lfloor #1 ,  #2 \rfloor \hspace{-0.8mm} \rfloor} }

\newcommand{\inner}[2]{\langle#1  , #2\rangle}
\newcommand{\iinner}[2]{\langle\!\langle#1\,,\,#2\rangle\!\rangle}

\newcommand{\ang}[1]{\langle #1 \rangle}
\newcommand{\Ang}[1]{\big\langle #1 \big\rangle }


\newcommand{\abs}[1]{ \big| #1 \big| }
\newcommand{\Abs}[1]{ \left| #1 \right| }
\renewcommand{\mag}[1]{ |\!| #1 |\!| }
\newcommand{\nrm}[1]{{\lvert #1 \rvert} }
\renewcommand{\det}[0]{ {\textrm{\textup{\footnotesize{det}}}\hspace{0.5mm}} }
\renewcommand{\dim}[0]{ {\textrm{\textup{\footnotesize{dim}}}\hspace{0.5mm}} }
\newcommand{\codim}[0]{ {\textrm{\textup{\footnotesize{codim}}}\hspace{0.5mm}} } %
\newcommand{\rnk}[0]{ {\textrm{\textup{\small{rnk}}}\hspace{0.5mm}} }
\newcommand{\sgn}[0]{ {\textrm{\textup{\small{sgn}}}\hspace{0.5mm}} }
\newcommand{\tr}[0]{ {\textrm{\textup{\small{tr}}}\hspace{0.5mm}} }
\renewcommand{\ker}[0]{ {\textrm{\textup{\footnotesize{ker}}}\hspace{0.5mm}} }
\newcommand{\img}[0]{ {\textrm{\textup{\footnotesize{im}}}\hspace{0.5mm}} }
\newcommand{\spn}[0]{ {\textrm{\textup{\footnotesize{span}}}\hspace{0.5mm}} }
\renewcommand{\div}[1][\,]{ {{\textrm{\textup{\footnotesize{div}}}}_{\!#1}}}

\newcommand{\snn}[1][]{\sin{{\mkern-2mu}#1}}
\newcommand{\csn}[1][]{\cos{{\mkern-2mu}#1}}

\newcommand{\crd}[2]{ \tensor*[^{\ss{#2}}]{[#1]}{} }
\newcommand{\crdl}[2]{ \tensor*[^{\ss{#2}}]{#1}{} }
\newcommand{\cordl}[2]{  #1_{\scriptscriptstyle{\!#2}}  } 
\newcommand{\cord}[2]{  #1_{\scriptscriptstyle{#2}}  }

\newcommand{\eval}[2]{\left.{#1}\right|_{#2}}

\newcommand{\pdup}[0]{\uppartial} %
\newcommand{\upd}[0]{{\rotatebox[origin=t]{10}{$\partial$} \mkern-2mu}}
\newcommand{\pd}[0]{{\rotatebox[origin=t]{15}{$\partial$} \mkern-1mu}}


\newcommand{\mypd}[0]{\rotatebox[origin=t]{10}{$\partial$}} 

\newcommand*{\ffsize}[1]{{\scalebox{0.7}{$#1$}}}

\newcommand*{\bpart}[1]{\bs{\mypd}_{\hspace{-.3mm}#1}}
\newcommand*{\bpartup}[1]{\bs{\mypd}^{#1}}
\newcommand*{\tbpart}[1]{\tilde{\bs{\mypd}}_{\hspace{-.3mm}#1}}
\newcommand*{\tbpartup}[1]{{\tilde{\bs{\mypd}}\vphantom{l}^{#1}}}
\newcommand*{\hbpart}[1]{\hbs{\pd}_{\hspace{-.3mm}#1}}
\newcommand*{\hbpartup}[1]{{\hbs{\pd}\vphantom{l}^{#1}}}

\newcommand*{\pdii}[1]{ \bs{\mypd}_{\ffsize{\!#1}}}
\newcommand*{\pdupii}[1]{ \bs{\mypd}^{\ffsize{#1}} }
\newcommand*{\pdiiup}[1]{ \bs{\mypd}^{\ffsize{#1}} }
\newcommand*{\tpdii}[1]{ \tilde{\bs{\mypd}}_{\hspace{-.5mm}\ffsize{#1}} }
\newcommand*{\tpdupii}[1]{ \tilde{\bs{\mypd}}{}^{\ffsize{#1}} }
\newcommand*{\hpdii}[1]{ \hbs{\pd}_{\ffsize{\!#1}} }
\newcommand*{\hpdiiup}[1]{ {\hbs{\pd}}^{\ffsize{#1}} }
\newcommand*{\barpdii}[1]{ \bar{\bs{\mypd}}_{\hspace{-.6mm}\ffsize{#1}} }
\newcommand*{\barpdiiup}[1]{ {\bs{\bar{\pd}}}^{\ffsize{#1}} }

\newcommand{\bpd}[1][]{ {\bs{\mypd}_{\ffsize{\!#1}}} }%
\newcommand{\bpdup}[1][]{ {\bs{\mypd}^{\ffsize{#1}}} }
\newcommand{\hbpd}[1][]{ {\hbs{\pd}_{\ffsize{\!#1}}} }%
\newcommand{\hbpdup}[1][]{ {\hbs{\pd}^{\ffsize{#1}}} }
\newcommand{\tbpd}[1][]{ {\tbs{\pd}_{\ffsize{\!#1}}} }%
\newcommand{\tbpdup}[1][]{ {\tbs{\pd}^{\ffsize{#1}}} }
\newcommand{\barbpd}[1][]{ {\barbs{\pd}_{\ffsize{\!#1}}} }%
\newcommand{\barbpdup}[1][]{ {\barbs{\pd}^{\ffsize{#1}}} }

\newcommand*{\bpdh}[1]{ {\bs{\mypd}_{\mkern-1mu\hat{#1}}} }
\newcommand*{\bpdhup}[1]{ {\bs{\mypd}^{\hat{#1}}} }
\newcommand*{\bpdt}[1]{ {\bs{\mypd}_{\mkern-1mu\til{#1}}} }
\newcommand*{\bpdtup}[1]{ {\bs{\mypd}^{\til{#1}}} }
\newcommand*{\bpdbar}[1]{ {\bs{\mypd}_{\mkern-1mu\bar{#1}}} }
\newcommand*{\bpdbarup}[1]{ {\bs{\mypd}^{\bar{#1}}} }
\newcommand*{\bdelh}[1]{ \bs{\delta}^{\hat{#1}}\vphantom{l} }
\newcommand{\bdelhdn}[1][]{ {\bs{\delta}_{\mkern-1mu\hat{#1}}\vphantom{l}} }
\newcommand*{\bdelt}[1]{ \bs{\delta}^{\til{#1}}\vphantom{l} }
\newcommand{\bdeltdn}[1][]{ {\bs{\delta}_{\mkern-1mu\til{#1}}\vphantom{l}} }
\newcommand*{\bdelbar}[1]{ \bs{\delta}^{\bar{#1}}\vphantom{l} }
\newcommand{\bdelbardn}[1][]{ {\bs{\delta}_{\mkern-1mu\bar{#1}}\vphantom{l}} }



\newcommand{\bdel}[1][]{ {\bs{\delta}^{\ffsize{#1}}} }%
\newcommand{\bdeldn}[1][]{ {\bs{\delta}_{\ffsize{\!#1}}} }
\newcommand{\hbdel}[1][]{ \hbs{\delta}^{\ffsize{#1}}\vphantom{l} }%
\newcommand{\hbdeldn}[1][]{ {\hbs{\delta}_\ffsize{\!#1}\vphantom{l}} } 
\newcommand{\tbdel}[1][]{ \tbs{\delta}^{\ffsize{#1}}\vphantom{l} }%
\newcommand{\tbdeldn}[1][]{ {\tbs{\delta}_\ffsize{\!#1}\vphantom{l}} }
\newcommand{\barbdel}[1][]{ \barbs{\delta}^{\ffsize{#1}}\vphantom{l} }%
\newcommand{\barbdeldn}[1][]{ {\barbs{\delta}_\ffsize{\!#1}\vphantom{l}} }

\newcommand{\bD}[1][]{{ \fnsz{\sfb{D}}_{\ffsize{\!#1}}} }
\newcommand{\bDel}[1][]{{ \fnsz{\bs{\Delta}}^{\mkern-1mu\ffsize{#1}}}}%
\newcommand{\bDeldn}[1][]{{ \fnsz{\bs{\Delta}}_{\ffsize{#1}}}}%

\newcommand{\ff}[2]{\sfb{#1}_{\ffsize{\!#2}}} 
\newcommand{\ffup}[2]{\sfb{#1}^{\ffsize{#2}}} 
\newcommand{\coff}[2]{\bs{#1}^{\ffsize{#2}}}
\newcommand{\coffdn}[2]{\bs{#1}_{\ffsize{\!#2}}}

\newcommand{\bi}[1][]{{\sfb{i}_{\ffsize{\!#1}}}}
\newcommand{\bio}[1][]{{\bs{\iota}^{\ffsize{#1}}} \vphantom{\iota}}
\newcommand{\hbi}[1][]{{\hsfb{\i}_{\ffsize{\!#1}}}} 
\newcommand{\hbio}[1][]{{\hbs{\iota}^{\ffsize{#1}}} \vphantom{\iota}}

\newcommand{\be}[1][]{{\sfb{e}_{\ffsize{\!#1}}}}
\newcommand{\beup}[1][]{{\sfb{e}^{\ffsize{#1}}}} 
\newcommand{\tbe}[1][]{ {\tsfb{e}_{\ffsize{\!#1}}} }
\newcommand{\hbe}[1][]{ {\hsfb{e}_{\ffsize{\!#1}}} }
\newcommand{\bep}[1][\!]{{\bs{\epsilon}^{\ffsize{#1}}}}
\newcommand{\bepdn}[1][]{{\bs{\epsilon}_{\ffsize{\!#1}}}}
\newcommand{\tbep}[1][\!]{\tbs{\epsilon}^{\ffsize{#1}} \vphantom{\epsilon}}
\newcommand{\hbep}[1][\!]{{\hbs{\epsilon}^{\ffsize{#1}}} \vphantom{\epsilon}}

\newcommand{\bt}[1][]{ \sfb{t}_{\ffsize{\!#1}} \vphantom{t}}
\newcommand{\tbt}[1][]{\tsfb{t}_{\ffsize{\!#1}} \vphantom{t}}
\newcommand{\hbt}[1][]{\hsfb{t}_{\ffsize{\!#1}} \vphantom{t}}
\newcommand{\btup}[1][\,]{ \sfb{t}^{\ffsize{#1}\!} \vphantom{t}}
\newcommand{\tbtup}[1][\,]{\tsfb{t}^{\ffsize{#1}\!} \vphantom{t}}
\newcommand{\hbtup}[1][\,]{\hsfb{t}^{\ffsize{#1}\!} \vphantom{t}}

\newcommand{\btau}[1][\!]{{\bm{\tau}^{\ffsize{#1}}}}
\newcommand{\tbtau}[1][\!]{{\tbm{\tau}^{\ffsize{#1}}}}
\newcommand{\hbtau}[1][\!]{{\hbm{\tau}^{\ffsize{#1}}}}
\newcommand{\btaudn}[1][]{{\bm{\tau}_{\ffsize{\!#1}}}}
\newcommand{\tbtaudn}[1][]{{\tbm{\tau}_{\ffsize{\!#1}}}}
\newcommand{\hbtaudn}[1][]{{\hbm{\tau}_{\ffsize{\!#1}}}}


\newcommand{\Zvar}{{ \text{\ooalign{\hidewidth\raisebox{0.25ex}{--}\hidewidth\cr$Z$\cr}}\vphantom{Z}} }

\newcommand{\zvar}{%
  \text{\ooalign{\hidewidth -\kern-.3em-\hidewidth\cr$z$\cr}} \vphantom{z}%
}

\renewcommand{\Zbar}{{\textnormal{\ooalign{\hidewidth\raisebox{0.3ex}{\footnotesize{--}}\hidewidth\cr$Z$\cr}} \vphantom{Z}} }

\newcommand{\zbar}{ \mkern0.5mu
  {\text{\ooalign{\hidewidth\raisebox{0.1ex}{\scriptsize{--}}\hidewidth\cr$\mkern-0.5muz$\cr}} \vphantom{z}}%
} 

\newcommand{\Jbar}{ \mkern2mu {\textnormal{\ooalign{\hidewidth\raisebox{0.4ex}{\scriptsize{--}}\hidewidth\cr$\mkern-2mu J$\cr}} \vphantom{J}} }

\newcommand{\Ibar}{{\textnormal{\ooalign{\hidewidth\raisebox{0.4ex}{\scriptsize{--}}\hidewidth\cr$I$\cr}} \vphantom{I}} }

\newcommand{\Sbar}{{\mkern0mu \textnormal{\ooalign{\hidewidth\raisebox{0.4ex}{\scriptsize{\textbf{---}}}\hidewidth\cr$\mkern-0mu S$\cr}} \vphantom{S}} }

\newcommand{\sbar}{ \mkern-0.5mu
  {\text{\ooalign{\hidewidth\raisebox{0.1ex}{\scriptsize{--}}\hidewidth\cr$\mkern0.5mu s$\cr}} \vphantom{s}}%
}

\newcommand{\Cbar}{{\mkern-7mu \textnormal{\ooalign{\hidewidth\raisebox{0.35ex}{\footnotesize{--}}\hidewidth\cr$\mkern7muC$\cr}} \vphantom{C}} }

\newcommand{\cbar}{%
  {\mkern-3mu\text{\ooalign{\hidewidth\raisebox{0.1ex}{\scriptsize{--}}\hidewidth\cr$\mkern3muc$\cr}} \vphantom{c}}%
}

\newcommand{\ldash}[0]{{\mkern-2mu \textit{\l}\mkern1mu}} 
\newcommand{\Ldash}[0]{{\textit{\L}}} 
\newcommand{\ldashb}[0]{{\itsb{\l}}}
\newcommand{\Ldashb}[0]{{\txsfb{\L}}} 
\newcommand{\Lbar}{ {\mkern-2mu \textnormal{\ooalign{\hidewidth\raisebox{0.35ex}{\scriptsize{--}}\hidewidth\cr$\mkern2mu L$\cr}} \vphantom{L}}  }

\newcommand{\lbar}{%
    {\textnormal{\ooalign{\hidewidth\raisebox{0.3ex}{-}\hidewidth\cr$l$\cr}} \vphantom{l}}%
} 

\newcommand{\blbar}{%
    {\textbf{\ooalign{\hidewidth\raisebox{0.4ex}{\scriptsize{--}}\hidewidth\cr$\bs{l}$\cr}} \vphantom{l}}%
} 

\newcommand{\kbar}{%
   {\mkern-1mu \text{\ooalign{\hidewidth\raisebox{0.8ex}{\scriptsize{--}}\hidewidth\cr$\mkern1mu k$\cr}} \vphantom{k} }%
}
\newcommand{\kbarup}{%
   {\mkern-4mu \textnormal{\ooalign{\hidewidth\raisebox{0.8ex}{\scriptsize{--}}\hidewidth\cr$\mkern4mu \mrm{k}$\cr}} \vphantom{k} }%
}
\newcommand{\txkbar}[0]{\kbarup}

\newcommand{\hvar}{%
   {\mkern-1mu \text{\ooalign{\hidewidth\raisebox{0.8ex}{\scriptsize{--}}\hidewidth\cr$\mkern1mu h$\cr}} \vphantom{h} }%
}
\newcommand{\hbarup}{%
   {\mkern-4mu \textnormal{\ooalign{\hidewidth\raisebox{0.8ex}{\scriptsize{--}}\hidewidth\cr$\mkern4mu \mrm{h}$\cr}} \vphantom{k} }%
}
\newcommand{\txhbar}[0]{\hbarup}

\newcommand{\bbar}{%
   {\mkern-1mu \text{\ooalign{\hidewidth\raisebox{0.8ex}{\scriptsize{--}}\hidewidth\cr$\mkern1mu b$\cr}} \vphantom{b} }%
}
\newcommand{\bbarup}{%
   {\mkern-4mu \textnormal{\ooalign{\hidewidth\raisebox{0.8ex}{\scriptsize{--}}\hidewidth\cr$\mkern4mu \mrm{b}$\cr}} \vphantom{b} }%
} 
\newcommand{\txbbar}[0]{\bbarup}

\newcommand{\dbar}{%
   {\mkern4.5mu \text{\ooalign{\hidewidth\raisebox{0.8ex}{\scriptsize{--}}\hidewidth\cr$\mkern-4.5mu d$\cr}} \vphantom{d} }%
}
\newcommand{\dbarup}{%
   {\mkern3mu \textnormal{\ooalign{\hidewidth\raisebox{0.8ex}{\scriptsize{--}}\hidewidth\cr$\mkern-3mu \mrm{d}$\cr}} \vphantom{d} }%
} 
\newcommand{\txdbar}[0]{\dbarup}

\newcommand{\Gambar}{%
   {\mkern-4mu \text{\ooalign{\hidewidth\raisebox{0.4ex}{\scriptsize{--}}\hidewidth\cr$\mkern4mu \Gamma$\cr}} \vphantom{h} }%
}


\renewcommand{\v}[0]{{\nu}} 
\newcommand{\txv}[0]{{\txi{v}}}
\newcommand{\altv}[0]{{\scriptstyle{\mathcal{V}}}}
\newcommand{\newv}[0]{{\scriptstyle{\mathscr{V}}}}

\newcommand{\one}[0]{{\ss{1}}}
\newcommand{\two}[0]{{\ss{2}}}
\newcommand{\three}[0]{{\ss{3}}}
\newcommand{\four}[0]{{\ss{4}}}
\newcommand{\six}[0]{{\ss{6}}}
\newcommand{\eight}[0]{{\ss{8}}}


\renewcommand{\a}[0]{{\alpha}}
\renewcommand{\b}[0]{{\ii{\beta}}}
\newcommand{\g}[0]{{\ss{\gamma}}}
\newcommand{\gam}[0]{{\gamma}}
\newcommand{\Gam}[0]{{\Gamma}}
\newcommand{\y}[0]{{\ss{\lambda}}}
\newcommand{\lam}[0]{{\lambda}}
\newcommand{\Lam}[0]{{\Lambda}}
\newcommand{\sig}[0]{{\sigma}}
\newcommand{\Sig}[0]{{\Sigma}}
\newcommand{\varsig}[0]{{\varsigma}}
\newcommand{\ep}[0]{{\epsilon}}
\newcommand*{\varep}[0]{{\varepsilon}}
\newcommand{\omg}[0]{{\omega}}
\newcommand{\Omg}[0]{{\Omega}}
\newcommand{\kap}[0]{{\kappa}}
\newcommand{\varkap}[0]{{\varkappa}}
\newcommand{\ro}[0]{{\varrho}}
\renewcommand{\th}[0]{{\theta}}
\newcommand{\Th}[0]{{\Theta}}
\newcommand{\varth}[0]{{\vartheta}}

\DeclareRobustCommand{\rchi}{{\mathpalette\irchi\relax}}
\newcommand{\irchi}[2]{\raisebox{\depth}{$#1\chi$}} 
\DeclareRobustCommand{\pup}{{\mathpalette\newp\relax}}
\newcommand{\newp}[2]{\raisebox{0.8\depth}{$#1p$}}
\DeclareRobustCommand{\jup}{{\mathpalette\newj\relax}}
\newcommand{\newj}[2]{\raisebox{0.8\depth}{$#1\jmath$}}

\newcommand{\upa}[0]{{\upalpha}}
\newcommand{\upb}[0]{{\upbeta}}
\newcommand{\upgam}[0]{{\upgamma}}
\newcommand{\upGam}[0]{{\upGamma}}
\newcommand{\Upgam}[0]{{\Upgamma}}
\newcommand{\uplam}[0]{{\uplambda}}
\newcommand{\upLam}[0]{{\upLambda}}
\newcommand{\Uplam}[0]{{\Uplambda}}
\newcommand{\updel}[0]{{\updelta}}
\newcommand{\upsig}[0]{{\upsigma}}
\newcommand{\upSig}[0]{{\upSigma}} 
\newcommand{\Upsig}[0]{{\Upsigma}}
\newcommand{\upomg}[0]{{\upomega}}
\newcommand{\upOmg}[0]{{\upOmega}}
\newcommand{\Upomg}[0]{{\Upomega}}
\newcommand{\upep}[0]{\upepsilon}
\newcommand{\upth}[0]{\uptheta}
\newcommand{\upvarth}[0]{\upvartheta}
\newcommand{\upTh}[0]{\upTheta}
\newcommand{\Upth}[0]{\Uptheta}
\newcommand{\upkap}[0]{\upkappa}
\newcommand{\upvark}[0]{\upvarkappa}
\newcommand{\upn}[0]{\upeta}
\newcommand{\upi}[0]{\upiota}

\newcommand{\aup}[0]{{\alphaup}}
\newcommand{\bup}[0]{{\betaup}}
\newcommand{\gamup}[0]{{\gammaup}}
\newcommand{\Gamup}[0]{{\Gammaup}}
\newcommand{\lamup}[0]{{\lambdaup}}
\newcommand{\Lamup}[0]{{\Lambdaup}}
\newcommand{\delup}[0]{{\deltaup}}
\newcommand{\sigup}[0]{{\sigmaup}}
\newcommand{\varsigup}[0]{{\varsigmaup}}
\newcommand{\Sigup}[0]{{\Sigmaup}}
\newcommand{\wup}[0]{{\omegaup}} 
\newcommand{\omgup}[0]{{\omegaup}}
\newcommand{\Omgup}[0]{{\Omegaup}}
\newcommand{\epup}[0]{\epsilonup}
\newcommand{\thup}[0]{\thetaup}
\newcommand{\varthup}[0]{\varthetaup}
\newcommand{\kup}[0]{\kappaup}
\newcommand{\kapup}[0]{\kappaup}
\newcommand{\varkapup}[0]{\varkappaup}
\newcommand{\varkup}[0]{\varkappaup}
\newcommand{\nup}[0]{\etaup}
\newcommand{\iup}[0]{\iotaup}

\newcommand{\txa}[0]{{\text{\textalpha}}}
\newcommand{\txb}[0]{{\text{\textbeta}}}
\newcommand{\txbet}[0]{{\text{\textbeta}}}
\newcommand{\txg}[0]{{\text{\textgamma}}}
\newcommand{\txgam}[0]{{\text{\textgamma}}}
\newcommand{\txGam}[0]{{\text{\textGamma}}}
\newcommand{\txy}[0]{{\text{\textlambda}}}
\newcommand{\txlam}[0]{{\text{\textlambda}}}
\newcommand{\txr}[0]{{\text{\textrho}}}
\newcommand{\txrho}[0]{{\text{\textrho}}}
\newcommand{\txsig}[0]{{\text{\textsigma}}}
\newcommand{\txSig}[0]{{\text{\textSigma}}}
\newcommand{\txw}[0]{{\text{\textomega}}}
\newcommand{\txomg}[0]{{\text{\textomega}}}
\newcommand{\txep}[0]{{\text{\textepsilon}}}
\newcommand{\txth}[0]{{\text{\texttheta}}}
\newcommand{\txvarth}[0]{{\text{\textvartheta}}}
\newcommand{\txphi}[0]{{\text{\textphi}}}
\newcommand{\txvarphi}[0]{{\text{\textvarphi}}}
\newcommand{\txpsi}[0]{{\text{\textpsi}}}
\newcommand{\txd}[0]{{\text{\textdelta}}}
\newcommand{\txdel}[0]{{\text{\textdelta}}}
\newcommand{\txiota}[0]{{\text{\textiota}}}
\newcommand{\txn}[0]{{\text{\texteta}}}
\newcommand{\txeta}[0]{{\text{\texteta}}}
\newcommand{\txxi}[0]{{\text{\textxi}}}
\newcommand{\txchi}[0]{{\text{\textchi}}}
\newcommand{\txpi}[0]{{\text{\textpi}}}
\newcommand{\txtau}[0]{{\text{\texttau}}}
\newcommand{\txz}[0]{{\text{\textzeta}}}
\newcommand{\txzeta}[0]{{\text{\textzeta}}}
\newcommand{\txnu}[0]{{\text{\textnu}}}
\newcommand{\txmu}[0]{{\text{\textmugreek}}} 
\newcommand{\txups}[0]{{\text{\textupsilon}}} 
\newcommand{\txk}[0]{{\text{\textkappa}}}
\newcommand{\txkap}[0]{{\text{\textkappa}}}

\newcommand{\tximu}[0]{ {\textit{\textmugreek}} } %
\newcommand{\txiu}[0]{ {\textit{\textupsilon}} } 
\newcommand{\txitau}[0]{ {\textit{\texttau}} } %
\newcommand{\txipi}[0]{ {\textit{\textpi}} }%
\newcommand{\txirho}[0]{ {\textit{\textrho}} } 
\newcommand{\txilam}[0]{ {\textit{\textlambda}} }

\newcommand{\wtxt}[0]{{\textomega}}
\newcommand{\ytxt}[0]{{\textlambda}}
\newcommand{\thtxt}[0]{{\texttheta}}
\newcommand{\btxt}[0]{{\textbeta}}
\newcommand{\atxt}[0]{{\textalpha}}
\newcommand{\gtxt}[0]{{\textgamma}}
\newcommand{\tautxt}[0]{{\texttau}}

\newcommand{\gvgr}[1]{\vec{\tbf{#1}}}
\newcommand{\uvgr}[1]{\hat{\tbf{#1}}}
\newcommand{\tvgr}[1]{\widetilde{\tbf{#1}}}
\newcommand{\bvgr}[1]{\bar{\tbf{#1}}}
\newcommand{\gvbb}[1]{\vec{\pmb{#1}}}
\newcommand{\uvbb}[1]{\hat{\pmb{#1}}}
\newcommand{\tvbb}[1]{\widetilde{\pmb{#1}}}
\newcommand{\bvbb}[1]{\bar{\pmb{#1}}}

\newcommand{\sigvec}[0]{\vec{\pmb{\upsigma}} }
 \newcommand{\sigtvec}[0]{\widetilde{\pmb{\upsigma}} }

\newcommand{\gvec}[1]{\vec{\bs{\mathrm{#1}}}}
\newcommand{\gv}[1]{\vec{\bm{#1}}}
\newcommand{\gvb}[1]{\vec{\bs{#1}}}
\newcommand*{\bsvec}[1]{\vec{\bs{#1}}}
\newcommand*{\sfvec}[1]{\,\vec{\!\sfb{#1}}}

\newcommand{\uvec}[1]{\hat{\bs{#1}}}
\newcommand{\uv}[1]{\hat{\bs{#1}}}
\newcommand{\ihat}[0]{\hat{\bs{\iota}}}
\newcommand{\ehat}[0]{\hat{\bs{e}}}
\newcommand{\bhat}[0]{\hat{\bs{b}}}
\newcommand{\shat}[0]{\hat{\bs{s}}}
\newcommand{\nhat}[0]{\uvec{n}}
\newcommand{\evec}[0]{\vec{\bs{e}}}
\newcommand{\epvec}[0]{\vec{\bs{\epsilon}}}

\newcommand{\uvecb}[1]{ {\hbs{#1}} }
\newcommand{\uvb}[1]{ {\hbs{#1}} }
\newcommand{\ihatb}[0]{\hat{\pmb{\upiota}}} 
\newcommand{\ehatb}[0]{\uvecb{e}} 
\newcommand{\bhatb}[0]{\!\!\uvecb{\;b}}
\newcommand{\hhatb}[0]{\!\!\uvecb{\;h}}
\newcommand{\shatb}[0]{\uvecb{s}}
\newcommand{\nhatb}[0]{\uvecb{n}}
\newcommand{\ohatb}[0]{\uvecb{o}}
\newcommand{\ahatb}[0]{\uvecb{a}}
\newcommand{\epvecb}[0]{\vec{\bs{\epsilon}}}
\newcommand{\evecb}[0]{ \hspace{0.3mm}\vec{\mbf{e}}\hspace{0.1mm} }

\newcommand{\tvb}[1]{\widetilde{\mbf{#1}}}
\newcommand{\tvecb}[1]{\bar{\bs{#1}}}
\newcommand{\bv}[1]{\bar{\bs{#1}}}
\newcommand{\bvb}[1]{\bar{\bs{\mathrm{#1}}}}


\newcommand{\rnote}[1]{\noindent{\footnotesize\textit{{\color{red}${\circ}$}  #1}}}
\newcommand{\red}[1]{{\color{red}#1}}

\newcommand{\bnote}[1]{\noindent{\footnotesize\textit{{\color{blue}${\circ}$} #1}}}

\newcommand{\gnote}[1]{\noindent{\footnotesize\textit{$\circ$  #1}}}

 \newcommand{\note}[1]{\noindent{\footnotesize\textit{{\color{darkgray} #1}}}}

\newcommand{\nsp}[0]{\!\!\!\!}    
\newcommand{\nquad}[0]{\hspace{-1em}} 
\newcommand{\nqquad}[0]{\hspace{-2em}} 


\theoremstyle{plain} 
\newtheorem{thrm}{Theorem}[section]
\newtheorem{defn}{Def.}[section]

\theoremstyle{plain} 
\newtheorem{remit}{Remark}[section]

\theoremstyle{definition} 
\newtheorem{remark}{Remark}[section]


\newtheoremstyle{remsans}
{8pt} 
{12pt} 
{\sffamily\slshape\small}
{}
{\rmfamily\small}
{.}
{.3em}
{}

\theoremstyle{remsans}
\newtheorem{remsf}{$\iisqr\,$\rmsb{remark}}[section]

\newtheoremstyle{remrmsmall}
{8pt} 
{12pt} 
{\rmfamily\small}
{}
{\rmfamily\small\slshape}
{.}
{0.3em}
{}

\theoremstyle{remrmsmall}
\newtheorem{remrm}{$\iisqr\,$Remark}[section]
\newtheorem{noether}[remark]{$\iisqr\,$\rmsb{Noether}}

\newtheoremstyle{remslant}
{8pt} 
{12pt} 
{\rmfamily\slshape\small}
{}
{\rmfamily\slshape\bfseries\small}
{.}
{.3em}
{}

\theoremstyle{remslant} 
\newtheorem{remsl}{$\iisqr\,$\rmsb{Remark}}[section]

\newtheoremstyle{remnopunc}
{8pt} 
{8pt} 
{\rmfamily\small}
{}
{\bfseries\small}
{}
{.2em}
{}
\theoremstyle{remnopunc} 
\newtheorem*{noteblt}{$\nmblt$}
\newtheorem*{notestr}{\raisebox{0.1ex}{$\star$}}
\newtheorem*{noteast}{$\bs{*}$}

\newenvironment{notesq}
    {\begin{small}
     \begin{itemize}[left=0pt,labelsep=0.3ex,topsep=10pt]
     \item[$\iisqr$] }
    {\end{itemize}
     \end{small} }

\newtheoremstyle{remslantnopunc}
{8pt} 
{8pt} 
{\slshape\rmfamily\small}
{}
{\slshape\bfseries\small}
{}
{.2em}
{}
\theoremstyle{remslantnopunc} 
\newtheorem*{notesl}{$\iisqr$}

\newtheoremstyle{remrmnonbold}
{} 
{} 
{\rmfamily\footnotesize}
{}
{\rmfamily\itshape\footnotesize}
{.}
{.2em}
{}
\theoremstyle{remrmnonbold}
\newtheorem*{notation}{Notation}



\newcommand{\eq}[1]{\text{$#1$}}

\newenvironment{eqn}
    {\begin{align} }
    {\end{align}}

\newenvironment{fleqn}
    {\begin{flalign} }
    {\end{flalign}}

\newenvironment{smeqn}
    {\begin{small}\begin{align} }
    {\end{align}\end{small} }

\newenvironment{flsmeqn}
    {\begin{small}\begin{flalign} }
    {\end{flalign}\end{small} }


%% file: Mysecs_prj/00_abstract.tex

 \begin{abstract}
    We introduce a method of regularizing and linearizing central-force Hamiltonian dynamics using projective point transformations extended to the momentum level. By considering a generalization of the usual projective decomposition, we obtain a family of canonical transformations which differ in their momentum coordinates. 
    From this family, a preferred canonical coordinate set is chosen that possesses a simple, intuitive connection to orbital reference frames and attitude dynamics. 
    Using this transformation, 
    closed-form phase space solutions and state transition matrices are readily obtained for inverse square and inverse cubic radial forces, or any superposition thereof (i.e., Kepler or Manev dynamics). 
    The classic $J_2$-perturbed two-body problem is formulated in the new coordinates and used for numerical verification. 
\end{abstract}


%% file: Mysecs_prj/00_notation.tex

\subsubsection*{Nomenclature}

\vspace{-1ex}
 \begin{small}
\begin{longtable}[htbp]{@{}p{0.13\textwidth} p{0.80\textwidth}@{}}
  $\imat_\en, \, \kd_{ij}, \, \lc_{i_1 \dots i_\en}$
  & $\en$-dim identity matrix, Kronecker delta \textit{symbol}, and Levi-Civita permutation \textit{symbol}, respectively. 
\\[1ex]
    $\ibase_i \in\mbb{R}^\en$
    & standard basis for $\mbb{R}^\en$ (columns of \eq{\imat_\en}, usually \eq{\en=3}). 
\\[1ex]
       $ \Somat{\en}, \, \somat{\en} $ &  Lie group of real special orthogonal matrices and its Lie algebra of antisymmetric matrices.
\\[1ex]
       $ \Spmat{2\en}, \, \spmat{2\en} $ &  Lie group of real symplectic matrices and its Lie algebra of Hamiltonian matrices.
\\[1ex]
       $ \jmat_{2\en} = \fnpmat{0 & \imat_\en \\ -\imat_\en & 0 } $
       &  Standard symplectic matrix on $\mbb{R}^{2\en}$. $\inv{J}=\trn{J}=-J\in\Spmat{2\en}$. 
\\[2ex]
    $(\tup{r},\tup{v})$ 
   &inertial cartesian position and velocity/momentum 
    coordinates, $(\tup{r},\tup{v})\in\mbb{R}^{6}$.  
\\[1ex]
      $(\bartup{q},\bartup{p})$
   & projective coordinates, $(\bartup{q},\bartup{p})=(\tup{q},u,\tup{p},p_\ss{u})\in\mbb{R}^{6+2}$.
   (some exceptions\footnote{This applies everywhere in this work other than section \ref{sec:Pxform} where \eq{(\bartup{q},\bartup{p})} denote arbitrary redundant phase space coordinates. Also in Appx.~\ref{app:Ham_mech_crd} where \eq{(\tup{q},\tup{p})} denote arbitrary phase space coordinates.}).
\\[1ex]
    $n,\,m$
    & ``knobs'' we may turn in the point transformation $\tup{r}=u^n q^m \tup{q}$. (we choose \eq{n=m=-1}). 
\\[1ex]
   $\tup{\slangup},\; \hdge{\tup{\slangup}}$ 
   & (specific) angular momentum coordinate vector $\tup{\slangup}\in\mbb{R}^3 $, and matrix $\hdge{\tup{\slangup}}\in \somat{3}$. 
\\[1ex]
    $w:= u^2 p_\ss{u}$ 
    &quasi-momenta coordinate in place of $p_\ss{u}$.
\\[2ex]
   $\mscr{K},\, \mscr{H}$ 
   &``cartesian coordinate Hamiltonian'' $\mscr{K}$, and ``projective coordinate Hamiltonian''  $\mscr{H}$.
\\[1ex]
    $V^0,\; V^1$
    & central-force potential $V^0(r)$, and arbitrary perturbing potential $V^1(\tup{r},t)$. 
\\[1ex]
    $\kconst_1,\,\kconst_2$
    & scalar constants for Manev-type potential $V^0=-\kconst_1/r - \tfrac{1}{2}\kconst_2/r^2$ (with $\kconst_2=0$ for Kepler). 
\\[1ex]
    $\tup{a}^\nc \in\mbb{R}^3$
    & inertial cartesian components of any/all nonconservative perturbing forces.
\\[1ex]
    $\tup{F}:=-\pderiv{V^1}{\tup{r}}+\tup{a}^\nc$
    & inertial cartesian components of total perturbing forces.
\\[1ex]
    $\bartup{\alpha}\in\mbb{R}^4$
    & generalized nonconservative perturbing forces for projective coordinates. Split as $\bartup{\alpha}=(\tup{\alphaup},\alpha_\ss{u})\in\mbb{R}^4$.
\\[1ex]
    $\bartup{f}:=-\pderiv{V^1}{\bartup{q}}+\bartup{\alphaup}$
    & generalized total perturbing forces for projective coordinates. Split as $\bartup{f}=(\tup{f},f_\ss{u})\in\mbb{R}^4$.
\\[2ex]
    $\dot{\square}:=\diff{\square}{t}$ 
    &derivative ``with respect to'' \eq{t} (time).
\\[1ex]
    $\pdt{\square}:=\diff{\square}{s}$ 
    &derivative ``with respect to'' \eq{s}, where $\mrm{d} t = r^2 \mrm{d}s = u^{2n} \mrm{d}s$.
\\[1ex]
    $\rng{\square}:=\diff{\square}{\tau}$
    &derivative ``with respect to'' \eq{\tau}, where $\mrm{d} t = (r^2/\slang)\mrm{d} \tau = (u^{2n}/\slang) \mrm{d} \tau$.
\\[2ex]
    $a:=b,\, b=:a$ 
    & $a$ is defined as $b$.
\\[1ex]
    $a\simeq b$
    & $a=b$ under certain conditions/simplifications (usually, $\mag{\tup{q}}=1$ and \eq{\htup{q}\cdot\tup{p}=0}). 
\\[2ex]
       ODE     & ordinary differential equation.
     \\[1ex]
      KS    &  Kustaanheimo-Stiefel (two people).
   \\[1ex]
      BF    &  Burdet-Ferrándiz (two people).
  \\[1ex]
     DEF   & Deprit-Elipe-Ferrer (three people).
    \\[1ex]
        LVLH     & local vertical, local horizontal (basis).
\end{longtable}
\end{small}

%% file: Mysecs_prj/00_INTRO2.tex
\section{Introduction}




In the Newtonian gravity model, the non-linear equation of motion for two bodies of mass \eq{m_a} and \eq{m_b} in Euclidean 3-space is given, after the usual reduction from 6-\sc{dof} to 3-\sc{dof}, as follows:
\begin{small}
\begin{flalign} \label{ddr_2bp}
     && \ddiff{}{t}\vecbs{r} \,=\, -\tfrac{\kconst_1}{r^2}\hbs{r} \,+\, \vecbs{F}
      \;\;,&& \kconst_1=G(m_a+m_b)
      \qquad
\end{flalign}
\end{small}
 where \eq{\vecbs{r}} is the displacement vector from \eq{m_a} to \eq{m_b} with norm \eq{r:=\mag{\vecbs{r}}} and unit vector \eq{\hbs{r}:=\vecbs{r}/r}, 
 where \eq{\ddiff{}{t}\vecbs{r}} is the acceleration, 
 \eq{G} is the gravitational constant, 
 and \eq{\vecbs{F}} is the perturbing force vector per unit reduced mass, \eq{m:=\tfrac{m_a m_b}{m_a+m_b}}, arising from any non-Keplerian forces (e.g., gravitational perturbations, thrust, drag, etc.).\footnote{Note in Eq.\eqref{ddr_2bp} that we have already reduced the initial \eq{6}-\sc{dof} problem of two particles to an equivalent 3-\sc{dof} problem, as usual. It is often further assumed that \eq{m_b<<m_a}  such that \eq{m_b} orbits around \eq{m_a}, which is then the approximately stationary ``central body''. In that case, \eq{\kconst_1\approx G m_a} is approximately the gravitational parameter of the central body (often denoted \eq{\mu}) and the reduced mass \eq{m:=\tfrac{m_a m_b}{m_a+m_b}\approx m_b} is approximately the mass of the orbiting body.} 
The above vector equation may be expressed using a variety of coordinates and considerable work in celestial mechanics has focused on developing coordinate representations that remove the singularity at \eq{r=0} or lead to other desirable features in the equations of motion (e.g., linearizing coordinate representations or slowly-varying non-singular orbit element sets) \cite{schumacher1987results,stiefel1975linear,deprit1994linearization,roa2017regularization,ferrer1983regularization,bau2013new,bau2015non,majji2020regular}. This process is broadly referred to as \textit{regularization} and often incorporates a coordinate transformation 
and/or a transformation of the evolution parameter to something other than time \cite{schumacher1987results,bond1981canonical}. See Eq.\eqref{tparams} for examples of such parameters. 


This work is not concerned with ``mere regularization'' of Kepler-type dynamics but, more specifically, \textit{linearization}, which can be seen as an ideal case of regularization.
Not only that, we are concerned specifically with linearization in a \textit{Hamiltonian} framework.
That is, \textit{we are interested in obtaining canonical transformations such that Eq.\eqref{ddr_2bp} can be represented as a Hamiltonian system of ODEs  which are fully linear in the unperturbed case \eq{\vecbs{F}=0}}. 

To the authors' knowledge, only two categories of such transformations exist:~\textit{quaternion}-based transformations and
\textit{projective} transformations.
The development of these transformations as canonical transformations in the Hamiltonian framework is a non-trivial problem due to the fact that both types utilize redundant coordinate representations as well as transformations of the evolution parameter. 
Yet, this has been accomplished and is now well-established for quaternionic transformations, with
the most successful example being the  \bemph{Kustaanheimo-Stiefel (KS) transformation} \cite{stiefel1975linear,stiefel1973linear,kurcheeva1977kustaanheimo,bond1974uniformKS}. Although originally formulated in terms of spinors \cite{KS1kustaanheimo1965perturbation}, several sources have since clarified the KS transformation's close connection to quaternions and orbital reference frames \cite{deprit1994linearization,saha2009interpreting,waldvogel2008quaternions}.
Chelnokov is particularly prolific on the subject (e.g., two separate series of papers with initial entries \cite{chelnokov2013quatTraj1} and \cite{chelnokov2017quat3BP1}).
It has seen use ranging from celestial mechanics to quantum mechanics, and has been well-studied from a mathematical perspective \cite{cordani2003kepler,ferrer2018altKSgeo,meer2021reduction}. 
In comparison, the use of \textit{projective} transformations for linear regularization of Kepler dynamics in a Hamiltonian framework is considerably less established; that will be our focus here.  For context, further details and references are given below.



\paragraph*{Projective Point Transformations for Linearized \textit{Newtonian} Kepler Dynamics.}
A seminal example of a linear and regular representation of Eq.\eqref{ddr_2bp}  \textit{not} employing Hamiltonian or symplectic methods is the so-called 
\textit{projective decomposition} notably used by Burdet and Vitins \cite{burdet1969mouvement,Burdet+1969+71+84,vitins1978keplerian}. 
This method uses a redundant set of four configuration coordinates, along with a transformation of the evolution parameter, to transform the \eq{3}-dim Kepler problem into a \eq{4}-dim linear oscillator. Though different formulations exist — and Schumacher has detailed several of them \cite{schumacher1987results} —  
the general process is summarized as follows:
\begin{small}
\begin{enumerate}
\item Let \eq{\tup{r}\in\mbb{R}^3}  denote inertial cartesian coordinates in some fixed orthonormal basis
such that the two-body dynamics in Eq.\eqref{ddr_2bp} are written simply as follows (where \eq{\dot{\square} :=\diff{\square}{t}}):
    \begin{align}  \label{ddr_2bp2}
          \ddot{\tup{r}} \,=\, -\tfrac{\kconst_1}{r^2}\htup{r} \,+\,  \fnsize{(perturbations)}
    \end{align}
\item Specify a point transformation \eq{\mbb{R}^4\ni(\tup{y},z)\mapsto \tup{r}\in\mbb{R}^3} with \eq{\tup{y}=\htup{r}} the three inertial cartesian components of the radial unit vector and \eq{z} is a single coordinate related to the radial distance by \eq{z = r^\ss{1/n}} for some real number \eq{n\neq0} 
(often, \eq{n=\pm 1}). 
This is commonly called the \textit{projective decomposition/transformation}:
\begin{align} \label{proj_pt}
    \tup{r} \,=\, z^n \tup{y}
    \qquad \leftrightarrow \qquad 
    \tup{y} = \htup{r} \;\;,\;\; z = r^\ss{1/n}
\end{align}
\item Transform the evolution parameter from time, \eq{t}, to a parameter \eq{s} related by some specified differential relation 
(examples\footnote{Some common evolution parameters used for regularization of the Kepler problem are given below, where \eq{\slang} is angular momentum magnitude: 
\begin{align}\label{tparams}
\begin{array}{lllll}
      \d t = r \d s  & \fnsize{(unnamed)}
\\[3pt]
    \d t = r^2 \d s &  \fnsize{(unnamed)}
\end{array}
&&
  \d t = \tfrac{r^2}{\slang} \d \tau \;\; \fnsize{(true anomaly)}
  &&
\begin{array}{lllll}
      \d t = r\sqrt{a/\kconst_1} \, \d \itsc{e} = \tfrac{r}{\sqrt{-2 \mscr{E}^\zr}} \d \itsc{e}  
      & \fnsize{(eccentric anomaly)}
\\[3pt]
    \d t = \sqrt{a^3/\kconst_1} \, \d \itsc{m} 
    & \fnsize{(mean anomaly)}
\end{array}
\end{align}
which defines \eq{\tau},  \eq{\itsc{e}}, and \eq{\itsc{m}} as the true, eccentric, and mean anomaly (up to an additive constant), and \eq{\mscr{E}^\zr=\v^2/2-\kconst_1/r} is the Keplerian energy. Note that, when combined with projective transformations of the form in Eq.\eqref{ddr_2bp2}, not all of the above lead to \textit{linear} regularization of the Kepler problem.}):
\begin{align} \label{t'_proj}
    \pdt{t} \,=\,  \diff{t}{s} \,=\, f(\tup{r},\dot{\tup{r}},t)
\end{align}
\item Substitute Eq.\eqref{proj_pt} into Eq.\eqref{ddr_2bp2} and convert all \eq{\diff{}{t}} terms to \eq{\diff{}{s}} terms using \eq{\diff{}{s}=\pdt{t}\diff{}{t} = f\diff{}{t}} and \eq{\ddiff{}{s} = f^2 \ddiff{}{t} + f \diff{f}{t}\diff{}{t} }. The result — for appropriate choices of \eq{n} and \eq{f} — is four equations for \eq{(\tup{y},z)} of the following general form (where \eq{\pdt{\square}:=\diff{\square}{s}}):
\begin{align} \label{ddru_proj}
     \begin{array}{ll}
         &  \pddt{\tup{y}} + \omega^2 \tup{y} \,=\,  \tup{c} \,+\, \fnsize{(perturbations)}
     \\[4pt]
         &  \pddt{z} + \omega^2 z \,=\, b \,+\, \fnsize{(perturbations)}
    \end{array}
\end{align}
where \eq{\omega}, \eq{b}, and \eq{\tup{c}} are all either numeric constants or integrals of motion for unperturbed Kepler dynamics. In that  case, the above describes a 4-dim linear harmonic oscillator with constant frequency \eq{\omega} and constant driving ``forces'' \eq{\tup{c}} and \eq{b}.
The particular meaning of  \eq{\omega}, \eq{b}, and \eq{\tup{c}} depend on the choice of \eq{n} in Eq.\eqref{proj_pt} and the choice of \eq{\diff{t}{s}=f(\tup{r},\dot{\tup{r}},t)} in Eq.\eqref{t'_proj}.
Some notable examples are given in the footnote\footnote{Burdet used \eq{n=-1} (that is,  \eq{\tup{y}=\htup{r}} and \eq{z=1/r}) and transformed the evolution parameter using  \eq{\diff{t}{s}=r^2}, leading to \eq{b=\kconst_1},  \eq{\tup{c}=0}, and \eq{\omega^2=\slang^2} — where \eq{\slang} is the (specific) angular momentum magnitude \cite{burdet1969mouvement,Burdet+1969+71+84}. Vitins used the same coordinates but with \eq{\diff{t}{s}=r^2/\slang} — making  \eq{s} the true anomaly up to some additive constant — which leads to  \eq{\omega^2=1}, \eq{b=\kconst_1/\slang^2}, and \eq{\tup{c}=0} \cite{vitins1978keplerian}.  Vitins also made a subsequent point transformation from projective coordinates to quaternionic attitude coordinates (Euler parameters) for the LVLH basis and used their closed-form Kepler solutions to develop a corresponding set of eight orbit elements.
Vitins's quaternionic configuration coordinates differ from the KS coordinates. See the authors' work \cite{peterson2023regularized} for a Hamiltonian formulation of Vitins's quaterion-based regularization. We lastly note that, like Burdet, Bond also used \eq{\diff{t}{s}=r^2} to achieve linearization, but via a trigonometric transformation of spherical coordinates \cite{bond1985transformation}. 
}.
\end{enumerate}
\end{small}
Broadly speaking, the goal of this work is to extend the above developments to a canonical transformation in the Hamiltonian framework. Previous work has investigated this problem, as described below.

\paragraph*{Projective Canonical Transformations for Linearized \textit{Hamiltonian} Kepler Dynamics.} 
The projective point transformations described above are not formulated in the Hamiltonian framework. A large part of this is due to the fact that they employ redundant coordinates, which are not easily compatible with classic Hamiltonian formalism. However, prior work has explored methods of overcoming this obstacle.
In particular, Ferrándiz et al.~considered the general problem of dimension-raising canonical transformations and used this to extend Burdet-type projective point transformations of the form \eq{\tup{r}=\tfrac{1}{z}\tup{y}} or \eq{\tup{r}=z\tup{y}} to a canonical transformation in the Hamiltonian framework \cite{ferrandiz1987general,ferrandiz1988extended,ferrandiz1994extended,ferrandiz1992increased}. 
We will refer to this as the \bemph{Burdet-Ferrándiz (BF) transformation}, which follows the naming convention of Deprit et al.~\cite{deprit1994linearization}, who re-visited Ferrándiz's BF transformation, offering comparisons to the KS transformation as well as their own formulation, 
which they dub the \bemph{Deprit-Elipe-Ferrer (DEF) transformation}, along with notable efforts to clarify the structure and derivations of the prior BF transformation.

We would be remiss not to also mention Moser's 1970 regularization of the Kepler problem which is based on a \textit{stereographic} projective transformation \cite{moser1970regularization}. This construction differs substantially from the projective transformations discussed above and is \textit{not} the focus of the present work.
Nevertheless, it should be noted that Moser's transformation was extended to the Hamiltonian framework and that, when combined with a transformation of the evolution parameter, likewise yields a linearization of Kepler-type dynamics. While Moser's regularization is mathematically transparent, it is inconvenienced by an explicit dependence on energy levels. We refer the reader to \cite{moser1970regularization,kummer1982regularization,heckman2012regularization} for details.




\paragraph*{The Present Work.}\hspace{-1ex}\footnote{The authors have previously documented parts of this work in earlier conference papers \cite{peterson2022nonminimal,peterson2023regularized}, the first author's dissertation \cite{peterson2025phdThesis}, and related preprints \cite{peterson2025prjGeomech,peterson2025prjElements}.}
An outline and detailed summary of this work is given in section \ref{sec:prj_sum}.
Like the BF transformation, 
we develop a canonical extension of the projective point transformation and use it to linearize Kepler (and Manev) dynamics in a Hamiltonian framework. Actually, we develop a family of such transformations (Appendix \ref{sec:2bp_gen}), with the BF transformation included as one case. From this family, we ultimately choose a transformation that differs subtly from the BF transformation at the configuration level, and significantly at the momentum level.
In the authors' opinion, many properties of the transformed system, as well as its derivation, are notably simplified in comparison to the previous works mentioned above.
Our coordinates are intimately and intuitively related to the orbiting particle's local vertical local horizontal (LVLH) basis — attitude dynamics and the angular momentum matrix 
play a central role in the transformed Hamiltonian system. 
Furthermore, we first consider arbitrary central-forces along with arbitrary perturbations and, notably, we find that our projective transformation fully linearizes not only Kepler-type dynamics, but also \textit{Manev}-type dynamics. The latter warrants clarification:

    The Manev potential is an augmented  Kepler potential that includes an additional inverse square term:
    \begin{small}
    \begin{align} \label{Manev_def_intro}    
    \begin{array}{cc}
             \fnsize{Kepler}  \\[-1pt]
              \fnsize{potential}
        \end{array}
        \!\!:\quad 
        V^0 = -\tfrac{\kconst_1}{r}
        \qquad\quad,\qquad\quad
        \begin{array}{cc}
             \fnsize{Manev}  \\[-1pt]
              \fnsize{potential}
        \end{array}
        \!\!:\quad 
         V^0 = -\tfrac{\kconst_1}{r} - \tfrac{1}{2} \tfrac{\kconst_2}{r^2}
    \end{align}
    \end{small}
    for scalars \eq{\kconst_1,\kconst_2\in\mbb{R}} (the negative signs and factor of \eq{\sfrac{1}{2}} are included for later convenience). 
     Originally introduced as a classical approximation to certain relativistic corrections, the Manev potential  maintains key features of the Kepler potential such as conic orbits and integrability, while simultaneously introducing key relativistic features like perihelion precession 
    (a hallmark prediction of general relativity).\footnote{The Manev potential is not a true physics-based model in the sense that it is not derived from general relativity. Rather, the additional term \textit{qualitatively} reproduces the leading-order relativistic correction to the Kepler problem predicted by the relativistic Schwarzschild solution, while still maintaining the integrability of the classic Kepler problem.} 
    However, other than considering the above mathematical form of the Manev potential, we provide little discussion of practical applications in regards to relativity-motivated corrections in orbital mechanics. We refer the reader to Marmo et al.~\cite{marmo2006manev1,marmo2007manev2symmetries,marmo2008manev3kepler}, Diacu et al.~\cite{diacu1995manev1,diacu1996manev2,diacu2000manevPhase3,diacu2004manev4}, and references therein. 

\paragraph*{Comparisons to the KS Transformation.}
Given the similar end result of the KS transformation and the transformation developed in this work, some comparative comments are in order. 
The KS transformation has received notably more interaction and use compared to the BF projective transformation 
or its subsequent modifications (including the present work\footnote{Indeed, when writing the first iteration of this work, the authors themselves were unaware of Ferrándiz's BF transformation (canonical extension of Burdet's projective point transformation) and mistakenly thought they were delving into unexplored territory. It was coincidental that our transformation differed from the BF transformation in a manner of which we are fond. It was only upon discovering Ferrándiz's earlier work (predating our own by decades) that the authors were motivated to explore the family of canonical projective transformation starting from a point transformation of the generalized form \eq{\tup{r}=u^n q^m\tup{q}}, as detailed in Appendix \ref{sec:2bp_gen}. The BF transformation corresponds to the case \eq{m=0} and \eq{n=\pm 1}. Our preferred transformation corresponds to the case \eq{m=n=-1}. This difference is shown to have surprising effects on the properties of the resulting Hamiltonian system and meaning of the momenta coordinates. }).
On this note, Deprit et al.~make the the following remark \cite{deprit1994linearization}: 
``We claim [the projective transformation] achieves equally well all the objectives of the KS transformation — linearization, regularization and canonicity [Hamiltonian/symplic structure] — although, we are inclined to believe, in a simpler and more intuitive way''.
The authors are inclined to agree with this comment and are of the opinion that their own canonical extension of the projective transformation presented here is perhaps even more simple and intuitive, though this is subjective.
As is the case for the BF transformation, the transformation we develop in this work has several features in common with the popular KS transformation:
\begin{small}
\begin{itemize}
    \item Both transform the 3-dim Kepler problem into a 4-dim harmonic oscillator (though the natural frequencies are different).  
    \item Both require the use of redundant coordinates (+1 configuration coordinates, +2 phase space coordinates). 
     \item Both involve a transformation of the evolution parameter (though the KS parameter differs from those used here). 
    \item Both can be extended to a canonical transformation compatible with Hamiltonian analytical dynamics (this is central to the goals of the present work). 
\end{itemize}
\end{small} 
Yet, there are several ways in which our transformation is quite different from the KS transformation:
\begin{small}
\begin{itemize}
    \item The KS transformation is quaternion-based whereas ours is projective (similar to the BF transformation). As such, 
    the meaning of our coordinates differs greatly from the KS coordinates.
    \item The KS transformation linearizes the dynamics for Kepler-type forces of the form \eq{-\tfrac{\kconst_1}{r^2}\htup{r}}.  Our projective transformation does the same but for the slightly larger class of Manev-type forces of the form \eq{-\tfrac{\kconst_1}{r^2}\htup{r} -  \tfrac{\kconst_2}{r^3}\htup{r} } for any \eq{\kconst_1,\kconst_2\in\mbb{R}} (with Kepler easily recovered by \eq{\kconst_2=0}). 
    Such Manev-type forces have physically meaningful applications, as mentioned at Eq.\eqref{Manev_def_intro}.
    \item Using our projective transformation, the rotational/angular motion and the radial motion are decoupled from one another. In addition, the rotational part of the dynamics (6 of the 8 phase space dimensions) is invariant and fully linear for \textit{any} arbitrary central-forces. These statements do not hold for the KS transformation. 
    \item The KS transformation, in its usual presentations, is only defined for motion in 2 or 3 dimensions.\footnote{Though methods of generalizing to other dimensions have been investigated \cite{cordani1989generalKS}.}
    The projective transformation, though only explicitly developed here in the classic 3-dim setting, easily generalizes to higher dimensions.\footnote{See \cite{peterson2025prjGeomech} or later chapters of \cite{peterson2025phdThesis} for a geometric formulation of the present work which, among other things, generalizes the developments of this work to arbitrary finite dimensional real inner product spaces.} 
    \item The harmonic oscillator solutions resulting from the KS transformation are valid for closed orbits (\eq{\mscr{E}^0<0}), whereas those obtained by the projective transformation in this work are valid for all energy levels and all types of orbits are treated the same.
\end{itemize}
\end{small} 

\paragraph*{Mathematical Conventions.}
This work is written in the coordinate-dominant manner of classic analytical dynamics with all transformations framed as ``passive'' coordinate transformations. Our mathematical arena is only ever regarded as some \eq{\en}-dim real coordinate vector space, \eq{\mbb{R}^\en}, equipped with the usual structures, operations, and abuses of notation. No mention of vector bundles, differential forms, Lie derivatives, or other geometric constructs will be made.
In addition to the list of notation given at the start of this work, 
some general notation and conventions are as follows: 
\begin{small}
\begin{itemize}
    \item Matrix notation is used frequently, along with occasional ``cartesian index notation''. For the latter, all indices appear as subscripts with summation implied over any repeated indices. 
    Usually, Latin indices such as \eq{i,j,k,l} range from \eq{1} to \eq{3}.\footnote{Though many developments easily generalize to case that \eq{3} is replaced by any finite dimension \eq{\en}.}
    \item Expressions such as \eq{\tup{u}\in \mbb{R}^\en} are often an abuse of notation indicating that \eq{\tup{u}} is an object that \textit{takes values in} \eq{\mbb{R}^\en}. For any such \eq{\tup{u}}, we make no distinction between \eq{\tup{u}} as an ordered \eq{\en}-tuple, \eq{(u_1,\dots,u_\en)}, and \eq{\tup{u}} as a column vector, \eq{\trn{[u_1 \dots u_\en]}}. 
    \item  The norm of some  \eq{\tup{u}\in \mbb{R}^\en} is denoted by 
    \eq{u:=\mag{\tup{u}}:= \sqrt{\tup{u}\cdot \tup{u}}}, and the normalized unit vector by \eq{\htup{u}:=\tfrac{1}{u}\tup{u}} such that \eq{\htup{u}\cdot\htup{u}=1}.
    \item \sloppy  The ``dot product'' denotes contraction in the usual sense on \eq{\mbb{R}^\en}. That is, for \eq{\tup{u},\tup{v}\in \mbb{R}^\en}, then \eq{\tup{u}\cdot\tup{v} = \kd_{ab}u_a \v_b = u_a \v_a} (for \eq{a,b= 1,\dots, \en}). Contraction with some  \eq{M \in\mbb{R}^{\en \times \en}} is denoted either with or without a dot. That is,   \eq{M\tup{u}=M\cdot\tup{u}=M_{ab}u_b \ibase_a} and \eq{\trn{M}\tup{u} = \trn{M}\cdot\tup{u}=\tup{u}\cdot M = M_{ab}u_a \ibase_b = M_{ba}u_b \ibase_a }.   
   \item  The tensor product, \eq{\otimes}, and exterior/wedge product, \eq{\wedge}, also have their usual interpretation on \eq{\mbb{R}^\en}. In particular, the tensor product of any \eq{\tup{u},\tup{v}\in \mbb{R}^\en} is equivalent to the matrix \eq{\tup{u}\otms\tup{v}=\tup{u}\trn{\tup{v}}\in\mbb{R}^{\en\times\en}}, and the exterior product is equivalent to the antisymmetric matrix  \eq{\tup{u}\wdg\tup{v} = \tup{u}\otms \tup{v}-\tup{v}\otms \tup{u} = \tup{u}\trn{\tup{v}}-\tup{v}\trn{\tup{u}} \in\somat{\en}\subset \mbb{R}^{\en\times\en}}.
    \item  \textit{Hodge dual on \eq{\mbb{R}^3}.} For any \eq{\tup{u}\in\mbb{R}^3}, the Hodge dual, \eq{\hdge{\tup{u}}\in\somat{3}}, is the  antisymmetric matrix with components given in terms of the 3-dim Levi-Civita permutation symbol by \eq{\hdge{u}_{ij} := \lc_{ijk}u_k}. For \eq{\mbb{R}^3}, this differs from the usual ``axial dual'', or ``cross product matrix'', \eq{\ax{\tup{u}}}, only by a sign: \eq{\hdge{\tup{u}} = -\ax{\tup{u}}}. We note the following relations (specific to \eq{\mbb{R}^3}):
     \begin{align} \label{hodge_cord}
      \begin{array}{rcll}
         \hdge{u}_{ij} := \lc_{ijk}u_k 
           &\leftrightarrow& u_i = \tfrac{1}{2}\lc_{ijk}\hdge{u}_{jk}
         \\[3pt]
           \hdge{(\tup{u}\tms \tup{v})} =  \tup{u} \wdg \tup{v} &\leftrightarrow &
           \hdge{(\tup{u} \wdg \tup{v})} = \tup{u} \tms \tup{v} 
     \end{array}
     &&,&&
     \begin{array}{lll}
          \hdge{\tup{u}}\cdot\tup{v} =-\tup{u}\tms \tup{v}  = \tup{v}\tms \tup{u} = -\hdge{\tup{v}}\cdot\tup{u}
     \\[3pt]
        \hdge{\tup{u}} \cdot \hdge{\tup{v}} \cdot \tup{w} = \tup{u}\tms (\tup{v}\tms  \tup{w}) = (\tup{v}\wdg\tup{w})\cdot\tup{u}
     \end{array}
     &&,&& 
     \begin{array}{lll}
         \hhdge{\tup{u}} =\tup{u}
     \\[3pt]
        \hdge{\tup{u}}\cdot\hdge{\tup{v}} = \tup{v}\otms \tup{u} - (\tup{u}\cdot\tup{v})\imat_3
     \end{array}
     \end{align}
     \item \textit{A note on ``\eq{y=y(x)}''.}  We (ab)use the notation \eq{y=y(x)}, or just \eq{y(x)}, to indicate there is something called \eq{y} that is a ``function of'' something called \eq{x}. Similarly, under some transformation \eq{x\leftrightarrow s}, the notation \eq{y(x)} and \eq{y(s)} should be interpreted, respectively, as ``\eq{y} expressed in terms of \eq{x}'' and ``\eq{y} expressed in terms of \eq{s}''. 
    \item \textit{A note on mass.} Throughout this work, the mass is scaled out of the Hamiltonian, kinetic/potential energy, forces, angular momentum, etc. All such quantities are given per unit mass (or, for two-body dynamics, per unit 
    \textit{reduced} mass, \eq{m:=\tfrac{m_a m_b}{m_a+m_b}}). For instance, when we refer to the angular momentum we always mean the 
    specific\footnote{The term ``specific angular momentum'' means the angular momentum per unit mass (or per unit reduced mass). The same applies to terms such as ``specific energy''. }
    angular momentum, even if not stated explicitly. 
    \item \textit{A note on inertial cartesian velocity vs.~momentum coordinates.} These are the same for the purposes of this work wherein the (reduced) mass is scaled out. As such, there is there is no difference, in practice, between \textit{inertial cartesian} velocity and momentum coordinates.\footnote{Though, in mathematical/geometric terms, there \textit{is} still a difference:~velocity coordinates are functions on velocity phase space (a tangent bundle) whereas momentum coordinates are functions on phase space (a cotangent bundle) and they therefore cannot be mathematically equal, even in the inertial cartesian case. We ignore such details in this work.}   
    \item \textit{A note on ``canonical transformation''.}   This work often (everywhere other than Appendix \ref{app:Ham_mech_crd}) 
    deals with dimension-raising coordinate transformations at the configuration level that are ``lifted'' or ``extended'' to dimension-raising coordinate transformations at the momentum level in such a way that Hamiltonian/symplectic structure is preserved in the higher-dimensional redundant coordinate description.  
    We often refer to such transformations as ``canonical extensions" or even ``canonical transformations''. Some may object to this language as such dimension-raising transformations do \textit{not} have a unique inverse mapping. Even though these transformations  preserve Hamiltonian/symplectic structure, their non-invertibility would traditionally disqualify them from being regarded as true canonical transformations.  
    Despite this, we often describe such transformations as ``canonical'', with the present note serving as the reader's only warning.\footnote{A more appropriate term for such dimension-raising transformations would perhaps be "symplectic submersion".}
    \item \textit{A note on ``projective''.} Although the term "projective" appears frequently in this work, it should not be interpreted too strictly. 
    In the mathematical context of projective spaces or projective geometry, terms like "projective transformation" and "projective coordinates" have definitions which do not always align exactly with our more relaxed use of such terms. 
    This clash of terminology arises because, in the celestial mechanics literature,  point transformations of the general form seen in Eq.\eqref{proj_pt} are commonly
    referred to  as a ``projective decomposition'' or a ``projective transformation'' (we use either interchangeably) and the resulting coordinates 
    are then referred to as ``projective coordinates'' (Deprit et al.~\cite{deprit1994linearization} attribute this terminology to Ferrándiz). While there is certainly a connection to be made with the mathematical construct of projective spaces, the terminology is not always in one-to-one agreement.
\end{itemize}
\end{small}

%% file: Mysecs_prj/new_prjSUM_noElems.tex
\section{Summary of this work} \label{sec:prj_sum}

 Here, we will give a summary of the canonically-extended  projective coordinate transformation and its application to regularization and linearization of central-force particle dynamics — in particular, Kepler-type and Manev-type dynamics. The derivations and details can be found within the relevant sections, which are outlined as follows:
\begin{small}
\begin{enumerate}
    \item[2.] For the reader's benefit, we summarize below in section \ref{sec:prj_sum} the main developments, dispensing with derivations and details.
    \item[3.] First, in section \ref{sec:Pxform}, before discussing anything to do with regularization or orbital dynamics, we address the general problem of transforming a given Hamiltonian system from a minimal coordinate representation to some new, non-minimal/redundant, coordinate representation. Our approach, based on Hamilton's principle, allows for time-dependent transformations and constraints. 
    \begin{itemize}[topsep=-20pt]
         \item In Appendix \ref{sec:2bp_gen}, the methods from section \ref{sec:Pxform} are used to develop a general family of BF-like transformation. We use a coordinate point transformation, \eq{\mbb{R}^4\ni (\tup{q},u)\mapsto \tup{r}\in\mbb{R}^3}, given by \eq{\tup{r} = u^n \mag{\tup{q}}^m \tup{q}} — for arbitrary scalars \eq{n,m\in\mbb{R}} — whose  ``canonical extension'' induces a corresponding family of canonical transformations, \eq{\mbb{R}^8\ni(\tup{q},u,\tup{p},p_u)\mapsto (\tup{r},\tup{v})\in\mbb{R}^6}.
        By choosing \eq{m=0} and \eq{n=-1}, one recovers the BF transformation \cite{ferrandiz1987general}. 
        We instead choose to focus on the case that \eq{m=-1} (and, eventually, \eq{n=-1}), leading to several different features. 
    \end{itemize} 
  \item[4.] In section \ref{sec:central force}, we focus on the canonical extension of a projective point transformation \eq{(\tup{q},u)\mapsto \tup{r} = u^n \htup{q}} (the case \eq{m=-1}). Before addressing specifically the Kepler problem, we consider arbitrary central-forces along with arbitrary conservative and nonconservative perturbations. It is shown that, by choosing \eq{n=-1},  our transformation fully linearizes particle dynamics for any potential of Kepler type or Manev-type seen in Eq.\eqref{Manev_def_intro}. This requires a transformation of the evolution parameter; two such parameters are given, both achieving linearization: (1) a parameter \eq{s} defined by \eq{\mrm{d} t = r^2 \mrm{d} s}, and (2) a parameter \eq{\tau} defined by \eq{\mrm{d} t = (r^2/\slang)\mrm{d} \tau} (making \eq{\tau} the true anomaly up to an additive constant).  
    \item[5.] In section \ref{sec:2BP}, we focus specifically on linear regularization of Kepler-type and Manev-type dynamics using the canonical extension of a projective point transformation \eq{(\tup{q},u)\mapsto \tup{r} = \tfrac{1}{u} \htup{q}} (the case \eq{n=m=-1}, our preferred transformation).   
    We present closed-form solutions and state transition matrices in terms of evolution parameters \eq{s} and \eq{\tau} (where \eq{\mrm{d} t = r^2 \mrm{d} s} and \eq{\mrm{d} t = (r^2/\slang) \mrm{d} \tau}).  
    The $J_2$-perturbed Kepler problem is formulated in the projective coordinates and used to numerically verify the equations of motion. 
     \item[\scriptsize{\textbullet}]  Other appendices: 
     Appendix \ref{app:Ham_mech_crd} contains a short review of Hamiltonian analytical mechanics with nonconservative forces, as well as with transformations of the evolution parameter using extended phase space; Appendix \ref{sec:J2_plots_prj} contains plots for numerical verification of the \eq{J_2}-perturbed Kepler dynamics in projective coordinates developed in section \ref{sec:2BP} (in particular, \ref{sec:j2}). 
\end{enumerate}
\end{small}
The impatient reader may read the following summary in section \ref{sec:prj_sum} and refer to relevant sections as desired. 
For the reader interested in following the details and derivations in logical sequence, the order of the paper would be: section \ref{sec:Pxform}, Appendix \ref{sec:2bp_gen}, section \ref{sec:central force}, section \ref{sec:2BP}.

 \begin{notesq}
\textsl{For brevity and simplicity,  the following summary will include only conservative central-forces arising from a potential function \eq{V^{0}(r)}.  Additional conservative and nonconservative perturbations of an arbitrary nature are included throughout the main sections of this work.} 
\end{notesq}

\subsection{Summary:~canonical transformations from minimal to redundant coordinates}

In section \ref{sec:Pxform}, we first  treat the general problem of constructing a canonical transformation that is induced by a point transformation between minimal and redundant coordinates. 
We consider some \eq{2\en}-dim phase space with canonical coordinates \eq{(\tup{x},\tup{\piup})\in\mathbb{R}^{2\en}} and known Hamiltonian \eq{\mscr{K}(\tup{x},\tup{\piup},t)} and show that, given some specified point transformation \eq{\tup{x}=\tup{x}(\bartup{q},t)} — where  \eq{\bartup{q}\in\mathbb{R}^{\en+\emm}} are some new set of \textit{redundant} coordinates — along with \eq{\emm} specified independent constraint relations, \eq{0=\tup{\varphiup}(\bartup{q},t)\in\mathbb{R}^\emm}, then we can construct a new set of redundant phase space coordinates \eq{(\bartup{q},\bartup{p})\in\mathbb{R}^{2(\en+\emm)}} with a new Hamiltonian \eq{\mscr{H}(\bartup{q},\bartup{p},t)}. In short, this is accomplished using Hamilton's principle and, in particular, by requiring that:
\begin{align}
    \tup{\piup}\cdot \mrm{d} \tup{x} \,-\, \mscr{K} \mrm{d} t \,=\, \bartup{p}\cdot \mrm{d} \bartup{q} 
    \,-\, \mscr{H} \mrm{d} t \,-\, \tup{\lambdaup}\cdot \mrm{d} \tup{\varphiup} 
\end{align}
for Lagrange multipliers \eq{\tup{\lambdaup}\in\mathbb{R}^\emm}. The momenta transformation is found to be
\begin{align} \label{Pxform_finaly1}
    \bartup{p} =
    \trn{B}
    \begin{pmatrix}
     \tup{\piup} \\ \tup{\lambdaup} 
    \end{pmatrix} 
\quad \leftrightarrow \quad 
    \begin{pmatrix}
     \tup{\piup} \\ \tup{\lambdaup} 
    \end{pmatrix}  = \invtrn{B} \bartup{p}
&&,&&
     B =
     \begin{pmatrix}
         \pderiv{\tup{x}}{\bartup{q}}
         \\[5pt] \pderiv{\tup{\varphiup}}{\bartup{q}}
    \end{pmatrix}
    \in \Glmat{\en+\emm} 
\end{align}
The new coordinates \eq{(\bartup{q},\bartup{p})} satisfy Hamilton's canonical equations for a Hamiltonian \eq{ \mscr{H}(\bartup{q},\bartup{p},t)} given by
\begin{align} 
        \mscr{H} \,=\, \mscr{K} \,-\, \tup{\lambdaup}\cdot\pderiv{\tup{\varphiup}}{t} \,-\, \tup{\piup}\cdot\pderiv{\tup{x}}{t}
\end{align}
with the right-hand-side expressed in terms of \eq{(\bartup{q},\bartup{p})} using the specified \eq{\tup{x}(\bartup{q},t)} and the induced \eq{\tup{\piup}(\bartup{q},\bartup{p},t)} and  \eq{\tup{\lambdaup}(\bartup{q},\bartup{p},t)} from Eq.\eqref{Pxform_finaly1}.

\subsection{Summary:~a family of canonical transformations for projective coordinates}

\begin{notation}
    \sloppy Above, \eq{(\bartup{q},\bartup{p})} were arbitrary redundant phase space coordinates.  We now use \eq{(\bartup{q},\bartup{p})=(\tup{q},u,\tup{p},p_\ss{u})\in\mbb{R}^8} to denote our projective coordinates, with \eq{\bartup{q}=(\tup{q},u)\in\mbb{R}^4} and \eq{\bartup{p}=(\tup{p},p_\ss{u})\in\mbb{R}^4}. 
    We denote by \eq{(\tup{r},\tup{v})\in\mbb{R}^6} inertial cartesian coordinates. 
\end{notation}

\noindent In Appendix \ref{sec:2bp_gen}, we detail a generalized version of the projective point transformation given by \eq{(\tup{q},u)\mapsto\tup{r}= u^n q^m \tup{q}} for arbitrary \eq{n,m\in\mbb{R}}, and then ``lift'', or ``extend'',  this point transformation to a canonical transformation, \eq{(\tup{q},u,\tup{p},p_\ss{u})\mapsto (\tup{r},\tup{v})}. This is summarized below.


\paragraph*{Original Cartesian Coordinate Formulation.}
We start with the Hamiltonian (per unit mass) and canonical equations of motion for a particle in Euclidean 3-space, subject to conservative forces corresponding to some central-force potential function (per unit mass), \eq{V^{0}(r)}, that depends only on the radial distance. 
Letting \eq{(\tup{r},\tup{v})\in\mbb{R}^6} be cartesian position and velocity
(or momentum\footnote{As mentioned,  \textit{inertial cartesian} velocity and momentum coordinates are the same for the purposes of this work wherein the (reduced) mass is scaled out.})
coordinates in an orthonormal inertial frame, the original Hamiltonian system is then given by:
\begin{align} \label{K0}
    &\mscr{K}(\tup{r},\tup{v}) \,=\, \tfrac{1}{2} \v^2 \,+\, V^{0}(r)
&& \begin{array}{ll}
      \dot{\tup{r}} \,=\, \pderiv{\mscr{K}}{\tup{v}} \,=\, \tup{v} 
\qquad,\qquad 
    \dot{\tup{v}} \,=\, -\pderiv{\mscr{K}}{\tup{r}} 
    \,=\, -\pderiv{V^{0}}{r}\htup{r}
\end{array}
\end{align}
where \eq{r :=\mag{\tup{r}}}, and \eq{\htup{r}:=\tup{r}/r}.
When the particular form of \eq{V^{0}(r)} is relevant (for many developments, it is not) then we will consider the Manev-type potential for scalars \eq{\kconst_1,\kconst_2\in\mbb{R}}:
    \begin{flalign} \label{Vmanev_sum}
     \begin{array}{cc}
     \fnsize{Manev-type}  \\
     \fnsize{potential} 
    \end{array} 
    &&
        V^{0} \,=\, -\tfrac{\kconst_1}{r} - \tfrac{1}{2} \tfrac{\kconst_2}{r^2}
        &&
    \end{flalign}
    where Kepler-type dynamics are always easily recovered with \eq{\kconst_2=0}.

\paragraph*{A Family of Canonical Transformations.}  Starting from the above, we then transform from the cartesian coordinates \eq{(\tup{r},\tup{v})\in\mbb{R}^6} to
new, redundant, ``canonical projective coordinates'' \eq{(\bartup{q},\bartup{p})=(\tup{q},u,\tup{p},p_\ss{u}) \in\mbb{R}^8}. 
The usual generalized projective point transformation would be given by \eq{\tup{r}=u^n\tup{q}} for some \eq{n\neq0} (as in Eq.\eqref{proj_pt}). We instead start with the following further-generalized family of projective point transformations: 
\begin{align} \label{PT}
     \gam:\mbb{R}^4\to\mbb{R}^3 \;,
     \qquad \tup{r} \,=\, \gam(\bartup{q}) \,=\, u^n q^m \tup{q}
     &&
     \fnsize{subject to} \;\; \varphi(\tup{q}) = q-1 =0
\end{align}
where \eq{q:=\mag{\tup{q}}\neq \mag{\bartup{q}}} and
where we may choose the values of \eq{n\neq0,m\in \mathbb{R}}.\footnote{We ultimately prefer \eq{n=m=-1}.}
Given that the constraint is equivalent to \eq{q=1}, 
the factor of \eq{q^m} in the above transformation may seem rather pointless.  
It turns out to be important.\footnote{Namely, different values of \eq{m} result in very different momenta coordinates \eq{\tup{p}} conjugate to \eq{\tup{q}}. Additionally, the value \eq{m} affects the  dependence of the radial distance, \eq{r=u^n q^{m+1}}, on the new coordinates \eq{(\tup{q},u)}. In particular, only for \eq{m=-1} is \eq{r=u^n} independent of \eq{\tup{q}}.   } 
The momenta coordinate transformation, along with a Lagrange multiplier \eq{\lambda}, is then obtained as in Eq.\eqref{Pxform_finaly1}, leading to a canonical transformation \eq{\Gam:\mbb{R}^8\to \mbb{R}^6} given 
by\footnote{The point transformation \eq{\gam:\mbb{R}^4\ni(\tup{q},u)\mapsto \tup{r}\in\mbb{R}^3}  is
a \textit{submersion}; it is  surjective with \eq{\text{rnk} \dif \gam \equiv \text{rnk} \pderiv{\tup{r}}{(\tup{q},u)}=3}. It is then lifted to a ``symplectic submersion'', \eq{\Gam:\mbb{R}^8\ni(\tup{q},u,\tup{p},p_\ss{u})\mapsto (\tup{r},\tup{v})\in\mbb{R}^6}, which is surjective with \eq{\text{rnk} \dif \Gam \equiv \text{rnk}\pderiv{(\tup{r},\tup{v})}{(\tup{q},u,\tup{p},p_\ss{u})} =6} (for \eq{n\neq 0}). The domain and codomains of these maps are not actually all of the indicated \eq{\mbb{R}^\en}; they are subsets excluding the case \eq{r=0}, which, in the projective coordinates for \eq{m=0} or \eq{m=-1}, means excluding the case \eq{u=\infty}.  }
\begin{flalign} \label{PT_0}
\qquad
\begin{array}{cc}
   \Gam: \mbb{R}^8\to \mbb{R}^6 
   \\[3pt]
   (\bartup{q},\bartup{p})\mapsto (\tup{r},\tup{v})
\end{array}
\; \left\{ \;\;
\begin{array}{lllll}
      \tup{r} \,=\, u^n q^m \tup{q} 
\\[5pt]
     \tup{v} \,=\,  \tfrac{1}{u^n q^m} \big(  (\imat_3 - \htup{q}\otms\htup{q})\cdot\tup{p} \,+\, \tfrac{u}{n q} p_\ss{u}\htup{q}\big)
\end{array}  \right.
&&,&&
    \begin{array}{rllll}
        &\varphi = q-1 = 0  
         \\[5pt]
          &\lambda=\tfrac{1}{q}(\tup{q}\cdot\tup{p} - \tfrac{m+1}{n} u p_\ss{u})
    \end{array}
\qquad
\end{flalign}
 where \eq{\htup{q}:= \tup{q}/q}, and where \eq{\rnk \dif \gam = 3} and \eq{\rnk \dif \Gam = 6}. 
Note we \textit{specify} the point transformation and constraint in Eq.\eqref{PT} which then \textit{induces} the above momenta transformation and Lagrange multiplier, \eq{\lambda}. This \eq{\lambda} is involved in obtaining an inverse transformation given below in Eq.\eqref{qp_rv_0_gen}-Eq.\eqref{qp_rv_0}.

\begin{notesq}
    \textit{Angular momentum.} It is worth noting that the (specific) 
    angular momentum, \eq{\tup{\slangup}}, takes the same form when expressed in terms of \eq{(\tup{r},\tup{v})} or \eq{(\bartup{q},\bartup{p})}. It is independent of \eq{(u,p_\ss{u})}, depending only on \eq{(\tup{q},
    \tup{p})}:
    \begin{align} \label{lang_qp_SUM}
        \hdge{\tup{\slangup}} 
    \,=\, \tup{r}\wdg \tup{v} \,=\, \tup{q}\wdg\tup{p}
    &&,&&
    \tup{\slangup} 
    \,=\, \tup{r}\tms \tup{v} \,=\, \tup{q}\tms \tup{p} 
    &&,&&
    \slang^2  \,=\, r^2 \v^2 - (\tup{r}\cdot\tup{v})^2 \,=\, q^2 p^2 - (\tup{q}\cdot\tup{p})^2
    \end{align}
    This holds for any choice of \eq{n} and \eq{m} in Eq.\eqref{PT_0}. A collection of useful angular momentum relations is given in Appendix \ref{sec:ang_momentum}. 
\end{notesq}

\noindent Now, the transformation in Eq.\eqref{PT_0} is not time-dependent. As such, it takes the cartesian coordinate Hamiltonian, \eq{\mscr{K}}, to a projective coordinate Hamiltonian, \eq{\mscr{H}}, given by direct substitution, \eq{\mscr{H}=\mscr{K}\circ \Gam}. This leads to:
\begin{align} \label{Ham0}
\mscr{K} = \tfrac{1}{2} \v^2 + V^{0}(r)
\qquad\Rightarrow \qquad
  \mscr{H} \,=\,  \tfrac{1}{u^{2n}q^{2m+2}} \tfrac{1}{2} \big( \slang^2 + \tfrac{1}{n^2} u^2 p_\ss{u}^2 \big)  \,+\, V^{0}(\bartup{q}) 
\end{align}
with \eq{\slang^2 = q^2 p^2 - (\tup{q}\cdot\tup{p})^2} the angular momentum as noted above, and where our (ab)use of notation  \eq{V^{0}(\bartup{q})} simply indicates the original \eq{V^{0}(r)} expressed in the new coordinates using \eq{\eq{\tup{r}=u^n q^m \tup{q}}} and, thus, \eq{r=u^n q^{m+1}}
    (see footnote for clarification\footnote{Note that \eq{\mscr{K}} and \eq{\mscr{H}} are two representations of the same function; the latter is just the former re-expressed in the new coordinates \eq{(\bartup{q},\bartup{p})} (this would not be true for a time-dependent coordinate transformation).
    We denote them by different symbols only for the sake of clearer distinction between the original cartesian coordinate description and the ``new'' projective coordinate description. We generally do \textit{not} make such notational distinctions for other functions or objects encountered in this work. For instance, some generic potential function \eq{V} will simply be denoted \eq{V} in any coordinate description. If greater clarity is needed, we may write, for example,  \eq{V(\tup{r})} or \eq{V(\bartup{q})} to indicate that a specific coordinate description is being used. }).


\paragraph*{An Inverse Transformation.}
Eq.\eqref{PT_0} is a (local) submersion, \eq{\Gam:\mbb{R}^8 \to \mbb{R}^6}, with no  global unique inverse. However, we show (Appendix \ref{sec:lambda_const}) that the new Hamiltonian \eq{\mscr{H}(\bartup{q},\bartup{p})} for the redundant canonical coordinates, \eq{(\bartup{q},\bartup{p})=(\tup{q},u,\tup{p},p_\ss{u})\in\mbb{R}^8}, permits two ``extra'' integrals of motion 
which allow for an inverse transformation. These are: 
\begin{flalign} \label{qlam_0}
     \begin{array}{cc}
     \fnsize{integrals}  \\
     \fnsize{of motion} 
\end{array} 
\left. \quad
\begin{array}{ll}
    q = \mag{\tup{q}} 
     \\[5pt]
    \lambda =  \tfrac{1}{q}(\tup{q}\cdot\tup{p} - \tfrac{m+1}{n} u p_\ss{u})   
\end{array} \right.
&&
\begin{array}{ll}
    \dot{q} = \pbrak{q}{\mscr{H}} + \pderiv{q}{\bartup{p}} \cdot \bartup{\alphaup}  = 0 
   \\[5pt]
    \dot{\lambda}  = \pbrak{\lambda}{\mscr{H}} + \pderiv{\lambda}{\bartup{p}} \cdot \bartup{\alphaup}  = 0
\end{array} 
&&  \xRightarrow{\text{choose}}
&&
\begin{array}{ll}
     q = q_\zr  = 1
     \\[5pt]
       \lambda = \lambda_\zr  = 0
\end{array}
\quad
\end{flalign}
where \eq{\bartup{\alphaup}=(\tup{\alphaup},\alpha_\ss{u})} are  generalized nonconservative forces (omitted from the present summary). 
That is, the above functions \eq{q} and \eq{\lambda} — our constraint and Lagrange multiplier — are
integrals of motion of the transformed Hamiltonian system \textit{even in the presence of arbitrary, perhaps nonconservative, forces}. 
When transforming/specifying the initial conditions, we are free to limit consideration to any chosen values \eq{0<q_0 \in\mbb{R}_\ii{+}}
and \eq{\lambda_0\in\mbb{R}}, and they will then remain constant at these values for all time. 
Identifying these two \textit{functions} with their constant \textit{values} allows us to invert Eq.\eqref{PT_0} as:
\begin{align}  \label{qp_rv_0_gen}
\begin{array}{cc}
     \fnsize{restrict to}  \\
     \fnsize{$q=q_0$} \\
      \fnsize{$\lambda=\lambda_0$} 
\end{array} 
\quad \Rightarrow \qquad
 (\tup{r},\tup{v}) \mapsto (\bartup{q},\bartup{p})
\quad \left\{\quad
\begin{array}{llll}
 \tup{q} \,=\,  q_0\htup{r} 
      &,\quad 
      \tup{p} \,=\, \tfrac{r}{q_0}( \imat_3 + m \htup{r}\otms\htup{r})\cdot \tup{v} + \lambda_0 \htup{r}
  \\[5pt]
    u \,=\, \big(\tfrac{r}{ q_0^{m+1}} \big)^\ss{1/n} 
   &, \quad 
   p_\ss{u} \,=\,  n \big(\tfrac{r}{q_0^{m+1}} \big)^\ss{-1/n} \tup{r}\cdot\tup{v}  
\end{array} \right.
\end{align}
The above is not a true inverse of Eq.\eqref{PT_0}; it is only unique up to some chosen values of \eq{q} and \eq{\lambda}.
We consider the specific values
\eq{q = 1} and \eq{\lambda=0} in which case the above 
     simplifies to:\footnote{\eq{q=1} ensures that \eq{\tup{q}=\htup{q}=\htup{r}} is the radial unit vector. The value \eq{\lambda=0} is chosen because we see no advantage in considering any other value and because  \eq{\lambda=0} leads to \eq{\tup{p}} having a ``nice'' interpretation.}
\begin{align}  \label{qp_rv_0}
\begin{array}{cc}
     \fnsize{restrict to}  \\
    \fnsize{$q=1$} \\
      \fnsize{$\lambda=0$}
\end{array} 
\quad \Rightarrow \qquad
 (\tup{r},\tup{v}) \mapsto (\bartup{q},\bartup{p})
\quad \left\{\quad
\begin{array}{llll}
     \tup{q} \,=\,  \htup{r} 
      &,\quad \tup{p}  \,=\,  r( \imat_3 + m \htup{r}\otms\htup{r})\cdot \tup{v}
  \\[5pt]
  u \,=\, r^\ss{1/n}  
   &, \quad 
   p_\ss{u} \,=\, 
      n r^\ss{-1/n } \tup{r}\cdot\tup{v} 
     \,=\,  n r^\ss{(n-1)/n}\,\dot{r}
\end{array}  \right.
\end{align}
which is our chosen default inverse transformation of Eq.\eqref{PT_0}.

\begin{notesq}
    We will often use ``\eq{\simeq}'' to indicate relations which have been simplified using the integrals of motion \eq{q=1} and \eq{\lambda=0}. We do not need to do anything to enforce such relations other than transform initial conditions \eq{(\tup{r}_\zr,\tup{v}_\zr) \mapsto (\bartup{q}_\zr,\bartup{p}_\zr)} using Eq.\eqref{qp_rv_0}. This places no restrictions on the cartesian coordinates. 
\end{notesq}

\subsection{Summary:~a preferred projective coordinate set}

The values of \eq{n\neq 0,m\in\mbb{R}}  in the initial point transformation \eq{\tup{r}= u^n q^m \tup{q}} of Eq.\eqref{PT},   are ``knobs'' that we may turn to construct different coordinate transformations:  
\textit{different values of \eq{m} give different  \eq{\tup{p}} while different values of \eq{n} give different \eq{u} and \eq{p_\ss{u}}}.\footnote{For instance, for \eq{m=0} vs.~\eq{m=-1}, or for \eq{n=1} vs.~\eq{n=-1}, then we see from Eq.\eqref{qp_rv_0} that: 
\begin{align} \nonumber
\begin{array}{rllllll}
      m=0:& \tup{q}=\htup{r} &,\quad
      \tup{p}= r\tup{v} 
     & = r \dot{\tup{r}}
 \\[4pt]
 m=-1:& \tup{q}=\htup{r}&,\quad
 \tup{p} = -\hdge{\tup{\slangup}}\cdot\htup{r} 
  & = r^2\dot{\htup{r}}
\end{array}
 &&,&&
 \begin{array}{rllllll}
     n=1:&  u=r &,\quad
     p_\ss{u} = \htup{r}\cdot\tup{v} 
      & = \dot{r}
 \\[4pt]
  n=-1:&  u= 1/r  &,\quad
      p_\ss{u} = -r^2 \htup{r}\cdot\tup{v}
       & = -r^2 \dot{r}
\end{array}
\end{align}
}
It would seem intuitive to choose \eq{m=0} 
such that the point transformation is simply the classic  \eq{\tup{r}= u^n\tup{q}}. Indeed, using \eq{m=0} along with \eq{n=-1} (i.e., \eq{\tup{r}=\tfrac{1}{u}\tup{q}}) is precisely the point transformation used by Ferrándiz to develop a canonical/Hamiltonian version of Burdet's projective transformation for regular and linear Kepler dynamics \cite{ferrandiz1987general}. 
However, we ultimately prefer a slightly different point transformation (which leads to significant differences at the momenta level):

\begin{notesq}
\textit{We prefer the projective transformation for the values \eq{n=m=-1}, that is, the canonical extension of a point transformation \eq{\tup{r}=\tfrac{1}{u q}\tup{q}}.} 
Some reasons for this are given below (with further discussion in Appendix \ref{sec:BF_vs_mine} and section \ref{sec:central force}). 
\end{notesq}

\noindent 
With \eq{n=m=-1}, the family of transformations \eq{(\bartup{q},\bartup{p})\leftrightarrow (\tup{r},\tup{v})} in Eq.\eqref{PT_0} and Eq.\eqref{qp_rv_0} then lead to the following:
\begin{align} \label{rv_qu_full}
\begin{array}{lllllllll}
 \Gam: (\bartup{q},\bartup{p})\mapsto (\tup{r},\tup{v})
    & \left\{ \;\; \boxed{ \begin{array}{llllll}
          \tup{r} \,=\, \tfrac{1}{u}\htup{q}  
         &,\quad 
         \tup{v} & =\,  u q  (\imat_3 - \htup{q}\otms\htup{q})\cdot\tup{p} - u^2 p_\ss{u}\htup{q}
         \\[4pt]
        &   & =\,
         - u \hdge{\tup{\slangup}}\cdot\htup{q} - u^2 p_\ss{u} \htup{q}
    \end{array} } \right. 
    \qquad\qquad\qquad
    \begin{array}{rllll}
         \fnsize{constaint:} &\varphi = q-1 = 0  
         \\[4pt]
          &\lambda=\htup{q}\cdot\tup{p}
    \end{array}
\\[14pt]
    \inv{\Gam} : (\tup{r},\tup{v}) \mapsto  (\bartup{q},\bartup{p})
    &\left\{ \;\; 
    \boxed{ \begin{array}{lllllll}
         \tup{q} \,=\,  \htup{r} 
        &,\quad 
          \tup{p}  \,=\,  r( \imat_3 - \htup{r}\otms\htup{r})\cdot \tup{v}
             \,=\, -\hdge{\tup{\slangup}}\cdot\htup{r}
    \\[5pt]
      u \,=\, 1/r
       &, \quad 
       p_\ss{u} \,=\, 
          - r^2\htup{r}\cdot\tup{v} 
  \end{array} } \right.
    \begin{array}{llll}
          =\, r^2 \dot{\htup{r}}
     \\[6pt]
           =\, -r^2 \dot{r}
    \end{array}
\end{array}
\end{align}
where \eq{\hdge{\tup{\slangup}} = \tup{r}\wdg\tup{v} = \tup{q}\wdg\tup{p}} and where the above inverse transformation for \eq{(\bartup{q},\bartup{p})= \inv{\Gam}(\tup{r},\tup{v})} follows from restricting consideration of the integrals of motion \eq{q} and \eq{\lambda = \htup{q}\cdot\tup{p}} to the values \eq{q=1} and \eq{\lambda =0}.
The Hamiltonian and time-parameterized equations of motion are then given in the above projective coordinates by 
\begin{align} \label{qpdot_SUM}
\begin{array}{ccccc}
\mscr{K} = \tfrac{1}{2} \v^2 + V^{0}(r)
\qquad\Rightarrow \qquad
\\[6pt]
\vphantom{ \dot{q}_i = \pderiv{\mscr{H}}{p_i}}
\\[6pt]
\vphantom{ \dot{q}_i = \pderiv{\mscr{H}}{p_i}}
\end{array}  
\begin{array}{ccccc}
    \boxed{\mscr{H}  \,=\, 
    \tfrac{1}{2} u^2 (\slang^2 + u^2 p_\ss{u}^2) + V^{0}(u) }
    \phantom{XXXXXXXXX}
\\[10pt] 
\begin{array}{lll}
     \dot{\tup{q}} = \pderiv{\mscr{H}}{\tup{p}} \,=\,  -u^2 \hdge{\tup{\slangup}}\cdot \tup{q}&,
\\[6pt]
    \dot{u}  = \pderiv{\mscr{H}}{p_u}  \,=\,  u^4 p_\ss{u} &,
 \end{array}
\quad
\begin{array}{ll}
    \dot{\tup{p}} = -\pderiv{\mscr{H}}{\tup{q}} \,=\,
    - u^2 \hdge{\tup{\slangup}}\cdot \tup{p}
\\[6pt]
     \dot{p}_\ss{u}  = -\pderiv{\mscr{H}}{u} \,=\,
    -u (\slang^2 +  2 u^2 p_\ss{u}^2 )   -  \pderiv{V^{0}}{u}  
\end{array}  
\end{array} 
\end{align}
where \eq{\slang^2 = q^2 p^2 - (\tup{q}\cdot\tup{p})^2} and \eq{\hdge{\tup{\slangup}}} are integrals of motion for any cental-force dynamics.
An evolution parameter transformation is still needed before the above yields what we are chasing (linearity). First, we note some features of the transformation in Eq.\eqref{rv_qu_full} (i.e., features that follow from choosing \eq{m=n=-1} in Eq.\eqref{PT_0}):
\begin{small}
\begin{itemize}
    \item \textit{Central-forces.} Only for \eq{m=-1} (and any \eq{n}) does \eq{r=\mag{\tup{r}} = u^n q^{m+1} = u^n} depend only on \eq{u}, but not \eq{\tup{q}}. 
    Thus, a central-force potential \eq{V^{0}(r)} is taken to some \eq{\tup{q}}-independent \eq{V^{0}(u)} — \textit{any central-forces are eliminated from the Hamiltonian dynamics for \eq{(\tup{q},\tup{p})}, and appear only in those for \eq{(u,p_\ss{u})}}.\footnote{We show that the same turns out to be true of nonconservative central-forces.}
    This does not hold for the BF transformation (\eq{m=0}).  
    \item \textit{Relation to LVLH basis.} For any \eq{n,m\in\mbb{R}} the angular momentum satisfies \eq{\tup{\slangup}=\tup{q}\tms\tup{p}} such that \eq{ \tup{q}\cdot\tup{\slangup} = \tup{p}\cdot\tup{\slangup} = 0}. But, only for \eq{m=-1} (and any \eq{n}) does it hold that \eq{\lambda = \htup{q}\cdot\tup{p}} such that \eq{\lambda=0} is equivalent to \eq{\tup{q}\cdot\tup{p}=0} and, thus, \eq{(\tup{q},\tup{p},\tup{\slangup})} are mutually orthogonal. In fact,  their normalization is precisely the inertial cartesian components of the LVLH basis. That is, the integrals of motion \eq{\lambda=0} and \eq{q=1} imply that the coordinates \eq{(\tup{q},\tup{p})} in Eq.\eqref{rv_qu_full} satisfy: 
    \begin{flalign} \label{qp_m1_props_sum}
       \fnsz{\begin{array}{cc}
            \fnsize{using:} \\[2pt]
             q=1  \\[2pt]
             \lambda=\htup{q}\cdot\tup{p}=0 
        \end{array}}
         \Rightarrow \;\; 
        \left\{
        \quad
        \begin{array}{llll}
            \tup{q} \simeq \htup{q} \simeq \htup{p}\tms \htup{\slangup} 
          \\[4pt]
            \tup{p} \simeq \tup{\slangup}\tms \htup{q} 
        \\[4pt]
              \slang^2 \simeq q^2 p^2 \simeq p^2  
        \end{array} \right.
        \qquad,\qquad
        \begin{array}{lll}
            \hdge{\tup{\slangup}}\cdot\tup{q} \simeq -q^2\tup{p} \simeq -\tup{p}
            \\[4pt]
            \hdge{\tup{\slangup}}\cdot\tup{p} \simeq p^2\tup{q} \simeq \slang^2 \tup{q} 
        \end{array}
         \qquad,\qquad
        \begin{array}{lll}
            \imat_3 \simeq \htup{q} \otms\htup{q} + \htup{p} \otms\htup{p} + \htup{\slangup} \otms\htup{\slangup}
            \\[4pt]
             \{\htup{q},\htup{p},\htup{\slangup}\} \simeq \{\htup{t}_r,\htup{t}_\tau, \htup{t}_\slang \}
             = \!\!
             \begin{array}{cc}
            \fnsize{LVLH}\\[-1pt]
            \fnsize{basis}  
            \end{array}
        \end{array}
        \quad
    \end{flalign}
    The relations in Eq.\eqref{qp_m1_props_sum} hold along any solution curve of Eq.\eqref{qpdot_SUM} for which \eq{q_\zr=1} and \eq{\lambda_\zr=0}. In that case, the ODEs themselves could be simplified using these 
    relations.\footnote{E.g., the \eq{(\tup{q},\tup{p})} dynamics in Eq.\eqref{qpdot_SUM} simplify to \eq{\dot{\tup{q}} \simeq  u^2 \tup{p}} and \eq{ \dot{\tup{p}} \simeq - u^2 p^2 \tup{q} \simeq  - u^2 \slang^2 \tup{q}}.}
    \begin{small}
    \begin{itemize}[nosep]
        \item In contrast, for the the BF transformation  (\eq{m=0} and \eq{n=-1}), then \eq{\lambda=0} would instead give \eq{\bartup{q}\cdot\bartup{p}=\tup{q}\cdot\tup{p}+u p_\ss{u} = 0} and the above 
    would not hold.
    \end{itemize}
    \end{small}
    \item \textit{Recovering cartesian coordinate solutions.}
    Given some solution curve \eq{(\bartup{q}_t,\bartup{p}_t)}, the transformation \eq{\Gam:(\bartup{q},\bartup{p})\mapsto(\tup{r},\tup{v})} from Eq.\eqref{rv_qu_full} can, of course, be used to recover the corresponding cartesian coordinate solution \eq{(\tup{r}_t,\tup{v}_t)}. 
    Yet, when converting numerical solutions, the transformation in Eq.\eqref{rv_qu_full} may be simplified considerably using the integrals of motion \eq{q=1} and \eq{\lambda=\htup{q}\cdot\tup{p}=0}, leading to:
    \begin{flalign} \label{rv_qu_simp}
    &&
    \begin{array}{lllllll}
          \tup{r} \,=\, \tfrac{1}{u}\htup{q} 
    \\[4pt]
          \tup{v} \,=\, -u(\tup{q}\wdg\tup{p})\cdot\htup{q} - u^2 p_\ss{u} \htup{q} 
    \end{array}
    \qquad
    \xRightarrow[\htup{q}\cdot\tup{p}\,=\,0]{q\,=\,1}
    \qquad 
     \boxed{ \begin{array}{llll}
        \tup{r}_t  \,\simeq\, \tfrac{1}{u_t}\tup{q}_t
    \\[4pt]
        \tup{v}_t 
        \,\simeq\, 
        u_t\tup{p}_t - w_t \tup{q}_t 
    \end{array} }
    &&
    w:=u^2 p_\ss{u}
    \qquad
    \end{flalign}
    where \eq{w= u^2 p_\ss{u}} is introduced for later convenience. 
    So long as \eq{(\bartup{q}_t , \bartup{p}_t)} is a solution curve starting with \eq{q_\zr=1} and \eq{\lambda_\zr=\htup{q}_\zr\cdot\tup{p}_\zr=0}, then \eq{(\tup{r}_t,\tup{v}_t)} solutions recovered using the above simplified relations will be \textit{numerically} equivalent to those recovered using the full projective transformation in Eq.\eqref{rv_qu_full}. 
    \textit{The simplified relations in Eq.\eqref{rv_qu_simp} above are \emph{not} the projective transformation}.\footnote{We stress that Eq.\eqref{rv_qu_full} is the actual map for the projective transformation that is used to construct the Hamiltonian, derive equations of motion, transform forces, etc. The simplified relations Eq.\eqref{rv_qu_simp} should only be used to transform numerical points/solutions (assuming \eq{q_\zr=1} and \eq{\lambda_\zr=\htup{q}_\zr\cdot\tup{p}_\zr =0}). }
\end{itemize}
\end{small}

\subsection{Summary:~time reparameterization, linearization, closed-form solutions}  \label{sec:SUM_linear_sol}

    The following is developed in sections \ref{sec:ext_cen}-\ref{sec:2BP_lin}, and restated concisely in section \ref{sec:2BP} in the context of Kepler and Manev dynamics. 

Continuing from the above, we use the extended phase space\footnote{See Appx.~\ref{sec:ext_coord} for  details on extended phase space.}
to transform the evolution parameter from the time, \eq{t}, to two new evolution parameters, \eq{s} and \eq{\tau}, which are related to the time through:\footnote{The evolution parameter \eq{s} is the same as used by Burdet \cite{burdet1969mouvement,Burdet+1969+71+84}, \eq{\tau} (true anomaly) is equivalent to the one used by Vitins \cite{vitins1978keplerian},  and Ferrándiz considered both \eq{s} and \eq{\tau} \cite{ferrandiz1987general}.} 
\begin{align} \label{dtds_0}
      \mrm{d} t \,=\, r^2 \mrm{d} s \,=\, u^\ss{-2} \mrm{d} s
\qquad,\qquad 
     \mrm{d} t \,=\, \tfrac{r^2}{\slang} \mrm{d} \tau 
     \,=\, \tfrac{1}{\slang u^2} \mrm{d} \tau 
\qquad,\qquad     
      \mrm{d} \tau \,=\, \slang \mrm{d} s
\end{align} 
where \eq{\slang^2 =q^2 p^2 -(\tup{q}\cdot\tup{p})^2 } and 
where \eq{\tau} is equivalent to the true anomaly up to an additive constant.
For any central-force dynamics, \eq{\tup{\slangup}=\tup{\slangup}_\zr} is conserved such that \eq{s} and \eq{\tau} are then related simply as follows (assuming \eq{\tau_\zr=s_\zr=0}):
\begin{align}
    \fnsize{if} \; \slang = \slang_\zr \qquad \Rightarrow\qquad  \tau = \slang s
\end{align}
Either \eq{s} or \eq{\tau} can be used, in conjunction with the projective coordinates, to linearize orbital dynamics. 
We first consider the parameter \eq{s}. 
The extended Hamiltonian \eq{\wtscr{H}(\bartup{q},\bartup{p},t,p_t)} for \eq{s} as the evolution parameter is then defined as \eq{\wtscr{H} := \diff{t}{s}(\mscr{H} + p_t )}, where \eq{\diff{t}{s}=1/u^2} and where \eq{p_t} is the momenta conjugate to the time, \eq{t}. 
The \eq{s}-parameterized Hamiltonian dynamics are then given in projective coordinates by:
\begin{gather} \label{dqp_s_0}
   \boxed{ \wtscr{H} 
    \,=\, \tfrac{1}{2}\big(\slang^2 + u^2 p_\ss{u}^2 \big)  \,+\, \wt{V}^\zr(u)  \,+\,  
 u^\ss{-2} p_t }
\\ \nonumber
\begin{array}{llll}
     \pdt{\tup{q}} \,=\, \pderiv{\wtscr{H}}{\tup{p}}   \,=\, 
      -\hdge{\tup{\slangup}}\cdot\tup{q} 
\\[5pt]
      \pdt{\tup{p}} \,=\, -\pderiv{\wtscr{H}}{\tup{q}} 
      \,=\,  -\hdge{\tup{\slangup}}\cdot\tup{p} 
\\[4pt]
    \vphantom{\pderiv{V^{0}}{u}}
\end{array}
\quad,\quad 
\begin{array}{llll}
    \pdt{u} \,=\, \pderiv{\wtscr{H}}{p_u}  &=\,   u^2 p_\ss{u}
\\[5pt]
      \pdt{p}_\ss{u} \,=\, -\pderiv{\wtscr{H}}{u}
      &=\, -u p_\ss{u}^2    -  \pderiv{\wt{V}^0}{u} +  \tfrac{2}{u^3} p_t 
\\[4pt]
      &=\, -\tfrac{1}{u}\big( \slang^2 + 2 u^2 p_\ss{u}^2 \big)   - \pdt{t}\pderiv{V^{0}}{u}
\end{array}
\quad,\quad 
\begin{array}{llll}
    \pdt{t}  \,=\, \pderiv{\wtscr{H}}{p_t} \,=\,  1/u^2
\\[5pt]
     \pdt{p}_t \,=\, -\pderiv{\wtscr{H}}{t}
     \,=\, 0 
 \\[4pt]
    \vphantom{\pderiv{V^{0}}{u}}
\end{array}
\end{gather}
where  \eq{\pdt{\square}:=\diff{\square}{s}= u^\ss{-2}\diff{\square}{t}} and where \eq{\wt{V}^0:=  u^\ss{-2} V^0} still depends only on \eq{u}.
The second expression above for \eq{\pdt{p}_\ss{u}} has used the well-known relation  \eq{p_t=-\mscr{H}} (cf.~Appx.~\ref{sec:ext_coord} or \cite{struckmeier2005hamiltonian,lanczos2012variational}) to eliminate \eq{p_t} from the \eq{\pdt{p}_\ss{u}} equation.\footnote{Differentiating \eq{\wtscr{H}} directly, and using \eq{\pderiv{\wt{V}}{u} = -\tfrac{2}{u^3}V + \tfrac{1}{u^2}\pderiv{V}{u}}, leads to: 
\begin{align} \nonumber
    \pdt{p}_\ss{u} = -\pderiv{\wtscr{H}}{u}  = -u p_\ss{u}^2    -  \pderiv{\wt{V}}{u} +  \tfrac{2}{u^3} p_t 
    \quad =\,  -u p_\ss{u}^2 - \tfrac{1}{u^2}\pderiv{V}{u} +  \tfrac{2}{u^3}(V + p_t )
    \,=\, 
    - \tfrac{1}{u}(\slang^2 +  2 u^2 p_\ss{u}^2 )   -  \tfrac{1}{u^2} \pderiv{V}{u} + \tfrac{2}{u} \wtscr{H}
    \quad=\, -\tfrac{1}{u}\big( \slang^2 + 2 u^2 p_\ss{u}^2 \big)   - \tfrac{1}{u^2}\pderiv{V}{u}
\end{align} 
where the last equality follows from substitution of \eq{p_t= -\mscr{H}= - \tfrac{1}{2}u^2 (\slang^2 + u^2 p_\ss{u}^2) - V }, leading to the equation for \eq{\pdt{p}_\ss{u}=\tfrac{1}{u^2}\dot{p}_\ss{u}} seen in Eq.\eqref{dqp_s_0}. 
}

We also consider the \eq{\tau}-parameterized dynamics, obtained in a similar manner as Eq.\eqref{dqp_s_0} but using an extended Hamiltonian \eq{\whscr{H}:=\diff{t}{\tau}(\mscr{H}+p_t)=\tfrac{1}{\slang}\wtscr{H}}, where \eq{\diff{t}{\tau}=\tfrac{1}{\slang u^2}}. This leads to the following (with \eq{\rng{\square}:=\diff{\square}{\tau}}):
\begin{gather} \label{dqp_TA_SUM}
    \whscr{H}\,=\,  \tfrac{1}{\slang}\wtscr{H} 
\qquad,\qquad 
\begin{array}{llll}
      \rng{\tup{q}}    \,=\, 
      -\hdge{\htup{\slangup}}\cdot\tup{q} 
\\[4pt]
     \rng{\tup{p}} 
      \,=\,  -\hdge{\htup{\slangup}}\cdot\tup{p} 
\end{array}
\quad,\quad
\begin{array}{lll}
     \rng{u} \,=\,  \tfrac{1}{\slang} u^2 p_\ss{u}
\\[4pt]
      \rng{p}_\ss{u}  
      \,=\, -\tfrac{1}{\slang u}\big( \slang^2 + 2 u^2 p_\ss{u}^2 \big)  - \rng{t}\pderiv{V^{0}}{u}
\end{array}
\quad,\quad
\begin{array}{lll}
      \rng{t} \,=\,  \tfrac{1}{\slang u^2}
\\[4pt]
     \rng{p}_t \,=\,  0 
\end{array}
\end{gather}
where the relation \eq{p_t=-\mscr{H}} has been used to rewrite the extended phase space ODEs in the form \eq{\rng{\square}=\tfrac{1}{\slang}\pdt{\square} = \tfrac{1}{\slang u^2}\dot{\square}} seen above.\footnote{See the developments leading to Eq.\eqref{qpdot_ext_alt_apx} in Appendix \ref{sec:ext_coord} for details.}

\begin{notesq}
    \textit{A quasi-momenta coordinate.} It can be convenient to exchange the conjugate momentum coordinate \eq{p_\ss{u}} for 
    a quasi-momentum coordinate \eq{w:=u^2 p_\ss{u}}.\footnote{That is:
        \eq{\quad  w  :=\,  u^2 p_\ss{u}  \,=\, \pdt{u} \,=\,  -\dot{r}
         \quad \leftrightarrow \quad
        p_\ss{u} \,=\, w/u^2  \,=\, r^2 w   \,=\,  -\pdt{r} } .  }
    For instance, replacing the pair \eq{(u,p_\ss{u})} with \eq{(u,w)}, the above \eq{s}- or \eq{\tau}-parameterized ODEs are equivalent to the following:
    \begin{flalign} \label{uw_eom_gen}
    w := u^2 p_\ss{u} 
     \quad \Rightarrow \quad\;\;
    \boxed{\begin{array}{llllll}
         \pdt{\tup{q}}  \,=\, 
      -\hdge{\tup{\slangup}}\cdot\tup{q} 
      &,\quad 
       \pdt{u} \,=\, w
    \\[4pt]
        \pdt{\tup{p}}  \,=\, 
      -\hdge{\tup{\slangup}}\cdot\tup{p} 
        &,\quad 
         \pdt{w} = -\slang^2 u - \pderiv{V^0}{u} 
    \end{array}
    \qquad\quad  \fnsize{or,}  \qquad\quad 
     \begin{array}{llllll}
         \rng{\tup{q}}  \,=\, 
      -\hdge{\htup{\slangup}}\cdot\tup{q} 
      &,\quad 
       \rng{u} \,=\,  w/\slang 
    \\[4pt]
       \rng{\tup{p}}  \,=\, 
      -\hdge{\htup{\slangup}}\cdot\tup{p} 
        &,\quad 
         \rng{w} = -\slang u - \tfrac{1}{\slang} \pderiv{V^0}{u} 
    \end{array} }
    \end{flalign}
    where \eq{\tup{\slangup}=\tup{\slangup}_\zr} is conserved such that the above systems are fully linear for certain forms of \eq{V^\zr} (e.g., Kepler or Manev). 
    Interestingly, the quasi-momentum coordinate \eq{w=u^2 p_\ss{u}} and the conformal factor \eq{\pdt{t}=1/u^2} cancel out in such a way that the above \eq{s}-parameterized ODEs for \eq{(u,w)} obey Hamilton's canonical equations in their usual form (without the scale factor):
    \begin{align} \label{uw_eom_gen_intereting}
         \mscr{H} = \tfrac{1}{2} u^2 ( \slang^2 +  u^2 p_\ss{u}^2 )  + V^{0} \,=\,  \tfrac{1}{2} ( u^2  \slang^2 +  w^2 )  + V^{0}
         &&,&&
         \begin{array}{llllll}
                \pdt{\tup{q}} = \pdt{t}\pderiv{\mscr{H}}{\tup{p}} 
                 &, \quad
                  \pdt{\tup{p}} = -\pdt{t}\pderiv{\mscr{H}}{\tup{q}} 
           \\[5pt]
                \pdt{u} = \pderiv{\mscr{H}}{w}
                &, \quad \pdt{w} = -\pderiv{\mscr{H}}{u} 
         \end{array}
    \end{align}
\end{notesq}

\begin{notesq}
     \textit{Simplifications.} 
     The dynamics in Eq.\eqref{dqp_s_0}-Eq.\eqref{uw_eom_gen} have \textit{not} been simplified with \eq{q=1} or \eq{\lambda=\htup{q}\cdot\tup{p}=0}, which imply the relations in Eq.\eqref{qp_m1_props_sum}. 
     We are free to make these simplifications such that, for example, the above ODEs for \eq{(\tup{q},\tup{p})} are equivalent to:
    \begin{flalign}
    &&
        \begin{array}{lllll}
        \pdt{\tup{q}} \simeq \tup{p} 
         \\[4pt]
         \pdt{\tup{p}} 
         \simeq -\slang^2\tup{q} 
    \end{array}
        \qquad\fnsize{or,}\qquad 
     \begin{array}{lllll}
         \rng{\tup{q}} \simeq \tfrac{1}{\slang}\tup{p} \simeq \htup{p}
         \\[4pt]
         \rng{\tup{p}} \simeq -\slang \tup{q} 
    \end{array}
    &&
        \slang^2 \simeq p^2 \qquad
    \end{flalign}
\end{notesq}

\paragraph*{Rotational Motion.} The rotational/angular motion is described by the coordinates \eq{(\tup{q},\tup{p})}. For central-force dynamics, \eq{\tup{\slangup}=\tup{\slangup}_\zr} is conserved such that the above \eq{s}- or \eq{\tau}-parameterized ODEs for \eq{(\tup{q},\tup{p})} in Eq.\eqref{dqp_s_0} or Eq.\eqref{dqp_TA_SUM}:
\begin{flalign} \label{qp_Frad_SUM}
&\begin{array}{cc}
   \fnsize{arbitrary}  \\
   \fnsize{central-force}\\
   \fnsize{dynamics} 
 \end{array}
 &&
 \begin{array}{lllll}
       \pdt{\tup{q}} = 
      -\hdge{\tup{\slangup}}\cdot\tup{q} 
\\[4pt]
      \pdt{\tup{p}} =  -\hdge{\tup{\slangup}}\cdot\tup{p} 
\end{array}
\qquad 
\fnsize{or,} \qquad 
\begin{array}{lllll}
       \rng{\tup{q}}    = 
      -\hdge{\htup{\slangup}}\cdot\tup{q} 
 \\[4pt]
     \rng{\tup{p}} =  -\hdge{\htup{\slangup}}\cdot\tup{p}
\end{array} 
 &&
\end{flalign}
are linear equations with constant coefficients \eq{\hdge{\tup{\slangup}}, \hdge{\htup{\slangup}} \in \somat{3}}. The solutions to either of  the above coincide and are given by the matrix exponential \eq{\mrm{e}^{-\hdge{\tup{\slangup}}s} = \mrm{e}^{-\hdge{\htup{\slangup}}\tau} =: R_\tau(\htup{\slangup}) \in \Somat{3}}, which is simply a special orthogonal rotation by \eq{\tau=\slang s} about \eq{\htup{\slangup}=\htup{q}\tms\htup{p}} (the orbit normal direction). Using the Rodrigues rotation formula, this leads to:
\begin{flalign} \label{dqp_ug_sol}
 &&
\begin{array}{lllllll}
     \tup{q}_\tau 
      \,=\,
      R_\tau(\htup{\slangup}) \cdot\tup{q}_\zr \;=\;
      \tup{q}_\zr \csn{\tau} -  \hdge{\htup{\slangup}} \cdot \tup{q}_\zr \snn{\tau} 
      &\;\;\simeq\, \tup{q}_\zr \csn{\tau} + \tfrac{1}{\slang} \tup{p}_\zr \snn{\tau} 
\\[4pt]
      \tup{p}_\tau 
      \,=\, 
      R_\tau(\htup{\slangup}) \cdot\tup{p}_\zr \;=\; 
      \tup{p}_\zr \csn{\tau} -  \hdge{\htup{\slangup}} \cdot \tup{p}_\zr \snn{\tau} 
      &\;\;\simeq\,  \tup{p}_\zr \csn{\tau} - \slang \tup{q}_\zr \snn{\tau} 
\end{array} 
&&
\tau = \slang s
\qquad
\end{flalign}
where  ``\eq{\simeq}'' indicates relations simplified using \eq{q=1} and \eq{\htup{q}\cdot\tup{p}=0} (and thus \eq{\slang^2\simeq p^2}). 

\paragraph*{Radial Motion.} The radial motion is described by the coordinates \eq{(u,p_\ss{u})} or, alternatively, by the non-conjugate pair \eq{(u,w:=u^2 p_\ss{u})}.
Unlike the rotational/angular dynamics, the transformed radial dynamics are not linear for any arbitrary central-forces. 
 Yet, they \textit{are} linear for Kepler-type or, more generally, Manev-type central-forces. That is, we now consider specifically a central-force potential of the form seen in Eq.\eqref{Vmanev_sum}:
\begin{flalign}
\begin{array}{cc}
     \fnsize{Manev-type}  \\
      \fnsize{potential} 
\end{array} 
&&
     V^{0} = -\tfrac{\kconst_1}{r} - \tfrac{1}{2} \tfrac{\kconst_2}{r^2} 
     \;=\; -\kconst_1 u - \tfrac{1}{2}\kconst_2 u^2
&&
\end{flalign}
 leading to the following projective coordinate Hamiltonian:
\begin{align}
   \mscr{H}  = 
    \tfrac{1}{2} u^2 ( \slang^2 + u^2 p_\ss{u}^2 ) - \kconst_1 u - \tfrac{1}{2} \kconst_2 u^2
    &&,&&
    \wtscr{H} = \tfrac{1}{2} (\slang^2 + u^2 p_\ss{u}^2 ) - \tfrac{\kconst_1}{u} -   \tfrac{1}{2} \kconst_2  +  \tfrac{1}{u^2} p_t
\end{align}
such that the \eq{s}- or \eq{\tau}-parameterized  dynamics for \eq{(u,p_\ss{u})} follow from Eq.\eqref{dqp_s_0} or Eq.\eqref{dqp_TA_SUM} as:
\begin{flalign}
\begin{array}{cc}
     \fnsize{Manev-type }  \\
      \fnsize{dynamics} 
\end{array} 
&&
\begin{array}{llll}
    \pdt{u} =   u^2 p_\ss{u}
\qquad,\qquad 
      \pdt{p}_\ss{u} = - u p_\ss{u}^2 - \tfrac{\kconst_1}{u^2} + \tfrac{2}{u^3}p_t
      \,=\, -\tfrac{1}{u}\big( \slang^2 + 2 u^2 p_\ss{u}^2 \big)   + \tfrac{\kconst_1}{u^2} + \tfrac{\kconst_2}{u} 
\qquad,\qquad 
    \rng{\square} = \tfrac{1}{\slang}\pdt{\square}
\end{array}
&&
\end{flalign}
As a \textit{first}-order system, the above ODEs for the conjugate pair \eq{(u,p_\ss{u})} are still nonlinear — though they are indeed equivalent to a linear \textit{second}-order ODE for \eq{u}.
Linear first-order dynamics are realized if we consider the non-conjugate pair \eq{(u,w)} from Eq.\eqref{uw_eom_gen}. The above is then equivalent to:
\begin{flalign} \label{uw_Fmanev_SUM}
\begin{array}{cc}
     \fnsize{Manev-type }  \\
      \fnsize{dynamics} 
\end{array} 
&&
 w := u^2 p_\ss{u} 
 \quad \Rightarrow \qquad 
\begin{array}{llll}
     \pdt{u}  =  w  
   \\[4pt]
       \pdt{w} = -\omega^2 u +\kconst_1 
\end{array}
\qquad \scrsize{or,} \qquad 
\begin{array}{llll}
    \rng{u}  =  w/\slang 
     \\[4pt] 
     \rng{w} =  -\tfrac{1}{\slang}\omega^2 u + \tfrac{\kconst_1}{\slang} 
\end{array} 
\quad,&&
\begin{array}{lllll}
     \omega^2 := \slang^2 - \kconst_2 
\end{array}
\quad
\end{flalign}
where we have implicitly assumed that \eq{\kconst_2<\slang^2} such that \eq{\omega} is real (for the Kepler case, \eq{\kconst_2=0} and \eq{\omega=\slang}).  
Since \eq{\tup{\slangup}=\tup{\slangup}_\zr} is conserved for any central-force dynamics, the above are linear equations with a constant ``driving force'' term. The solutions to either of the above coincide:
\begin{flalign}
&&
\begin{array}{llll}
      u_{\tau} \,=\, (u_\zr-\tfrac{\kconst_1}{\omega^2})\csn{\varpi \tau}
    \,+\, \tfrac{1}{\omega}w_\zr\sin{\varpi\tau} \,+\, \tfrac{\kconst_1}{\omega^2}
\\[5pt]
    w_{\tau} \,=\, -\slang(u_\zr-\tfrac{\kconst_1}{\omega^2})\snn{\varpi \tau} \,+\, w_\zr\csn{\varpi \tau}
\end{array} 
   &&
\begin{array}{llll}
    \varpi := \omega/\slang
    \\[3pt]
      \varpi \tau = \omega s
\end{array} 
\qquad 
\end{flalign}
and the solutions for the conjugate pair \eq{(u,p_\ss{u}=w/u^2)} are easily recovered from the above. Kepler-type dynamics are also easily recovered with \eq{\kconst_2=0} such that \eq{\omega=\slang} and \eq{\varpi=1}.

\paragraph*{Second-Order Dynamics.}
It is shown that, even \textit{without} making use of \eq{q=1} or \eq{\lambda=\htup{q}\cdot\tup{p}=0}, the first-order dynamics for \eq{(\tup{q},u,\tup{p},p_\ss{u})} in Eq.\eqref{dqp_s_0} or Eq.\eqref{dqp_TA_SUM} are equivalent to the following second-order dynamics for \eq{(\tup{q},u)}:

\begin{flalign}
\;\;\begin{array}{rllllllllll}
   \fnsize{arbitrary $V^0$}:  &\qquad\qquad 
      \pddt{\tup{q}} + \slang^2 \tup{q} = 0  
    &,\quad
         \pddt{u} + \slang^2 u  = -\pderiv{V^0}{u}  
     & \quad\qquad \fnsize{or,} \qquad\qquad
           \rrng{\tup{q}} + \tup{q}   =  0
    &,\quad
         \rrng{u} +  u =  -\tfrac{1}{\slang^2}\pderiv{V^0}{u} 
\end{array} 
&&
\end{flalign}
where the particular form of the central-force potential, \eq{V^0}, is only relevant for the ``\eq{u}-part'' of the dynamics. For the case of a Manev-type \eq{V^0 = -\kconst_1/r- \tfrac{1}{2}\kconst_2/r^2}, or Kepler-type (\eq{\kconst_2=0}), then the above leads to:
\begin{flalign}
\;\;\begin{array}{rllllllllll}
   \fnsize{Manev-type $V^0$}: &\qquad\qquad 
   \pddt{\tup{q}} + \slang^2\tup{q} = 0
    &,\quad
         \pddt{u} + \omega^2 u  = \kconst_1
     & \quad\qquad \fnsize{or,} \qquad\qquad
           \rrng{\tup{q}} + \tup{q}   =  0
    &,\quad
         \rrng{u} + \varpi^2 u =  \kconst_1/\slang^2
\\[8pt]
  \fnsize{Kepler-type $V^0$}:  &\qquad\qquad 
  \pddt{\tup{q}} + \slang^2\tup{q} = 0
    &,\quad
         \pddt{u} + \slang^2 u  = \kconst_1
     & \quad\qquad \fnsize{or,} \qquad\qquad
           \rrng{\tup{q}} + \tup{q}   =  0
    &,\quad
         \rrng{u} + u =  \kconst_1/\slang^2
\end{array} 
&&
\end{flalign}
with \eq{\slang}, \eq{\omega^2=\slang^2-\kconst_2}, and \eq{\varpi=\omega/\slang} integrals of motion.  The above describe a \eq{4}-dim linear harmonic oscillator.

\paragraph*{Kepler Dynamics \& Solutions.}
Kepler-type dynamics are easily recovered from the preceding developments by setting \eq{\kconst_2=0} such that \eq{\omega=\slang} and \eq{\varpi=1}.
Using the transformation in Eq.\eqref{rv_qu_full}, 
the Kepler Hamiltonian in cartesian coordinates, \eq{\mscr{K}}, 
is transformed to the projective coordinate Hamiltonian, \eq{\mscr{H}}:
\begin{align}
    \mscr{K} = \tfrac{1}{2} \v^2 - \tfrac{\kconst_1}{r}
\qquad \Rightarrow \qquad 
    \mscr{H} \,=\, \tfrac{1}{2} u^2 ( \slang^2 + u^2 p_\ss{u}^2) - \kconst_1 u
    \qquad,\qquad 
    \wtscr{H} = \tfrac{1}{2} ( \slang^2 + u^2 p_\ss{u}^2) - \tfrac{\kconst_1}{u}  + \tfrac{1}{u^2}p_t
\end{align}
with \eq{ \slang^2 = q^2 p^2 - (\tup{q}\cdot\tup{p})^2}.
For \eq{s} and \eq{\tau} defined in Eq.\eqref{dtds_0}, 
the Kepler dynamics follow from Eq.\eqref{uw_eom_gen}-Eq.\eqref{uw_eom_gen_intereting} as:
\begin{flalign}
\begin{array}{cc}
     \fnsize{Kepler-type }  \\
      \fnsize{dynamics} 
\end{array} 
&&
\begin{array}{llll}
       \pdt{\tup{q}}
      \,=\,  -\hdge{\tup{\slangup}}\cdot\tup{q} 
      \;\simeq\, \tup{p} 
      &,\quad  \pdt{u} 
      \,=\,  w
\\[5pt]
      \pdt{\tup{p}} 
      \,=\,  -\hdge{\tup{\slangup}}\cdot\tup{p} 
     \;\simeq\,  -\slang^2 \tup{q} 
     &, \quad
     \pdt{w} = -\slang^2 u + \kconst_1 
\end{array}
&&\fnsize{or,}&&
 \begin{array}{llll}
      \rng{\tup{q}}
      \,=\,  -\hdge{\htup{\slangup}}\cdot\tup{q} 
      \;\simeq\,  \tfrac{1}{\slang}\tup{p}  
      &,\quad 
       \rng{u} 
      \,=\, \tfrac{1}{\slang} w
\\[5pt]
     \rng{\tup{p}} 
      \,=\,  -\hdge{\htup{\slangup}}\cdot\tup{p} 
     \;\simeq\,  -\slang \tup{q}  
     &, \quad 
       \rng{w} = -\slang u + \tfrac{\kconst_1}{\slang}
\end{array} 
\quad
\end{flalign}
with \eq{w:=u^2 p_\ss{u}} and where \eq{\tup{\slangup}=\tup{\slangup}_\zr} is an integral of motion.
Either of  the above leads to closed-form solutions in terms of \eq{s} or \eq{\tau=\slang s}: 
\begin{align} \label{kep_sol_SUM}
    \begin{array}{rllll}
         \tup{q}_\tau  \,=\,  \tup{q}_\zr \csn{\tau} -  \hdge{\htup{\slangup}} \cdot \tup{q}_\zr \snn{\tau} 
        & \simeq\,
         \tup{q}_\zr\csn{\tau}  \,+\,  \tfrac{1}{\slang}\tup{p}_\zr\snn{\tau}
\\[6pt]
       \tup{p}_\tau \,=\, \tup{p}_\zr \csn{\tau} - \hdge{\htup{\slangup}} \cdot \tup{p}_\zr \snn{\tau} 
        & \simeq\,
       \tup{p}_\zr\csn{\tau} - \slang\tup{q}_\zr\snn{\tau}
\end{array} 
&&,&&
 \begin{array}{rllll}
       u_{\tau} & =\, (u_\zr-\tfrac{\kconst_1}{\slang^2})\csn{\tau}
    \,+\, \tfrac{1}{\slang}w_\zr\snn{\tau} \,+\, \tfrac{\kconst_1}{\slang^2}
\\[5pt]
     w_{\tau} & =\, -\slang(u_\zr-\tfrac{\kconst_1}{\slang^2})\snn{\tau} \,+\, w_\zr\csn{\tau}
\end{array}
\end{align}
Kepler solutions for cartesian coordinates \eq{(\tup{r},\tup{v})} are recovered from the above using Eq.\eqref{rv_qu_simp} — the result is given explicitly in closed-form in Eq.\eqref{rv_kepsol_nosimp}-Eq.\eqref{rv_kepsol} of this work.  
Closed-form Kepler solutions for \eq{t(\tau)} are given in Eq.\eqref{TA_t_sol}. 
In section \ref{sec:prj_STM_new}, the above solutions are also used to obtain closed-form Kepler state transition matrices. 

%% file: Mysecs_prj/new_0_CTnonmin.tex
\section{Point transformations with redundant coordinates extended to canonical transformations} \label{sec:Pxform}

\begin{notation}
    Here in section \ref{sec:Pxform}, the coordinates \eq{(\tup{x},\tup{\piup})} and \eq{(\bartup{q},\bartup{p})} have no particular significance. Elsewhere, we use the same notation for specific coordinate sets but, here, they are arbitrary. 
\end{notation}

\noindent 
In order to extend projective point transformations to a canonical transformation compatible with Hamiltonian dynamics, we must first address the general problem of constructing a canonical transformation that is induced by a point transformation from minimal to redundant configuration coordinates.
Though the result of our method agrees in the end with with Ferrándiz's treatment of the same problem \cite{ferrandiz1994extended}, we take a different approach (based on Hamilton's principle) that we feel is more transparent and easily understood from an analytical dynamics perspective. We also expand upon Ferrándiz's treatment by allowing for time-dependent transformations and constraints.


 Consider a point transformation (possibly time-dependent) from some minimal set of \eq{\en} generalized configuration coordinates, \eq{\tup{x}=(x_1,\dots,x_\en)}, to some new set of \eq{\en+\emm} \textit{redundant} coordinates, \eq{\bartup{q}}:
\begin{small}
\begin{flalign} \label{Qq_nonmin}
&&
    \tup{x} \,=\, \tup{x}(\bartup{q},t)
    &&
    \fnsize{with:} \;\; 
         \tup{x}\in\mbb{R}^\en 
         \;\;,\;\; \bartup{q}\in\mbb{R}^{\en+\emm} 
    \quad
\end{flalign}
\end{small}
Since \eq{\bartup{q}} are over-parameterized by \eq{\emm} degrees, there exists some \eq{\emm} independent constraint functions, \eq{\tup{\varphiup}(\bartup{q},t)\in\mbb{R}^\emm} — which we shall take to be holonomic and rheonomic — of the new coordinates such that
\begin{small}
\begin{align}  \label{const_1}
    \tup{\varphiup}(\bartup{q},t) \,=\, 0 
    \qquad,\qquad
     \mrm{d} \tup{\varphiup} \,=\, \pderiv{\tup{\varphiup}}{\bartup{q}} \cdot \mrm{d} \bartup{q} 
    \,+\, \pderiv{\tup{\varphiup}}{t}\mrm{d} t \,=\, 0 
\end{align}
\end{small}
We wish to find the canonical transformation induced by the point transformation of Eq.\eqref{Qq_nonmin}, subject to the \eq{\emm} constraints \eq{\tup{\varphiup}(\bartup{q},t) = 0 }. That is, the general problem we wish to solve is as follows:

\begin{notesq} 
    \sloppy \textsl{Suppose we are given the Hamiltonian \eq{\mscr{K}(\tup{x},\tup{\piup},t)} for known canonical coordinates \eq{(\tup{x},\tup{\piup})\in\mbb{R}^{2\en}}. If we then specify a point transformation, \eq{\tup{x}=\tup{x}(\bartup{q},t)} — for redundant coordinates \eq{\bartup{q}\in\mbb{R}^{\en+\emm}} subject to \eq{\emm} constraints \eq{\tup{\varphiup}(\bartup{q},t)\in\mbb{R}^\emm} — what then is the momenta transformation, \eq{\tup{\piup}=\tup{\piup}(\bartup{q},\bartup{p},t)}, and the new Hamiltonian, \eq{\mscr{H}(\bartup{q},\bartup{p},t)}, such that the redundant coordinate set, \eq{(\bartup{q},\bartup{p})\in\mbb{R}^{2(\en+\emm)}}, is canonical (in the sense they obey Hamilton's equations for \eq{\mscr{H}})? }
\end{notesq}

\noindent There are several approaches one may take to answer this question.  We will start with Hamilton's Principle which says that, out of all possible paths between two points, the actual path taken is the one for which the action is stationary:\footnote{In the case that nonconservative forces are present,  the work from these forces would also be added to the integrand. See Appx.~\ref{app:non_conserv}. }
\begin{small}
\begin{align}
\begin{array}{lll}
   \delta I \,=\, \delta \int_{t_{\zr}}^{t_f}\mscr{L}_{\sscr{K}}\,\mrm{d} t \,=\, \delta \int_{t_{\zr}}^{t_f}(\tup{\piup}\cdot\dottup{x}-\mscr{K})\,\mrm{d} t \,=\, 0
\end{array}
\end{align}
\end{small}
For any canonical transformation, the relation between the original and new phase spaces can be found from the condition that Hamilton's principle in the above form
must hold in both phase spaces:
\begin{small}
\begin{align} \label{HamP_nonMin}
\begin{array}{lll}
   \delta I\,=\, \delta\int_{t_\zr}^{t_f} (\tup{\piup}\cdot\dottup{x} - \mscr{K})\, \mrm{d} t \;\,=\,\; 0 \;\,=\,\;
     \delta\int_{t_\zr}^{t_f} (\bartup{p}\cdot\dot{\bartup{q}} 
     -\mscr{H})\, \mrm{d} t
\end{array}
\end{align}
\end{small}
The above will certainly hold if the integrands are equal:\footnote{Note that this is \textit{not} the most general relation between the two phase spaces for a canonical transformation. Hamilton's principle still holds if the time-derivative of some arbitrary generating function is added to either side and/or if either side is multiplied by some scalar constant \cite{lanczos2012variational}. Such an approach leads to a more general family of conical transformations. However, we are only concerned with the canonical point transformation following from some specified \eq{\tup{x}(\bartup{q},t)}. For this, the generating function approach is unnecessary.} 
\begin{small}
\begin{align} \label{diff_rel_0}
    \tup{\piup}\cdot \mrm{d} \tup{x} \,-\, \mscr{K} \mrm{d} t \,=\, \bartup{p}\cdot \mrm{d} \bartup{q} 
    \,-\, \mscr{H} \mrm{d} t 
\end{align}
\end{small}
However, because the coordinates are redundant, the \eq{\mrm{d} \bartup{q}} are  not independent but, rather, are related through the \eq{\emm} constraints in Eq.\eqref{const_1}. We can include these \eq{\emm} constraints in the above expression by adding an additional term, \eq{\tup{\lambdaup}\cdot\mrm{d} \tup{\varphiup}} — where \eq{\tup{\lambdaup}\in\mbb{R}^\emm} are Lagrange multipliers — to the right-hand-side of the above, leading to:
\begin{small}
\begin{align} \label{diff_rel1}
    &\tup{\piup}\cdot \mrm{d} \tup{x} \,-\, \mscr{K} \mrm{d} t \,=\, \bartup{p}\cdot \mrm{d} \bartup{q} 
    \,-\, \mscr{H} \mrm{d} t \,-\, \tup{\lambdaup}\cdot \mrm{d} \tup{\varphiup} 
\end{align}
\end{small}
Since \eq{\mrm{d} \tup{\varphiup}=0}, the above relation still ensures that Hamilton's principle, as given by equation  Eq.\eqref{HamP_nonMin}, holds in both phase spaces. Then, using the relations \eq{\tup{x}=\tup{x}(\bartup{q},t)} and \eq{\tup{\varphiup}=\tup{\varphiup}(\bartup{q},t)} — which are presumed to be specified — to express \eq{\mrm{d} \tup{x}} and \eq{\mrm{d} \tup{\varphiup}} in terms of \eq{\mrm{d} \bartup{q}} and \eq{ \mrm{d} t}, the above leads to
\begin{small}
\begin{align} \label{turky}
\begin{array}{ll}
    &\quad \tup{\piup}\cdot (\pderiv{\tup{x}}{\bartup{q}} \cdot \mrm{d} \bartup{q}  \,+\, \pderiv{\tup{x}}{t} \mrm{d} t) \,-\, \mscr{K} \mrm{d} t \;\,=\,\; 
    \bartup{p}\cdot \mrm{d} \bartup{q} \,-\, \mscr{H} \mrm{d} t \,-\, \tup{\lambdaup}\cdot(\pderiv{\tup{\varphiup}}{\bartup{q}} \cdot \mrm{d} \bartup{q} \,+\,  \pderiv{\tup{\varphiup}}{t} \mrm{d} t)
\\[8pt]
   \Rightarrow  &\quad (\bartup{p} - \trn{\pderiv{\tup{x}}{\bartup{q}}} \cdot \tup{\piup} \,-\, \trn{\pderiv{\tup{\varphiup}}{\bartup{q}}} \cdot \tup{\lambdaup}  )\cdot \mrm{d} \bartup{q} 
    \,-\, (\mscr{H} -\mscr{K} \,+\, \tup{\lambdaup}\cdot\pderiv{\tup{\varphiup}}{t} \,+\, \tup{\piup}\cdot\pderiv{\tup{x}}{t}) \mrm{d} t \,=\, 0
\end{array}
\end{align}
\end{small}
With the Lagrange multiplier allowing the differentials to be treated as independent, the above condition requires that the coefficients of \eq{\mrm{d} \bartup{q}} and \eq{\mrm{d} t} must vanish. 
Setting the coefficient of \eq{\mrm{d} \bartup{q}} to zero defines the new momenta, \eq{\bartup{p}}:
\begin{small}
\begin{align} \label{qpH_new}
    &\bartup{p} \,=\, \trn{\pderiv{\tup{x}}{\bartup{q}}} \cdot \tup{\piup} \,+\, \trn{\pderiv{\tup{\varphiup}}{\bartup{q}}} \cdot \tup{\lambdaup} 
\end{align}
\end{small}
This may be combined into a matrix equation for the momenta transformation:
\begin{small}
\begin{align} \label{Pxform_finaly}
    \bartup{p}
    \,=\,
    \trn{B}\cdot
    \begin{pmatrix}
     \tup{\piup} \\ \tup{\lambdaup} 
    \end{pmatrix} 
\qquad \leftrightarrow \qquad
   \begin{pmatrix}
     \tup{\piup} \\ \tup{\lambdaup}
    \end{pmatrix} 
   \,=\, \invtrn{B}\cdot\bartup{p}
&&,&&
     B(\bartup{q},t) :=\, 
     \begin{pmatrix}
         \pderiv{\tup{x}}{\bartup{q}}
         \\[5pt] \pderiv{\tup{\varphiup}}{\bartup{q}}
    \end{pmatrix} \in\mbb{R}^{(\en+\emm)\times (\en+\emm)}
\end{align}
\end{small}
Where \eq{\inv{B}} exists provided that the constraint functions are independent.  
The above gives \eq{\tup{\piup}} and \eq{\tup{\lambdaup}} as functions of \eq{(\bartup{q},\bartup{p},t)} as desired, but gives \eq{\bartup{p}} as a function of \eq{(\bartup{q},\tup{\piup},\tup{\lambdaup},t)}. In order to obtain the desired expression for \eq{\bartup{p}(\tup{x},\tup{\piup},\tup{\lambdaup},t)}, 
Eq.\eqref{Qq_nonmin} and Eq.\eqref{const_1} must  be used together in order to re-express \eq{B} in terms of \eq{\tup{x}} rather than \eq{\bartup{q}}.

Lastly, setting the coefficient of \eq{ \mrm{d} t} equal to zero in Eq.\eqref{turky} yields the relation between the  Hamiltonian for the redundant coordinates, \eq{\mscr{H}(\bartup{q},\bartup{p},t)}, and the Hamiltonian for the original coordinates, \eq{\mscr{K}(\tup{x},\tup{\piup},t)}:
\begin{small}
\begin{align} \label{Ham_xform}
        \mscr{H} \,=\, \mscr{K} \,-\, \tup{\lambdaup}\cdot\pderiv{\tup{\varphiup}}{t} \,-\, \tup{\piup}\cdot\pderiv{\tup{x}}{t}
\end{align}
\end{small}
with the right-hand-side expressed in terms of \eq{(\bartup{q},\bartup{p})} using the specified \eq{\tup{x}(\bartup{q},t)} from  Eq.\eqref{Qq_nonmin} and the induced \eq{\tup{\piup}(\bartup{q},\bartup{p},t)} and  \eq{\tup{\lambdaup}(\bartup{q},\bartup{p},t)} from Eq.\eqref{Pxform_finaly}.

We have now answered the problem we set out to solve in this section: given the Hamiltonian \eq{\mscr{K}(\tup{x},\tup{\piup},t)} for known minimal canonical coordinates \eq{(\tup{x},\tup{\piup})\in\mbb{R}^{2\en}}, along with a specified point transformation \eq{\tup{x}(\bartup{q},t)} subject to \eq{\emm} constraints \eq{\tup{\varphiup}(\bartup{q},t)=0}, then the momenta transformation, \eq{\tup{\piup}(\bartup{q},\bartup{p},t)}, is given by Eq.\eqref{Pxform_finaly} and the new Hamiltonian, \eq{\mscr{H}(\bartup{q},\bartup{p},t)}, is given by Eq.\eqref{Ham_xform}. 


\begin{small}
\begin{itemize} 
    \item  \textit{Dimensions and coordinate partitioning.} 
    The redundant phase space coordinates are over-parameterized by some integer number \eq{\emm\geq 1} such that \eq{\bartup{q}\in\mbb{R}^{\en+\emm}} and \eq{\bartup{p}\in\mbb{R}^{\en+\emm}}, where \eq{\en} is the degrees of freedom of the system (minimal number of coordinates needed).
    It can be helpful to partition \eq{\bartup{q}} and \eq{\bartup{p}} into \eq{\en}-dimensional and \eq{\emm}-dimensional parts as: 
    \begin{small}
    \begin{align} \label{qu_part}
         \bartup{q} = (\tup{q},\tup{u})  \in\mbb{R}^{\en+\emm}
    \qquad, \qquad 
        \bartup{p} = (\tup{p},\tup{p}_u) \in\mbb{R}^{\en+\emm}  
    \end{align}
    \end{small}
    where \eq{\tup{q}} and \eq{\tup{p}} each contain the minimal number of \eq{\en} coordinates, while  \eq{\tup{u}} and \eq{\tup{p}_u} each contain \eq{\emm} ``extra''  coordinates.\footnote{In this formulation, though we have defined \eq{\tup{q}} to have the minimal number of coordinates, these \eq{\en} coordinates generally do \textit{not} constitute a set of minimal coordinates themselves. I.e., \eq{\tup{q}\in\mbb{R}^\en} alone do not fully define the configuration; the full, redundant set \eq{\bartup{q}=(\tup{q},\tup{u})\in\mbb{R}^{\en+\emm}} is still needed. }
    The dimensions of various tuples appearing previously are then as follows:
    \begin{small}
    \begin{align}
       \left. \begin{array}{lcl}
             \tup{x} , 
             \tup{\piup},  
             \tup{q} ,
             \tup{p} \in \mbb{R}^\en
        \end{array}\;\right.
       \qquad, \qquad
       \left. \begin{array}{lcl}
              \tup{\varphiup} , 
             \tup{\lambdaup}, 
              \tup{u} ,
             \tup{p}_u \in \mbb{R}^\emm
        \end{array}\;\right. 
    \end{align}
    \end{small}
    and Eq.\eqref{Pxform_finaly} is written as
    \begin{small}
    \begin{align} \label{Pxform_div}
            \begin{pmatrix}
         \tup{p} \\ \tup{p}_u 
        \end{pmatrix}
        \,=\,
        \trn{B}
        \begin{pmatrix}
         \tup{\piup} \\ \tup{\lambdaup} 
        \end{pmatrix} 
        \qquad\leftrightarrow\qquad
          \begin{pmatrix}
         \tup{\piup} \\ \tup{\lambdaup}
        \end{pmatrix}
        \,=\,
        \invtrn{B}
        \begin{pmatrix}
         \tup{p} \\ \tup{p}_u 
        \end{pmatrix} 
    \qquad,\qquad
     B :=
     \begin{pmatrix}
         \pderiv{\tup{x}}{\bartup{q}}
         \\[5pt] \pderiv{\tup{\varphiup}}{\bartup{q}}
    \end{pmatrix}
    \,=\, 
        \begin{pmatrix}
             \big(\pderiv{\tup{x}}{\tup{q}}\big)_{\en \times \en} & \left(\pderiv{\tup{x}}{\tup{u}}\right)_{\en \times \emm}
             \\[5pt] 
             \left(\pderiv{\tup{\varphiup}}{\tup{q}}\right)_{\emm \times \en} 
             & \left(\pderiv{\tup{\varphiup}}{\tup{u}}\right)_{\emm \times \emm}
        \end{pmatrix} 
    \end{align}
    \end{small}
    \item  \textit{Generating function.}
    Though this work makes no use of generating functions, readers familiar with this approach to canonical transformations may care to note that the above developments may also be obtained using the following "type-2" generating function, \eq{G_2(\bartup{q},\tup{\piup},t)}, satisfying the usual relations:
    \begin{small}
    \begin{align}
        G_2(\bartup{q},\tup{\piup},t) \,=\, \tup{x}(\bartup{q},t)\cdot\tup{\piup} \,+\, \tup{\varphiup}(\bartup{q},t)\cdot \tup{\lambdaup}
    \qquad,\qquad 
        \tup{p} = \pderiv{G_2}{\tup{q}}
    \qquad,\qquad 
        \tup{x} =\pderiv{G_2}{\tup{\piup}}
    \end{align}
    \end{small}
\end{itemize}
\end{small}

%% file: Mysecs_prj/new_2_prjFrad.tex
\section{Central-force Hamiltonian dynamics in projective coordinates} \label{sec:central force}


We now turn our attention to the classic Hamiltonian dynamics of a particle moving in a central-force field (along with other arbitrary perturbations). In particular, how such a Hamiltonian system is transformed when represented in certain projective coordinates (defined below in section \ref{sec:sum_m=-1}). 
Then, in section \ref{sec:ext_cen}, we will use the extended phase space to introduce a transformation of the evolution parameter to something other than time. This will then be used in section \ref{sec:2BP_lin} to obtain equations of motion that are fully linear for central-forces arising from any potential of the form in Eq.\eqref{MANEV_cartesian} (i.e., for Manev-type or Kepler-type forces). 





\paragraph*{Original Cartesian Coordinate Formulation.}
Let \eq{(\tup{r},\tup{v})\in\mbb{R}^6} denote cartesian position and velocity coordinates for an orthonormal inertial frame with origin at \eq{r:=\mag{\tup{r}}=0}. 
 We consider a system of an unconstrained particle moving in Euclidean 3-space.
subject to conservative forces modeled by some potential function \eq{V(\tup{r},t)}, as well as arbitrary nonconservative forces with cartesian coordinates \eq{\tup{a}^\nc\in\mbb{R}^3}.\footnote{E.g., \eq{\tup{a}^\nc} might be the cartesian coordinate vector for thrust or drag forces.}  
See  Appx.~\ref{app:non_conserv} for details on Hamiltonian dynamics with nonconservative forces.
We assume that the potential decomposes as  \eq{V = V^{0}(r) + V^{1}(\tup{r},t)}, where \eq{V^{0}(r)} accounts for all conservative central-forces and where \eq{V^{1}(\tup{r},t)} is arbitrary, accounting for all other conservative perturbations\footnote{ E.g., in the context of orbital dynamics, \eq{V^1} could be the \eq{J_2} potential or higher order terms, a third body potential, etc.}.
The Hamiltonian dynamics for inertial cartesian coordinates \eq{(\tup{r},\tup{v})} are then:\footnote{As before, \eq{\tup{r}} and \eq{\tup{v}} are cartesian position and velocity coordinates for an orthonormal inertial frame, with origin at \eq{r=0}.}
\begin{small}
\begin{align} \label{K_cen}
    \mscr{K} \,=\, \tfrac{1}{2}\v^2 \,+\, V^{0}(r) \,+\,V^{1}(\tup{r},t)
      &&
\begin{array}{ll}
      \dot{\tup{r}} \,=\, \pderiv{\mscr{K}}{\tup{v}} \,=\, \tup{v} 
\\[5pt] 
    \dot{\tup{v}} \,=\, -\pderiv{\mscr{K}}{\tup{r}} \,+\, \tup{a}^\nc \,=\,
    -\pderiv{V^0}{r}\htup{r}  + \tup{F}  
\end{array}
&&
\tup{F} := -\pderiv{V^1}{\tup{r}} + \tup{a}^\nc
\end{align}
\end{small}
where \eq{\tup{F}\in\mbb{R}^3} denotes the cartesian components of the total perturbing forces from both conservative sources (\eq{V^1}) and nonconservative sources (\eq{\tup{a}^\nc}). 
The particular form of \eq{V^{0}(r)} will not matter until later on. When it does matter, we consider so-called \textit{Manev}-type potentials (which include Kepler-type for the special case \eq{\kconst_2=0}):
\begin{small}
\begin{align} \label{MANEV_cartesian}
      V^{0} = -\tfrac{\kconst_1}{r} - \tfrac{1}{2}\tfrac{\kconst_2}{r^2}
      \qquad,\qquad 
      -\pderiv{V^0}{\tup{r}} = -\tfrac{\kconst_1}{r^2}\htup{r} - \tfrac{\kconst_2}{r^3}\htup{r}
      \qquad,\qquad 
      \kconst_1,\kconst_2\in\mbb{R}
\end{align}
\end{small}

\subsection{A modified BF-like projective transformation} \label{sec:sum_m=-1}

We now use the methods of section \ref{sec:Pxform} to construct the canonical extension of a projective point transformation given by \eq{\tup{r} = u^n \htup{q}}, and transform the above Hamiltonian system to a new redundant-coordinate representation. 
Many developments in this section follow as special cases of those in Appendix \ref{sec:2bp_gen} where further details and derivations may be found. In particular, Appendix  \ref{sec:2bp_gen} develops the canonical extension of a family of projective point transformations of the more general form \eq{\tup{r}=u^n q^m \tup{q}} for any \eq{n,m\in\mbb{R}}. Here, we consider only the \eq{m=-1} case (and, eventually, \eq{n=-1}) and remark upon some properties in comparison to any \eq{m\neq 1} case such as the BF transformation.   

Starting with the inertial cartesian coordinate Hamiltonian dynamics in Eq.\eqref{K_cen}, we specify a point transformation from \eq{\tup{r}\in\mbb{R}^3} to redundant coordinates \eq{\bartup{q}=(\tup{q},u)\in\mbb{R}^4} given by \eq{\tup{r} = \tfrac{u^n}{q}\tup{q}= u^n \htup{q}}, subject to the constraint \eq{\mag{\tup{q}}=1}.  
Using the developments of section \ref{sec:Pxform} we obtain the associated momenta coordinate transformation \eq{\tup{v}(\bartup{q},\bartup{p})},
as in Eq.\eqref{Pxform_finaly} (see Appendix \ref{sec:2bp_gen} for more explicit details).
This leads to the full ``canonically-extended'' projective transformation for some yet-to-be-chosen \eq{0\neq n\in\mbb{R}}:\footnote{where \eq{\htup{q}=\tfrac{1}{q}\tup{q}},
 and \eq{\imat_3 - \htup{q}\otms\htup{q} = -\hdge{\htup{q}}\cdot\hdge{\htup{q}}} with \eq{\hdge{q}_{ij} = \lc_{ijk}q_k}. }
\begin{small}
\begin{align} \label{sum_rv_m} 
    (\bartup{q},\bartup{p}) \mapsto (\tup{r},\tup{v})
\; \left\{ \;\;
\boxed{ \begin{array}{lllll}
      \tup{r} \,=\, u^n \htup{q}
  \\[5pt]
     \tup{v} \,=\,  \tfrac{q}{u^n} (\imat_3 - \htup{q}\otms\htup{q})\cdot\tup{p}
         \,+\, \tfrac{1}{n} u^{1-n} p_\ss{u} \htup{q}
\end{array} } \right. 
&&
\begin{array}{llll}
   \fnsize{constraint:}  & \varphi = q-1 = 0 
\\[5pt]
     & \lambda  \,=\, \htup{q}\cdot\tup{p} 
\end{array}
\end{align}
\end{small}
where the constraint \eq{\varphi=0} (equivalent to \eq{q=1}) is built into the derivation of the momenta transformation, which includes the Lagrange multiplier \eq{\lambda=\htup{q}\cdot\tup{p}} associated with \eq{\varphi}.
As mentioned in Appendix \ref{sec:ang_momentum}, a notable property of the above is that
the specific angular momentum, \eq{\tup{\slangup}}, is given solely in terms of \eq{(\tup{q},\tup{p})} in exactly the same manner as in terms of \eq{(\tup{r},\tup{v})}:
\begin{small}
\begin{align} \label{sum_ell_m}
    \hdge{\tup{\slangup}} 
    \,=\, \tup{r}\wdg \tup{v} \,=\, \tup{q}\wdg\tup{p}
&&,&&
    \tup{\slangup} \,=\, \hdge{\tup{v}}\cdot\tup{r} 
    \,=\, \hdge{\tup{p}}\cdot\tup{q} 
&&,&&
    \slang^2  \,=\, r^2 \v^2 - (\tup{r}\cdot\tup{v})^2 \,=\, q^2 p^2 - (\tup{q}\cdot\tup{p})^2
\end{align}
\end{small}
Since the transformation is time-independent, the new Hamiltonian \eq{\mscr{H}} is found simply by substitution of  Eq.\eqref{sum_rv_m} into the cartesian coordinate  Hamiltonian \eq{\mscr{K}} from Eq.\eqref{K_cen}. 
This leads to the following (with \eq{\slang^2} as above):\footnote{The terms \eq{V^{0}(u)} and \eq{V^{1}(\bartup{q},t)} appearing in \eq{\mscr{H}} are an abuse of notation indicating the original \eq{V^{0}(r)} and \eq{V^{1}(\tup{r},t)} rewritten in terms of \eq{\bartup{q}} using \eq{\tup{r}=u^n\htup{q}}.}
\begin{align} \label{sum_H_m}
    \mscr{K} = \tfrac{1}{2}\v^2 + V^{0}(r) + V^{1}(\tup{r},t) 
    \qquad \Rightarrow \qquad
    \boxed{ \mscr{H} \,=\, 
    \tfrac{1}{u^{2n}} \tfrac{1}{2} \big( \slang^2 + \tfrac{1}{n^2}u^2 p_\ss{u}^2 \big) \,+\, V^{0}(u) \,+\, V^{1}(\bartup{q},t) }
\end{align}

\begin{notesq} 
        Note from Eq.\eqref{sum_rv_m} that \eq{r=u^n}  is independent of \eq{\tup{q}}; a central-force potential \eq{V^{0}(r)} is taken to some \eq{\tup{q}}-independent \eq{V^{0}(u)}:
    \begin{small}
    \begin{align} \label{f(r)_m}
      \tup{r}=u^n\htup{q} \qquad \Rightarrow \qquad
        r=u^n \;\;, \quad 
         \pderiv{}{\tup{q}} V^{0} =  0 
        \quad \forall \, V^{0}(r) 
    \end{align}
     \end{small}
    As a consequence, all conservative central-forces will drop out of the dynamics for \eq{(\tup{q},\tup{p})} and appear only in those for \eq{(u,p_\ss{u})}.
    The same turns out to be true for \textit{non}conservative central-forces as well (cf.~Eq.\eqref{totalForce}).\footnote{This property does not hold for \eq{m\neq -1} in \eq{\tup{r}=u^n q^m \tup{q}}, leading to \eq{r=u^n q^{m+1}} and \eq{\pderiv{}{\tup{q}}V^{0} \neq 0 } for any \eq{m\neq -1},. }.
\end{notesq}

\noindent Now, Eq.\eqref{sum_rv_m} has no unique global inverse. However, as shown in Appendix \ref{sec:lambda_const}, \eq{q=\mag{\tup{q}}} and \eq{\lambda =\htup{q}\cdot\tup{p}}  are  integrals of motion for the new Hamiltonian \eq{\mscr{H}} — this holds even in the presence of arbitrary nonconservative forces.  That is, our eight coordinates always automatically satisfy two ``extra'' integrals of motion:
\begin{small}
\begin{flalign} \label{qp=0}
     \begin{array}{cc}
     \fnsize{integrals}  \\
     \fnsize{of motion} 
\end{array}\!\!:
\quad
\begin{array}{ll}
    q = \mag{\tup{q}} 
     \\[4pt]
    \lambda =  \htup{q}\cdot\tup{p}  
\end{array} 
&&
\begin{array}{ll}
    \dot{q} = \pbrak{q}{\mscr{H}} + \pderiv{q}{\bartup{p}} \cdot \bartup{\alphaup}  = 0 
   \\[4pt]
   \dot{\lambda}  = \pbrak{\lambda}{\mscr{H}} + \pderiv{\lambda}{\bartup{p}} \cdot \bartup{\alphaup}  = 0
\end{array} 
\qquad  \xRightarrow{\text{choose}} \qquad
\begin{array}{llll}
     q = q_\zr  = 1
     \\[4pt]
       \lambda = \lambda_\zr  = 0 \;=\htup{q}\cdot\tup{p}
\end{array}
 && 
\end{flalign}
\end{small}
where \eq{\bartup{\alphaup}=(\tup{\alphaup},\alpha_\ss{u})} account for any non-conservative perturbations (discussed soon). As per Appendix \ref{sec:2bp_gen}, we are free to limit consideration of the above integrals to any values \eq{q_\zr>0,\lambda_0\in\mbb{R}}. This then allows us to find the inverse (rather \textit{an} inverse) of the projective 
transformation.\footnote{See Eq.\eqref{prj_nmGen}-Eq.\eqref{qu_k1}.} 
In particular, by restricting consideration to \eq{q=1} and \eq{\lambda=0}, 
the inverse transformation (for these chosen values) of Eq.\eqref{sum_rv_m} is then obtained as:
\begin{small}
\begin{align}  \label{sum_qp_m}
\begin{array}{cc}
     \fnsize{choose}  \\
     \fnsize{$q=1$ } \\
     \fnsize{$\lambda=0$}
\end{array} 
\qquad \Rightarrow \qquad
(\tup{r},\tup{v}) \mapsto (\bartup{q},\bartup{p})
\quad \left\{ \quad 
\boxed{ \begin{array}{llll} 
      \tup{q} \,=\,  \htup{r}  
      &,\quad 
       \tup{p} \,=\,  r(\imat_3-\htup{r}\otms\htup{r})\cdot\tup{v}
     \,=\, -\hdge{\tup{\slangup}}\cdot\htup{r} 
 \\[4pt]
      u \,=\, r^\ss{1/n} 
      &,\quad
        p_\ss{u} \,=\,
    n r^\ss{\frac{1}{n}(n-1)} \htup{r}\cdot\tup{v} 
\end{array} }  \right.
\begin{array}{lllll}
      &=\;
     r\dot{\tup{r}} - \dot{r}\tup{r}  \,=\, r^2\dot{\htup{r}}
\\[4pt]
      &=\;   n r^\ss{\frac{1}{n}(n-1)}\,\dot{r}
\end{array}
\end{align}
\end{small}
which also gives the meaning of new canonical coordinates in terms of cartesian coordinates \eq{(\tup{r},\tup{v})}. 
Note that the common choices of \eq{n \pm 1} in the initial point transformation lead to the following \eq{(u,p_\ss{u})}:
\begin{small}
\begin{align}
n=1\;\left\{\quad
    \begin{array}{ll}
         u\,=\, r  \\[4pt]
         p_\ss{u} \,=\, \htup{r}\cdot\tup{v} \,=\, \dot{r} 
    \end{array}\right.
    &&
    n=-1\;\left\{\quad
     \begin{array}{ll}
         u \,=\, 1/r  \\[4pt]
         p_\ss{u} \,=\, -r^2\htup{r}\cdot\tup{v} \,=\, -r^2\dot{r}
    \end{array}\right.
\end{align}
\end{small}
Though we eventually choose \eq{n=-1}, we continue to assume arbitrary \eq{n\neq0} for now.
Before continuing to the dynamics, let us remark on some other features of the above transformation:

\begin{remrm} \label{rem:qp=0}
   One does not need to do anything to enforce the ``constraints'' \eq{q=1} and \eq{\lambda=0} other than to transform initial conditions \eq{(\tup{r}_\zr,\tup{v}_\zr) \mapsto (\bartup{q}_\zr,\bartup{p}_\zr)} using Eq.\eqref{sum_qp_m}. This places no restrictions on \eq{(\tup{r}_\zr,\tup{v}_\zr)} and it guarantees that  \eq{q_\zr=1} and \eq{\lambda_\zr=0}. Since \eq{q} and \eq{\lambda} are integrals of motion, it also guarantees that  \eq{q=1} and \eq{\lambda=0} hold for all time along any solution curve.
    We often use   ``\eq{\simeq}'' to denote relations which have been simplified using \eq{q=1} and \eq{\lambda=\htup{q}\cdot\tup{p}=0}. 
    To reiterate, such relations hold along any
    solution curve \eq{(\bartup{q}_t,\bartup{p}_t)} with initial conditions \eq{(\bartup{q}_\zr,\bartup{p}_\zr)} satisfying \eq{q_\zr=1} and \eq{\lambda_\zr =0}. Equivalently,  they hold along any solution curve with initial conditions \eq{(\bartup{q}_\zr,\bartup{p}_\zr)} satisfying Eq.\eqref{sum_qp_m} for some given \eq{(\tup{r}_\zr,\tup{v}_\zr)}.
\end{remrm}

\begin{remrm}[\textit{Relation to LVLH basis}]  \label{rem:LVLH}
    By choosing \eq{m=-1} (i.e., using \eq{\tup{r}=u^n\htup{q}}), we have that \eq{\lambda=\htup{q}\cdot\tup{p}}
    such that \eq{\lambda=0} directly implies \eq{\tup{q}\cdot\tup{p}=0} and, thus, \eq{(\tup{q},\tup{p},\tup{\slangup})} are mutually orthogonal. Along with \eq{q=1}, this means  \eq{(\tup{q},\tup{p})} have the following convenient properties:
    \begin{small}
    \begin{flalign} \label{qpl_rels}
    \quad\;\;
    \begin{array}{lll}
          \tup{\slangup}=\tup{q}\tms \tup{p}
        \\[4pt]
          \tup{q}\cdot\tup{\slangup} = \tup{p}\cdot\tup{\slangup}  = 0
       \\[4pt]
        ` \slang^2 = q^2 p^2-(\tup{q}\cdot\tup{p})^2
    \end{array}
    &&,&&
    \begin{array}{cc}
        \fnsize{using:} \\
         q=1  \\
         \tup{q}\cdot\tup{p}=0 
    \end{array}
     \Rightarrow \;\; 
    \left\{ \;\; 
    \begin{array}{cccc}
        \begin{array}{ll}
            \tup{q} \simeq \htup{q} \,\simeq\, \htup{p}\tms \htup{\slangup} 
            &,\qquad 
            \hdge{\tup{\slangup}}\cdot\tup{q} 
            \,\simeq\, -q^2\tup{p} \,\simeq\, -\tup{p}
          \\[4pt]
            \tup{p} \,\simeq\, \tup{\slangup}\tms \htup{q} 
             &,\qquad 
              \hdge{\tup{\slangup}}\cdot\tup{p}
              \,\simeq\, p^2\tup{q} \,\simeq\, \slang^2 \tup{q} 
        \\[4pt]
              \slang^2 \,\simeq\, p^2  
           &,\qquad 
            \imat_3 \,\simeq\, \htup{q} \otms\htup{q} + \htup{p} \otms\htup{p} + \htup{\slangup} \otms\htup{\slangup}
        \end{array} 
    \end{array}
    \right.
     \qquad
    \end{flalign}
    \end{small}
    In particular, \eq{(\tup{q},\tup{p})} define not just the angular momentum, but the inertial cartesian components of the particle's LVLH basis:
    \begin{small}
    \begin{align} \label{qpl_rels2}
     \begin{array}{llll}
             \begin{array}{cc}
                 \fnsize{orthonormal}  \\
                 \fnsize{LVLH basis}
                \end{array}
                \!=\, \{\htup{t}_r,\htup{t}_\tau, \htup{t}_\slang \}
                =  \{\htup{q},-\hdge{\htup{\slangup}}\cdot\htup{q},\htup{\slangup}\} \simeq \{\htup{q},\htup{p},\htup{\slangup}\}
            \qquad\qquad 
            \fnsize{with:} \;\; q^2 \simeq 1\;,\;\; \slang^2 \simeq p^2
    \end{array}
\end{align}
\end{small}
\end{remrm}


\subsection{Hamilton's equations for perturbed central-force dynamics}

Although the primary advantage (linear orbital dynamics) of the projective transformation described above is not realized until it is combined with a transformation of the evolution parameter (section \ref{sec:ext_cen}), we note here the time-parametrized dynamics and generalized forces. 
The transformed Hamiltonian is given by Eq.\eqref{sum_H_m}, leading to the following \eq{t}-parametrized 
canonical equations of motion:\footnote{Note that these ODEs are simplified considerably compared to what they would be for any \eq{m\neq-1}, as given by Eq.\eqref{sum_qpdot} or Eq.\eqref{sum_qpdot_k1}.}
 \begin{small}\begin{gather} \label{eoms_fterm} 
  \mscr{H} \,=\, 
    \tfrac{1}{u^{2n}} \tfrac{1}{2}\big( \slang^2 + \tfrac{1}{n^2}u^2 p_\ss{u}^2 \big) \,+\, V^{0}(u) \,+\, V^{1}
    \qquad\qquad \qquad\qquad 
    \phantom{F = VQ +A}
\\ \nonumber
\begin{array}{llll}
     \dot{\tup{q}} \,=\, \pderiv{\mscr{H}}{\tup{p}} 
     \,=\,-\tfrac{1}{u^{2n}} \hdge{\tup{\slangup}}\cdot \tup{q} 
\\[6pt] 
    \dot{u} \,=\, \pderiv{\mscr{H}}{p_u} \,=\,\tfrac{1}{u^{2n}}\tfrac{1}{n^2}u^2 p_\ss{u} 
\end{array}
 \quad,\quad
\begin{array}{llll}
    \dot{\tup{p}} \,=\, -\pderiv{\mscr{H}}{\tup{q}} \,+\, \tup{\alphaup} 
     \,=\, -\tfrac{1}{u^{2n}} \hdge{\tup{\slangup}}\cdot \tup{p} \,+\, \tup{f} 
\\[6pt]
    \dot{p}_\ss{u} \,=\, -\pderiv{\mscr{H}}{u} + \alpha_\ss{u} \,=\,
     \tfrac{n}{u^{2n+1}}\big( \slang^2 + \tfrac{n-1}{n^3}u^2 p_\ss{u}^2 ) -  \pderiv{V^{0}}{u} + f_\ss{u}
\end{array} 
\qquad\quad 
\left| \quad\begin{array}{llll}
 \tup{f} :=-\pderiv{V^{1}}{\tup{q}} + \tup{\alphaup}
\\[6pt]
  f_\ss{u} := -\pderiv{V^{1}}{u} + \alpha_\ss{u} 
\end{array}  \right.
\end{gather} \end{small}
where  \eq{\slang^2 = q^2 p^2  - (\tup{q}\cdot\tup{p})^2}, where \eq{V^1=V^1(\bartup{q},t)} is arbitrary, and where \eq{(\tup{\alphaup},\alpha_\ss{u})} are generalized \textit{non}conservative forces (as detailed soon). Above, these have been combined with the conservative perturbations from \eq{V^{1}} into the terms \eq{(\tup{f},f_\ss{u})} — the details are given in Eq.\eqref{sum_dV_m}-Eq.\eqref{totalForce} below. 
First, we note:


\begin{notesq}
    \textit{Angular momentum.}
    The angular momentum plays a central role in this work and it is convenient to obtain the ODEs governing its evolution.
    From the above, along with Eq.\eqref{sum_ell_m}, we obtain the governing equations for the angular momentum in terms of the projective coordinates
    as follows (confirming angular momentum is conserved for any
     central-force dynamics):\footnote{Note that \eq{\dot{\tup{\slangup}}= \dot{\tup{q}}\tms \tup{p} + \tup{q}\tms \dot{\tup{p}} = \tup{q}\tms  \tup{f} = -\hdge{\tup{q}}\cdot\tup{f}}.
        And also that \eq{ \dot{\slang} = \tfrac{q^2}{\slang}\tup{p}\cdot\dot{\tup{p}} =  \tfrac{q^2}{\slang}\tup{p} \cdot\tup{f}}, and that \eq{ \dot{p} = \htup{p} \cdot \dot{\tup{p}} = \htup{p} \cdot\tup{f} =  u^n \htup{p} \cdot \tup{F}}.
    }
    \begin{align} \label{ldot_cen}
    \begin{array}{rlll}
         &\hdge{\dot{\tup{\slangup}}} = \tup{q}\wdg\tup{f}
         \,=\, u^n \htup{q} \wdg \tup{F}
    \end{array}
    &&,&&
     \dot{\tup{\slangup}} \,=\, -\hdge{\tup{q}}\cdot\tup{f}
        \,=\, - u^n \hdge{\htup{q}}\cdot\tup{F}
    &&,&&
     \dot{\slang} \,=\,  \tfrac{q^2}{\slang}\tup{p} \cdot\tup{f}  
     \,=\, u^n \htup{q}\cdot\hdge{\htup{\slangup}} \cdot\tup{F}
     &&,&&
     \dot{p} \,=\, \htup{p} \cdot\tup{f} 
    \end{align}
    where we used the relation between \eq{\tup{f}} and \eq{\tup{F}} given below in Eq.\eqref{totalForce}.
    Note that the above is equivalent to \eq{\hdge{\dot{\tup{\slangup}}} = \tup{r} \wdg \tup{F}}, and \eq{\dot{\tup{\slangup}} = -\hdge{\tup{r}} \cdot \tup{F}= \tup{r}\tms \tup{F}}, and \eq{\dot{\slang} = \tup{r}\cdot\hdge{\htup{\slangup}} \cdot\tup{F}}.
\end{notesq}

\begin{notesq}
    \textit{Simplified dynamics.} We may use the integrals of motion \eq{q=1} and \eq{\lambda=\htup{q}\cdot\tup{p}=0}  (which imply \eq{\slang^2\simeq p^2}) to simplify the dynamics in Eq.\eqref{eoms_fterm}. In particular, the ODEs for \eq{(\tup{q},\tup{p})} simplify 
    to:\footnote{The \eq{\dot{p}_\ss{u}} equation could also be written using \eq{\slang^2 \simeq p^2} as \eq{\dot{p}_\ss{u} =
     \tfrac{n}{u^{2n+1}}\big( p^2 + \tfrac{n-1}{n^3}u^2 p_\ss{u}^2 ) -  \pderiv{V^{0}}{u} + f_\ss{u}} .}
    \begin{small}
    \begin{align} \label{sum_qpdot_m1}
    \fnsz{\begin{array}{cc}
             q=1  \\
             \tup{q}\cdot\tup{p}=0 
        \end{array}}
        \quad \Rightarrow \qquad 
    \begin{array}{ll}
         \dot{\tup{q}} \,\simeq\,\tfrac{1}{u^{2n}}  \tup{p}  
    \\[4pt] 
        \dot{\tup{p}} \,\simeq\,
        -\tfrac{1}{u^{2n}}p^2\tup{q} + \tup{f} 
        \,\simeq\,
        -\tfrac{1}{u^{2n}}\slang^2\tup{q} + \tup{f} 
     \end{array}
     &&
     \fnsize{also,}\quad 
     \begin{array}{lll}
         \slang^2 \simeq p^2
         \\[4pt]
         \dot{\slang} \simeq \dot{p} = \htup{p} \cdot\tup{f}  \simeq  u^n \htup{p} \cdot \tup{F}
     \end{array}
    \end{align}
    \end{small} 
    where \eq{q=1} and \eq{\tup{q}\cdot\tup{p}=0} are used only \textit{after} differentiating \eq{\mscr{H}} in Eq.\eqref{eoms_fterm}. None of the main developments in this work require simplifying the dynamics as above. 
\end{notesq}

\begin{notesq}
    \textit{Simplified Hamiltonian.} 
    Note the above simplified dynamics are the same as what we would obtain if we simplified \eq{\mscr{H}} itself by eliminating the term \eq{(\tup{q}\cdot\tup{p})^2 = 0^2 } (appearing in the \eq{\slang^2} term). That is, we could, if we wished, use a simplified  
    Hamiltonian:\footnote{Whether we simplify using \eq{\tup{q}\cdot\tup{p}=0} before or after finding the equations of motion makes no difference in this case. That is, differentiating the full \eq{\mscr{H}} in Eq.\eqref{eoms_fterm} and then simplifying the result with \eq{\tup{q}\cdot\tup{p}=0} leads to the same thing as differentiating \eq{ \mscr{H}_\ss{\mrm{simplified}}} (which has been pre-simplified using \eq{\tup{q}\cdot\tup{p}=0}).
    This is \textit{not} true for the integral \eq{q=1}.}
    \begin{align} \label{sum_H_m2}
       \mscr{H}_\ss{\mrm{simplified}} \,=\,  \tfrac{1}{u^{2n}} \tfrac{1}{2}\big( q^2 p^2  + \tfrac{1}{n^2}u^2 p_\ss{u}^2 \big) \,+\, V^{0}(u) + V^{1}
    \end{align}
    The above generates the same (simplified) dynamics as the full Hamiltonian from Eq.\eqref{eoms_fterm}. Yet, it obscures the important role of the angular momentum and we make no use of \eq{\mscr{H}_\ss{\mrm{simplified}}} in this work. 
\end{notesq}

\paragraph*{Perturbation Terms.} 
We now clarify the perturbation terms \eq{\bartup{f}=(\tup{f},f_\ss{u})} and \eq{\tup{F}} appearing in the above dynamics. 
Recall that \eq{\tup{F}=(F_{r_1},F_{r_2},F_{r_3})}  denotes the cartesian components of the \textit{total} — conservative and nonconservative — perturbing forces:
\begin{small}
\begin{align} \label{a1_def}
 \tup{F} \,:=\, -\pderiv{V^{1}}{\tup{r}} \,+\,\tup{a}^\nc 
    \,=\, 
\left( \begin{array}{cc}
     \fnsize{cartesian components of total}
 \\
     \fnsize{perturbations (per unit mass)}  
\end{array} \right) 
\end{align}
\end{small}
For the transformed dynamics in projective coordinates, the terms  \eq{(-\pderiv{V^{1}}{\tup{q}},-\pderiv{V^{1}}{u})} and \eq{\bartup{\alphaup}=(\tup{\alphaup},\alpha_\ss{u})}  account for all conservative and nonconservative perturbations, respectively. 
As shown in section \ref{sec:gen_burdet_eom}, these may be expressed in terms of the cartesian components, \eq{-\pderiv{V^1}{\tup{r}}} and \eq{\tup{a}^{\ii{\mrm{nc}}}}, as: 
\begin{small}
\begin{align} \label{sum_dV_m}
 \begin{array}{cc}
        \fnsize{generalized} \\
         \fnsize{conservative:}  
\end{array}  \left. 
 \begin{array}{llll}
      \pderiv{V^{1}}{\tup{q}} \,=\, \trn{\pderiv{\tup{r}}{\tup{q}}} \cdot \pderiv{V^{1}}{\tup{r}} 
\\[4pt]
     \pderiv{V^{1}}{u} \,=\, \pderiv{\tup{r}}{u}\cdot\pderiv{V^{1}}{\tup{r}} 
  \end{array} \right.
&&,&&
 \begin{array}{cc}
       \fnsize{generalized} \\
         \fnsize{nonconservative:}  
\end{array}   \left.
\begin{array}{llll}
  \tup{\alphaup} \,=\, \trn{\pderiv{\tup{r}}{\tup{q}}} \cdot  \tup{a}^\nc 
\\[4pt]
  \alpha_\ss{u} \,=\, \pderiv{\tup{r}}{u}\cdot\tup{a}^\nc 
\end{array} \right.
&&,&&
\begin{array}{cc}
        \fnsize{generalized} \\
         \fnsize{total:} 
\end{array}   \left. 
\begin{array}{llll}
  \tup{f} :=\, -\pderiv{V^{1}}{\tup{q}} + \tup{\alphaup}  \,=\, \trn{\pderiv{\tup{r}}{\tup{q}}} \cdot  \tup{F}
\\[4pt]
  f_\ss{u} :=\, -\pderiv{V^{1}}{u} + \alpha_\ss{u}  \,=\, \pderiv{\tup{r}}{u}\cdot\tup{F} 
\end{array} \right.
\end{align}
\end{small}
where we have combined the conservative and nonconservative terms into the total \eq{\bartup{f}=(\tup{f},f_\ss{u})}. 
Plugging the appropriate \eq{\pderiv{\tup{r}}{\tup{q}}} and \eq{\pderiv{\tup{r}}{u}}
(for \eq{\tup{r}=u^n \htup{q}}) into Eq.\eqref{sum_dV_m}, we obtain the following relations between \eq{(\tup{f},f_\ss{u})} and \eq{\tup{F}}:\footnote{Note that \eq{\tup{F}\cdot\htup{q} = \tup{F}\cdot\htup{r} = F_r } is simply the radial component of the total force.} 
\begin{small}
\begin{align} \label{totalForce}
\begin{array}{lllll} 
   \tup{f} :=\, -\pderiv{V^{1}}{\tup{q}} + \tup{\alphaup} 
   & =\, \trn{\pderiv{\tup{r}}{\tup{q}}} \cdot  \tup{F} 
    & =\,  \tfrac{u^n}{q} (\imat_3  - \htup{q} \otms \htup{q} ) \cdot \tup{F}
\\[6pt]
    f_\ss{u} :=\, -\pderiv{V^{1}}{u} + \alpha_\ss{u} 
    & =\, \pderiv{\tup{r}}{u} \cdot  \tup{F}
    & =\,  n u^\ss{n-1} \htup{q}\cdot \tup{F}
\end{array}   
&&,&&
\begin{array}{llll}
     \tup{q}\cdot\tup{f} = 0 
\\[4pt]
     \tup{p} \cdot \tup{f} 
     = \tfrac{1}{q} \htup{q}\cdot\hdge{\tup{\slangup}} \cdot \tup{f}
     = \tfrac{u^n}{q^2}\htup{q}\cdot\hdge{\tup{\slangup}}\cdot\tup{F}
     \,\simeq u^n \tup{p}\cdot\tup{F}
\\[4pt]
    \tup{\slangup}\cdot\tup{f} = \tfrac{u^n}{q} \tup{\slangup}\cdot \tup{F}
\end{array}
\end{align}
\end{small} 
The above relations also hold individually for the conservative and nonconservative parts  of \eq{\bartup{f} = -\pderiv{V^{1}}{\bartup{q}} + \bartup{\alphaup}} and \eq{\tup{F} = -\pderiv{V^{1}}{\tup{r}} +\tup{a}^\nc }. 
The relation \eq{\tup{p} \cdot \tup{f}  \simeq u^n \tup{p}\cdot\tup{F}} makes use of the integrals \eq{q=1} and \eq{\tup{q}\cdot\tup{p}=0}. 

\begin{remrm} \label{rem:Frad}
    If \eq{\tup{F}} is purely radial — i.e., if \eq{\tup{F} = F_r \htup{r} = F_r \htup{q}} — then \eq{\tup{f}=0}. As such, any central-forces of any kind, conservative or nonconservative, are absent from the ODEs for \eq{(\tup{q},\tup{p})}, appearing only in those for \eq{(u,p_\ss{u})}. 
\end{remrm}

\subsection{Transformation of the evolution parameter} \label{sec:ext_cen}

So far, we have only developed the Hamiltonian which generates equations of motion with respect to time as given by Eq.\eqref{eoms_fterm}.  These equations are non-linear even in the case of unperturbed (\eq{\tup{F}= 0 }) motion, regardless of the form of the central-force potential, \eq{V^{0}}. However, we can obtain equations that are fully linear in the unperturbed case if we transform the evolution parameter to something other than time.  We will now consider a transformation of the evolution parameter from the time, \eq{t}, to two new evolution parameters,  \eq{s} and \eq{\tau} which we define through the following differential relations: 
\begin{small}
\begin{gather} \label{dtds_m1}
\boxed{\begin{array}{llll}
    \pdt{t}\,=\,  \diff{t}{s} \,=\, r^2 \,=\, u^{2n} 
&,\qquad 
     \rng{t} \,=\,  \diff{t}{\tau}\,=\, r^2/\slang  \,=\, u^{2n}/\slang 
\end{array}}
\\ \nonumber
\dot{\square} := \diff{\square}{t}
\quad,\quad 
\pdt{\square} := \diff{\square}{s} = u^{2n} \dot{\square}
\quad,\quad 
\rng{\square} := \diff{\square}{\tau} = \tfrac{1}{\slang}\pdt{\square} =  \tfrac{u^{2n}}{\slang} \dot{\square}
\end{gather}
\end{small}
where \eq{\slang^2 = q^2 p^2 -(\tup{q}\cdot\tup{p})^2}, where
a prime will denote differentiation with respect to \eq{s} and a ring will denote differentiation with respect to \eq{\tau}. For orbital motion, \eq{\tau} corresponds to the true anomaly (up to some additive constant). 
In the case of pure central-force dynamics, then the angular momentum magnitude is conserved, \eq{\slang=\slang_\zr}, and the parameters \eq{s} and \eq{\tau} defined above are then related simply by \eq{ \tau = \slang s}  (assuming \eq{\tau_\zr=s_\zr=0}).

\paragraph*{Extended Phase Space Dynamics (General).}
As shown in Appx.~\ref{sec:ext_coord},
we may transform the independent variable from time, \eq{t}, to any other evolution parameter, \eq{\varep}, and still maintain the usual Hamiltonian structure by making use of the \textit{extended phase space} where time is treated as an additional configuration coordinate with conjugate momenta \eq{p_t}. Given a specified differential relation \eq{\diff{t}{\varep}} — which is some function of the canonical coordinates as in Eq.\eqref{dtds_m1}  — the \textit{extended Hamiltonian}, denoted \eq{\wtscr{H}(\bartup{q},\bartup{p},t,p_t)}, is then given by Eq.\eqref{Hext_app} in Appx.~\ref{sec:ext_coord} as
\begin{small}
\begin{align} \label{Hext_def}
     \fnsize{extended Hamiltonian: } & \quad \wtscr{H} \,=\, \diff{t}{\varep}(\mscr{H}+p_t) 
\end{align}
\end{small}
As seen in Eq.\eqref{qp'_coord} of Appx.~\ref{sec:ext_coord}, the canonical equations of motion with \eq{\varep} as the evolution parameter are then obtained from \eq{\wtscr{H}} in the familiar manner:
\begin{align} \label{qp_eom_ext_def}
\begin{array}{ll}
     & \diff{\bartup{q}}{\varep} \,=\, \pderiv{\wtscr{H}}{\bartup{p}} 
     \\[4pt]
     & \diff{t}{\varep} \,=\, \pderiv{\wtscr{H}}{p_t}
\end{array}
\qquad\qquad
\begin{array}{ll}
     & \diff{\bartup{p}}{\varep} \,=\, -\pderiv{\wtscr{H}}{\bartup{q}} \,+\,\diff{t}{\varep}\bartup{\alphaup}
     \\[4pt]
     & \diff{p_t}{\varep} \,=\, -\pderiv{\wtscr{H}}{t} \,+\,\diff{t}{\varep}\alpha_t
\end{array}
\end{align}
Where \eq{\bartup{\alphaup}:=\tup{a}^\nc\cdot\pderiv{\tup{r}}{\bartup{q}}} are the usual generalized nonconservative forces (per unit mass) and where \eq{\alpha_t} has units of power (per unit mass) and is given by Eq.\eqref{gt_qp_def} in  Appx.~\ref{sec:ext_coord} as
\begin{small}
\begin{align}  \label{g_t}
    \alpha_t \,=\, -\bartup{\alphaup}\cdot\dot{\bartup{q}} 
    \,=\, -\tup{a}^\nc\cdot\pderiv{\tup{r}}{\bartup{q}} \cdot \dot{\bartup{q}}
    \,=\, -\tup{a}^\nc\cdot\dot{\tup{r}} 
&&,&& 
    \diff{t}{\varep}\alpha_t \,=\, -\bartup{\alphaup}\cdot\diff{\bartup{q}}{\varep} \,=\, -\tup{a}^\nc\cdot\diff{\tup{r}}{\varep}
\end{align}
\end{small}
where \eq{\pderiv{\tup{r}}{\bartup{q}}\cdot \dot{\bartup{q}} =\dot{\tup{r}}} holds in the case that \eq{\tup{r}(\bartup{q}}) is not an explicit function of time (true in our case). 
Thus, for the problem at hand, \eq{\alpha_t}  is the negative of the rate that nonconservative forces do work (per unit mass) on the system. 
This makes sense when we 
recall\footnote{See Eq.\eqref{pt_deff} of Appx.~\ref{sec:ext_coord}.}
that, along any solution curve of Eq.\eqref{qp_eom_ext_def},  it holds that \eq{p_t} is related the the \textit{value} of the non-extended Hamiltonian along the curve simply by \eq{ p_t  = -\mscr{H} }. 
In many cases, including the present, \eq{\mscr{H}} is the total (specific) energy of the system/particle.  

When differentiating the extended Hamiltonian, \eq{\wtscr{H}}, \eq{p_t} should be treated as an independent coordinate, \textit{not} as a function of the other coordinates. However, after obtaining the equations of motion, we may use the fact that  \eq{ p_t  = -\mscr{H} } holds along any solution curve. That is, the relation \eq{ p_t  = -\mscr{H} } may be substituted into the resulting ODEs themselves (although this is not necessary). Upon doing so, one  finds that the general extended phase space ODEs for \eq{(\bartup{q},\bartup{p})} in Eq.\eqref{qp_eom_ext_def} are equivalent to a conformal scaling by \eq{\diff{t}{\varep}}:
\begin{small}
\begin{align} \label{qpdot_ext_altyalt}
\begin{array}{llll}
      \diff{\bartup{q}}{\varep} \,=\, \pderiv{\wtscr{H}}{\bartup{p}}
      \;=\; \diff{t}{\varep} \pderiv{\mscr{H}}{\bartup{p}} \,=\, \diff{t}{\varep} \dot{\bartup{q}}
\qquad,\qquad 
      \diff{\bartup{p}}{\varep} \,=\, -\pderiv{\wtscr{H}}{\bartup{q}} +\diff{t}{\varep}\bartup{\alphaup}
    \;=\;  \diff{t}{\varep}(-\pderiv{\mscr{H}}{\bartup{q}} + \bartup{\alphaup} )
    \,=\, \diff{t}{\varep} \dot{\bartup{p}}
\end{array}
\end{align}
\end{small}
See the developments leading to Eq.\eqref{qpdot_ext_alt_apx} in Appendix \ref{sec:ext_coord} for details.



\paragraph*{The $s$- and $\tau$-Parameterized Projective Coordinate Dynamics.}
We now return to the specific case at hand. 
With \eq{\diff{t}{s}=u^{2n}}, 
the extended Hamiltonian \eq{\wtscr{H}:=\diff{t}{s}(\mscr{H}+p_t)} for the \eq{s}-parameterized dynamics is then as 
follows:\footnote{As per Eq.\eqref{f(r)_m}, the fact that \eq{\wt{V}^\zr(u)} is still a function only of \eq{u}, but not \eq{\tup{q}}, holds only for the the point transformation of the form \eq{\tup{r}(\bartup{q})=u^n\htup{q}\,} (\eq{m=-1}).  If we had instead used the point transformation \eq{\tup{r}(\bartup{q})=u^n\tup{q}\,} (\eq{m=0}), then we would have \eq{\pdt{t}=r^2=u^{2n}\mag{\tup{q}}^2} and \eq{\wt{V}^\zr(\bartup{q})}.} 
\begin{small}
\begin{align} \label{Hs_m1}
    \begin{array}{lr}
        \wtscr{H} \,=\, \frac{1}{2}\big( \slang^2 + \tfrac{1}{n^2}u^2 p_\ss{u}^2 \big)  \,+\, \wt{V}^\zr(u)  \,+\, \wt{V}^\ss{1}(\bartup{q},t)
    \,+\, u^{2n}p_t  
    \end{array}
    &&
    \wt{V} := \pdt{t}V = u^{2n} V 
\end{align}
\end{small}
where \eq{\slang^2 = q^2 p^2 -(\tup{q}\cdot\tup{p})^2}. 
The \eq{s}-parameterized  dynamics are then found from \eq{\wtscr{H}} as in Eq.\eqref{qp_eom_ext_def}:
\begin{small}
\begin{align} \label{qp_eom_s_m1}
\begin{array}{llll}
     \pdt{\tup{q}} \,=\, \pderiv{\wtscr{H}}{\tup{p}} \,=\,  -\hdge{\tup{\slangup}}\cdot\tup{q}   &,
 \\[6pt]
       \pdt{u} \,=\, \pderiv{\wtscr{H}}{p_u}  \,=\,  \tfrac{u^2}{n^2}p_\ss{u}  &,
 \\[6pt]
     \pdt{t} \,=\, \pderiv{\wtscr{H}}{p_t} \,=\,  u^{2n} &,
\end{array}
 \qquad
\begin{array}{llll}
      \pdt{\tup{p}} \,=\, -\pderiv{\wtscr{H}}{\tup{q}} +\pdt{t}\tup{\alphaup} 
      & =\, 
      -\hdge{\tup{\slangup}}\cdot\tup{p}  \,+\, \pdt{t} \tup{f}
 \\[6pt]
       \pdt{p}_\ss{u} \,=\, -\pderiv{\wtscr{H}}{u} +\pdt{t}\alpha_\ss{u} 
       & =\,   -\tfrac{1}{n^2}u p_\ss{u}^2     -  \pderiv{\wt{V}^\zr}{u}  -  \pderiv{\wt{V}^\ss{1}}{u} -   2n u^{2n-1} p_t + \pdt{t}\alpha_\ss{u}
 \\[6pt]
     \pdt{p}_t \,=\, -\pderiv{\wtscr{H}}{t} + \pdt{t}\alpha_t
     & =\,   \pdt{t}(-\pderiv{V^{1}}{t} + \alpha_t) \,=\, -(\pdt{t}\pderiv{V^1}{t} + \tup{\alphaup}\cdot\pdt{\tup{q}} +  \alpha_\ss{u}\pdt{u})
\end{array} 
\end{align}
\end{small}
Note we may use the relation \eq{p_t=-\mscr{H}}, with \eq{\mscr{H}} given in Eq.\eqref{eoms_fterm}, 
to eliminate \eq{p_t} from the  above \eq{\pdt{p}_\ss{u}}  
equation:\footnote{Differentiating \eq{\wtscr{H}} directly gives: \\
\eq{\qquad\qquad \pdt{p}_\ss{u} \,=\,  -\tfrac{u}{n^2}p_\ss{u}^2 - \pderiv{\wt{V}^\zr}{u}   -  2nu^{2n-1}(p_t +V^{1})  + u^{2n}(-\pderiv{V^1}{u} +\alpha_\ss{u}) 
\quad = \tfrac{n}{u} ( \slang^2  +\tfrac{n-1}{n^3}u^2 p_\ss{u}^2 ) 
    - u^{2n} \pderiv{V^{0}}{u} - \tfrac{2n}{u}\wtscr{H}
    + u^{2n}(-\pderiv{V^1}{u} +\alpha_\ss{u}) } .\\
Using \eq{p_t=-\mscr{H} = -\tfrac{1}{2 u^{2n}}(\slang^2 + \tfrac{1}{n^2}u^2 p_\ss{u}^2) -V^0 -V^1} (i.e., \eq{\wtscr{H}=0}) and \eq{f_\ss{u}=-\pderiv{V^1}{u} +\alpha_\ss{u}}, we then obtain the expression in Eq.\eqref{dpds_m1}.     }
\begin{small}
\begin{align} \label{dpds_m1} 
 p_t  = -\mscr{H} 
\qquad \Rightarrow \qquad
     \pdt{p}_\ss{u} \,=\,
    \tfrac{n}{u} ( \slang^2  +\tfrac{n-1}{n^3}u^2 p_\ss{u}^2 ) 
    - \pdt{t} \pderiv{V^{0}}{u}
    + \pdt{t} f_\ss{u}
    \;=  \pdt{t}\dot{p}_\ss{u}
\end{align}
\end{small}
where, as usual, we collect the perturbation terms into \eq{\tup{f}:=-\pderiv{V^1}{\tup{q}}+\tup{\alphaup}} and  \eq{f_\ss{u}:=-\pderiv{V^1}{u}+ \alpha_\ss{u}}.
We also consider the evolution parameter \eq{\tau} defined by \eq{\diff{t}{\tau} = u^{2n}/\slang}. The associated extended Hamilton is given by \eq{\whscr{H}:= \diff{t}{\tau}(\mscr{H}+p_t) = \tfrac{1}{\slang}\wtscr{H} }, with \eq{\wtscr{H}} given above in Eq.\eqref{Hs_m1}. The \eq{\tau}-parameterized dynamics follow from Eq.\eqref{qp_eom_ext_def}:
\begin{small}
\begin{gather} \label{dqp_TA_full}
    \whscr{H}\,=\,  \tfrac{1}{\slang}\wtscr{H} 
\qquad,\qquad 
\begin{array}{lllll}
      \rng{\tup{q}} \,=\, \pderiv{\whscr{H}}{\tup{p}}   \,=\, 
      -(1 - \tfrac{\whscr{H}}{\slang})  \hdge{\htup{\slangup}}\cdot\tup{q} 
      &,
\\[5pt]
      \rng{u} \,=\, \pderiv{\whscr{H}}{p_u}  \,=\,   \tfrac{u^2}{n^2 \slang}p_\ss{u}
      &,
\\[5pt]  
    \rng{t}  \,=\, \pderiv{\whscr{H}}{p_t} \,=\,  u^{2n}/\slang
    &,
\end{array}
\qquad 
\begin{array}{llllll}
    \rng{\tup{p}} \,=\, -\pderiv{\whscr{H}}{\tup{q}} +\rng{t}\tup{\alphaup} 
      & =\,  -(1- \tfrac{\whscr{H}}{\slang} ) \hdge{\htup{\slangup}}\cdot\tup{p}  
      +\rng{t}\tup{f} 
\\[5pt]
     \rng{p}_\ss{u} \,=\, -\pderiv{\whscr{H}}{u} +\rng{t} \alpha_\ss{u}
      & =\, \tfrac{n}{u \slang} ( \slang^2  +\tfrac{n-1}{n^3}u^2 p_\ss{u}^2 )   -  \rng{t} \pderiv{V^{0}}{u} - \tfrac{2n}{u} \whscr{H}  +  \rng{t} f_\ss{u}
\\[5pt]
     \rng{p}_t \,=\, -\pderiv{\whscr{H}}{t} + \rng{t} \alpha_t  
      & =\, 
      -( \rng{t}\pderiv{V^1}{t} + \tup{\alphaup}\cdot\rng{\tup{q}} + \alpha_\ss{u} \rng{u} )
\end{array}
\end{gather}
\end{small}
Yet, as per Eq.\eqref{qpdot_ext_altyalt}, we may use the relation \eq{p_t=-\mscr{H}} (i.e., \eq{\wtscr{H}=0=\whscr{H}}) to rewrite the above ODEs in form \eq{\rng{\square} = \rng{t}\dot{\square} = \tfrac{1}{\slang}\pdt{\square}}:
\begin{small}
\begin{align} \label{qp_eom_f_m1} 
p_t = -\mscr{H} 
\quad \Rightarrow \quad
\left\{ \qquad \begin{array}{ll}
     \rng{\tup{q}}  =    -\hdge{\htup{\slangup}}\cdot\tup{q}    &, 
     \\[7pt]
       \rng{u}    =   \tfrac{u^2}{n^2 \slang}p_\ss{u}  &, 
     \\[7pt]
     \rng{t}  =   {u^{2n}}/{\slang} &, 
\end{array} \right.
\qquad
\begin{array}{llll}
      \rng{\tup{p}}    =    -\hdge{\htup{\slangup}}\cdot\tup{p}   +  \rng{t} \tup{f}
 \\[6pt]
       \rng{p}_\ss{u}  =   \tfrac{n}{u \slang} ( \slang^2  +\tfrac{n-1}{n^3}u^2 p_\ss{u}^2 )  
    -  \rng{t}\pderiv{V^{0}}{u}
    + \rng{t} f_\ss{u}
 \\[6pt]
     \rng{p}_t   
     =  -( \rng{t}\pderiv{V^1}{t} + \tup{\alphaup}\cdot\rng{\tup{q}} + \alpha_\ss{u} \rng{u} )
\end{array}
\end{align}
\end{small}
For any of the above ODEs, the governing equations for the angular momentum functions are given as in Eq.\eqref{ldot_cen}, along with \eq{\pdt{\square} = \pdt{t}\dot{\square}} or \eq{\rng{\square} = \rng{t}\dot{\square} = \tfrac{1}{\slang}\pdt{\square}}.



\begin{notesq} 
\textit{Simplifications.}
    The above \eq{s}- and \eq{\tau}-parameterized dynamics have \textit{not} been simplified using the integrals of motion \eq{q=1} or \eq{\lambda=\htup{q}\cdot\tup{p}=0}, which imply \eq{\slang \simeq p} and the other relations in Eq.\eqref{qpl_rels}. As usual, we are free to simplify the dynamics use said relations from Eq.\eqref{qpl_rels} such that, for example, the above ODEs for \eq{(\tup{q},\tup{p})} simplify to:
    \begin{small}
    \begin{flalign} \label{Hs_m17_simp}
    && \fnsz{\begin{array}{cc}
         q=1  \\
         \htup{q}\cdot\tup{p}=0 
    \end{array}}
    \Rightarrow \quad
    \left\{\qquad\begin{array}{llll}
         \pdt{\tup{q}} \simeq  q^2\tup{p}   \,\simeq \tup{p}
         \\[4pt]
         \pdt{\tup{p}} 
         \simeq  -p^2\tup{q} + \pdt{t} \tup{f} 
          \,\simeq
          -\slang^2 \tup{q}  + \pdt{t}\tup{f} 
    \end{array} \right.
    \qquad\fnsize{or,}\qquad
    \begin{array}{llll}
         \rng{\tup{q}} \simeq \tfrac{1}{\slang} \tup{p} 
            \,\simeq \htup{p} 
         \\[4pt]
         \rng{\tup{p}} \simeq -\slang \tup{q} + \rng{t}\tup{f}
    \end{array}
    &&,&&
    \slang^2 \simeq p^2 
    \quad
    \end{flalign}
    \end{small}
where the relation \eq{\slang^2\simeq p^2} could also be used
to rewrite the the ODEs for 
    \eq{(u,p_\ss{u})}.\footnote{Recall that \eq{\htup{q}\cdot\tup{p}=0} could, if we wished, also be used to simplify the Hamiltonian itself. We could then regard the \textit{extended} Hamiltonian (for \eq{s}) as \eq{ \wtscr{H}_\ss{\mrm{simplified}} = \tfrac{1}{2}\big( q^2 p^2  + \tfrac{1}{n^2}u^2 p_\ss{u}^2 \big)  + \wt{V}^\zr  + \wt{V}^\ss{1}
    + u^{2n}p_t  }. This results in the same simplified \eq{s}-parameterized ODEs (and likewise for \eq{\whscr{H}_\ss{\mrm{simplified}}}).  Yet, unlike the BF transformation, we find little use for the simplified Hamiltonian.}
\end{notesq}

\subsection{Linear dynamics \& closed-form solutions}  \label{sec:2BP_lin}

The dynamics of \eq{(\tup{q},\tup{p})} correspond to the angular/rotational motion, while the dynamics of \eq{(u,p_\ss{u})} correspond to the radial motion. 
We will address \eq{(\tup{q},\tup{p})} and  \eq{(u,p_\ss{u})} separately as the particular form of the central-force potential \eq{V^0} only matters for the latter, having no impact whatsoever on the ``\eq{(\tup{q},\tup{p})}-part'' of the dynamics.

In the following, we will often present the \eq{s}- and \eq{\tau}-parameterized dynamics in tandem with one another (with \eq{s} and \eq{\tau} defined in Eq.\eqref{dtds_m1} such that \eq{\mrm{d} \tau = \slang \mrm{d} s}).
When considering pure central-force dynamics — that is, the unperturbed case  \eq{\tup{F}=0} or even the radially-perturbed case \eq{\tup{F}=F_r\htup{r}} — then recall from remark \ref{rem:Frad} that
this leads to \eq{\tup{f}=0} such that angular momentum is preserved 
and, as such, the parameters \eq{s} and \eq{\tau} are then related simply by \eq{\tau = \slang s} (assuming \eq{\tau_\zr=s_\zr=0}):\footnote{The relation \eq{\tau = \slang s} follows from \eq{\mrm{d} \tau = \slang \mrm{d} s} along with preservation of \eq{\slang}.}
\begin{small}
\begin{flalign} \label{Frad_f0}
\begin{array}{cc}
   \fnsize{arbitrary}  \\
   \fnsize{central-force}\\
   \fnsize{dynamics} 
 \end{array} 
&&
    \fnsize{if }  \; \tup{F} = F_r \htup{r}
    \qquad \Rightarrow \qquad 
    \tup{f} = 0 
    \quad,\quad 
    \eq{\tup{\slangup}=\tup{\slangup}_\zr}
    \quad,\quad 
    \tau = \slang s
&&
\end{flalign}
\end{small}
In this case, it makes little difference whether one wishes to express solutions in terms of \eq{s} or \eq{\tau} as reparameterization between them is trivial. For brevity, we will often express solutions in terms of \eq{\tau}. The equivalent expressions in terms of \eq{s} are then easily obtained using \eq{\tau = \slang s}.


\subsubsection{The rotational/angular motion (for any central-force dynamics)} 

The \eq{s}- or \eq{\tau}-parameterized dynamics for \eq{(\tup{q},\tup{p})} are given in terms of the antisymmetric angular momentum matrix, \eq{\hdge{\tup{\slangup}}= \tup{q} \wdg \tup{p}},  as follows (from Eq.\eqref{qp_eom_s_m1} and Eq.\eqref{qp_eom_f_m1}):
\begin{small}
\begin{align} \label{dqdp_nonsimp}
 \begin{array}{lllll}
      \pdt{\tup{q}}  \,=\,  -\hdge{\tup{\slangup}}\cdot\tup{q}  
 \\[4pt]
     \pdt{\tup{p}} \,=\, -\hdge{\tup{\slangup}}\cdot\tup{p}  \,+\, \pdt{t} \tup{f}
 \end{array}
\qquad\fnsize{or,}\qquad
 \begin{array}{lllll}
      \rng{\tup{q}}  \,=\,  -\hdge{\htup{\slangup}}\cdot\tup{q}  
 \\[4pt]
     \rng{\tup{p}} \,=\, -\hdge{\htup{\slangup}}\cdot\tup{p}  \,+\, \rng{t} \tup{f}
 \end{array}
\end{align}
\end{small}
where  \eq{\pdt{t}=u^{2n}}, \eq{\rng{t}=u^{2n}/\slang}, and \eq{\hdge{\htup{\slangup}}=\tfrac{1}{\slang}\hdge{\tup{\slangup}}}. 
In the unperturbed case \eq{\tup{F}=0} — or even the radially-perturbed case \eq{\tup{F}=F_r \htup{r}} —  then \eq{\tup{f}=0} and \eq{\hdge{\tup{\slangup}}=\hdge{\tup{\slangup}}_\zr} is  constant such that the above can be seen as a linear first-order ODEs with solutions given by the matrix exponential:
\begin{small}
\begin{flalign} \label{E3sol0_s_prj}
\begin{array}{cc}
   \fnsize{arbitrary}  \\
   \fnsize{central-force}\\
   \fnsize{dynamics} 
 \end{array}
\quad \Rightarrow \qquad 
\tup{f} = 0 
&& \Rightarrow &&
\begin{array}{lcllllll}
     \tup{q}_s \,=\,  \mrm{e}^\fnsz{-\hdge{\tup{\slangup}} s} \cdot \tup{q}_\zr  
     &=&
       \mrm{e}^\fnsz{-\hdge{\htup{\slangup}} \tau} \cdot \tup{q}_\zr \,=\,  \tup{q}_\tau 
\\[4pt]
      \tup{p}_s \,=\,  \mrm{e}^\fnsz{-\hdge{\tup{\slangup}} s} \cdot \tup{p}_\zr 
      &=&  \mrm{e}^\fnsz{-\hdge{\htup{\slangup}} \tau} \cdot \tup{p}_\zr \,=\,  \tup{p}_\tau 
\end{array}
&&,&&
 \begin{array}{lllllll}
        \hdge{\tup{\slangup}} =  \tup{q} \wdg \tup{p} = \hdge{\tup{\slangup}}_\zr
        \,\in \somat{3}
     \\[4pt]
         \mrm{e}^\fnsz{-\hdge{\tup{\slangup}} s} = \mrm{e}^\fnsz{-\hdge{\htup{\slangup}} \tau}
         \in \Somat{3}
\end{array}
\quad
\end{flalign}
\end{small}
where \eq{\tup{\slangup}=\tup{\slangup}_\zr} is identified with its constant value along the solution curve such that \eq{\tau = \slang s}. 
Note the \eq{s}- or \eq{\tau}-parameterized solutions coincide; the above matrix exponentials are equivalent since \eq{\tau = \slang s} implies \eq{\hdge{\tup{\slangup}} s = \hdge{\htup{\slangup}}\tau  }. 
These exponentials are conveniently expressed using 
the Rodrigues rotation formula\footnote{For any \eq{\tup{\rhoup}\in\mbb{R}^3} with \eq{\hdge{\tup{\rhoup}}\in\somat{3}}, and some parameter \eq{\varep}, then the Rodrigues formula gives \eq{\mrm{e}^{\hdge{\tup{\rhoup}}\varep}\in \Somat{3}} as follows (with \eq{\rho=\mag{\tup{\rhoup}}}):
\begin{gather} \nonumber
 \mrm{e}^{\hdge{\tup{\rhoup}}\varep}  =   \imat_3 + \tfrac{1}{\rho}( \sin{\rho \varep} )\hdge{\tup{\rhoup}} + \tfrac{1}{\rho^2}(1-\cos{\rho \varep}) \hdge{\tup{\rhoup}} \cdot \hdge{\tup{\rhoup}}  
 \,=\,
 \imat_3 + ( \sin{\rho \varep} )\hdge{\htup{\rhoup}} + (1-\cos{\rho \varep}) \hdge{\htup{\rhoup}} \cdot \hdge{\htup{\rhoup}}  
 \;\;\;=\;
  \csn{\rho \varep}\,\imat_3 + ( \sin{\rho \varep} )\hdge{\htup{\rhoup}} + (1-\cos{\rho \varep}) \htup{\rhoup} \otms \htup{\rhoup}
\\ \nonumber 
 \inv{(\mrm{e}^{\hdge{\tup{\rhoup}}\varep})} = \eq{\mrm{e}^{-\hdge{\tup{\rhoup}}\varep}  =  \trn{(\mrm{e}^{\hdge{\tup{\rhoup}}\varep})} }
\end{gather} 
Note that if \eq{\sigma:=\rho\varep} then \eq{ \mrm{e}^{\hdge{\htup{\rhoup}}\sigma}= \mrm{e}^{\hdge{\tup{\rhoup}}\varep}}.  
},
leading to
\begin{small}
\begin{flalign}
&&
    \mrm{e}^\fnsz{-\hdge{\tup{\slangup}} s} \,=\,  \mrm{e}^\fnsz{-\hdge{\htup{\slangup}} \tau}
    \,=\,   \imat_3 -  \snn{\tau}\,\hdge{\htup{\slangup}} + (1-\csn{\tau}) \hdge{\htup{\slangup}} \cdot \hdge{\htup{\slangup}}
    \;=\;
   \csn{\tau}\, \imat_3 -  \snn{\tau}\,\hdge{\htup{\slangup}} + (1-\csn{\tau}) {\htup{\slangup}} \otms {\htup{\slangup}}
    \,=: R_\tau(\htup{\slangup}) \, \in \Somat{3}
    &&,&&
    \tau = \slang s
\end{flalign}
\end{small}
Using the above — along with relations from Eq.\eqref{angmoment_rels_crd} — the solutions in Eq.\eqref{E3sol0_s_prj} are expressed in  
closed-form as:\footnote{The solutions are inverted using \eq{\inv{R}_\tau(\htup{\slangup})  = R_{-\tau}(\htup{\slangup}) = \trn{R}_\tau(\htup{\slangup})}: 
\eq{\quad 
 \tup{q}_\zr =  \inv{R}_\tau(\htup{\slangup})  \cdot \tup{q}_\tau  
     =
      \tup{q}_\tau \csn{\tau} +  \hdge{\htup{\slangup}} \cdot \tup{q}_\tau \snn{\tau} 
\quad,\quad 
      \tup{p}_\zr =  \inv{R}_\tau(\htup{\slangup}) \cdot \tup{p}_\tau  
      =  \tup{p}_\tau \csn{\tau} + \hdge{\htup{\slangup}} \cdot \tup{p}_\tau \snn{\tau}\;. }
}
\begin{small}
\begin{flalign} \label{E3sol_s_prj}
\begin{array}{cc}
      \fnsize{arbitrary}  \\
       \fnsize{central-force}\\
       \fnsize{dynamics} 
 \end{array}
&&
\boxed{\begin{array}{llllll}
     \tup{q}_\tau \,=\, R_\tau(\htup{\slangup})  \cdot \tup{q}_\zr  
      &  =\,  \tup{q}_\zr \csn{\tau} - \hdge{\htup{\slangup}} \cdot \tup{q}_\zr \snn{\tau} 
\\[4pt]
       \tup{p}_\tau \,=\, R_\tau(\htup{\slangup})  \cdot \tup{p}_\zr  
      &  =\,   \tup{p}_\zr \csn{\tau} - \hdge{\htup{\slangup}} \cdot \tup{p}_\zr \snn{\tau} 
\end{array}} 
&&
\left.\begin{array}{lllll}
\tup{\slangup}=\tup{\slangup}_\zr
\\[3pt]
    \tup{q}_\tau \cdot \tup{q}_\tau  = \tup{q}_\zr\cdot \tup{q}_\zr 
\\[3pt]
    \tup{q}_\tau \cdot \tup{p}_\tau = \tup{q}_\zr\cdot \tup{p}_\zr 
\\[3pt]
     \tup{p}_\tau \cdot \tup{p}_\tau = \tup{p}_\zr\cdot \tup{p}_\zr 
\end{array}\right. 
\quad 
\end{flalign}
\end{small}
The above is a special orthogonal rotation through angle \eq{\tau=\slang s} about the fixed axis \eq{\htup{\slangup}=\htup{q}\tms \htup{p}=\htup{\slangup}_\zr} (orbit normal). Therefore, as indicated above, the magnitudes of, and inner product of, \eq{\tup{q}} and \eq{\tup{p}} are all preserved (this is quick to verify from the above solutions). 
As per Eq.\eqref{qp=0}, preservation of  \eq{\mag{\tup{q}}} and \eq{\tup{q} \cdot \tup{p} = q\lambda} is guaranteed for any and all dynamics (not just for central-forces). Preservation of \eq{\tup{\slangup}} and \eq{\mag{\tup{p}}} is specific to any central-force dynamics. 


\begin{notesq}
    \textit{Simplifications.} 
    Nothing so far has been simplified using the  integrals of motion \eq{q=1} or \eq{\lambda=\htup{q}\cdot\tup{p}=0}. As discussed, we are free to make these simplifications such that the above unperturbed dynamics for \eq{(\tup{q},\tup{p})} are then equivalent to:
    \begin{small}
    \begin{flalign} \label{qp_sol_f_m1}
    \;\;
    \fnsz{\begin{array}{cc}
             q=1  \\
             \tup{q}\cdot\tup{p}=0 
        \end{array}}
        \;\Rightarrow \quad
    \left\{\;\;
     \begin{array}{lllll}
         \pdt{\tup{q}} \simeq \tup{p} 
         \\[4pt]
         \pdt{\tup{p}} \simeq -\slang^2 \tup{q} 
    \end{array}
     \right.
        \;\;\fnsize{ or,}\quad 
     \begin{array}{lllll}
         \rng{\tup{q}} \simeq \tfrac{1}{\slang}\tup{p} 
         \\[4pt]
         \rng{\tup{p}} \simeq -\slang \tup{q} 
    \end{array}
    &&\Rightarrow &&
    \begin{array}{llllll}
           \tup{q}_{\tau}  
           \, \simeq\,  \tup{q}_{\zr}\csn{\tau}  +  \tfrac{1}{\slang}\tup{p}_{\zr}\snn{\tau}
           & \simeq\,  \tup{q}_{\zr}\csn{\tau}  +  \htup{p}_{\zr}\snn{\tau}
    \\[5pt] 
           \tup{p}_{\tau} 
           \,\simeq\, \tup{p}_{\zr}\csn{\tau}   -\slang\tup{q}_{\zr}\snn{\tau} 
            & \simeq\, \tup{p}_{\zr}\csn{\tau}   -p_\zr\tup{q}_{\zr}\snn{\tau} 
    \end{array}
    &&
    \end{flalign}
    \end{small}
    which also follow immediately from the solutions in Eq.\eqref{E3sol_s_prj} using \eq{q=1} and \eq{\tup{q}\cdot\tup{p}=0} (which imply Eq.\eqref{qpl_rels}).
    Note that, unlike the unsimplified solutions in Eq.\eqref{E3sol_s_prj}, the above do not \textit{directly} verify preservation of \eq{\mag{\tup{q}}}, \eq{\tup{q}\cdot\tup{p}}, or \eq{\mag{\tup{p}}}; they are preserved \textit{iff} we limit initial conditions  to \eq{q_\zr=1} and \eq{\tup{q}_\zr\cdot \tup{p}_\zr = 0} (which is what the above already assumes), in which case the above leads to: 
    \begin{small}
    \begin{flalign} \label{qp_sol_constants_simp}
        &&
          \fnsz{\begin{array}{cc}
             q_0=1  \\
             \tup{q}_0\cdot\tup{p}_0=0 
        \end{array}}
        \;\Rightarrow \qquad
          \mag{\tup{q}_\tau} \simeq \mag{\tup{q}_\zr}\simeq 1
          \quad,\quad 
          \tup{q}_\tau \cdot \tup{p}_\tau \simeq \tup{q}_\zr\cdot \tup{p}_\zr \simeq 0
          \quad,\quad 
           \mag{\tup{p}_\tau} \simeq \mag{\tup{p}_\zr} \simeq \slang_\zr
           &&
    \end{flalign}
    \end{small}
\end{notesq}

\paragraph*{As a Linear Hamiltonian System.}
The above developments may also be combined into a single linear ODE for \eq{(\tup{q},\tup{p})\in\mbb{R}^6}:
\begin{small}
\begin{flalign} \label{dz_ds_E3}
   &&
 \begin{array}{llll}
   \diff{}{s} \fnpmat{
        \tup{q} \\
        \tup{p} }
    = L \cdot \fnpmat{
        \tup{q}\\
        \tup{p} }
    \qquad \fnsize{or,} \qquad 
      \diff{}{\tau} \fnpmat{
        \tup{q} \\
        \tup{p} }
    = \hat{L} \cdot \fnpmat{
        \tup{q}\\
        \tup{p} }    
\end{array}
&& \Rightarrow &&
     \fnpmat{
        \tup{q}_s \\
        \tup{p}_s } 
        = \mrm{e}^\fnsz{L s} \cdot  \fnpmat{
        \tup{q}_\zr \\
        \tup{p}_\zr } 
        \;\;=\;\;
        \mrm{e}^\fnsz{\hat{L} \tau} \cdot  \fnpmat{
        \tup{q}_\zr \\
        \tup{p}_\zr } 
         =
        \fnpmat{
        \tup{q}_\tau \\
        \tup{p}_\tau }
        &&
\end{flalign}
\end{small}
where the matrices are given in terms of \eq{\hdge{\tup{\slangup}} =\hdge{\tup{\slangup}}_\zr } and \eq{\tau =\slang s} by:
\begin{small}
\begin{align} \label{Lmat_spso}
       L\, ,\;  \hat{L}:=\tfrac{1}{\slang} L = \smpmat{
        -\hdge{\htup{\slangup}}  & 0 \\  0 & -\hdge{\htup{\slangup}}  }  \in \spmat{6} \cap \somat{6} 
   \qquad,\qquad 
    \mrm{e}^\fnsz{L s} = \mrm{e}^\fnsz{\hat{L} \tau} = \smpmat{
         R_\tau(\htup{\slangup}) & 0 \\
          0 &   R_\tau(\htup{\slangup})  }
    \in \Spmat{6} \cap \Somat{6} 
\end{align}
\end{small}
 Interestingly, one may also express the above as a different, but equivalent, matrix ODE.
Using the relations \eq{-\hdge{\tup{\slangup}}\cdot\tup{q} = (q^2 \imat_3 - \tup{q}\otms\tup{q})\cdot\tup{p}} and \eq{-\hdge{\tup{\slangup}}\cdot\tup{p} = -(p^2 \imat_3 - \tup{p}\otms\tup{p})\cdot\tup{q}}, the above dynamics for \eq{(\tup{q},\tup{p})} are equivalent to:
\begin{small}
\begin{align} \label{dz_ds_E3_alt}
   \diff{}{s}  \fnpmat{ \tup{q} \\ \tup{p}}
     =  A \cdot   \fnpmat{ \tup{q} \\ \tup{p}}
 \qquad \fnsize{or,} \qquad 
    \diff{}{\tau} 
        \fnpmat{ \tup{q} \\ \tup{p}}
     =  \tfrac{1}{\slang} A \cdot  \fnpmat{ \tup{q} \\ \tup{p}}
    &&\Rightarrow&&
    \fnpmat{ \tup{q}_s \\ \tup{p}_s }
     = \mrm{e}^{A s} \cdot 
    \fnpmat{ \tup{q}_\zr \\ \tup{p}_\zr }
    \;\;=\;\;
    \mrm{e}^{ \frac{1}{\slang} A \tau} \cdot 
     \fnpmat{ \tup{q}_\zr \\ \tup{p}_\zr } 
     =
    \fnpmat{ \tup{q}_\tau\\ \tup{p}_\tau }
\end{align}
\end{small}
where \eq{A=A_\zr \in \spmat{6}} is again a matrix of integrals of motion (for any central-force dynamics), given as follows:
\begin{small}
\begin{align} \label{dz_ds_E3_alt2}
     A := \smpmat{ -(\tup{q}\cdot\tup{p})\imat_3 & q^2 \imat_3 \\  -p^2 \imat_3 & (\tup{q}\cdot\tup{p})\imat_3 }
     \in \spmat{6}
     \qquad,\qquad 
    \mrm{e}^{A s}  =  \mrm{e}^{ \frac{1}{\slang}A \tau} = 
  \smpmat{
    \big(\csn{\tau} - \tfrac{1}{\slang}(\tup{q} \cdot \tup{p}) \snn{\tau} \big) \imat_3 
    & \big(\tfrac{1}{\slang}q^2 \snn{\tau} \big) \imat_3 \\
    -\big( \tfrac{1}{\slang}p^2 \snn{\tau} \big) \imat_3 & 
      \big( \csn{\tau} +  \tfrac{1}{\slang}(\tup{q} \cdot \tup{p}) \snn{\tau} \big) \imat_3
    }
    \in \Spmat{6}
\end{align}
\end{small}
One can verify that any of the above matrix ODEs and solutions agree with those in Eq.\eqref{E3sol_s_prj}.
\begin{small}
\begin{itemize}
    \item \sloppy \textit{Eigenvalues.}
    We note that \eq{L\in\spmat{6}} in Eq.\eqref{dz_ds_E3} has eigenvalues \eq{(0,0,\mrm{i}\slang,\mrm{i}\slang,-\mrm{i}\slang,-\mrm{i}\slang)}, 
    whereas \eq{A\in\spmat{6}} in Eq.\eqref{dz_ds_E3_alt} has eigenvalues \eq{(\mrm{i}\slang,\mrm{i}\slang,\mrm{i}\slang,-\mrm{i}\slang,-\mrm{i}\slang,-\mrm{i}\slang)}. For the \eq{\tau}-parameterized dynamics, \eq{\slang} drops out and these reduce to \eq{0} or \eq{\pm{\mrm{i}}}.  
     \item \textit{Simplifications.}  If we use \eq{q=1} and \eq{\lambda=\htup{q}\cdot\tup{p}=0} to simplify the central-force dynamics for \eq{(\tup{q},\tup{p})}, then the above dynamics in the form of Eq.\eqref{dz_ds_E3_alt}-Eq.\eqref{dz_ds_E3_alt2} simplify with:
    \begin{small}
    \begin{align} \label{dqp_Fcen_simp}
           A \,\simeq\, \fnpmat{ 0 &  \imat_3 \\[2pt]
           -\slang^2 \imat_3 & 0 }
         \in \spmat{6}
         \qquad,\qquad 
         \mrm{e}^{ A s} = \mrm{e}^{ \frac{1}{\slang}A \tau} \,\simeq\,
      \fnpmat{
       \csn{\tau}\,\imat_3 
        & \tfrac{1}{\slang}\snn{\tau}\,\imat_3 \\[2pt]
        - \slang \snn{\tau}\,\imat_3 & 
           \csn{\tau}\,\imat_3
        }
        \in \Spmat{6}
    \end{align}
    \end{small}
    The above leads to the same simplified solutions in Eq.\eqref{qp_sol_f_m1}. 
\end{itemize}
\end{small}

\vspace{1ex} 
\noindent  For Eq.\eqref{dqdp_nonsimp}-Eq.\eqref{dqp_Fcen_simp} above, note the following:
\begin{small}
\begin{itemize}
     \item  All of the equations are valid for any arbitrary
     central-forces of any kind (conservative or otherwise). Such forces never show up in the dynamics for \eq{(\tup{q},\tup{p})}. 
     \item The central-force dynamics and solutions for \eq{(\tup{q},\tup{p})} are fully decoupled from the radial motion coordinates \eq{(u,p_\ss{u})}. 
    \item We assumed \eq{s_{\zr}=\tau_{\zr}=0} for convenience. The solutions are equally valid with \eq{s} and \eq{\tau} replaced by \eq{\Delta s=s-s_{\zr}} and \eq{\Delta\tau=\tau-\tau_{\zr}}.
\end{itemize}
\end{small}

\paragraph*{Second-Order ODEs for $\tup{q}$ as a Linear Oscillator.} 
Even \textit{without} simplifying using \eq{q=1} and \eq{\tup{q}\cdot\tup{p}=0}, it is relatively straightforward to show that, for any arbitrary central-force potential \eq{V^{0}(r)},  the second-order equations for \eq{\pddt{\tup{q}}} or \eq{\rrng{\tup{q}}} are those of a perturbed harmonic oscillator given as follows (derived below, starting at Eq.\eqref{some_ds_rels}): 
\begin{small}
\begin{align} \label{sotired_nonsimp}
\boxed{\begin{array}{rllllll}
      \pddt{\tup{q}} \,+\, \slang^2 \tup{q}  &  =\,  q^2 u^{2n}\tup{f}
  \\[4pt]
       &  =\, q u^{3n}(\imat_3 - \htup{q} \otms  \htup{q}) \cdot \tup{F}
\end{array} 
\qquad\quad\fnsize{or,}\qquad\quad
\begin{array}{rllllll}
    \rrng{\tup{q}} \,+\, \tup{q}  
       &  =\,
         \tfrac{1}{\slang^2} q^2 u^{2n} (\imat_3 - \tfrac{1}{\slang} \rng{\tup{q}}\otms \tup{p}) \cdot \tup{f}
\\[4pt]
      &  =\, 
      \tfrac{1}{\slang^2} q u^{3n} ( \imat_3 - \htup{q}\otms\htup{q} -  \tfrac{1}{q^2}\rng{\tup{q}} \otms \rng{\tup{q}}
    ) \cdot \tup{F}
\end{array} }
\end{align}
\end{small}
where \eq{\rng{\tup{q}} =-\hdge{\htup{\slangup}}\cdot\tup{q}} and
where the right-hand-sides — which account for all conservative and nonconservative perturbations — may be expressed in terms of either \eq{\tup{f}=-\pderiv{V^{1}}{\tup{q}}+\tup{\alphaup}} or in terms of the  cartesian components \eq{\tup{F}=-\pderiv{V^{1}}{\tup{r}}+\tup{a}^\nc} (see Eq.\eqref{totalForce}). In the case these perturbations vanish — or even in the case that \eq{\tup{F}=F_r\htup{r}} such that \eq{\eq{\tup{f}=0}} — then \eq{\tup{\slangup}} are integrals of motion and \eq{\pddt{\tup{q}}+ \slang^2 \tup{q}=0} and \eq{\rrng{\tup{q}} + \tup{q} =0 }  
 describe a linear oscillator with natural frequency \eq{\slang=\slang_\zr} and \eq{1}, respectively. Either is easily solved for: 
 \begin{small}
 \begin{flalign} \label{qp_sol_m1}
 \begin{array}{cc}
      \fnsize{arbitrary}  \\
       \fnsize{central-force}\\
       \fnsize{dynamics} 
 \end{array}
 \qquad
\tup{f} = 0 
&& \Rightarrow &&
 \begin{array}{lclll}
         \tup{q}_{s}  \,=\,   \tup{q}_{\zr}\cos{\slang s}  +  \tfrac{1}{\slang}\pdt{\tup{q}}_{\zr}\sin{\slang s}
          &\quad =
          &\quad    \tup{q}_{\tau}  \,=\,  \tup{q}_{\zr}\csn{\tau}  +  \rng{\tup{q}}_{\zr}\snn{\tau}
\\[5pt]
     \pdt{\tup{q}}_{s}  \,=\,  -\slang\tup{q}_{\zr}\sin{\slang s}
            +  \pdt{\tup{q}}_{\zr}\cos{\slang s}    
             &\quad\neq &\quad
             \rng{\tup{q}}_{\tau}  \,=\,  -\tup{q}_{\zr}\snn{\tau}
            + \rng{\tup{q}}_{\zr}\csn{\tau} 
\end{array}
 &&
 \slang s = \tau
\quad
\end{flalign}
\end{small}
which agrees with the solutions given previously for \eq{(\tup{q},\tup{p})} in Eq.\eqref{E3sol_s_prj}-Eq.\eqref{qp_sol_f_m1}.

\begin{notesq}
    \textit{Simplifications.}
      Eq.\eqref{sotired_nonsimp} does \textit{not} require any simplifications using the integrals of motion \eq{q=1} or \eq{\tup{q}\cdot\tup{p}=0}. Yet, We are indeed free to make these simplifications, leading to \eq{\slang \simeq p}, \eq{\pdt{\tup{q}} \simeq \tup{p}}, and \eq{\rng{\tup{q}}\simeq \htup{p}}, such that the right-hand-sides are equivalent to:
     \begin{align} \label{sotired}
        \fnsz{\begin{array}{cc}
             q=1  \\
             \tup{q}\cdot\tup{p}=0 
        \end{array}}
        \;\; \Rightarrow \qquad
        \left\{ \qquad 
      \begin{array}{rllllll}
      \pddt{\tup{q}} \,+\, \slang^2 \tup{q}  &  \simeq\,   u^{2n}\tup{f}
        \\[4pt]
            &  \simeq\,
        u^{3n}(\imat_3 - \htup{q} \otms  \htup{q}) \cdot \tup{F}
     \\[6pt]
    \end{array} \right.
    \qquad\quad \fnsize{or,} \qquad\quad 
    \begin{array}{rllllll}
          \rrng{\tup{q}} \,+\, \tup{q}  
           &  \simeq\,
             \tfrac{u^{2n}}{\slang^2} (\imat_3 -  \rng{\tup{q}}\otms \rng{\tup{q}} ) \cdot \tup{f}
    \\[4pt]
         &  \simeq\, 
          \tfrac{ u^{3n}}{\slang^2}  ( \htup{\slangup}\otms\htup{\slangup} ) \cdot \tup{F}
    \end{array} 
    \end{align}
     where \eq{\htup{\slangup}\otms\htup{\slangup} \simeq   \imat_3 - \htup{q}\otms\htup{q} - \htup{p} \otms \htup{p} \simeq  \imat_3 - \htup{q}\otms\htup{q} - \rng{\tup{q}} \otms \rng{\tup{q}} }. We still write \eq{\slang} (rather than \eq{p}) simply to keep in mind this is the angular momentum.\footnote{The convenient relation that \eq{\slang=p} follows as a result of choosing \eq{m=-1} back in section \ref{sec:sum_m=-1}. It does not hold for other choices of \eq{m}.}
\end{notesq}

\begin{footnotesize}
\begin{itemize}
    \item[] \textit{Derivation of $\pddt{\tup{q}}$ in Eq.\eqref{sotired_nonsimp}.} 
    Consider the following relations shown previously  (Eq.\eqref{qp_eom_s_m1}, Eq.\eqref{ldot_cen}, and Eq.\eqref{totalForce}):
    \begin{align} \label{some_ds_rels}
    \pdt{\tup{q}} = -\hdge{\tup{\slangup}} \cdot \tup{q}
    &&,&&
        \hdge{\pdt{\tup{\slangup}}} = \pdt{t} \tup{q}\wdg\tup{f} 
    &&,&&
       \pdt{t}=u^{2n}
    &&,&&
        \htup{q}\cdot\tup{f} = 0
    \end{align}
    Differentiating the first of the above, and making subsequent use of the other relations, we obtain:
    \begin{align} \nonumber 
        \pddt{q}_i = -\slang_{ij}\pdt{q}_j - \pdt{\slang}_{ij}q_j 
        \,=\,  \slang_{ij}\slang_{jk}q_k - \pdt{t}(q_i f_j - q_j f_i) q_j 
        \,=\, -\slang^2 q_i +  q^2 u^{2n}(\kd_{ij} - \hat{q}_i\hat{q}_j)f_j
        \,=\,
         -\slang^2 q_i +  q^2 u^{2n} f_i 
    \end{align}
    where the last equality makes use of \eq{\hat{q}_j f_j = 0}, which is easily verified by the relation \eq{f_i = \tfrac{u^n}{q}(\kd_{ij} -\hat{q}_i\hat{q}_j)F_j}. 
    \item[] \textit{Derivation of $\rrng{\tup{q}}$ in Eq.\eqref{sotired_nonsimp}.} 
     Consider the following relations shown previously (Eq.\eqref{qp_eom_f_m1}, Eq.\eqref{ldot_cen}, and Eq.\eqref{totalForce}):
     \begin{align} \label{some_df_rels}
        \rng{\tup{q}} = -\tfrac{1}{\slang}\hdge{\tup{\slangup}} \cdot \tup{q}  &&,&&
            \hdge{\rng{\tup{\slangup}}} = \rng{t} \tup{q}\wdg\tup{f} 
         &&,&&
           \rng{\slang}= \rng{t}\tfrac{q^2}{\slang}\tup{p}\cdot\tup{f}
        &&,&&
           \rng{t}= u^{2n}/\slang
        &&,&&
            \htup{q}\cdot\tup{f} = 0
    \end{align}
     Differentiating the first of the above, and making subsequent use of the other relations, we obtain:
    \begin{align}
    \begin{array}{rllllll}
          \rrng{q}_i = -\tfrac{1}{\slang}\slang_{ij}\rng{q}_j - \tfrac{1}{\slang}\rng{\slang}_{ij}q_j 
        \,+\, \tfrac{1}{\slang^2}\rng{\slang} \slang_{ij}q_j
            &\,=\, 
            \tfrac{1}{\slang^2}\slang_{ij}\slang_{jk}q_k -  \tfrac{1}{\slang} \rng{t} (q_i f_j - q_j f_i) q_j  \,+\, \tfrac{1}{\slang^2}\rng{t} \tfrac{q^2}{\slang} p_k f_k \slang_{ij}q_j
            &\qquad (\rng{t} = u^{2n}/\slang)
        \\[3pt]
            &=\;
            -q_i + \tfrac{q^2 u^{2n}}{\slang^2} (\kd_{ij} - \hat{q}_i\hat{q}_j)f_j \,+\, \tfrac{q^2 u^{2n}}{\slang^4}\slang_{ij}q_j p_k f_k 
            &\qquad (\hat{q}_j f_j =0)
         \\[3pt]
            &=\;
            -q_i + \tfrac{q^2 u^{2n}}{\slang^2} (\kd_{ij} + \tfrac{1}{\slang^2}\slang_{ik}q_k p_j) f_j
            \,=\,
            -q_i + \tfrac{q^2 u^{2n}}{\slang^2} (\kd_{ij} - \tfrac{1}{\slang} \rng{q}_i p_j) f_j
    \end{array}
    \end{align}
    where \eq{\rng{q}_i = -\tfrac{1}{\slang}\slang_{ik}q_k}. To express the above using the cartesian components \eq{F_i}, we simply use \eq{f_i = \tfrac{u^n}{q}(\kd_{ij}-\hat{q}_i\hat{q}_j)F_j}.  
    Details in the footnote\footnote{It was shown that \eq{\rrng{\tup{q}} + \tup{q} = \tfrac{1}{\slang^2} q^2 u^{2n} (\imat_3 - \tfrac{1}{\slang} \rng{\tup{q}}\otms \tup{p}) \cdot \tup{f}}, where \eq{\rng{\tup{q}} =-\tfrac{1}{\slang}\hdge{\tup{\slangup}}\cdot\tup{q}}. Substituting \eq{\tup{f}  =  \tfrac{u^n}{q} (\imat_3  - \htup{q} \otms \htup{q} ) \cdot \tup{F}} leads to:
    \begin{align} \nonumber
    \begin{array}{rllllll}
          \rrng{\tup{q}} \,+\, \tup{q} 
       &\,=\,  \tfrac{1}{\slang^2} q^2 u^{2n} (\imat_3 - \tfrac{1}{\slang} \rng{\tup{q}}\otms \tup{p}) \cdot\tup{f} 
        \,=\,
        \tfrac{1}{\slang^2} q u^{3n} \big( \imat_3 - \htup{q}\otms\htup{q} -  \tfrac{1}{\slang} \rng{\tup{q}}\otms \tup{p} +  \tfrac{1}{\slang}(\htup{q}\cdot\tup{p}) \rng{\tup{q}}\otms \htup{q}
        \big) \cdot \tup{F}
        \\[3pt]
        &\,=\,  \tfrac{1}{\slang^2} q u^{3n} \big( \imat_3 - \htup{q}\otms\htup{q} -  \tfrac{1}{\slang}\rng{\tup{q}} \otms [ \tup{p} - (\htup{q}\cdot\tup{p})\htup{q} ]
        \big) \cdot \tup{F}
        \;\;=\,  \tfrac{1}{\slang^2} q u^{3n} \big( \imat_3 - \htup{q}\otms\htup{q} -  \tfrac{1}{\slang q^2}\rng{\tup{q}} \otms [ q^2 \tup{p} - (\tup{q}\cdot\tup{p})\tup{q} ]
        \big) \cdot \tup{F}
         \\[3pt]
        &\,=\,  \tfrac{1}{\slang^2} q u^{3n} \big( \imat_3 - \htup{q}\otms\htup{q} -  \tfrac{1}{\slang q^2}\rng{\tup{q}} \otms [ -\hdge{\tup{\slangup}}\cdot\tup{q} ]
        \big) \cdot \tup{F}
        \;\,=\,  \tfrac{1}{\slang^2} q u^{3n} \big( \imat_3 - \htup{q}\otms\htup{q} -  \tfrac{1}{q^2}\rng{\tup{q}} \otms \rng{\tup{q}}
        \big) \cdot \tup{F}
    \end{array}
    \end{align}
    } 
 \end{itemize}
\end{footnotesize}

\subsubsection{The radial motion (for Kepler and Manev dynamics)} 


The dynamics for \eq{(u,p_\ss{u})} are not so obviously linear as those of \eq{(\tup{q},\tup{p})} discussed above. In fact, the first-order dynamics:
\begin{small}
\begin{align} \label{upu_eom_again}
\begin{array}{llll}
      \pdt{u} \,=\,  \tfrac{u^2}{n^2}p_\ss{u}
\\[5pt]
        \pdt{p}_\ss{u} \,=\,
    \tfrac{n}{u} ( \slang^2  +\tfrac{n-1}{n^3}u^2 p_\ss{u}^2 )  
    - \pdt{t} \pderiv{V^{0}}{u}
    + \pdt{t} f_\ss{u}
\end{array}
    &&\fnsize{or,} &&
\begin{array}{llll}
      \rng{u} \,=\,  \tfrac{u^2}{n^2 \slang}p_\ss{u}
\\[5pt]
        \rng{p}_\ss{u} \,=\,
    \tfrac{n}{u\slang} ( \slang^2  +\tfrac{n-1}{n^3}u^2 p_\ss{u}^2 )  
    - \rng{t} \pderiv{V^{0}}{u}
    + \rng{t} f_\ss{u}
\end{array}   
\end{align}
\end{small}
(with \eq{\pdt{t}=u^{2n}} and \eq{\rng{t}=u^{2n}/\slang})
are not linear at all, even in the unperturbed case. 
However, for certain forms of \eq{V^0}, the above \textit{are} equivalent to perturbed linear \textit{second}-order ODEs for \eq{\pddt{u}} or \eq{\rrng{u}}. This will require finally choosing the value of \eq{n} to be \eq{n=-1}.

\paragraph*{Second-Order ODEs for $u$ as a Linear Oscillator.} 

We will show that, for \eq{V^0} the Kepler or Manev potential, the canonical equations of motion for \eq{(u,p_\ss{u})} are equivalent to second-order equations for \eq{u} as a harmonic oscillator with a constant driving ``force'' term. Yet, obtaining linear equations for \eq{u} is not quite as simple as it was for \eq{\tup{q}}. Unlike the equations of motion for \eq{\tup{q}}, the equation of motion for \eq{u} does depend on the central-force potential, \eq{V^{0}}.
From Eq.\eqref{upu_eom_again}, we see that:
\begin{small}
\begin{align}
    & \pddt{u} \,=\, \tfrac{2}{n^2}u\pdt{u}p_\ss{u} \,+\, \tfrac{1}{n^2}u^2\pdt{p}_\ss{u} \,=\, \tfrac{2}{n^4}u^3 p_\ss{u}^2 \,+\, \tfrac{1}{n^2}u^2\pdt{p}_\ss{u}
\\
   & \rrng{u} \,=\, \tfrac{2}{\slang n^2}u\rng{u}p_\ss{u} \,+\, \tfrac{u^2}{\slang n^2}\rng{p}_\ss{u} \,-\, \tfrac{u^2 p_\ss{u}}{n^2}\tfrac{\rng{\slang}}{\slang^2} \,=\, \tfrac{2u^3}{\slang^2 n^4}p_\ss{u}^2 \,+\, \tfrac{u^2}{\slang n^2}\rng{p}_\ss{u} \,-\, \tfrac{u^2 p_\ss{u}}{n^2}\tfrac{\rng{\slang}}{\slang^2}
\end{align}
\end{small}
Substitution of \eq{\pdt{p}_\ss{u}} and \eq{\rng{p}_\ss{u}} from Eq.\eqref{upu_eom_again},  and \eq{\rng{\slang}= \rng{t} \tfrac{q^2}{\slang}\tup{p} \cdot\tup{f}  } from Eq.\eqref{ldot_cen}, into the above 
leads to:\footnote{These relations may also be found from one another using \eq{ \diff{}{\tau} = \tfrac{1}{\slang}\diff{}{s}} and \eq{ \ddiff{}{\tau} =\tfrac{1}{\slang^2}(\ddiff{}{s}- \tfrac{\pdt{\slang}}{\slang} \diff{}{s})}. }
\begin{small}
\begin{align} \label{u_dds_long}
     \pddt{u} &\,=\, 
    \tfrac{1}{n}\big( \slang^2 + \tfrac{n+1}{n^3}u^2 p_\ss{u}^2  \,-\, \tfrac{1}{n}u^{2n+1}\pderiv{V^{0}}{u}\big)u \,+\,  \tfrac{u^{2n+2}}{n^2} f_\ss{u}
\\ \label{ddu_f_long}
   \rrng{u}
    &\,=\,  \tfrac{1}{n}\big( 1  +  \tfrac{n+1}{n^3} \tfrac{1}{\slang^2} u^2 p_\ss{u}^2 - \tfrac{u^{2n+1}}{n\slang^2} \pderiv{V^{0}}{u}   \big) u 
    \,+\, \tfrac{u^{2n+2}}{n^2 \slang^2} \big( f_\ss{u}
    - \tfrac{q^2}{\slang^2} p_\ss{u} \tup{p}\cdot\tup{f} \big)
\end{align}
\end{small}
Looking at Eq.\eqref{u_dds_long} and Eq.\eqref{ddu_f_long}, we note that neither is  the equation of a perturbed oscillator, regardless of the central-force potential \eq{V^{0}}. The first term in parenthesis — those not involving the perturbations  —  is not, in general, constant for unperturbed motion on account of  \eq{p_\ss{u}} (which is related to \eq{\dot{r}}).

\begin{notesq}
    For the above equations to have any hope of being linear in the unperturbed case requires the elimination of the \eq{p_\ss{u}} term in the coefficient of \eq{u}. This, in turn, requires \eq{n=-1} such that \eq{u=1/r}.
    Eq.\eqref{u_dds_long}  and Eq.\eqref{ddu_f_long} then simplify to
    \begin{small}
    \begin{align} \label{ddu_sf}
     \smsz{n=-1}\;\;  \left\{ \;\;
        \begin{array}{rllll}
              \pddt{u}   +  \slang^2 u   \,+\, \pderiv{V^{0}}{u} 
              &=\, f_\ss{u} 
              &=\; -\tfrac{1}{u^2}\htup{q}\cdot\tup{F}
     \\[6pt]
        \rrng{u} + u \,+\,  \tfrac{1}{\slang^2}\pderiv{V^{0}}{u}    
         &=\, \tfrac{1}{\slang^2}\big(f_\ss{u} -  \tfrac{q^2}{\slang u^2} \rng{u} \tup{p}\cdot \tup{f} \big)
       &=\;  -\tfrac{1}{q u^3 \slang^2} ( u\tup{q} +  \rng{u} \rng{\tup{q}} ) \cdot\tup{F}   
        \end{array} \right.
    \end{align}
    \end{small}
    where \eq{\rng{u} = u^2 p_\ss{u}/\slang} and \eq{\rng{\tup{q}}=-\tfrac{1}{\slang}\hdge{\tup{\slangup}}\cdot\tup{q}}, and where we have used \eq{\tup{p}\cdot\tup{f} = \tfrac{1}{u q^3}(\tup{q}\cdot\hdge{\tup{\slangup}})\cdot\tup{F} } from Eq.\eqref{totalForce}. 
\end{notesq}


\noindent Now, looking at the above, it can be seen that these equations simplify to those of a perturbed oscillator for the following cases of the central-force potential:
\begin{small}
\begin{align} \label{Vo_cases}
\begin{array}{rllllll}
     \fnsize{case 1 (Kepler):}  
     &  V^{0} = -\tfrac{\kconst_1}{r}  
     &=\;   -\kconst_1 u 
\\[6pt]
     \fnsize{case 2 (Manev):}
     &   V^{0} =  -\tfrac{\kconst_1}{r} - \tfrac{1}{2}\tfrac{\kconst_2}{r^2}  
     &=\;  -\kconst_1 u - \tfrac{1}{2}\kconst_2 u^2   
\end{array}
\end{align}
\end{small}
 for any scalar constants \eq{\kconst_1,\kconst_2 \in\mbb{R}}. The negative signs seen above are simply included such that the forces (per unit mass) resulting from the above potentials are attractive (directed towards the origin at \eq{\tup{r}= 0})  if \eq{\kconst_i>0} and are repulsive if   \eq{\kconst_i<0}.
 For the following developments, we shall consider the more general Manev potential (it includes Kepler as a special case).  Eq.\eqref{ddu_sf} then leads to:
 \begin{small}
\begin{flalign} \label{Vo_cases_s}
     \fnsz{\begin{array}{cc}
          \text{for } n=-1 \text{ and,}  \\[3pt]
          V^{0} = -\kconst_1u - \tfrac{1}{2}\kconst_2 u^2
     \end{array} }
 &&
\boxed{\begin{array}{ll}
            \pddt{u} \,+\, \omega^2 u \,-\, \kconst_1  \,=\, f_\ss{u} 
            &\,=\, -\tfrac{1}{u^2}\htup{q} \cdot \tup{F} 
    \\[6pt]
          \rrng{u}  \,+\, \varpi^2 u  -  \tfrac{\kconst_1}{\slang^2} 
            \,=\,  \tfrac{1}{\slang^2}\big(f_\ss{u} -  \tfrac{q^2}{\slang u^2} \rng{u} \tup{p}\cdot \tup{f} \big)
          &\,=\, -\tfrac{1}{q u^3 \slang^2} ( u\tup{q} +  \rng{u} \rng{\tup{q}} ) \cdot\tup{F}  
 \end{array} } 
&&
\begin{array}{llll}
           \omega^2 :=\, \slang^2-\kconst_2
 \\[6pt]
          \varpi^2 :=\,  \omega^2/\slang^2
\end{array} 
  \qquad
\end{flalign}
\end{small}
where the Kepler case corresponds to \eq{\kconst_2=0} and thus \eq{\omega = \slang} and \eq{\varpi =1}. 
We implicitly assume that \eq{\slang^2-\kconst_2>0} such that \eq{\omega} and \eq{\varpi} are real. For unperturbed motion, the right-hand-side of the above vanishes and \eq{\slang} is constant (thus, so too are \eq{\omega} and \eq{\varpi}).  For general perturbed motion, \eq{\slang} is no longer constant but evolves according to 
\eq{ \pdt{\slang} = \pdt{t} \tfrac{q^2}{\slang}\tup{p} \cdot\tup{f} = \tfrac{u^2 q^2}{\slang}\tup{p} \cdot\tup{f} }, as per Eq.\eqref{ldot_cen}. 


\begin{notesq}
    \textit{Simplifications.}
   Nothing in Eq.\eqref{upu_eom_again}-Eq.\eqref{Vo_cases_s} has required using \eq{q=1} or \eq{\lambda=\htup{q}\cdot\tup{p}=0}. We are, as usual, free to make these simplifications, which lead to  \eq{\slang^2 \simeq p^2} and \eq{\rng{\tup{q}} \simeq \tfrac{1}{\slang}\tup{p} \simeq \htup{p}}.  
    This does not lead to anything particularly illuminating for the ``\eq{u}-part'' of the dynamics (other than the obvious \eq{q=1} on the right-hand-side of the above ODEs). Still, we note, for example, that this leads to \eq{ \rrng{u} + \varpi^2 u  -  \tfrac{\kconst_1}{\slang^2}  \simeq  \tfrac{1}{u^2\slang^2} (u^2 f_\ss{u} -   \rng{u} \rng{\tup{q}}\cdot \tup{f} )}. 
\end{notesq}


\noindent 
Now, consider the \textit{un}perturbed Manev problem given by Eq.\eqref{Vo_cases_s} with \eq{\tup{F}=0=(\tup{f},f_\ss{u})}:  
\begin{small}
\begin{flalign} \label{ddu_0}
\begin{array}{cc}
      \fnsize{Manev-type}  \\
       \fnsize{dynamics} 
 \end{array}
&&
\begin{array}{lllll}
          \pddt{u}    \,+\, \omega^2u  \,-\, \kconst_1
          \,=\, 0 
 \qquad\;\;\fnsize{or,}\;\;\qquad 
          \rrng{u}  \,+\, \varpi^2u \,-\, \tfrac{\kconst_1}{\slang^2}
          \,=\, 0 
\end{array} 
&&
\end{flalign}
\end{small}
As is the case for any central-force dynamics, the angular momentum is preserved, \eq{\slang=\slang_\zr}, such that \eq{\omega^2:=\slang^2-\kconst_2} and \eq{\varpi:=\omega/\slang} are as well, leading to: 
\begin{small}
\begin{align}
       \slang = \slang_\zr
      \qquad \Rightarrow \qquad
      \omega=\omega_\zr \quad,\quad \varpi=\varpi_\zr 
      \quad,\quad 
      \tau = \slang s 
      \quad,\quad 
      \varpi \tau  =  \omega s 
\end{align}
\end{small}
where  we assume \eq{s_\zr=\tau_\zr=0} for convenience. 
The ODEs in Eq.\eqref{ddu_0} is then those of a linear harmonic oscillator with a constant driving ``force'' term. 
The solutions are easily found in terms of initial conditions and \eq{s} or \eq{\tau = \slang s}:
\begin{small}
\begin{flalign} \label{bokky}
\qquad
 \begin{array}{lclll}
      u_{s} \,=\, (u_{\zr}-\tfrac{\kconst_1}{\omega^2})\cos{\omega s}
    + \tfrac{1}{\omega}\pdt{u}_{\zr}\sin{\omega s} + \tfrac{\kconst_1}{\omega^2}
    &=&
     u_{\tau} \,=\, (u_{\zr}-\tfrac{\kconst_1}{\omega^2})\cos{\varpi \tau}
    + \tfrac{1}{\varpi}\rng{u}_{\zr}\sin{\varpi\tau} + \tfrac{\kconst_1}{\omega^2}
\\[5pt]
    \pdt{u}_{s} \,=\, -\omega(u_{\zr}-\tfrac{\kconst_1}{\omega^2})\sin{\omega s} + \pdt{u}_{\zr}\cos{\omega s}
    &\neq &
    \rng{u}_{\tau} \,=\, -\varpi(u_{\zr}-\tfrac{\kconst_1}{\omega^2})\sin{\varpi \tau} + \rng{u}_{\zr}\cos{\varpi \tau}
\end{array} 
&&,&&
\begin{array}{lll}
     \tau = \slang s  
     \\[2pt]
     \varpi \tau = \omega s 
\end{array}
\quad
\end{flalign}
\end{small}
Then, using \eq{\pdt{u} = u^2 p_\ss{u}} or \eq{\rng{u} = \pdt{u}/\slang = u^2 p_\ss{u}/\slang}, we obtain the solutions for \eq{(u_s,p_{u_s}) = (u_\tau,p_{u_\tau}) }:
\begin{small}
\begin{align} \label{up_sol_f}
\begin{array}{llllll}
  u_{\tau} \,=\, (u_{\zr}-\tfrac{\kconst_1}{\omega^2})\cos{\varpi \tau}
    \,+\, \tfrac{1}{\omega}u_{\zr}^2 p_{u_{\zr}}\sin{\varpi \tau} \,+\, \tfrac{\kconst_1}{\omega^2}
\qquad,\qquad 
      p_{u_\tau}  \,=\, \dfrac{ -\omega(u_{\zr}- \frac{\kconst_1}{\omega^2})\sin{\varpi \tau} \,+\, u_{\zr}^2 p_{u_{\zr}}\cos{\varpi \tau} }
           {\big[ (u_{\zr}-\frac{\kconst_1}{\omega^2})\cos{\varpi \tau}
    \,+\, \tfrac{1}{\omega}u_{\zr}^2 p_{u_{\zr}}\sin{\varpi \tau} \,+\, \frac{\kconst_1}{\omega^2} \big]^2}
\end{array}  
\end{align}
\end{small}
where Kepler-type dynamics are recovered with \eq{\kconst_2=0} such that \eq{\omega=\slang} and \eq{\varpi = 1}. 

%% file: Mysecs_prj/new_3_2BP.tex
\section{Kepler \& Manev dynamics in projective coordinates}  \label{sec:2BP}

Much of the following is a direct result of setting 
\eq{n=-1} in section \ref{sec:central force} (equivalently, 
\eq{m=n=-1} in Appendix \ref{sec:2bp_gen}) together with considering a particular central-force potential of the Manev-type, \eq{V^{0}(r)=-\kconst_1/r - \tfrac{1}{2}\kconst_2/r^2}, where Kepler-type problems are easily recovered with \eq{\kconst_2=0}. As such, some the following will be a bit repetitive. However, we have so far been building our projective coordinate transformation in a piece-wise fashion, and have yet to explicitly and concisely state the full regularizing transformation — for chosen values of \eq{n=m=-1} — that we prefer for both Kepler-type and Manev-type dynamics. That is the purpose of this section. 
We also present the \eq{J_2}-perturbed Kepler problem (i.e., the ``main'' satellite problem) formulated using our preferred projective transformation.

\subsection{A preferred projective transformation}

\paragraph*{The Transformation.}
It was shown in section \ref{sec:2BP_lin} that, in order to obtain linear equations for \eq{u=r^\ss{1/n}} using the evolution parameters \eq{s} or \eq{\tau} defined in Eq.\eqref{dtds_m1}, it is required to choose \eq{n=-1} in the point transformation given by Eq.\eqref{sum_rv_m}. Therefore, we consider a canonically-extended projective transformation starting with a point transformation
\eq{\tup{r} =  \tfrac{1}{u q}\tup{q}=  \tfrac{1}{u}\htup{q}}, subject to \eq{\varphi(\tup{q}) = q - 1 = 0}. This is extended to the momentum level for a canonical transformation \eq{\Gam:\mbb{R}^8\to\mbb{R}^6} given by:
\begin{small}
\begin{align}\label{sum_PT_n}
   \Gam: (\bartup{q},\bartup{p}) 
    \mapsto (\tup{r},\tup{v})
\;\; \left\{ \;\;
\boxed{\begin{array}{lllll}
      \tup{r} 
      \,=\, \tfrac{1}{u} \htup{q}
  \\[5pt]
     \tup{v} 
     \,=\,  u q ( \imat_3 - \htup{q}\otms\htup{q})\cdot\tup{p} - u^2 p_\ss{u}\htup{q} 
     \;=\,  -u (\tup{q}\wdg\tup{p})\cdot\htup{q} - u^2 p_\ss{u}\htup{q} 
\end{array} } \right. 
&& 
\begin{array}{lllll}
        \varphi(\tup{q}) = q-1 = 0
\\[5pt]
    \lambda  \,=\, \htup{q}\cdot\tup{p}
\end{array}
\end{align}
\end{small}
with \eq{\hdge{\tup{\slangup}}=\tup{q}\wdg\tup{p}} and
where, as per section \ref{sec:Pxform} and Appendix \ref{sec:2bp_gen}, the constraint \eq{\varphi} and associated Lagrange multiplier \eq{\lambda} are built into the momenta coordinate transformation.
We recall the following:
\begin{small}
\begin{itemize}
    \item \textit{Angular momentum.} The transformation in Eq.\eqref{sum_PT_n} preserves the form of the  angular momentum:
    \begin{small}
    \begin{flalign} \label{angMom_million}
    &&
            \hdge{\tup{\slangup}} 
        \,=\, \tup{r}\wdg\tup{v} \,=\, \tup{q}\wdg\tup{p}
     &&,&&  
        \tup{\slangup} 
        \,=\, \tup{r}\tms \tup{v} \,=\, \tup{q}\tms \tup{p}  
     &&,&& 
     \slang^2  \,=\, r^2 \v^2 - (\tup{r}\cdot\tup{v})^2 \,=\, q^2 p^2 - (\tup{q}\cdot\tup{p})^2
     &&
    \end{flalign}
    \end{small}
    where \eq{\slang_{ij} = \lc_{ijk} \slang_k \leftrightarrow \slang_i =  \tfrac{1}{2} \lc_{ijk}\slang_{jk} } are the Hodge duals of one another.\footnote{See Eq.\eqref{hodge_cord} for the Hodge dual on \eq{\mbb{R}^3}. In components, Eq.\eqref{angMom_million} reads: 
    \eq{\;\; \slang_{ij}  = r_i \v_j - \v_i r_j = q_i p_j - p_i q_j\;} and \eq{\;\slang_i = \lc_{ijk}r_j \v_k = \lc_{ijk}q_j p_k}.}
    \item \textit{The Hamiltonian.} A cartesian coordinate mechanical Hamiltonian, \eq{\mscr{K}=\tfrac{1}{2}\v^2 + V} — where we assume the potential, \eq{V}, has the general form \eq{V=V^0(r)+V^1(\tup{r},t)} — is taken by the transformation in Eq.\eqref{sum_PT_n} to a projective coordinate Hamiltonian given by \eq{\mscr{H}:= \mscr{K}\circ\Gam}. Direct substitution leads to:
    \begin{small}
    \begin{align} \label{Hxform_again}
        \mscr{K} = \tfrac{1}{2}\v^2 + V^0(r)+ V^1(\tup{r},t)
        \qquad \Rightarrow \qquad 
         \mscr{H} = \tfrac{1}{2}u^2\big( \slang^2  + u^2 p_u^2 \big) + V^0(u) + V^1(\bartup{q},t)
    \end{align}
    \end{small}
    where \eq{V^0(u)} and \eq{V^1(\bartup{q},t)} simply indicate the original  \eq{V^0(r)} and \eq{V^1(\tup{r},t)}
    rewritten using \eq{\tup{r} = \tfrac{1}{u} \htup{q}} and, thus, \eq{r = 1/u}. 
\end{itemize}
\end{small}

\noindent Now, in order to invert the transformation in Eq.\eqref{sum_PT_n} — which is a submersion with \eq{\rnk \pderiv{(\tup{r},\tup{v})}{(\bartup{q},\bartup{p})}=6} and therefore has no unique inverse — we use the fact that the transformed Hamiltonian, \eq{\mscr{H}=\mscr{K}\circ\Gam}, admits two ``extra'' integrals of motion, 
 \eq{q=\mag{\tup{q}}} and \eq{\lambda =  \htup{q}\cdot\tup{p}} (this was discussed in  section \ref{sec:central force} and Appendix \ref{sec:2bp_gen}): 
 \begin{small}
\begin{flalign} \label{qlam_2BP}
     \begin{array}{cc}
     \fnsize{integrals}  \\
     \fnsize{of motion} 
\end{array}\!\!:
&&
\begin{array}{llll}
    \dot{q} = \pbrak{q}{\mscr{H}} + \pderiv{q}{\bartup{p}} \cdot \bartup{\alphaup}  = 0 
   \\[5pt]
     \dot{\lambda}  = \pbrak{\lambda}{\mscr{H}} + \pderiv{\lambda}{\bartup{p}} \cdot \bartup{\alphaup}  = 0
\end{array} 
&&  \xRightarrow{\text{choose}} &&
\begin{array}{llll}
     q = q_\zr  = 1
     \\[5pt]
       \lambda = \lambda_\zr  = 0
       \;= \htup{q}\cdot\tup{p}
\end{array}
&&
\end{flalign}
\end{small}
The above holds for any potential function (conservative forces) and any generalized nonconservative forces \eq{\bartup{\alphaup}=(\tup{\alphaup},\alpha_\ss{u})}. 
We may limit consideration to the values  \eq{q=1} and \eq{\lambda=\htup{q}\cdot\tup{p}=0} to obtain an inverse of Eq.\eqref{sum_PT_n}:
\begin{small}
\begin{align}  \label{sum_PT_n_inv}
\begin{array}{cc}
     \fnsize{restrict to}  \\
     \fnsize{$q=1$} \\
     \fnsize{$\lambda=\htup{q}\cdot\tup{p} =0$}
\end{array} 
\quad \Rightarrow \qquad
\inv{\Gam}:(\tup{r},\tup{v}) \mapsto (\bartup{q},\bartup{p})
\;\; \left\{ \;\;
\boxed{ \begin{array}{ll} 
      \tup{q} \,=\,  \htup{r}  &,
 \\[5pt]
      u \,=\, 1/r &,
\end{array}  
\;\;
\begin{array}{ll}
     \tup{p}  \,=\,  r( \imat_3 - \htup{r}\otms \htup{r}) \cdot \tup{v} \,=\, -\hdge{\tup{\slangup}} \cdot \htup{r} 
 \\[5pt]
      p_u \,=\, - r^2 \htup{r} \cdot\tup{v} 
\end{array} } \right.
\begin{array}{llllll}
     \,=\, r^2\dot{\htup{r}} = \pdt{\htup{r}}
\\[5pt]
     =\, 
    -r^2\dot{r} = -\pdt{r}
\end{array}
\end{align}
\end{small}
with \eq{\hdge{\tup{\slangup}}=\tup{r}\wdg\tup{v}}. 
So long as \eq{(\bartup{q}_\zr,\bartup{p}_\zr)} are obtained from given initial conditions \eq{(\tup{r}_\zr,\tup{v}_\zr)} using the above, then it automatically holds along any phase space solution that \eq{q=q_\zr=1} and \eq{\lambda=\lambda_\zr =0=\htup{q}\cdot\tup{p}}.
This places no restrictions on the cartesian coordinates \eq{(\tup{r},\tup{v})}.

\begin{small}
\begin{itemize}
    \item  We often use ``\eq{\simeq}'' to indicate relations which have been simplified using the integrals of motion \eq{q=1} and \eq{\lambda=\htup{q}\cdot\tup{p}=0}; this leads to the relations given in Eq.\eqref{qpl_rels}-Eq.\eqref{qpl_rels2}.  
    \item \textit{Relation to LVLH basis.} 
    As per Remark \ref{rem:LVLH}, the  coordinates \eq{(\tup{q},\tup{p})} directly define inertial cartesian components of the orthonormal LVLH basis, \eq{\htup{t}_i}. 
    Specifically, \eq{ \{\htup{t}_r,\htup{t}_\tau, \htup{t}_\slang \} =  \{\htup{q},-\hdge{\htup{\slangup}}\cdot\htup{q},\htup{\slangup}\} \simeq \{\tup{q},\htup{p},\htup{\slangup}\}}, with \eq{\htup{\slangup}=\htup{q}\tms\htup{p}},  \eq{q^2\simeq 1}, and \eq{\slang^2\simeq p^2}.  
    \item \textit{Recovering cartesian coordinate numerical solutions.}
    For any solution curve expression in projective coordinates, 
     the corresponding solution curve for the inertial cartesian coordinates, \eq{(\tup{r},\tup{v})}, may always recovered using   Eq.\eqref{sum_PT_n}. Yet, when recovering numerical solutions, we are free to simplify these relations using the integrals of motion \eq{q=1} and \eq{\lambda=\htup{q}\cdot\tup{p}=0} such that Eq.\eqref{sum_PT_n} simplifies to:
    \begin{small}
    \begin{flalign} \label{sum_rv_sol}
    \text{Eq.\eqref{sum_PT_n}}
    \quad 
    \xRightarrow[\htup{q}\cdot\tup{p}\,=\,0]{q\,=\,1}
    \qquad 
       \tup{r} \,\simeq\, \tfrac{1}{u}\tup{q} 
       \quad, \quad
        \tup{v} \,\simeq\,  u\tup{p} - u^2p_u\tup{q}
    \end{flalign}
    \end{small}
     It should be stressed that the above is \textit{not} the projective transformation that is used to construct the Hamiltonian, transform the forces, etc. 
     The actual  projective transformation, \eq{\Gam:(\bartup{q},\bartup{p})\mapsto (\tup{r},\tup{v})}, is given in Eq.\eqref{sum_rv_m}. 
     Yet, when restricting \eq{(\bartup{q},\bartup{p})} solutions to those satisfying \eq{q=q_\zr=1} and \eq{\htup{q}\cdot\htup{p}=\htup{q}_\zr\cdot\htup{p}_\zr=0} — which we are free to do as per Remark \ref{rem:qp=0} —  then the above is \textit{numerically} equivalent to the full transformation in  Eq.\eqref{sum_rv_m}.
\end{itemize}
\end{small}

\paragraph*{Kepler-Manev Orbital Dynamics.}
For the Hamiltonian given in Eq.\eqref{Hxform_again}, we consider the specific case that the central-force potential, \eq{V^0}, is a Manev-type potential (as per section \ref{sec:2BP_lin}) for some scalars \eq{\kconst_1,\kconst_2 \in\mbb{R}}:
\begin{small}
\begin{flalign}
\begin{array}{cc}
     \fnsize{Manev-type} \\
     \fnsize{Potential} 
\end{array}
&&
    V^{0} \,=\,   -{\kconst_1}/{r} - \tfrac{1}{2}{\kconst_2}/{r^2} 
   \;\;= \;\;   -\kconst_1 u - \tfrac{1}{2} \kconst_2 u^2 
&&
\end{flalign}
\end{small}
where the Kepler case is recovered with \eq{\kconst_2=0}.
For the above potential, the Hamiltonian and equations of motion for the projective coordinates \eq{(\bartup{q},\bartup{p})} is then given by:
\begin{small}\begin{gather} \label{Hqp_manev}
    \boxed{ \mscr{H} \,=\,
     \tfrac{1}{2}u^2\big( \slang^2  + u^2 p_u^2 \big) 
     - \kconst_1 u - \tfrac{1}{2}\kconst_2 u^2 +\, V^{1} }
 \\ \nonumber 
 \begin{array}{llll}
     \dot{\tup{q}} = \pderiv{\mscr{H}}{\tup{p}} \,=\, - u^2 \hdge{\tup{\slangup}}\cdot\tup{q} 
     &,\qquad
      \dot{\tup{p}} = - \pderiv{\mscr{H}}{\tup{q}} + \tup{\alphaup} \,=\, - u^2 \hdge{\tup{\slangup}}\cdot\tup{p} + \tup{f}
\\[4pt] 
    \dot{u} = \pderiv{\mscr{H}}{p_u} \,=\, u^4p_u 
    &,\qquad
    \dot{p}_u  = - \pderiv{\mscr{H}}{u} + \alpha_\ss{u} \,=\,
     -u\big( \slang^2 + 2u^2 p_u^2 )  
    +  \kconst_1 + \kconst_2 u + f_\ss{u}
\end{array}
\end{gather}\end{small}
where \eq{V^{1}= V^{1}(\bartup{q},t)} accounts for arbitrary conservative perturbations, and
where \eq{(\tup{\alphaup},\alpha_\ss{u})} are generalized nonconservative perturbing forces. In the above ODEs, note we have combined the conservative and nonconservative perturbations into a single total generalized perturbing force,
\eq{(\tup{f},f_\ss{u})}, which may be written in terms of the cartesian components of the total perturbing force \eq{\tup{F} := -\pderiv{V^1}{\tup{r}} + \tup{a}^\nc} 
as follows:\footnote{This was shown around Eq.\eqref{totalForce}.}
\begin{small}
\begin{align} \label{sum_dV_2bp}
   \boxed{ 
\begin{array}{lllll}
         \tup{f} :=\, -\pderiv{V^{1}}{\tup{q}} + \tup{\alphaup} 
    & =\,
    \trn{\pderiv{\tup{r}}{\tup{q}}} \cdot  \tup{F}
       \,=\,  
       \tfrac{1}{u q} ( \imat_3 -  \htup{q}\otms\htup{q} ) \cdot \tup{F}
\\[4pt]
          f_\ss{u} :=\, -\pderiv{V^{1}}{u} + \alpha_\ss{u}
           & =\,
    \pderiv{\tup{r}}{u} \cdot  \tup{F}
    \,=\,  -\tfrac{1}{u^2}\htup{q}\cdot\tup{F}
\end{array}
 }
&&
\begin{array}{llll}
     \tup{q}\cdot\tup{f} = 0 
\\[3pt]
     \tup{p} \cdot \tup{f} 
     =
    \tfrac{1}{q}\htup{q}\cdot\hdge{\tup{\slangup}}\cdot\tup{f}
     = \tfrac{1}{u q^2}\htup{q}\cdot\hdge{\tup{\slangup}} \cdot \tup{F}
     \,\simeq  \tfrac{1}{u}\tup{p}\cdot\tup{F}
 \\[3pt]
     \tup{\slangup} \cdot \tup{f} = \tfrac{1}{uq}\tup{\slangup}\cdot\tup{F}
\end{array}
\end{align}
\end{small}

\noindent We also note that the ODEs  governing the evolution of the angular momentum functions were given in Eq.\eqref{ldot_cen} which, now with \eq{n=-1}, leads to:
\begin{align} \label{ldot_cen_again}
\begin{array}{llll}
     \hdge{\dot{\tup{\slangup}}} = \tup{q}\wdg\tup{f}
     \,=\, \tfrac{1}{u} \htup{q} \wdg \tup{F}
     &\;\;(=\tup{r}\wdg \tup{F})
 \\[4pt]
     \dot{\tup{\slangup}} \,=\, \tup{q}\tms\tup{f}
    \,=\, \tfrac{1}{u} \htup{q}\tms\tup{F}
    &\;\;(= \tup{r}\tms\tup{F})
\end{array}
&&,&&
\begin{array}{lllllll}
     \dot{\slang} &=\,  \tfrac{q^2}{\slang}\tup{p} \cdot\tup{f}  \,=\, \tfrac{1}{u} \htup{q}\cdot\hdge{\htup{\slangup}} \cdot\tup{F}
      &\;\;(= \tup{r}\cdot\hdge{\htup{\slangup}} \cdot\tup{F})
 \\[3pt]
    &\simeq\, \dot{p} \,=\, \htup{p} \cdot\tup{f} \simeq  \tfrac{1}{u} \htup{p} \cdot \tup{F}
\end{array}
\end{align}
where only the relations \eq{\dot{\slang} \simeq \dot{p}\simeq  \tfrac{1}{u} \htup{p} \cdot \tup{F}} make use of \eq{q=1} and \eq{\lambda=\htup{q}\cdot\tup{p}=0}. 


\paragraph*{Transformation of the Evolution Parameter.}
As detailed in section \ref{sec:ext_cen}, 
we then transform the evolution parameter from time, \eq{t}, to two new parameters, \eq{s} and \eq{\tau}, with respective extended Hamiltonians \eq{\wtscr{H}} and \eq{\whscr{H}}:
\begin{small}
\begin{align}
   \mrm{d} t \,=\, r^2  \mrm{d} s \,=\,  \tfrac{1}{u^2}  \mrm{d} s
  \quad,\quad  \wtscr{H} = \tfrac{1}{u^2}(\mscr{H} + p_t) 
&&,&&
   \mrm{d} t \,=\, \tfrac{r^2}{\slang} \mrm{d} \tau \,=\, \tfrac{1}{u^2\slang} \mrm{d} \tau
    \quad,\quad \whscr{H} = \tfrac{1}{\slang}\wtscr{H}
\end{align}
\end{small}
where we define the short-hand \eq{\pdt{\square}:=\diff{}{s}} and \eq{\rng{\square}:=\diff{}{\tau}}. Note \eq{\tau} is the true anomaly up to an additive constant. We begin with the \eq{s}-parameterized dynamics;  this was given in Eq.\eqref{qp_eom_s_m1} which, now with \eq{n=-1}, leads to:
\begin{small}\begin{gather} \nonumber
\wtscr{H}  \,=\,  \tfrac{1}{2}\big( \slang^2  + u^2 p_u^2 \big) 
     - \tfrac{\kconst_1}{u} - \tfrac{1}{2}\kconst_2  +  \tfrac{1}{u^2} V^{1} + \tfrac{1}{u^2} p_t 
\\ \label{qp_eom_s_sum}
\begin{array}{llll}
     \pdt{\tup{q}} \,=\, -  \hdge{\tup{\slangup}}\cdot\tup{q} 
     &,\qquad
     \pdt{u}   \,=\,  u^2 p_u &,
 \\[4pt]
      \pdt{\tup{p}}  \,=\, -\hdge{\tup{\slangup}}\cdot\tup{p}  + \tfrac{1}{u^2} \tup{f}  
      &,\qquad
       \pdt{p}_u \,=\, -\tfrac{1}{u}\big( \slang^2  + 2u^2 p_u^2 )  
        +  \tfrac{\kconst_1}{u^2} + \tfrac{\kconst_2}{u} + \tfrac{1}{u^2} f_\ss{u} &,
\end{array}
   \qquad
\begin{array}{llll}
      \pdt{t} \,=\, 1/u^2 
 \\[4pt]
     \pdt{p}_t    \,=\,  -( \pdt{t}\pderiv{V^{1}}{t} + \tup{\alphaup}\cdot\pdt{\tup{q}} + \alpha_\ss{u}\pdt{u})
\end{array} 
\end{gather}\end{small}
where the above equation for \eq{\pdt{p}_\ss{u}} has made use of the relation \eq{p_t=-\mscr{H}}.\footnote{Direct differentiation of \eq{\wtscr{H}} leads to \eq{\pdt{p}_\ss{u}} given as follows (which is equivalent to \eq{\pdt{p}_\ss{u}} in Eq.\eqref{qp_eom_s_sum} using \eq{p_t=-\mscr{H}}):
\begin{align} \nonumber
   \pdt{p}_\ss{u} \,=\, -\pd_\ss{u} \wtscr{H} + \pdt{t}\alpha_\ss{u} \,=\,  -u p_\ss{u}^2     -  \tfrac{\kconst_1}{u^2} -   \pd_\ss{u} (\tfrac{1}{u^2}V^{1})   +   \tfrac{2}{u^3}p_t + \tfrac{1}{u^2}\alpha_\ss{u}
   \;=\;
    -u p_\ss{u}^2     -  \tfrac{\kconst_1}{u^2} +   \tfrac{2}{u^3}(V^{1}  +  p_t) + \tfrac{1}{u^2}f_\ss{u}   
    \;=\;
     - \tfrac{1}{u}(\slang^2 +  2 u^2 p_\ss{u}^2 ) 
     + \tfrac{\kconst_1}{u^2} + \tfrac{\kconst_2}{u}
      + \tfrac{2}{u} \wtscr{H}  + \tfrac{1}{u^2} f_\ss{u}
\end{align}
 }

The \eq{\tau}-parameterized dynamics (where \eq{\mrm{d}\tau = \slang \mrm{d} s}) are obtained similarly, but using an extended Hamiltonian \eq{\whscr{H}=\tfrac{1}{\slang}\wtscr{H}}. 
The result is equivalent to a scaling by \eq{1/\slang} of the \eq{s}-parameterized ODEs in Eq.\eqref{qp_eom_s_sum} (i.e., \eq{ \rng{\square} = \tfrac{1}{\slang} \pdt{\square}}):
\begin{small}
\begin{align} \label{qp_eom_ta_sum}
    \begin{array}{lllll}
        \rng{\tup{q}} = -\hdge{\htup{\slangup}}\cdot\tup{q} &,
        \\[3pt]
        \rng{\tup{p}} = -\hdge{\htup{\slangup}}\cdot\tup{p} + \tfrac{1}{\slang u^2}\tup{f} &,
    \end{array}
    \qquad 
    \begin{array}{llll}
         \rng{u}   \,=\,  u^2 p_u/\slang &,
    \\[3pt]
     \rng{p}_\ss{u} \,=\, \pdt{p}_\ss{u} / \slang &, 
    \end{array}
    \qquad 
    \begin{array}{llll}
         \rng{t}   \,=\,  1/(\slang u^2)
    \\[3pt]
     \rng{p}_t \,=\, \pdt{p}_t / \slang 
    \end{array}
\end{align}
\end{small}
where the above form of the ODEs follow from making use of \eq{p_t=-\mscr{H}}.\footnote{Direct differentiation of the extended Hamiltonian \eq{\whscr{H}=\tfrac{1}{\slang}\wtscr{H}} leads to \eq{\tau}-parameterized Hamiltonian ODEs given as follows (cf.~Eq.\eqref{dqp_TA_full}):  
    \begin{gather}  \nonumber 
    \begin{array}{lllll}
          \rng{\tup{q}} 
          \,=\, 
          -(1 - \tfrac{\whscr{H}}{\slang})  \hdge{\htup{\slangup}}\cdot\tup{q} 
          &,
    \\[3pt]
        \rng{\tup{p}} 
          \,=\,  -(1- \tfrac{\whscr{H}}{\slang} ) \hdge{\htup{\slangup}}\cdot\tup{p}  
          +\rng{t}\tup{f} &,
    \end{array}
    \qquad 
    \begin{array}{lllll}
        \rng{u} 
          \,=\,   \tfrac{u^2}{\slang}p_\ss{u}
          &,
    \\[3pt]
         \rng{p}_\ss{u} 
          \,=\, -\tfrac{1}{\slang u} ( \slang^2  + 2u^2 p_\ss{u}^2 )   +  \rng{t} (\kconst_1 + \kconst_2 u)   + \tfrac{2}{u} \whscr{H}  +  \rng{t} f_\ss{u} &,
    \end{array}
     \qquad 
    \begin{array}{lllll}
          \rng{t}  
        \,=\, 1/(\slang u^2)
    \\[3pt]
         \rng{p}_t 
          \,=\, -( \rng{t}\pderiv{V^1}{t} + \tup{\alphaup}\cdot\rng{\tup{q}} + \alpha_\ss{u} \rng{u} )
    \end{array}
    \end{gather}
    which, after use of \eq{p_t=-\mscr{H}} (i.e., \eq{\whscr{H}=\wtscr{H}=0}), is equivalent to Eq.\eqref{qp_eom_ta_sum}. 
    }

\begin{small}
\begin{itemize}
      \item \textit{Simplifications.} Dynamics in our projective coordinates can always be simplified with \eq{q=1} and \eq{\lambda=\htup{q}\cdot\tup{p}=0} (implying the relations in Eq.\eqref{qpl_rels}). For instance, the above  dynamics for \eq{(\tup{q},\tup{p})} simplify to:
    \begin{small}
    \begin{flalign} \label{qp_eom_s_simp_sum}  
    && \fnsz{\begin{array}{cc}
         q=1  \\
         \htup{q}\cdot\tup{p}=0 
    \end{array}}
    \Rightarrow \quad
    \left\{\qquad\begin{array}{llll}
         \pdt{\tup{q}} \simeq  \tup{p}  
         \\[4pt]
         \pdt{\tup{p}}  
         \simeq -\slang^2 \tup{q}  + \tfrac{1}{u^2} \tup{f} 
    \end{array} \right.
    \qquad\fnsize{or,}\qquad
    \begin{array}{llll}
         \rng{\tup{q}} \simeq \tfrac{1}{\slang} \tup{p} 
            \,\simeq \htup{p} 
         \\[4pt]
         \rng{\tup{p}} \simeq -\slang \tup{q} + \tfrac{1}{\slang u^2}\tup{f}
    \end{array}
    &&,&&
    \slang^2 \simeq p^2 
    \quad
    \end{flalign}
    \end{small}
    \item For pure central-force dynamics (i.e., the case \eq{\tup{f}=0}) then the \eq{s} or \eq{\tau}-parameterized ODEs for \eq{(\tup{q},\tup{p})} are already linear whether or not one chooses to simply with \eq{q=1} and \eq{\lambda=\htup{q}\cdot\tup{p}=0}. 
\end{itemize}
\end{small}



\subsection{Linear Manev dynamics \& closed-form solutions} \label{sec:2bp_sols}


\noindent We are considering Manev-type dynamics (including Kepler as a special case),  which are described by a  cartesian coordinate Hamiltonian, \eq{\mscr{K}}, that is transformed via Eq.\eqref{sum_PT_n} to a projective coordinate Hamiltonian \eq{\mscr{H}}: 
\begin{small}
\begin{align}
     \mscr{K} \,=\, \tfrac{1}{2} \v^2 - \tfrac{\kconst_1}{r} - \tfrac{1}{2} \tfrac{\kconst_2}{r^2} + V^1
     \qquad \Rightarrow  \qquad
     \mscr{H} \,=\,
     \tfrac{1}{2}u^2\big( \slang^2  + u^2 p_u^2 \big) 
     - \kconst_1 u - \tfrac{1}{2}\kconst_2 u^2 + V^{1} 
\end{align}
\end{small}
where \eq{\slang^2=q^2 p^2 - (\tup{q}\cdot\tup{p})^2} is the angular momentum and \eq{V^{1}} is arbitrary. 
We have noted that the \eq{s}- or \eq{\tau}-parameterized dynamics for \eq{(\tup{q},\tup{p})} are linear for \textit{any} arbitrary central-forces — this was shown in section \ref{sec:2BP_lin} and is also again apparent from Eq.\eqref{qp_eom_s_sum} or Eq.\eqref{qp_eom_s_simp_sum} above. 
In contrast, while the \eq{s}- or \eq{\tau}-parameterized dynamics for \eq{(u,p_\ss{u})} are linear 
(in a sense\footnote{Specifically, they are equivalent to a linear \textit{second}-order ODE for \eq{u}. Yet, the \textit{first}-order dynamics for \eq{(u,p_\ss{u})} are not linear, even for pure Kepler-type or Manev-type dynamics.})
for the \textit{specific} cases of the Kepler or Manev potential, this linear nature is  not obvious from direct inspection of the first-order ODEs for \eq{(u,p_\ss{u})} in Eq.\eqref{qp_eom_s_sum}. It is more easily realized if, in place of \eq{p_\ss{u}}, we introduce a new quasi-momenta coordinate, \eq{w:=  u^2 p_\ss{u}}. That is,
\begin{small}
\begin{align} \label{wdef_2bp}
   \boxed{ w \,:=\,  u^2 p_\ss{u} \,=\, \pdt{u} 
    \,=\,  -\dot{r}
    \qquad \leftrightarrow \qquad
    p_\ss{u} \,=\, w/u^2  \,=\, -r^2 \dot{r}  \,=\,  -\pdt{r} }
\end{align}
\end{small}
While \eq{(u,w)} is \textit{not} a conjugate pair, many relations are ``nicer'' when expressed using \eq{w} in place of \eq{p_\ss{u}}.
E.g., the \eq{s}-  or \eq{\tau}-parameterized perturbed Manev dynamics in Eq.\eqref{qp_eom_s_sum} or Eq.\eqref{qp_eom_ta_sum} then lead to the below perturbed linear system:
\begin{small}
\begin{align} \label{qpw_eom_s_sum}
\begin{array}{lllll}
     \pdt{\tup{q}} \,=\, -  \hdge{\tup{\slangup}}\cdot\tup{q} 
     &,
 \\[4pt]
     \pdt{\tup{p}}  \,=\, -\hdge{\tup{\slangup}}\cdot\tup{p}  + \tfrac{1}{u^2} \tup{f} 
     &,
\end{array}
\;\;
\begin{array}{llll}
     \pdt{u}   \,=\,  w  
 \\[4pt]
       \pdt{w} \,=\, -\omega^2 u + \kconst_1 + f_\ss{u}
\end{array} 
&& \Bigg| &&
\begin{array}{lllll}
     \rng{\tup{q}} \,=\, -  \hdge{\htup{\slangup}}\cdot\tup{q} 
     &,
 \\[4pt]
     \rng{\tup{p}}  \,=\, -\hdge{\htup{\slangup}}\cdot\tup{p}  + \tfrac{1}{\slang u^2} \tup{f} 
     &,
\end{array}
\;\;
\begin{array}{llll}
     \rng{u}   \,=\,  w /\slang
 \\[4pt]
       \rng{w} \,=\, \tfrac{1}{\slang} (-\omega^2 u + \kconst_1 + f_\ss{u} )
\end{array} 
\end{align}
\end{small}
where \eq{\omega^2:= \slang^2-\kconst_2} (we implicitly assume \eq{\slang^2-\kconst_2>0} such that \eq{\omega} is real). 
In the \textit{unperturbed} case (\eq{\tup{F}= 0 =(\tup{f},f_\ss{u})}), the angular momentum is conserved and the solutions to either of the above systems are readily obtained: 
 \begin{small}
 \begin{flalign} \label{qupw_sols_sum}
\quad (\tup{f},f_\ss{u}) = 0
\quad \Rightarrow &&
 \begin{array}{lllllll}
     \tup{q}_s = \tup{q}_\tau 
      & =\,  \tup{q}_\zr \csn{\tau} - \tfrac{1}{\slang} \hdge{\tup{\slangup}} \cdot \tup{q}_\zr \snn{\tau} 
        \quad \simeq\, 
      \tup{q}_\zr \csn{\tau} + \tfrac{1}{\slang} \tup{p}_\zr \snn{\tau} 
\\[4pt]
      \tup{p}_s = \tup{p}_\tau 
      & =\,   \tup{p}_\zr \csn{\tau} - \tfrac{1}{\slang} \hdge{\tup{\slangup}} \cdot \tup{p}_\zr \snn{\tau} 
        \quad \simeq\; 
      -\slang\tup{q}_{\zr}\snn{\tau}  +  \tup{p}_{\zr}\csn{\tau}
\\[4pt]
      u_s = u_\tau 
      & =\,   (u_{\zr}-\tfrac{\kconst_1}{\omega^2})\cos{\varpi \tau} \,+\, \tfrac{1}{\omega}w_{\zr}\sin{\varpi\tau} \,+\, \tfrac{\kconst_1}{\omega^2}
\\[4pt]
      w_s = w_\tau  
      & =\,   -\omega(u_{\zr}-\tfrac{\kconst_1}{\omega^2})\sin{\varpi \tau} \,+\, w_{\zr}\cos{\varpi \tau}
\end{array} 
&&
\left|\;\;\begin{array}{lll}
     \tau = \slang s
\\[3pt]
   \omega^2:= \slang^2-\kconst_2
 \\[3pt]
     \varpi^2 := \omega^2/ \slang^2 
 \\[3pt]
       \varpi \tau = \omega s
\end{array}\right.
\end{flalign}
\end{small}
where \eq{\tup{\slangup}=\tup{\slangup}_\zr} is conserved such that solutions are easily expressed in terms of either \eq{s} or \eq{\tau} as indicated above, and where ``\eq{\simeq}'' indicates relations simplified using the integrals of motion \eq{q=1} and \eq{\lambda=\htup{q}\cdot\tup{p}=0}. 
The solutions for the conjugate pair \eq{(u,p_\ss{u}=w/u^2)} are easily recovered from the above solutions for \eq{(u,w=u^2 p_\ss{u})}, as in Eq.\eqref{up_sol_f}. 
The inertial cartesian coordinate solutions may be recovered from the above projective coordinate solutions using either Eq.\eqref{sum_PT_n} or Eq.\eqref{sum_rv_sol}. This is given explicitly for the Kepler problem (\eq{\kconst_2=0}) in section \ref{sec:kep_sol}.


\paragraph*{Second-Order Dynamics:~Equivalence to a Linear Oscillator.}
For perturbed Manev-type dynamics,  the  \eq{s}-parameterized second-order equations for \eq{\tup{q}} and \eq{u} were found in section \ref{sec:2BP_lin} as follows (now with \eq{n=-1}): 
\begin{small}
\begin{flalign} \label{ddqu1}
&&
\begin{array}{rllll}
      \pddt{\tup{q}} \,+\, \slang^2\tup{q} & =\, \tfrac{q^2}{u^2}\tup{f}
     &\quad=\, \tfrac{q}{u^3} ( \imat_3 - \htup{q}\otms\htup{q} )\cdot \tup{F}
\\[6pt]
     \pddt{u} + \omega^2 u - \kconst_1  & =\, f_\ss{u}
     &\quad=\, -\tfrac{1}{u^2}\htup{q} \cdot \tup{F}
\end{array} 
&&
\omega^2:=\slang^2-\kconst_2
\quad
\end{flalign}
\end{small}
The \eq{\tau}-parameterized second-order equations were also found in section \ref{sec:2BP_lin} (now with \eq{n=-1}):
\begin{small}
\begin{flalign} \label{ddqu2}
&&
\begin{array}{rlllll}
       \rrng{\tup{q}} \,+\, \tup{q} 
       & =\, \tfrac{q^2}{u^2 \slang^2} ( \imat_3 - \tfrac{1}{\slang}\rng{\tup{q}}\otms\tup{p} ) \cdot \tup{f}
       &\quad=\, \tfrac{1}{q u^3 \slang^2} ( q^2\imat_3 - \tup{q}\otms\tup{q} -  \rng{\tup{q}} \otms \rng{\tup{q}} ) \cdot \tup{F}
\\[6pt]
     \rrng{u} + \varpi^2 u  -  \tfrac{\kconst_1}{\slang^2} 
      & =\, \tfrac{1}{\slang^2} ( f_\ss{u} - \tfrac{q^2}{\slang u^2} \rng{u} \tup{p} \cdot \tup{f} )
    &\quad=\, -\tfrac{1}{q u^3 \slang^2} ( u \tup{q} + \rng{u} \rng{\tup{q}} ) \cdot\tup{F}
\end{array} 
&&
\varpi := \omega/\slang
\quad
\end{flalign}
\end{small}
\begin{small}
\begin{itemize}
    \item \textit{Simplifications.} As usual, the right-hand-side of the above second-order dynamics can be simplified with with \eq{q=1} and \eq{\lambda=\htup{q}\cdot\tup{p}=0}  such that \eq{\rng{\tup{q}} \simeq \tup{p}/\slang \simeq \htup{p}}. For instance, the above \eq{\rrng{\tup{q}}} and \eq{\rrng{u}} ODEs are then equivalent to:
    \begin{small}
    \begin{align} \label{ddqu2_simp}
         \fnsz{\begin{array}{cc}
         q=1  \\
         \tup{q}\cdot\tup{p}=0 
    \end{array}}
    \Rightarrow \qquad
        \rrng{\tup{q}} + \tup{q} 
        \,\simeq\, \tfrac{1}{u^2 \slang^2} ( \imat_3 - \rng{\tup{q}}\otms\rng{\tup{q}} ) \cdot \tup{f}
        \;\;\simeq\,
             \tfrac{1}{u^3\slang^2}(\htup{\slangup} \otms \htup{\slangup}) \cdot \tup{F}
        \qquad,\qquad 
         \rrng{u} + \varpi^2 u  -  \tfrac{\kconst_1}{\slang^2}  \,\simeq\,  \tfrac{1}{u^2\slang^2} (u^2 f_\ss{u} -   \rng{u} \rng{\tup{q}}\cdot \tup{f} )
    \end{align}
    \end{small}
    where we have used \eq{ \htup{\slangup} \otms \htup{\slangup} \simeq \imat_3 -\htup{q}\otms\htup{q} - \htup{p}\otms\htup{p} \simeq \imat_3 - \tup{q}\otms\tup{q} - \rng{\tup{q}}\otms\rng{\tup{q}}}.\footnote{Recall that \eq{ \{\htup{t}_r,\htup{t}_\tau, \htup{t}_\slang \} =  \{\htup{q},-\hdge{\htup{\slangup}}\cdot\htup{q},\htup{\slangup}\} \simeq \{\tup{q},\htup{p},\htup{\slangup}\}} are the inertial cartesian components of the LVLH basis. with \eq{\htup{\slangup}=\htup{q}\tms\htup{p}},  \eq{q^2\simeq 1}, and \eq{\slang^2\simeq p^2}. As such:\\
    \eq{\qquad\qquad\qquad
    \imat_3 \,=\,  
     \htup{q}\otms\htup{q} + (\hdge{\htup{\slangup}}\cdot\htup{q})\otms (\hdge{\htup{\slangup}}\cdot\htup{q}) +  \htup{\slangup} \otms \htup{\slangup}
    \,\simeq\, \htup{q}\otms\htup{q} + \htup{p}\otms\htup{p} + \htup{\slangup}\otms\htup{\slangup} \,\simeq\,  \tup{q}\otms\tup{q} + \rng{\tup{q}}\otms\rng{\tup{q}} + \htup{\slangup}\otms\htup{\slangup}  }.}
\end{itemize}
\end{small}

\vspace{1ex}
\noindent 
Now, for \textit{unperturbed} motion (\eq{\tup{F}= 0 =\bartup{f}}), the right-hand-side of the above second-order ODEs vanish and \eq{\slang=\slang_\zr} is constant such that the second-order ODEs for \eq{(\tup{q},u)} in Eq.\eqref{ddqu1} and Eq.\eqref{ddqu2} are those of an inhomogeneous linear oscillator with solutions as follows (assuming \eq{\slang^2>\kconst_2} such that \eq{\omega} and \eq{\varpi} are real):
 \begin{small}
 \begin{align}
 \begin{array}{llll}
         \tup{q}_{s}  \,=\,   \tup{q}_{\zr}\cos{\slang s}  +  \tfrac{1}{\slang}\pdt{\tup{q}}_{\zr}\sin{\slang s}
\\[4pt]
     \pdt{\tup{q}}_{s}  \,=\,  -\slang\tup{q}_{\zr}\sin{\slang s} + \pdt{\tup{q}}_{\zr}\cos{\slang s}  
\\[4pt]
     u_{s} \,=\, (u_{\zr}-\tfrac{\kconst_1}{\omega^2})\cos{\omega s}
    \,+\, \tfrac{1}{\omega}\pdt{u}_{\zr}\sin{\omega s} \,+\, \tfrac{\kconst_1}{\omega^2}
\\[4pt]
    \pdt{u}_{s} \,=\, -\omega(u_{\zr}-\tfrac{\kconst_1}{\omega^2})\sin{\omega s} \,+\, \pdt{u}_{\zr}\cos{\omega s}
\end{array}
\qquad
\begin{array}{llll}
       = &\qquad    \tup{q}_{\tau}  \,=\,  \tup{q}_{\zr}\csn{\tau}  +  \rng{\tup{q}}_{\zr}\snn{\tau} 
\\[4pt]
       \neq &\qquad      \rng{\tup{q}}_{\tau}  \,=\,  -\tup{q}_{\zr}\snn{\tau}
            + \rng{\tup{q}}_{\zr}\csn{\tau} 
\\[4pt]
      = &\qquad 
      u_{\tau} \,=\, (u_{\zr}-\tfrac{\kconst_1}{\omega^2})\cos{\varpi \tau}
    \,+\, \tfrac{1}{\varpi}\rng{u}_{\zr}\sin{\varpi\tau} \,+\, \tfrac{\kconst_1}{\omega^2}
\\[4pt]
    \neq &\qquad   
    \rng{u}_{\tau} \,=\, -\varpi(u_{\zr}-\tfrac{\kconst_1}{\omega^2})\sin{\varpi \tau} \,+\, \rng{u}_{\zr}\cos{\varpi \tau}
\end{array} 
\end{align}
\end{small}
The Kepler case is recovered with \eq{\kconst_2=0} such that \eq{\omega=\slang} and \eq{\varpi=1}.

\subsection{Linear Kepler dynamics \& closed-form solutions} \label{sec:kep_sol}


Kepler-type dynamics and solutions are easily recovered from the preceding developments by setting \eq{\kconst_2=0} (thus, \eq{\omega=\slang} and \eq{\varpi=1}).  Regardless, we will go through the developments, adding many further details and useful relations beyond what was included above for the Manev case.
We will consider only the unperturbed Kepler case, \eq{(\tup{f},f_\ss{u})=0}, as the perturbation terms are arbitrary and no different from what was just given for perturbed Manev dynamics.  

The \eq{s}- or \eq{\tau}-parameterized unperturbed
Kepler dynamics are given by Eq.\eqref{qp_eom_s_sum} and Eq.\eqref{qp_eom_ta_sum}, respectively, with \eq{\kconst_2=0} and \eq{(\tup{f},f_\ss{u})=0}.  Exchanging \eq{p_\ss{u}} for the quasi-momenta coordinate \eq{w:=u^2 p_\ss{u}}, the reparameterized Kepler dynamics are given by the first-order ODEs:
\begin{small}
\begin{flalign} \label{qpw_kep_EOM}
\begin{array}{cc}
     \fnsize{Kepler}  \\
     \fnsize{dynamics} 
\end{array}
&&
\begin{array}{llll}
     \pdt{\tup{q}} \,=\, -  \hdge{\tup{\slangup}}\cdot\tup{q} 
     &,
 \\[4pt]
     \pdt{\tup{p}}  \,=\, -\hdge{\tup{\slangup}}\cdot\tup{p}  
     &,
\end{array}
\;
\begin{array}{llll}
     \pdt{u}   \,=\,  w  
 \\[4pt]
       \pdt{w} \,=\, -\slang^2 u + \kconst_1 
\end{array} 
&& \Bigg| &&
\begin{array}{llll}
     \rng{\tup{q}} \,=\, -  \hdge{\htup{\slangup}}\cdot\tup{q} 
     &,
 \\[4pt]
     \rng{\tup{p}}  \,=\, -\hdge{\htup{\slangup}}\cdot\tup{p} 
     &,
\end{array}
\;
\begin{array}{llll}
     \rng{u}   \,=\,  w/\slang
 \\[4pt]
       \rng{w} \,=\, -\slang u + \kconst_1/\slang
\end{array}
&&
\end{flalign}
\end{small}
Which were shown to be equivalent to the following  second-order ODEs (this does not use \eq{q=1} or \eq{\htup{q}\cdot\tup{p}=0}):
\begin{small}
\begin{flalign}
\begin{array}{cc}
     \fnsize{Kepler}  \\
     \fnsize{dynamics} 
\end{array}
&&
\begin{array}{llllll}
       \pddt{\tup{q}} + \slang^2\tup{q} \,=\, 0
\\[4pt]
     \pddt{u} + \slang^2 u - \kconst_1  \,=\, 0
\end{array}  
&& \Bigg| &&
\begin{array}{llll}
    \rrng{\tup{q}} + \tup{q} \,=\, 0
\\[4pt]
     \rrng{u} + u  -  \tfrac{\kconst_1}{\slang^2} \,=\, 0
\end{array}
&&
\end{flalign}
\end{small}
We will focus on the first-order dynamics in Eq.\eqref{qpw_kep_EOM}. 
The solutions to either system in Eq.\eqref{qpw_kep_EOM} coincide using the fact that \eq{\tau = \slang s} for any central-force dynamics. This is just a special case of what was already given in Eq.\eqref{qupw_sols_sum}:
\begin{small}
\begin{flalign} \label{qusol_s_2bp}
&&
\boxed{\begin{array}{rllllll}
     \tup{q}_\tau 
      & =\,  \tup{q}_\zr \csn{\tau} -  \hdge{\htup{\slangup}} \cdot \tup{q}_\zr \snn{\tau} 
\\[4pt]
      \tup{p}_\tau 
      & =\,   \tup{p}_\zr \csn{\tau} -  \hdge{\htup{\slangup}} \cdot \tup{p}_\zr \snn{\tau} 
\\[4pt]
      u_\tau 
      & =\, (u_{\zr}-\tfrac{\kconst_1}{\slang^2})\csn{\tau} \,+\, \tfrac{1}{\slang}w_{\zr}\snn{\tau} \,+\, \tfrac{\kconst_1}{\slang^2}
\\[4pt]
      w_\tau  
      & =\, -\slang(u_{\zr}-\tfrac{\kconst_1}{\slang^2})\snn{\tau} \,+\, w_{\zr}\csn{\tau}
\end{array} }
&&
\begin{array}{llll}
     \tau = \slang s
 \\[3pt]
     \tup{\slangup}=\tup{\slangup}_0
 \\[3pt]
     \tup{q}_\tau \cdot \tup{p}_\tau = \tup{q}_\zr\cdot \tup{p}_\zr 
\\[3pt]
    \mag{\tup{q}_\tau} = \mag{\tup{q}_\zr}
\\[3pt]
     \mag{\tup{p}_\tau} = \mag{\tup{p}_\zr}
\end{array}
 &&
\end{flalign}
\end{small}
As noted in Eq.\eqref{E3sol_s_prj}, the above \eq{(\tup{q}_\tau,\tup{p}_\tau)} solutions directly verify preservation of \eq{\tup{\slangup}=\tup{q}\tms\tup{p}}, \eq{\mag{\tup{p}}}, \eq{\mag{\tup{q}}}, and \eq{\tup{q}\cdot\tup{p}}. 
\begin{small}
\begin{itemize}
    \item \textit{Solution for \eq{(u,p_\ss{u})}.} The solution for the conjugate pair \eq{(u,p_\ss{u})} is easily recovered from the  above solution for \eq{(u,w)} using \eq{p_\ss{u} = w/u^2 \leftrightarrow w =u^2 p_\ss{u}}.  
    This is given in the footnote\footnote{The solution for the conjugate pair \eq{(u,p_\ss{u}=w/u^2)} is obtained by direct substitution:
    \begin{align} \label{upu_KEPsols_sum}
       u_{\tau} \,=\,  (u_{\zr}-\tfrac{\kconst_1}{\slang^2})\csn{\tau}
        \,+\, \tfrac{1}{\slang}u_{\zr}^2 p_{u_0}\snn{\tau} \,+\, \tfrac{\kconst_1}{\slang^2}
    \qquad,\qquad 
           p_{u_\tau} \,=\,  \frac{ -\slang(u_{\zr}-\frac{\kconst_1}{\slang^2})\snn{\tau} \,+\, u_{\zr}^2 p_{u_0}\csn{\tau} }
               {\big[ (u_{\zr}-\frac{\kconst_1}{\slang^2})\csn{\tau}
        \,+\, \frac{1}{\slang}u_{\zr}^2 p_{u_0}\snn{\tau} \,+\, \frac{\kconst_1}{\slang^2} \big]^2}
    \end{align} }.
     \item  \textit{Simplifications.} 
    Nothing so far has been simplified using the  integrals of motion \eq{q=1} or \eq{\lambda=\htup{q}\cdot\tup{p}=0} (which imply the relations in Eq.\eqref{qpl_rels}). As discussed, we are free to make these simplifications such that, for example, the above already-linear dynamics for \eq{(\tup{q},\tup{p})} in Eq.\eqref{qpw_kep_EOM} are equivalent to the following still-linear system:  
    \begin{small}
    \begin{flalign} \label{qp_KEPsol_f_m1}
    \;\;
    \fnsz{\begin{array}{cc}
             q=1  \\
             \tup{q}\cdot\tup{p}=0 
        \end{array}}
        \;\Rightarrow \quad
    \left\{\qquad
    \begin{array}{lllll}
         \pdt{\tup{q}} \simeq \tup{p} 
         \\[4pt]
         \pdt{\tup{p}} 
         \simeq -\slang^2\tup{q} 
    \end{array} \right.
        \qquad \fnsize{or,}\qquad
     \begin{array}{lllll}
         \rng{\tup{q}} \simeq \tfrac{1}{\slang}\tup{p} 
         \simeq \htup{p}
         \\[4pt]
         \rng{\tup{p}} \simeq -\slang \tup{q} 
    \end{array}
    &&
    \begin{array}{lll}
          \slang = \slang_\zr \simeq  p_\zr = p
    \end{array} 
       \qquad
    \end{flalign}
    \end{small}
    The relation \eq{\slang \simeq p} could also be used to rewrite the ODEs for \eq{(u,w)} in Eq.\eqref{qpw_kep_EOM}. 
    The solutions to the simplified Kepler dynamics are then given by
    \begin{small}
    \begin{flalign} \label{qusol_s_2bp_simp} 
    \;\;
    \fnsz{\begin{array}{cc}
             q=1  \\
             \tup{q}\cdot\tup{p}=0 
        \end{array}}
        \Rightarrow \;\;
    \left\{\;\;
   \begin{array}{lllllll}
         \tup{q}_\tau 
            \simeq\, 
          \tup{q}_\zr \csn{\tau} + \tfrac{1}{\slang} \tup{p}_\zr \snn{\tau} 
            &  \simeq\, 
            \tup{q}_\zr \csn{\tau} +  \htup{p}_\zr \snn{\tau} 
    \\[4pt]
          \tup{p}_\tau 
          \simeq\, 
          \tup{p}_{\zr}\csn{\tau} -\slang\tup{q}_{\zr}\snn{\tau} 
          &  \simeq\; 
          \tup{p}_{\zr}\csn{\tau} -p_\zr\tup{q}_{\zr}\snn{\tau} 
    \\[4pt]
          u_\tau 
          =\,  (u_{\zr}-\tfrac{\kconst_1}{\slang^2})\csn{\tau} + \tfrac{1}{\slang}w_{\zr}\snn{\tau} + \tfrac{\kconst_1}{\slang^2}
           &  \simeq 
            (u_{\zr}-\tfrac{\kconst_1}{p_0^2})\csn{\tau} + \tfrac{1}{p_0}w_{\zr}\snn{\tau} + \tfrac{\kconst_1}{p_0^2}
    \\[4pt]
          w_\tau  
           =\,  w_{\zr}\csn{\tau} - \slang(u_{\zr}-\tfrac{\kconst_1}{\slang^2})\snn{\tau} 
           &  \simeq\,  w_{\zr}\csn{\tau} - p_\zr(u_{\zr}-\tfrac{\kconst_1}{p_0^2})\snn{\tau}  
    \end{array}   \right.
    \quad,&&
    \begin{array}{llll}
         \tau = \slang s \simeq p s
     \\[3pt]
         \tup{\slangup} = \tup{\slangup}_\zr
    \end{array}
    \end{flalign}
    \end{small}
    which also follows directly from simplifying the Kepler solutions in Eq.\eqref{qusol_s_2bp} using the relations in Eq.\eqref{qpl_rels}. 
\end{itemize}
\end{small}

\paragraph*{Recovering Cartesian Coordinate Kepler Solutions.}
The Kepler solutions for cartesian position and velocity coordinates, \eq{(\tup{r},\tup{v})}, may now be recovered from the above solutions for \eq{(\tup{q},u,\tup{p},w=u^2 p_\ss{u})} using either the full transformation given in Eq.\eqref{sum_PT_n}, or using the simplified relations in Eq.\eqref{sum_rv_sol}:
\begin{small}
\begin{align} \label{rv_qu_simp_yetagain}
\begin{array}{rlllllll}
     \begin{array}{cc}
     \fnsize{unsimplified}
    \end{array} \!\!: 
    &\quad  \tup{r}_\tau &=\, \tfrac{1}{u_\tau}\htup{q}_\tau 
&\quad,\qquad 
      \tup{v}_\tau &=\, -u_\tau(\tup{q}_\tau\wdg\tup{p}_\tau) \cdot \htup{q}_\tau - w_\tau\htup{q}_\tau 
   \;=\, -u_\tau \hdge{\tup{\slangup}}_\tau \cdot \htup{q}_\tau - w_\tau\htup{q}_\tau 
\\
\begin{array}{cc}
     \fnsize{simplified with}  \\
     \fnsz{q=1, \; \htup{q}\cdot\tup{p}=0}  
\end{array} \!\!: 
&\quad
    &\simeq\, \tfrac{1}{u_\tau}\tup{q}_\tau
&\quad,\qquad 
    &\simeq\, 
    u_\tau \tup{p}_\tau - w_\tau \tup{q}_t 
\end{array}
\end{align}
\end{small}
We first consider the unsimplified relations. For the case at hand (Kepler dynamics), we note these are equivalent to \eq{\tup{r}_\tau =\tfrac{1}{q_0} \tfrac{1}{u_\tau}\tup{q}_\tau} and \eq{\tup{v}_\tau =  -\tfrac{1}{q_0} ( u_\tau \hdge{\tup{\slangup}}_\zr \cdot \tup{q}_\tau + w_\tau \tup{q}_\tau )} — because \eq{q} is always an integral of motion, and \eq{\hdge{\tup{\slangup}}=\tup{q}\wdg\tup{p}} is an integral of motion for any central-force dynamics. 
Substitution of the unsimplified Kepler solutions for \eq{(\tup{q}_\tau,u_\tau,\tup{p}_\tau,w_\tau)} from Eq.\eqref{qusol_s_2bp} leads to:
\begin{small}
\begin{flalign} \label{rv_kepsol_nosimp}
\begin{array}{llllll}
     \tup{r}_\tau = \dfrac{\htup{q}_\tau}{u_\tau}
  \,=\, 
   \dfrac{ \htup{q}_\zr \csn{\tau} -  \hdge{\htup{\slangup}} \cdot \htup{q}_\zr \snn{\tau} }{ (u_\zr-\frac{\kconst_1}{\slang^2})\csn{\tau} + \frac{1}{\slang}w_\zr\snn{\tau} + \frac{\kconst_1}{\slang^2} }
   \\[4pt]
   \;\;
\end{array} 
&&,&&
\begin{array}{llllll}
     \tup{v}_\tau  = - u_\tau \, \hdge{\tup{\slangup}}_\tau \cdot \htup{q}_\tau - w_\tau \htup{q}_\tau 
  & = \,   \tup{v}_\zr  - \tfrac{\kconst_1}{\slang^2} \hdge{\tup{\slangup}} \cdot \big[\htup{q}_\zr( \csn{\tau}-1) - \hdge{\htup{\slangup}}\cdot\htup{q}_\zr \snn{\tau} \big]  
\\[4pt]
      & =\, \tup{v}_\zr - \tfrac{\kconst_1}{\slang^2} \big[ \hdge{\tup{\slangup}}\cdot \htup{q}_\zr(\csn{\tau}-1) +  \slang\htup{q}_\zr\snn{\tau}  \big]
\\[4pt]
    & = \,  \tup{v}_\zr  - \tfrac{\kconst_1}{\slang^2} \hdge{\tup{\slangup}} \cdot (\htup{q}_\tau - \htup{q}_\zr) 
\end{array} 
\end{flalign}
\end{small}
with \eq{ \tup{v}_\zr = - u_\zr \, \hdge{\tup{\slangup}} \cdot \htup{q}_\zr - w_\zr \htup{q}_\zr}
the velocity at \eq{\tau_\zr=0}.\footnote{If \eq{\tau=0} coincides with periapsis, then \eq{\tau} is the true anomaly and \eq{\mag{\tup{v}_\tau} \leq \mag{\tup{v}_\zr}}.} 
The above may be written in terms of \eq{(\tup{r}_\zr,\tup{v}_\zr)} using:
\begin{small}
\begin{align} \label{qp_rv_init_rels}
\begin{array}{llll}
     \htup{q}_\zr = \htup{r}_0
 \\[4pt]
      u_\zr = 1/r_0
\\[4pt]
    w_\zr = -\htup{r}_\zr \cdot \tup{v}_\zr = - \dot{r}_0
\end{array}
 &&,&&
\begin{array}{llll}
        \hdge{\tup{\slangup}} =  \hdge{\tup{\slangup}}_\zr = \tup{q}_\zr\wdg \tup{p}_\zr =  \tup{r}_\zr\wdg \tup{v}_\zr
\\[4pt]
     \tup{\slangup}=\tup{\slangup}_0 = \hdge{\tup{p}}_\zr\cdot\tup{q}_\zr = \hdge{\tup{v}}_\zr\cdot\tup{r}_\zr
 \\[4pt]
      \slang^2 = \slang_\zr^2 = q_\zr^2 p_\zr^2 - (\tup{q}_\zr \cdot\tup{p}_\zr)^2 = r_\zr^2 \v_\zr^2 - (\tup{r}_\zr \cdot\tup{v}_\zr)^2
\end{array}
 &&,&&
\begin{array}{rllll}
      \tup{v}_\zr &=\, - u_\zr  \hdge{\tup{\slangup}} \cdot \htup{q}_\zr - w_\zr \htup{q}_\zr
\\[3pt]
     &\simeq\,  u_\zr\tup{p}_\zr-w_\zr\tup{q}_\zr
\end{array}
\end{align}
\end{small}
none of which require the simplifications \eq{q=1} or \eq{\htup{q}\cdot\tup{p}=0}. 
Eq.\eqref{rv_kepsol_nosimp} then leads to
\begin{small}
\begin{flalign} \label{rv_kepsol_rv0}
\;\;
      \tup{r}_\tau  
      \,=\,  \dfrac{ \htup{r}_\zr\csn{\tau}  \,-\,  \hdge{\htup{\slangup}} \cdot\htup{r}_\zr \snn{\tau}}{ (1/r_\zr-\frac{\kconst_1}{\slang^2})\csn{\tau} \,-\, \frac{1}{\slang}\htup{r}_\zr \cdot \tup{v}_\zr\snn{\tau} \,+\, \frac{\kconst_1}{\slang^2}   }  
      \,=\, \dfrac{\htup{r}_\tau}{1/r_\tau}
&&,&&
\begin{array}{llllll}
     \tup{v}_\tau  & =\, \tup{v}_\zr  - \tfrac{\kconst_1}{\slang^2} \hdge{\tup{\slangup}} \cdot \big[\htup{r}_\zr( \csn{\tau}-1) - \hdge{\htup{\slangup}}\cdot\htup{r}_\zr \snn{\tau} \big] 
\\[4pt]
      & =\, \tup{v}_\zr - \tfrac{\kconst_1}{\slang^2} \big[ \hdge{\tup{\slangup}}\cdot \htup{r}_\zr(\csn{\tau}-1) +  \slang\htup{r}_\zr\snn{\tau}  \big]
\\[4pt]
   & =\,  \tup{v}_\zr - \tfrac{\kconst_1}{\slang^2} \hdge{\tup{\slangup}}\cdot(\htup{r}_\tau - \htup{r}_\zr )
\end{array}
\;\;
\end{flalign}
\end{small}
with \eq{\tup{\slangup}=\tup{\slangup}_\zr} and \eq{\slang=\slang_\zr} given in terms of \eq{(\tup{r}_\zr,\tup{v}_\zr)} as in Eq.\eqref{qp_rv_init_rels}.

\begin{small}
\begin{itemize}
    \item \textit{Simplifications.} Consider the simplified relations in Eq.\eqref{rv_qu_simp_yetagain}. Substitution of the simplified Kepler solutions for \eq{(\tup{q}_\tau,u_\tau,\tup{p}_\tau,w_\tau)} from Eq.\eqref{qusol_s_2bp_simp} then leads to
    \begin{small}
    \begin{flalign} \label{rv_kepsol} 
    \qquad
    \begin{array}{llllll}
          \tup{r}_\tau \simeq \dfrac{\tup{q}_\tau}{u_\tau}  
          & \simeq\,  \dfrac{\tup{q}_\zr\csn{\tau}  \,+\,  \frac{1}{\slang} \tup{p}_\zr\snn{\tau}}{ (u_\zr-\frac{\kconst_1}{\slang^2})\csn{\tau} + \frac{1}{\slang}w_\zr\snn{\tau} + \frac{\kconst_1}{\slang^2}   }  
    \end{array} 
    &&,&&
    \begin{array}{llllll}
         \tup{v}_\tau  \simeq u_\tau\tup{p}_\tau - w_\tau\tup{q}_\tau 
        & \simeq\,  \tup{v}_\zr + \tfrac{\kconst_1}{\slang^2} \big[ \tup{p}_\zr(\csn{\tau}-1) -  \slang\tup{q}_\zr\snn{\tau}  \big]
        \\[4pt]
        & \simeq\,  \tup{v}_\zr + \tfrac{\kconst_1}{\slang^2} (\tup{p}_\tau - \tup{p}_\zr )
    \end{array}
    \;\;
    \end{flalign}
    \end{small}
    where  \eq{\tup{v}_\zr \simeq  u_\zr\tup{p}_\zr-w_\zr\tup{q}_\zr}. The above agrees with Eq.\eqref{rv_kepsol_nosimp} and Eq.\eqref{rv_kepsol_rv0} using the following relations (cf.~Eq.\eqref{qpl_rels}):
    \begin{small}
    \begin{align}
          \tup{q}\simeq \htup{q} = \htup{r}
          \qquad,\qquad 
          \tup{p} \simeq -\hdge{\tup{\slangup}}\cdot \tup{q} \simeq -\hdge{\tup{\slangup}}\cdot \htup{q} =  -\hdge{\tup{\slangup}}\cdot \htup{r} 
    \end{align}
    \end{small}
\end{itemize}
\end{small}

\paragraph{Time Solution.}
The Kepler solution for \eq{t(\tau)}
is obtained using that for \eq{u(\tau)} in Eq.\eqref{qusol_s_2bp} as follows:
\begin{small}
\begin{align} \label{t_sol_int_gen}
    \mrm{d} t = \tfrac{1}{\slang u^2} \mrm{d} \tau
    \;\;,\;\; \slang=\slang_\zr 
    && \Rightarrow && 
    \begin{array}{llll}
    t-t_{\zr} \,=\, \int_{0}^{\tau} \tfrac{1}{\slang} u_\tau^{-2} \mrm{d} \tau
    \,=\, \tfrac{1}{\slang} \int_{0}^{\tau}
      \Big( (u_{\zr}-\tfrac{\kconst_1 }{\slang^2})\csn{\tau}
    + \tfrac{1}{\slang}w_\zr \snn{\tau} + \tfrac{\kconst_1 }{\slang^2} \Big)^{-2}   \mrm{d} \tau
    \end{array}
\end{align}
\end{small}
were we have assumed \eq{\tau_\zr=0} at \eq{t_\zr}. 
Though the above has a closed-form solution, it is rather cumbersome. 
For simplicity, let us limit consideration to the natural special case that \eq{\tau_\zr=0} coincides with periapsis such that \eq{\tau} is the true anomaly. Then \eq{\dot{r}_\zr=0} and \eq{w_\zr = \pdt{u}_\zr = -\dot{r}_\zr = 0 } vanishes and the solutions for \eq{(u,w)} in Eq.\eqref{qusol_s_2bp} simplify to:
\begin{small}
\begin{flalign}
\quad
\begin{array}{cccc}
     \fnsize{if  $\tau_\zr=0$} \\
     \fnsize{is periapsis}
\end{array}  
\!\!:
&&
w_\zr = -\dot{r}_\zr = 0  
\qquad \Rightarrow \qquad 
\begin{array}{llll}
      u_\tau 
      &=\,
      (u_{\zr}-\tfrac{\kconst_1}{\slang^2})\csn{\tau}  \,+\, \tfrac{\kconst_1}{\slang^2}
\quad\;\;,\quad\;\; 
      w_\tau  
      & =\, -\slang(u_{\zr}-\tfrac{\kconst_1}{\slang^2})\snn{\tau} 
\end{array}
&&
\end{flalign}
\end{small}
In this case (when \eq{\tau} is the actual true anomaly) then we further find that \eq{ u_\zr -  \tfrac{\kconst_1}{\slang^2} =  \tfrac{\kconst_1}{\slang^2} e} where \eq{e} is the classic eccentricity (cf.~Eq.\eqref{Pslr_qp0_periapsis} below), and the above solution for \eq{u_\tau} then leads to
\begin{small}
\begin{flalign}
   \quad
\begin{array}{cccc}
     \fnsize{if  $\tau_\zr=0$} \\
     \fnsize{is periapsis}
\end{array}  
 \!\!:
&&
\begin{array}{llll}
      u_\tau 
      \,=\, 
      \tfrac{\kconst_1}{\slang^2} e \csn{\tau} \,+\, \tfrac{\kconst_1}{\slang^2}
      \,=\, 
      \tfrac{\kconst_1}{\slang^2} ( 1 + e \csn{\tau})
       \,=\, \tfrac{1}{P_{\mrm{slr}}}( 1 + e \csn{\tau})
       \,=\, 1/r_\tau
\end{array}
&&
\end{flalign}
\end{small}
thereby recovering the well-known conic section formula \eq{r_\tau = P_\ss{\mrm{slr}}/(1+e\csn{\tau})}, where  \eq{P_\ss{\mrm{slr}} = \slang^2/\kconst_1} is the \textit{semilatus rectum}. 
The integral in Eq.\eqref{t_sol_int_gen} for \eq{t(\tau)} then simplifies to:
\begin{small}
\begin{flalign}
    \quad
    \begin{array}{cccc}
         \fnsize{if  $\tau_\zr=0$} \\
         \fnsize{is periapsis}
    \end{array}  
     \!\!:
&&
    \begin{array}{llll}
        t - t_\zr \,=\, \int_{0}^{\tau} \tfrac{1}{\slang} u_\tau^{-2} \mrm{d} \tau
        &\,=\, \tfrac{\slang^3}{\kconst_1^2} \int_{0}^{\tau}  (1+e \csn{\tau})^{-2} \mrm{d} \tau
    \end{array}
&&
\end{flalign}
\end{small}
which leads to the following solutions for the indicated orbit types:
\begin{small}
\begin{flalign} \label{TA_t_sol}
      \quad
    \begin{array}{cccc}
         \fnsize{if  $\tau_\zr=0$} \\
         \fnsize{is periapsis}
    \end{array}  
     \!\!:
&&
     \begin{array}{llll}
       \tfrac{\kconst_1^2}{\slang^3}(t - t_\zr) 
        &\,=\, -\tfrac{2}{(e^2-1)^{3/2}} \mrm{arctanh} \Big( (\tfrac{e-1}{e+1})^{\frac{1}{2}} \tan{\!\tfrac{\tau}{2}} \Big) 
         \,+\, \tfrac{e\snn{\tau}}{(e^2-1)(1+e\csn{\tau})}
         &\qquad (e>1)
     \\[6pt]
         &\,=\, \tfrac{2}{(1-e^2)^{3/2}} \arctan{ \Big( (\tfrac{1-e}{1+e})^{\frac{1}{2}} \tan{\!\tfrac{\tau}{2}} \Big) }
         \,-\, \tfrac{e\snn{\tau}}{(1-e^2)(1+e\csn{\tau})}
         &\qquad (e<1)
     \\[6pt]
           &\,=\, \tfrac{1}{2} \big( \tan{\!\tfrac{\tau}{2}} + \tfrac{1}{3}\tan^3{\!\tfrac{\tau}{2}} \big)
            &\qquad (e=1)
    \\[6pt]
           &\,=\, \tau  
            &\qquad (e=0)
    \end{array}
\end{flalign}
\end{small}
where the equations for the \eq{e=1} and \eq{e=0} cases both follow from that for the \eq{e<1} case.\footnote{The equation for the \eq{e=0} case follows immediately as a simplification of the \eq{e<1} case (the later holds more specifically for \eq{0\leq e <1}). The equation for the \eq{e=1} case follows from taking the limit as \eq{e\to 1} of the equation for the \eq{e<1} case. It was also noted by Elipe et al.~in \cite{elipe2020explicitTA}}
The above relations, which were obtained using a symbolic processor in \textsc{matlab},  were also presented in more detail in \cite{elipe2020explicitTA}.

%% file: Mysecs_prj/new_3.3_prjPerifocal.tex
\paragraph*{Relations with the perifocal basis.} 
The close relation between the projective coordinates, \eq{(\tup{q},u,\tup{p},p_\ss{u})} and the LVLH basis was discussed in Eq.\eqref{qpl_rels}-Eq.\eqref{qpl_rels2}. 
We now consider another reference frame often encountered in orbital mechanics — the \textit{perifocal basis}. 
This is defined from a set of three  mutually orthogonal vectors commonly studied in Kepler-type dynamics; the Laplace-Runge-Lenz (LRL) vector/eccentricity vector, \eq{\tup{e}}, and the Hamilton vector, \eq{\tup{h}}, along with the usual specific angular momentum vector, \eq{{\tup{\slangup}}}.
Their cartesian components are given in terms of the inertial cartesian coordinates \eq{(\tup{r},\tup{v})} 
 as:\footnote{The Laplace-Runge-Lenz (LRL) vector and the Hamilton vector are often defined as \eq{\tup{E}=\tup{v}\tms {\tup{\slangup}}-\kconst_1 \htup{r}} and \eq{\tup{H}=\tfrac{1}{\slang^2}{\tup{\slangup}}\tms \tup{E}}, respectively.   The eccentricity vector \eq{\tup{e}} and Hamilton vector \eq{\tup{h}} seen in Eq.\eqref{LRL_e} are simply a scaling by \eq{\kconst_1}, that is, \eq{\tup{e}=\tfrac{1}{\kconst_1}\tup{E}} and \eq{\tup{h}=\tfrac{1}{\kconst_1}\tup{H}}. We should also not that, in the present context, all of these ``vectors'' are really just the inertial cartesian components (\eq{\mbb{R}^3}-valued functions). }
\begin{small}
\begin{flalign} \label{LRL_e_prj}
\quad
    \fnsize{for}\; V^0 = -\kconst_1/r\, :
&&
    \kconst_1  \tup{e} 
    \,=\,
    \hdge{\tup{\slangup}}\cdot\tup{v} - \kconst_1 \htup{r} 
\qquad,\qquad 
   \kconst_1  \tup{h} = -\tfrac{\kconst_1 }{\slang^2} \hdge{\tup{\slangup}}\cdot \tup{e}
   \,=\, \tup{v} + \tfrac{\kconst_1 }{\slang^2} \hdge{\tup{\slangup}}\cdot \htup{r} 
\qquad,\qquad 
    {\tup{\slangup}} = \hdge{\tup{v}}\cdot \tup{r} 
&&
\end{flalign}
\end{small}
It is well-known that the above \eq{\mbb{R}^3}-valued functions are all conserved for pure Kepler dynamics.
Their normalization defines the orthonormal \textit{perifocal basis}, \eq{\{\htup{e},\htup{h},\htup{\slangup}\}=:\{\htup{o}_e,\htup{o}_h,\htup{o}_\slang\}}.\footnote{For  Keplerian motion, this basis is fixed/constant, with \eq{\htup{o}_e=\htup{e}} and \eq{\htup{o}_h = \htup{h}} defining the orbital plane;  \eq{\htup{e}} directed towards periapsis,  \eq{\htup{h}} directed along the velocity at periapsis (tangent to the orbit), and \eq{\htup{o}_{\slang}=\htup{\slangup}} directed normal to the orbit plane, completing the right-handed triad.  The magnitude \eq{e=\mag{\tup{e}}} is the usual dimensionless eccentricity of the orbit. }
To above are then given in terms of the projective coordinates (using \eq{w:= u^2 p_\ss{u})} by substitution of the full transformation in Eq.\eqref{sum_PT_n}, or the simplified relations in Eq.\eqref{sum_rv_sol},
leading to:\footnote{It was already shown that \eq{\tup{\slangup} = \tup{r}\tms \tup{v} = \tup{q}\tms \tup{p}}. The expressions for \eq{ \kconst_1 \tup{e}} and \eq{\kconst_1  \tup{h}} in terms of \eq{(\tup{q},u,\tup{p},w=u^2 p_\ss{u})} are seen as follows, using Eq.\eqref{sum_PT_n} or Eq.\eqref{sum_rv_sol}:\\
        \eq{\qquad\qquad\qquad \kconst_1 \tup{e} \,=\
        \hdge{\tup{\slangup}}\cdot\tup{v} - \kconst_1 \htup{r} 
        \,=\, 
         \hdge{\tup{\slangup}}\cdot (- u \hdge{\tup{\slangup}}\cdot \htup{q} - w \htup{q}) 
         - \kconst_1 \htup{q} 
         \,=\, 
         \slang^2 u \htup{q} - w \hdge{\tup{\slangup}}\cdot \htup{q}   - \kconst_1 \htup{q}
         \,=\, 
         \slang^2(u -\tfrac{\kconst_1 }{\slang^2} )\htup{q} -   w \hdge{\tup{\slangup}}\cdot\htup{q}
          \,\simeq \,
          \slang^2(u -\tfrac{\kconst_1 }{\slang^2} )\tup{q} +   w \tup{p}
        } \\
        \eq{\qquad\qquad\qquad 
         \kconst_1  \tup{h} 
        \,=\, \tup{v} + \tfrac{\kconst_1 }{\slang^2} \hdge{\tup{\slangup}}\cdot \htup{r} 
        \,=\, 
        - u \hdge{\tup{\slangup}}\cdot \htup{q} - w \htup{q} + + \tfrac{\kconst_1 }{\slang^2} \hdge{\tup{\slangup}}\cdot \htup{q}
        \,=\,
        -(u - \tfrac{\kconst_1 }{\slang^2})\hdge{\tup{\slangup}}\cdot \htup{q} - w \htup{q}
        \,\simeq\,
        (u - \tfrac{\kconst_1 }{\slang^2})\tup{p} - w \tup{q}
        } . 
    }
\begin{small}
\begin{flalign} \label{LRL_e}
&&
\begin{array}{rllllll}
      \kconst_1 \tup{e}
      & =\, \hdge{\tup{\slangup}}\cdot\tup{v}  - \kconst_1 \htup{r} 
      &=\, 
       \slang^2(u -\tfrac{\kconst_1 }{\slang^2} )\htup{q} \,-\,   w \hdge{\tup{\slangup}}\cdot\htup{q}
      & \simeq\,  \slang^2(u -\tfrac{\kconst_1 }{\slang^2} )\tup{q} \,+\,   w \tup{p} 
\\[4pt]
    \kconst_1\tup{h} 
    & =\, \tup{v} + \tfrac{\kconst_1}{\slang^2} \hdge{\tup{\slangup}}\cdot \htup{r}
    & =\, 
     -(u -\tfrac{\kconst_1 }{\slang^2} )\hdge{\tup{\slangup}}\cdot\htup{q} \,-\,   w \htup{q} 
     & \simeq\, (u - \tfrac{\kconst_1 }{\slang^2})\tup{p} \,-\,  w \tup{q} 
\\[4pt]
    \tup{\slangup} & =\, \hdge{\tup{v}}\cdot\tup{r} 
     & =\,  \hdge{\tup{p}}\cdot\tup{q}
\end{array}
&&
w:=u^2 p_\ss{u}
\qquad
\end{flalign}
\end{small}
where ``\eq{\simeq}'' denotes expressions simplified using \eq{q=1} and \eq{\lambda=\htup{q}\cdot\tup{p}=0}.
It is well-known that the above \eq{\mbb{R}^3}-valued functions are all conserved for pure Kepler dynamics.
Their normalization defines an orthonormal basis, \eq{:\{\htup{o}_e,\htup{o}_h,\htup{o}_\slang\} := \{\htup{e},\htup{h},\htup{\slangup}\}}, often called the \textit{perifocal basis}. For Keplerian motion, this basis is 
constant/stationary.\footnote{with  \eq{\htup{e}} and \eq{\htup{h}} defining the orbital plane;  \eq{\htup{e}} directed towards periapsis,  \eq{\htup{h}} directed along the velocity at periapsis (direction tangent to the orbit), and \eq{\htup{\slangup}} directed normal to the orbit plane, completing the right-handed triad. The magnitude \eq{e=\mag{\tup{e}}} is the usual dimensionless eccentricity of the orbit. }
We may verify this by substituting into the above the (simplified) Kepler solutions from Eq.\eqref{qusol_s_2bp_simp}, leading to:
\begin{small}
\begin{align} \label{LRL_e_qp0}
\begin{array}{rllllll}
      \kconst_1 \tup{e}_\tau
      & \simeq\, 
      \slang^2_\tau (u_\tau -\tfrac{\kconst_1 }{\slang^2_\tau} )\tup{q}_\tau \,+\,   w_\tau \tup{p}_\tau 
      &\simeq\, \slang^2_\zr (u_\zr - \tfrac{\kconst_1}{\slang^2_0}) \tup{q}_\zr \,+\, w_\zr \tup{p}_\zr
\\[4pt]
    \kconst_1\tup{h}_\tau 
     & \simeq\, (u_\tau - \tfrac{\kconst_1 }{\slang^2_\tau})\tup{p}_\tau \,-\,  w_\tau \tup{q}_\tau 
     & \simeq\, (u_\zr - \tfrac{\kconst_1 }{\slang^2_0})\tup{p}_\zr \,-\,  w_\zr \tup{q}_\zr 
\\[4pt]
    \tup{\slangup}_\tau 
     & =\,  \hdge{\tup{p}}_\tau\cdot\tup{q}_\tau 
     & =\, \hdge{\tup{p}}_\zr \cdot\tup{q}_\zr 
\end{array}
\end{align}
\end{small}
verifying that \eq{(\tup{e},\tup{h},\tup{\slangup})=(\tup{e}_\zr,\tup{h}_\zr,\tup{\slangup}_\zr)} are indeed conserved for pure Kepler motion. 
Furthermore, if \eq{\tau_\zr=0} coincides with periapsis such that \eq{\tau} is the actual true anomaly, then \eq{\dot{r}_\zr=0} and \eq{w_\zr = \pdt{u}_\zr = -\dot{r}_\zr = 0 } vanishes and two of the above relations simplify to the following (with \eq{\slang=\slang_\zr}):
\begin{small}
\begin{flalign} \label{LRL_e_qp0_periapsis}
\quad
\begin{array}{cccc}
     \fnsize{if  $\tau_\zr=0$} \\
     \fnsize{is periapsis}
\end{array}  
\!\!:
&&
w_\zr = -\dot{r}_\zr = 0  
\qquad \Rightarrow \qquad 
\begin{array}{rllllll}
      \kconst_1 \tup{e}
       & \simeq\, \slang^2 (u_\zr - \tfrac{\kconst_1}{\slang^2}) \tup{q}_\zr 
\\[4pt]
    \kconst_1\tup{h} 
      & \simeq\, (u_\zr - \tfrac{\kconst_1 }{\slang^2})\tup{p}_\zr 
\end{array}
&&
\end{flalign}
\end{small}
Note that \eq{\kconst_1/\slang^2}, which appears frequently, is the inverse of the well-known \textit{semilatus rectum}, \eq{P_\ss{\mrm{slr}} = \slang^2/\kconst_1}. Still assuming \eq{\tau_\zr=0} coincides with periapsis (i.e., \eq{r_\zr = r_\ss{\mrm{min}}}), then we may use the conic section formula, \eq{r_\tau = P_\ss{\mrm{slr}}/(1+e\csn{\tau})}, to obtain:
\begin{small}
\begin{flalign} \label{Pslr_qp0_periapsis}
\quad
\begin{array}{cccc}
     \fnsize{if  $\tau_\zr=0$} \\
     \fnsize{is periapsis}
\end{array}  
\!\!:
&&
\begin{array}{llll}
    u_\tau -  \tfrac{\kconst_1}{\slang^2}  =  \tfrac{1}{r_\tau} - \tfrac{1}{P_{\mrm{slr}}} 
     =  \tfrac{\kconst_1}{\slang^2} e \csn{\tau}
&\quad,\qquad 
  u_\zr -  \tfrac{\kconst_1}{\slang^2}  =  \tfrac{1}{r_{\mrm{p}}} - \tfrac{1}{P_{\mrm{slr}}} 
   =  \tfrac{\kconst_1}{\slang^2} e 
   \,=\, \tfrac{e}{(1+e) r_{\mrm{p}} }
   \,=\, \tfrac{e}{(1+e)  } u_\ss{\mrm{p}}
\end{array}
&&
\end{flalign}
\end{small}
where \eq{e=\mag{\tup{e}}} is the classic eccentricity, and  where \eq{(\cdot)_\zr = (\cdot)_\ss{\mrm{p}}} corresponds to periapsis — that is, \eq{r_\zr = r_\ss{\mrm{p}}=r_\ss{\mrm{min}} = P_\ss{\mrm{slr}}/(1+e)} and \eq{u_\zr = u_\ss{\mrm{p}}=u_\ss{\mrm{max}} =(1+e)/P_\ss{\mrm{slr}}} correspond to the closest approach. 
We may also consider the above within specific cases (e.g., elliptic, parabolic, etc.):
\begin{small}
\begin{flalign}
\quad
\begin{array}{cccc}
     \fnsize{if  $\tau_\zr=0$} \\
     \fnsize{is periapsis}
\end{array}  
\!\!:
&&
\begin{array}{rlll}
     \fnsize{hyperbolic } (e>1): &\quad
  u_\zr -  \tfrac{\kconst_1}{\slang^2} =
  \tfrac{\kconst_1}{\slang^2} e 
   \,=\, \tfrac{e}{1+e} u_\ss{\mrm{p}}
   \;=\;
   \fnsize{same as Eq.\eqref{Pslr_qp0_periapsis}.}
\\[4pt]
     \fnsize{elliptic } (e < 1): &\quad
   u_\zr -  \tfrac{\kconst_1}{\slang^2} = \tfrac{\kconst_1}{\slang^2} e  =  \tfrac{1}{2} \tfrac{r_{\mrm{a}} - r_{\mrm{p}} }{ r_{\mrm{a}}r_{\mrm{p}}}
    =  \tfrac{1}{2}( \tfrac{1}{r_{\mrm{p}}} - \tfrac{1}{r_{\mrm{a}}})
    \;= 
   \tfrac{1}{2}( u_\ss{\mrm{p}} - u_\ss{\mrm{a}})
\\[4pt]
  \fnsize{parabolic } (e =1): &\quad
  u_\zr -  \tfrac{\kconst_1}{\slang^2} = \tfrac{1}{2 r_{\mrm{p}}}  \,=\, \tfrac{1}{2} u_\ss{\mrm{p}}
  \,=\, \tfrac{1}{2} u_\zr
\\[4pt]
     \fnsize{circular } (e =0): &\quad
  u_\zr -  \tfrac{\kconst_1}{\slang^2} =  u_\tau -  \tfrac{\kconst_1}{\slang^2} \,=\, 0
\end{array}
&&
\end{flalign}
\end{small}
The equalities indicated above for the elliptic  case (\eq{e<1}) have used the well-known elliptic relations \eq{P_\ss{\mrm{slr}} = \slang^2/\kconst_1 =2 \tfrac{ r_{\mrm{a}}r_{\mrm{p}}}{r_{\mrm{a}} + r_{\mrm{p}} }  } and \eq{e = \tfrac{ r_{\mrm{a}} - r_{\mrm{p}}}{r_{\mrm{a}} + r_{\mrm{p}} } }, where \eq{r_\ss{\mrm{a}}=r_\ss{\mrm{max}}} and \eq{u_\ss{\mrm{a}} = 1/r_\ss{\mrm{a}}=u_\ss{\mrm{min}}} correspond to apoapsis (farthest approach in an elliptic/closed orbit).
Although apoapsis only exists in truth for elliptic orbits (\eq{e<1}), note the above elliptic relations also hold in the parabolic case (\eq{e=1}) using the limits \eq{r_\ss{\mrm{a}}\mapsto \infty} and, thus, \eq{u_\ss{\mrm{a}}=1/r_\ss{\mrm{a}}\mapsto 0}.

%% file: Mysecs_prj/new_3.4_prjSTM.tex
\subsection{Matrix relations \& the Kepler state transition matrix} \label{sec:prj_STM_new}


Kepler dynamics and solutions in projective coordinates were given in section \ref{sec:kep_sol}. We now re-frame those results with more emphasis on their matrix structure as a linear ODE and present the corresponding Kepler state transition matrix (STM). 
First, we recall some basics of linear systems on real coordinate space, \eq{\mbb{R}^\en}.

\paragraph*{Some Review.}
Consider some \eq{\en}-dim linear autonomous and \textit{in}homogeneous ODE, \eq{\diff{}{\varep}\tup{x} = M\cdot\tup{x} + \tup{k}}, where \eq{M\in\mbb{R}^{\en\times\en}} and \eq{\tup{k}\in\mbb{R}^\en} are constant and were \eq{\varep} is some evolution parameter (independent variable). The solution to any such system may be expressed in terms of initial conditions \eq{\tup{x}_{\zr}} at \eq{\varep=0} as:
\begin{small}
\begin{flalign} \label{lin_system_gen}
\qquad
   \diff{}{\varep}\tup{x} = \tup{X}(\tup{x}) \,=\, M\cdot\tup{x} + \tup{k}
   \quad\;\; \Rightarrow \quad\;\; 
    \tup{x}_{\varep} = \phi_\varep(\tup{x}_{\zr}) \,=\,\Sigma_\varep  \cdot \tup{x}_{\zr} + \tup{\sigup}_\varep
    &&,&&
    \begin{array}{llll}
        \Sigma_\varep := \mrm{e}^{M \varep} \in\Glmat{\en}
         \;\;,\quad 
         \inv{\Sigma_\varep} =  \Sigma_{-\varep} 
     \\[4pt] 
     \tup{\sigup}_{\varep} := \int_{\zr}^{\varep}  \inv{\Sigma_\varep} \cdot \tup{k}\, \mrm{d} \varep 
    \end{array}
\end{flalign}
\end{small}
where \eq{\tup{X}} is the linear inhomogeneous vector field\footnote{In this work, a vector field \eq{\tup{X}} on \eq{\mbb{R}^\en} may be regarded simply as a smooth map \eq{\tup{X}:\mbb{R}^\en\to\mbb{R}^\en}. }
for \eq{M} and \eq{\tup{k}},
with \eq{\varep}-parameterized solution flow \eq{\phi_\varep:\mbb{R}^\en\to\mbb{R}^\en} given in terms of \eq{\Sigma_\varep} and \eq{\tup{\sigup}_\varep} defined as above. 
The matrix \eq{\Sigma_\varep =\mrm{e}^{M \varep} } is sometimes referred to as the STM. Yet, this term is only truly accurate in the linear \textit{homogeneous} case, \eq{ \diff{}{\varep}\tup{x} = M\cdot\tup{x}} (when \eq{\tup{k}=0=\tup{\sigup}_\varep}), which has solution \eq{\tup{x}_{\varep} = \Sigma_\varep   \cdot \tup{x}_{\zr}}. The term STM is also used more generally to refer to
a matrix \eq{\Phi_\varep:=\pderiv{\phi_\varep}{\tup{x}}\big|_{\tup{x}_0}} which, by abuse of notation, one might write as \eq{\pderiv{\tup{x}_{\varep}}{\tup{x}_0}}. 
This matrix itself satisfies the following ODE:
\begin{small}
\begin{flalign} \label{STM_ODE}
  \Phi_\varep:=\pderiv{\phi_\varep}{\tup{x}}\big|_{\tup{x}_0} \equiv   \pderiv{\tup{x}_{\varep}}{\tup{x}_0} \in\Glmat{\en}
\qquad,\qquad
  \diff{}{\varep} \Phi_\varep \,=\, \pderiv{\tup{X}}{\tup{x}} \cdot \Phi_\varep
\qquad,\qquad
  \Phi_\zr \,=\, \imat_\en
\end{flalign}
\end{small}
We use the term STM to refer to \eq{\Phi_\varep =\pderiv{\tup{x}_{\varep}}{\tup{x}_0}}, not \eq{\Sigma_\varep =\mrm{e}^{M \varep}}. 
In the case that \eq{M} and \eq{\tup{k}} in Eq.\eqref{lin_system_gen} are truly numeric constants, then one has \eq{\Phi_\varep = \Sigma_\varep } and \eq{\pderiv{\tup{X}}{\tup{x}} = M}. However, in the following, \eq{M} will instead be a matrix of integrals of motion. While this still allows for \eq{M} to be treated as constant and to use Eq.\eqref{lin_system_gen}, it means that  \eq{\Phi_\varep \neq \Sigma_\varep } and \eq{\pderiv{\tup{X}}{\tup{x}} \neq M}. 

\paragraph*{Keplerian Matrix ODEs \& Solutions (Unsimplified).}
Returning to the problem at hand, 
we note the \eq{s}- or \eq{\tau}-parameterized unperturbed
Kepler dynamics from Eq.\eqref{qpw_kep_EOM} can indeed be posed in the above general form of  Eq.\eqref{lin_system_gen} — though the matrix \eq{M} in Eq.\eqref{lin_system_gen} will be a matrix of integrals of motion, not literal real numbers. 
We already examined this for the ``\eq{(\tup{q},\tup{p})}-part'' of the dynamics in Eq.\eqref{dqdp_nonsimp}-Eq.\eqref{qp_sol_constants_simp}.
On that note, it will now be convenient to re-order our eight  phase space coordinates by splitting them into a ``\eq{(\tup{q},\tup{p})}-part'' (rotational motion) and ``\eq{(u,p_\ss{u})}-part'' (radial motion), where the latter may also be replaced by the pair \eq{(u,w:=u^2 p_\ss{u})}.

\begin{notesq}
    \textit{Coordinate ordering.}
    Rather than the standard configuration-momentum split ordering, we will adopt a modified ordering:
\begin{small}
\begin{align}
\begin{array}{cc}
      \fnsize{standard}  \\
     \fnsize{ordering}
\end{array}
\left\{\quad 
\begin{array}{llll}
      \tup{z} =  (\tup{q},u,\tup{p},p_\ss{u}) =(\bartup{q},\bartup{p}) 
\\[4pt]
      \tup{x} = (\tup{q},u,\tup{p},w) 
      \simeq (\bartup{q},\pdt{\bartup{q}})
\end{array}\right.
&&,&&
\begin{array}{cc}
      \fnsize{modified}  \\
     \fnsize{ordering}
\end{array}
\left\{\quad 
\begin{array}{llll}
      \tup{z} =  (\tup{q},\tup{p},u,p_\ss{u}) 
\\[4pt]
       \tup{x} = (\tup{q},\tup{p},u,w) 
      \simeq (\tup{q},\pdt{\tup{q}},u,\pdt{u})
\end{array}\right.
\end{align}
\end{small}
where \eq{\tup{z}\in\mbb{R}^8} denotes the set of ``canonical projective coordinates'' and \eq{\tup{x}\in\mbb{R}^8} denotes said coordinates with \eq{p_\ss{u}} replaced by the quasi-momentum coordinate \eq{w=u^2 p_\ss{u}}. 
The following developments will mostly use \eq{ \tup{x} = (\tup{q},\tup{p},u,w)}.
\end{notesq}

\noindent Recall that \textit{first}-order Kepler dynamics for \eq{\tup{z}} are not fully linear due to nonlinearity in the \eq{(u,p_\ss{u})} dynamics, but that those of \eq{\tup{x}} are indeed fully linear. They are fully linear whether or not one simplifies the ODEs using the integrals of motion \eq{q=1} and \eq{\lambda=\htup{q}\cdot\tup{p}=0}. We consider first the unsimplified case, then the simplified case.



The ODEs in Eq.\eqref{qpw_kep_EOM} and solutions in Eq.\eqref{qusol_s_2bp} describe unperturbed (and unsimplified) Kepler dynamics using the coordinates \eq{\tup{x}=(\tup{q},\tup{p},u,w)}. 
The ODEs from Eq.\eqref{qpw_kep_EOM} were given by:
\begin{small}
\begin{flalign} \label{uw_eom_kep_yetagain}
\begin{array}{cc}
     \fnsize{Kepler-type }  \\
      \fnsize{dynamics} 
\end{array} 
&&
\begin{array}{llllll}
     \pdt{\tup{q}}  \,=\, 
  -\hdge{\tup{\slangup}}\cdot\tup{q} 
  &,\quad 
  \pdt{\tup{p}}  \,=\, 
  -\hdge{\tup{\slangup}}\cdot\tup{p} 
\\[4pt]
   \pdt{u} \,=\, w
    &,\quad 
     \pdt{w} = -\slang^2 u + \kconst_1  
\end{array}
\qquad\quad\fnsize{or,}  \qquad\quad
 \begin{array}{llllll}
     \rng{\tup{q}}  \,=\, 
  -\hdge{\htup{\slangup}}\cdot\tup{q} 
  &,\quad 
  \rng{\tup{p}}  \,=\, 
  -\hdge{\htup{\slangup}}\cdot\tup{p} 
\\[4pt]
   \rng{u} \,=\,  w/\slang 
    &,\quad 
     \rng{w} = -\slang u + {\kconst_1}/{\slang} 
\end{array} 
&&
\end{flalign}
\end{small}
with \eq{\tup{\slangup}=\tup{\slangup}_\zr} an integral of motion. 
The above, and their solutions, are equivalent to the following: 
\begin{small}
\begin{align} \label{dqp_mat}
   \diff{}{s} \fnpmat{
        \tup{q} \\[2pt]
        \tup{p} }
    = L \cdot \fnpmat{
        \tup{q}\\[2pt]
        \tup{p} }
    \quad &\fnsize{or,} \quad 
      \diff{}{\tau} \fnpmat{
        \tup{q} \\[2pt]
        \tup{p} }
    = \tfrac{1}{\slang} L \cdot \fnpmat{
        \tup{q}\\[2pt]
        \tup{p} }    
&& \Rightarrow &&
     \fnpmat{
        \tup{q}_s \\[2pt]
        \tup{p}_s } 
        = \mrm{e}^{L s} \cdot  \fnpmat{
        \tup{q}_\zr \\[2pt]
        \tup{p}_\zr } 
        \;\;=\;\;
        \mrm{e}^{\frac{1}{\slang} L \tau} \cdot  \fnpmat{
        \tup{q}_\zr \\[2pt]
        \tup{p}_\zr } 
         =
        \fnpmat{
        \tup{q}_\tau \\[2pt]
        \tup{p}_\tau }
\\[6pt] \nonumber
      \diff{}{s} \fnpmat{u \\ w}
    = \Lambda_2 \cdot \fnpmat{u \\ w} +  \fnpmat{0 \\ \kconst_1}
    \quad &\fnsize{or,} \quad
     \diff{}{\tau} \fnpmat{u \\ w}
    = \tfrac{1}{\slang} \Lambda_2 \cdot \fnpmat{u \\ w} +  \tfrac{1}{\slang}\fnpmat{0 \\ \kconst_1}
   && \Rightarrow &&
    \fnpmat{u_s \\[2pt] w_s} = \mrm{e}^{\Lambda_2 s} \cdot \fnpmat{u_\zr \\[2pt] w_\zr } 
    + \fnpmat{ \sig^u_s \\[2pt] \sig^w_s}
  \;\;=\;\;
   \mrm{e}^{\frac{1}{\slang}\Lambda_2 \tau} \cdot \fnpmat{u_\zr \\[2pt] w_\zr } 
    + \fnpmat{ \sig^u_\tau \\[2pt] \sig^w_\tau }
    =   \fnpmat{u_\tau \\[2pt] w_\tau} 
\end{align}
\end{small}
where \eq{\tau=\slang s} (assuming \eq{\tau_\zr=s_\zr=0}) and where \eq{L} and \eq{\Lambda_2} are matrices of integrals of motion that depend only on \eq{\tup{\slangup}=\tup{\slangup}_\zr}. They are given as follows, along with the other terms appearing above:
\begin{small}
\begin{align} \label{Lmat_spso_again}
\begin{array}{llllll}
      L := \smpmat{
        -\hdge{\tup{\slangup}}  & 0 \\  0 & -\hdge{\tup{\slangup}}  }  \in \spmat{6} \cap \somat{6} 
&,\qquad 
     \mrm{e}^{L s} = \mrm{e}^{\frac{1}{\slang}L \tau} = \smpmat{
          R_{\tau}(\htup{\slangup}) & 0 \\
          0 &    R_{\tau}(\htup{\slangup})  }
    \in \Spmat{6} \cap \Somat{6}  
    &,\qquad
     R_{\tau}(\htup{\slangup}):= \mrm{e}^{-\hdge{\htup{\slangup}} \tau} 
     \in \Somat{3}
\\[12pt]
         \Lambda_2 := \smpmat{ 0 & 1 \\  -\slang^2  & 0 } \in\spmat{2}
 &,\qquad 
        \mrm{e}^{\Lambda_2 s} =  \mrm{e}^{ \frac{1}{\slang} \Lambda_2 \tau} = 
        \smpmat{
            \csn{\tau}  & \tfrac{1}{\slang}\snn{\tau}  \\ 
            -\slang \snn{\tau}  &  \csn{\tau} } 
            \in \Spmat{2}
&,\qquad 
      \smpmat{ \sig^u_\tau \\[2pt] \sig^w_\tau }
     = \smpmat{ 
        \tfrac{\kconst_1}{\slang^2}(1-\csn{\tau}) \\
          \tfrac{\kconst_1}{\slang}\snn{\tau}  }
\end{array}
\end{align}
\end{small}
where the inhomogeneous term \eq{(\sig^u ,\sig^w)} is defined as usual 
(see footnote\footnote{As per the general formulas in Eq.\eqref{lin_system_gen}, \eq{(\sig^u,\sig^w)} is defined by:
     \eq{\quad  \fnpmat{ \sig^u_s \\[2pt] \sig^w_s } := \int_0^s \big[ \mrm{e}^{-\Lambda_2 s} \cdot \fnpmat{0 \\ \kconst_1} \big] \mrm{d} s
     \,=\,  \fnpmat{ \sig^u_\tau \\[2pt] \sig^w_\tau } = \int_0^\tau \big[ \mrm{e}^{-\Lambda_2 \tau /\slang} \cdot \fnpmat{0 \\ \kconst_1/\slang} \big] \mrm{d} \tau } .
}),
and where we have specified the dimension, 2, on \eq{\Lambda_2} because this matrix reappears several times with different dimensions
 (see footnote\footnote{For any finite dimension \eq{\en}, we define \eq{\Lambda_{2\en}\in\spmat{2\en}} and \eq{\mrm{e}^{\Lambda_{2\en}\tau/\slang} \in \Spmat{2\en} } as follows: 
\begin{align} \label{ThisMat_Ngen}
   \Lambda_{2\en} := \fnpmat{0 & \imat_\en \\ -\slang^2 \imat_\en & 0 }
     \in \spmat{2\en}
    \qquad,\qquad 
    \mrm{e}^{\Lambda_{2\en} s} = \mrm{e}^{\Lambda_{2\en}\tau/\slang}  
    = \fnpmat{ \csn{\tau} \,\imat_\en & \tfrac{1}{\slang}\snn{\tau}\,\imat_\en \\ 
            -\slang \snn{\tau}\,\imat_\en &  \csn{\tau}\,\imat_\en  }
            \in \Spmat{2\en}
\end{align}
where \eq{\Lambda_{2\en}} has eigenvalues \eq{\pm\mrm{i}\slang} (\eq{\en} of each sign)
}). 
The above \eq{L}, \eq{\mrm{e}^{Ls}=\mrm{e}^{\frac{1}{\slang} L \tau}}, and \eq{R_{\tau}(\htup{\slangup})}  were discussed in more detail around Eq.\eqref{Lmat_spso}. 
The dynamics in Eq.\eqref{dqp_mat} can then be combined as a single linear inhomogeneous  ODE:
\begin{small}
\begin{flalign} \label{qusol_MAT_2bp_nonsimp}
\quad
\tup{x}=(\tup{q},\tup{p},u,w)
&&
   \boxed{
\begin{array}{rlll}
        &\pdt{\tup{x}} \,=\, M\cdot\tup{x} + \tup{k} 
\\[5pt] 
   \fnsize{or,} &\rng{\tup{x}} \,=\, \tfrac{1}{\slang}M \cdot \tup{x} \,+\, \tfrac{1}{\slang}\tup{k} 
\end{array} }
    \quad  \Rightarrow \qquad 
     \tup{x}_\tau \,=\, \Sig_\tau \cdot \tup{x}_\zr \,+\, \tup{\sigup}_\tau 
&&
\tau = \slang s
\quad
\end{flalign} 
\end{small}
with \eq{\tup{k}\in\mbb{R}^8} a numeric constant, with \eq{M = M_\zr\in \mbb{R}^{8 \times 8}} a matrix of integrals of motion given below, with \eq{\Sig_\tau := \mrm{e}^{\frac{1}{\slang}M \tau} = \mrm{e}^{M s}},  and with \eq{\tup{\sigup}_\tau := \int_{\zr}^{\tau} \inv{\Sig_\tau} \cdot \tup{k}\, \mrm{d} \tau  \in \mbb{R}^8}  defined as in  Eq.\eqref{lin_system_gen}. These are given for the ordering \eq{\tup{x}=(\tup{q},\tup{p},u,w)} 
as:\footnote{Using the submatrices defined in Eq.\eqref{Lmat_spso_again}, \eq{M} and \eq{\Sig_\tau} are given explicitly by: 
\begin{align} \nonumber
\begin{array}{lllll}
        M 
=
\fnpmat{
        \fnpmat{
        -\hdge{\tup{\slangup}}  & 0_{3\times 3} \\  0_{3\times 3} & -\hdge{\tup{\slangup}}  }
      & 0_{6\times 2}
  \\
        0_{2\times 6}  
        & \fnpmat{ 0 & 1 \\  -\slang^2  & 0 }
    }     
        \in \spmat{6;2}
   \qquad,\qquad 
       \Sig_\tau := \mrm{e}^{\frac{1}{\slang}M \tau}
 =
 \fnpmat{ 
         \fnpmat{
          R_{\tau}(\htup{\slangup}) & 0_{3\times 3} \\
          0_{3\times 3} &    R_{\tau}(\htup{\slangup})  }
             & 0_{6\times 2}
     \\
        0_{2\times 6}
        & \fnpmat{
            \csn{\tau}  & \tfrac{1}{\slang}\snn{\tau}  \\ 
            -\slang \snn{\tau}  &  \csn{\tau} }
     }
         \in \Spmat{6;2}
\end{array}
\end{align}
}
\begin{small}
\begin{align} \label{KepMats_nonsimp}
\begin{array}{lllll}
        M :=  \smpmat{ L & 0 \\ 0 & \Lambda_2 }   
        \in \spmat{6;2}
    &,\qquad 
      \Sig_\tau := \mrm{e}^{\frac{1}{\slang}M \tau}
     =  \smpmat{
           \mrm{e}^{\frac{1}{\slang}L \tau}   & 0 \\
           0 & \mrm{e}^{\frac{1}{\slang}\Lambda_2 \tau} 
           }
            \in \Spmat{6;2}
    &,\qquad
    \tup{k} := \smpmat{  \tup{0} \\  \tup{0} \\ 0  \\ \kconst_1}
    &,\qquad
     \tup{\sigup}_\tau = \smpmat{ \tup{0} \\
        \tup{0} \\
         \tfrac{\kconst_1}{\slang^2}(1-\csn{\tau}) \\
          \tfrac{\kconst_1}{\slang}\snn{\tau}  }
\end{array}
\end{align}
\end{small}
where \eq{\Spmat{6;2}} and \eq{\spmat{6;2}} denote, respectively, symplectic and Hamiltonian matrices with regards to \eq{\fnpmat{ \jmat_6 &0 \\ 0 & \jmat_2} } 
    (see footnote\footnote{We define \eq{\Spmat{6;2}} and \eq{\spmat{6;2}} as follows,  where \eq{\jmat_{2\en}\in \Spmat{2\en}} denotes the standard symplectic matrix on \eq{\mbb{R}^{2\en}}: 
    \begin{align} \label{Sp_alt_def}
       \Spmat{6;2} :=\; \big\{\, S \in \mbb{R}^{8\times 8}    \;\big|\;  \trn{S}\jmat_\ss{6;\!2} S = \jmat_\ss{6;\!2} \, \big\}
    \qquad,\qquad 
       \spmat{6;2}  :=\; \big\{\, H \in \mbb{R}^{8\times 8}    \;\big|\;  \jmat_\ss{6;\!2} H + \trn{H} \jmat_\ss{6;\!2} = 0 \,\big\}
        \qquad,\qquad 
        \jmat_\ss{6;\!2} := \jmat_6\oplus \jmat_2 = \fnpmat{
        \jmat_6 &0 \\ 0 & \jmat_2}
    \end{align} }).
Note that \eq{M} has eigenvalues \eq{(0,0,\mrm{i}\slang,\mrm{i}\slang,-\mrm{i}\slang,-\mrm{i}\slang)} (from \eq{L}) and \eq{(\mrm{i}\slang,-\mrm{i}\slang)} (from \eq{\Lambda_2}). 
It can be verified that the solution in Eq.\eqref{qusol_MAT_2bp_nonsimp}, when written out explicitly, agrees with that given previously in Eq.\eqref{qusol_s_2bp}: 
\begin{small}
\begin{flalign}  \label{qusol_MAT_2bp_inv_nonsimp}
\qquad
\begin{array}{ccccc}
    \boxed{ \tup{x}_\tau = \phi_\tau (\tup{x}_\zr) = \Sig_\tau \cdot \tup{x}_\zr + \tup{\sigup}_\tau 
    \qquad \leftrightarrow \qquad 
    \tup{x}_\zr = \inv{\phi_\tau} (\tup{x}_\tau) = \inv{\Sig_\tau} \cdot \tup{x}_\tau + \tup{\varsigup}_\tau }
\\[10pt]
\begin{array}{lllllll}
    \tup{q}_\tau
       \,=\, 
      \tup{q}_\zr \csn{\tau} - \hdge{\htup{\slangup}}_\zr\cdot\tup{q}_\zr \snn{\tau} 
\\[4pt]
       \tup{p}_\tau
       \,=\, 
      \tup{p}_\zr \csn{\tau} - \hdge{\htup{\slangup}}_\zr \cdot\tup{p}_\zr \snn{\tau} 
\\[4pt]
      u_\tau
      \,=\,  u_{\zr} \csn{\tau} + \tfrac{1}{\slang_0}w_{\zr}\snn{\tau} + \tfrac{\kconst_1}{\slang_0^2}(1-\csn{\tau})
\\[4pt]
       w_\tau 
      \,=\, w_{\zr}\csn{\tau} -\slang_\zr u_{\zr} \snn{\tau} + \tfrac{\kconst_1}{\slang_0}\snn{\tau}
\end{array}
 \qquad \leftrightarrow \qquad 
\begin{array}{llll}
     \tup{q}_\zr \,=\, \tup{q}_\tau\csn{\tau} + \hdge{\htup{\slangup}}_\tau\cdot\tup{q}_\tau \snn{\tau} 
\\[4pt]
      \tup{p}_\zr \,=\, \tup{p}_\tau\csn{\tau} + \hdge{\htup{\slangup}}_\tau\cdot\tup{p}_\tau \snn{\tau}
\\[4pt]
     u_\zr \,=\,  u_\tau \csn{\tau} - \tfrac{1}{\slang_\tau}w_\tau\snn{\tau} + \tfrac{\kconst_1}{\slang_\tau^2}(1-\csn{\tau})
\\[4pt]
      w_\zr \,=\,    w_\tau\csn{\tau} + \slang_\tau u_\tau \snn{\tau}  - \tfrac{\kconst_1}{\slang_\tau}\snn{\tau}
\end{array}
\end{array}
&& \tup{\slangup}_\tau =\tup{\slangup}_\zr
\end{flalign}
\end{small}
where we have denoted by \eq{\phi_\tau:\mbb{R}^8\to\mbb{R}^8} the (unsimplified) Kepler flow in coordinates \eq{\tup{x}}, and where we note:
\begin{small}
\begin{align} \label{qusol_MAT_invTerm}
    \inv{\phi_\tau} =  \phi_{-\tau}
\qquad,\qquad 
    \inv{\Sig_\tau} =  \Sig_{-\tau}
\qquad,\qquad 
     \tup{\varsigup}_\tau := -\inv{\Sig_\tau} \cdot \tup{\sigup}_\tau  
           = \tup{\sigup}_{-\tau}
\end{align}
\end{small}

\begin{small}
\begin{notesq}
   \textit{Angular momentum dependence.} Though not explicit in our notation, many terms in the above Kepler dynamics depend on the angular momentum, \eq{\tup{\slangup}}. In particular:
    \begin{small}
    \begin{align} \label{KepODE_ang_depend1}
         M = M(\tup{\slangup}) 
         \qquad,\qquad 
         \Sig_\tau = \Sig_\tau(\tup{\slangup}) 
         \qquad,\qquad 
         \tup{\sigup}_\tau = \tup{\sigup}_\tau(\tup{\slangup})
    \end{align}
    \end{small}
    with \eq{\tup{\slangup} = \tup{\slangup}(\tup{x}) = \tup{q}\tms\tup{p}} a function of the coordinates. 
    Yet, since \eq{\tup{\slangup}=\tup{\slangup}_\zr} is conserved for any central-force dynamics, we have treated it as a constant when solving the Kepler ODEs. Still, it should be noted that Eq.\eqref{qusol_MAT_2bp_nonsimp} is really a family of linear ODEs parameterized by angular momentum values (set by initial conditions).
    For instance, for a given \eq{\tup{\slangup}_\zr}, it would be more accurate to write the ODE and the solution in  Eq.\eqref{qusol_MAT_2bp_nonsimp} as:
    \begin{small}
    \begin{align} \label{KepODE_ang_depend2}
        \rng{\tup{x}} \,=\, \tfrac{1}{\slang_0} M(\tup{\slangup}_\zr) \cdot \tup{x} \,+\, \tfrac{1}{\slang_0}\tup{k} 
    \qquad  \Rightarrow \qquad 
         \tup{x}_\tau \,=\, \phi_\tau (\tup{x}_\zr) \,=\, \Sig_\tau(\tup{\slangup}_\zr) \cdot \tup{x}_\zr \,+\, \tup{\sigup}_\tau(\tup{\slangup}_\zr)
         &&,&&
         \pderiv{\tup{x}_\tau}{\tup{x}_0} \neq \Sig_\tau
    \end{align}
    \end{small}
    This dependence on \eq{\tup{\slangup}} of the solution flow, \eq{\phi_\tau}, is  reflected in the way we have written the explicit solutions in Eq.\eqref{qusol_MAT_2bp_inv_nonsimp}.
\end{notesq}
\end{small}

\paragraph*{Kepler STM in Projective Coordinates.}
As indicated above, the matrices \eq{\Sig_\tau} and \eq{\pderiv{\tup{x}_\tau}{\tup{x}_0}} — both of which are sometimes referred to as the STM — do \textit{not} coincide. When perturbations are considered and \eq{\tup{\slangup}} is no longer conserved, then the object of interest is the more general version of the STM given by \eq{\Phi_\tau := \pderiv{\tup{x}_\tau}{\tup{x}_0}} (one might write this more precisely as \eq{ \pderiv{\phi_\tau}{\tup{x}}\big|_{\tup{x}_0}}), which we note is partitioned as:
\begin{small}
\begin{flalign} 
\qquad
\begin{array}{llll}
      \tup{x} =  (\tup{q},\tup{p},u,w)
\end{array}
&&
 \Phi_\tau := \pderiv{\phi_\tau}{\tup{x}}\big|_{\tup{x}_0} \equiv  \pderiv{\tup{x}_\tau}{\tup{x}_0}
     \,=
 \pmat{ 
     \pderiv{(\tup{q}_\tau,\tup{p}_\tau)}{(\tup{q}_0,\tup{p}_0)}
     &    \pderiv{(\tup{q}_\tau,\tup{p}_\tau)}{(u_0,w_0)} 
     \\[6pt]
     \pderiv{(u_\tau,w_\tau)}{(\tup{q}_0,\tup{p}_0)}
     &
      \pderiv{(u_\tau,w_\tau)}{(u_0,w_0)}
     }
\;=\;
 \smpmat{ 
     \smpmat{ 
        \pderiv{\tup{q}_\tau}{\tup{q}_0}
     & \pderiv{\tup{q}_\tau}{\tup{p}_0}
     \\[4pt]
      \pderiv{\tup{p}_\tau}{\tup{q}_0}
     & \pderiv{\tup{p}_\tau}{\tup{p}_0}    
     }_{6\times 6}
     & 
     \smpmat{
     \tup{0} &  \tup{0}  \\[4pt] \tup{0} &  \tup{0}
     }_{6\times 2}
 \\[14pt]
    \smpmat{ 
       \trn{\pderiv{u_\tau}{\tup{q}_0}}
     &  \trn{\pderiv{u_\tau}{\tup{p}_0}}
     \\[4pt]
       \trn{\pderiv{w_\tau}{\tup{q}_0}}
     &  \trn{\pderiv{w_\tau}{\tup{p}_0}}
    }_{2\times 6}
 &
    \smpmat{ 
     \pderiv{u_\tau}{u_0}
     &\pderiv{u_\tau}{w_0} 
     \\[4pt]
     \pderiv{w_\tau}{u_0}
     &\pderiv{w_\tau}{w_0}
    }_{2\times 2}
 }  
&&   
\end{flalign}
\end{small}
with \eq{\phi_\tau} the Kepler flow in Eq.\eqref{qusol_MAT_2bp_inv_nonsimp} and
where \eq{(u_\tau,w_\tau)} depend on \eq{(\tup{q}_\zr,\tup{p}_\zr)} only through \eq{\slang=\slang_\zr}. 
Using the solutions in Eq.\eqref{qusol_MAT_2bp_inv_nonsimp} — along with angular momentum relations from Eq.\eqref{l_qp_general}-Eq.\eqref{angmoment_rels_crd2} of Appx.~\ref{sec:ang_momentum} — the above leads to:
\begin{small}
\begin{align} \label{STM_prj_nosimp}
  \Phi_\tau = \pderiv{\tup{x}_\tau}{\tup{x}_0}
   &\,=\,\left.\fnpmat{ 
    \fnpmat{ 
    \imat_3 \csn{\tau} - \tfrac{1}{\slang} \big( \hdge{\tup{\slangup}} - \hdge{\tup{q}}\cdot\hdge{\tup{p}} -   (\hdge{\htup{\slangup}}\cdot\tup{q})\otms (\hdge{\htup{\slangup}}\cdot\tup{p}) \big) \snn{\tau} 
    &  \tfrac{q^2}{\slang}  \htup{\slangup}\otms\htup{\slangup}  \snn{\tau} 
    \\[3pt]
    -\tfrac{p^2}{\slang}  \htup{\slangup}\otms\htup{\slangup}  \snn{\tau} 
    &  \imat_3 \csn{\tau} -   \tfrac{1}{\slang} \big( \hdge{\tup{\slangup}} +  \hdge{\tup{p}}\cdot\hdge{\tup{q}} +   (\hdge{\htup{\slangup}}\cdot\tup{p})\otms (\hdge{\htup{\slangup}}\cdot\tup{q}) \big) \snn{\tau}  
    }
    &
    \fnpmat{ \tup{0} & \tup{0} \\[3pt] \tup{0} & \tup{0} }
\\[18pt]
    \fnpmat{ 
         \tfrac{1}{\slang}\big( \tfrac{2\kconst_1}{\slang^2}(1-\csn{\tau}) + \tfrac{1}{\slang}w \snn{\tau} \big)\tup{p}\cdot\hdge{\htup{\slangup}}
         &\;\;  -\tfrac{1}{\slang}\big( \tfrac{2\kconst_1}{\slang^2} (1-\csn{\tau}) + \tfrac{1}{\slang}w \snn{\tau} \big) \tup{q}\cdot\hdge{\htup{\slangup}}
     \\[3pt]
          (u + \tfrac{\kconst_1}{\slang^2})\snn{\tau}\,   \tup{p}\cdot\hdge{\htup{\slangup}}
         &   -(u + \tfrac{\kconst_1}{\slang^2})\snn{\tau}\, \tup{q}\cdot\hdge{\htup{\slangup}}
     }
    &
    \fnpmat{
    \csn{\tau} & \tfrac{1}{\slang}\snn{\tau}
    \\[3pt]
    -\slang \snn{\tau} & \csn{\tau}
    }
} \right|_{\tup{x}=\tup{x}_0}
\\[8pt] \nonumber
   \inv{\Phi_\tau} =  \pderiv{\tup{x}_0}{\tup{x}_\tau} 
&\,=\,
\left.\fnpmat{ 
    \fnpmat{ 
    \imat_3 \csn{\tau} + \tfrac{1}{\slang} \big( \hdge{\tup{\slangup}} - \hdge{\tup{q}}\cdot\hdge{\tup{p}} -   (\hdge{\htup{\slangup}}\cdot\tup{q})\otms (\hdge{\htup{\slangup}}\cdot\tup{p}) \big) \snn{\tau} 
    &  -\tfrac{q^2}{\slang}  \htup{\slangup}\otms\htup{\slangup}  \snn{\tau}  
    \\[3pt]
     \tfrac{p^2}{\slang}  \htup{\slangup}\otms\htup{\slangup}  \snn{\tau}  
    &  \imat_3 \csn{\tau} +   \tfrac{1}{\slang} \big( \hdge{\tup{\slangup}} +  \hdge{\tup{p}}\cdot\hdge{\tup{q}} +   (\hdge{\htup{\slangup}}\cdot\tup{p})\otms (\hdge{\htup{\slangup}}\cdot\tup{q}) \big) \snn{\tau}  
    } 
    &
    \fnpmat{ \tup{0} & \tup{0} \\[3pt] \tup{0} & \tup{0} }
\\[18pt]
    \fnpmat{ 
         \tfrac{1}{\slang}\big( \tfrac{2\kconst_1}{\slang^2}(1-\csn{\tau}) - \tfrac{1}{\slang}w \snn{\tau} \big) \tup{p}\cdot\hdge{\htup{\slangup}}
         &\;\;  -\tfrac{1}{\slang}\big( \tfrac{2\kconst_1}{\slang^2} (1-\csn{\tau}) - \tfrac{1}{\slang}w \snn{\tau} \big) \tup{q}\cdot\hdge{\htup{\slangup}}
     \\[3pt]
          - (u + \tfrac{\kconst_1}{\slang^2})\snn{\tau}\,   \tup{p}\cdot\hdge{\htup{\slangup}}
         &   (u + \tfrac{\kconst_1}{\slang^2})\snn{\tau}\,  \tup{q}\cdot\hdge{\htup{\slangup}}
     }
    &
    \fnpmat{
    \csn{\tau} & -\tfrac{1}{\slang}\snn{\tau}
    \\[3pt]
    \slang \snn{\tau} & \csn{\tau}
    }  
} \right|_{\tup{x}=\tup{x}_\tau}
=\, \Phi_{-\tau}\big|_{\tup{x}=\tup{x}_\tau}
\end{align}
\end{small}
The above follow from the unsimplified Kepler flow, \eq{\phi_\tau}, in Eq.\eqref{qusol_MAT_2bp_inv_nonsimp}.
If we now use the integrals of motion \eq{q=1} and  and \eq{\lambda=\htup{q}\cdot\tup{p}=0} — implying the relations in Eq.\eqref{qpl_rels_apx}-Eq.\eqref{angmoment_rels_simplified} of Appx.~\ref{sec:ang_momentum} —  to simplify the above, we obtain:
\begin{small}
\begin{align} \label{STM_prj_somesimp}
    \Phi_\tau 
    &\,\simeq\,
    \left.\fnpmat{ 
    \fnpmat{ 
    \imat_3 \csn{\tau} -  \hdge{\htup{\slangup}} \snn{\tau} 
    & \tfrac{1}{\slang} \hdge{\htup{\slangup}}\otms \hdge{\htup{\slangup}} \snn{\tau} 
    \\[3pt]
    -\slang \hdge{\htup{\slangup}}\otms \hdge{\htup{\slangup}} \snn{\tau} 
    &     \imat_3 \csn{\tau} -  \hdge{\htup{\slangup}} \snn{\tau} 
    } 
    &
    \fnpmat{ \tup{0} & \tup{0} \\[3pt] \tup{0} & \tup{0} }
\\[18pt]
   \fnpmat{ 
         -\big( \tfrac{2\kconst_1}{\slang^2}(1-\csn{\tau}) + \tfrac{1}{\slang}w \snn{\tau} \big) \trn{\tup{q}}
         &\;\;  -\tfrac{1}{\slang}\big( \tfrac{2\kconst_1}{\slang^2}(1-\csn{\tau}) + \tfrac{1}{\slang}w \snn{\tau} \big) \trn{\htup{p}}
     \\[3pt]
           -\slang (u + \tfrac{\kconst_1}{\slang^2})\snn{\tau}\, \trn{\tup{q}}
         &  - (u + \tfrac{\kconst_1}{\slang^2})\snn{\tau}\, \trn{\htup{p}}
     }
    &
    \fnpmat{
    \csn{\tau} & \tfrac{1}{\slang}\snn{\tau}
    \\[3pt]
    -\slang \snn{\tau} & \csn{\tau}
    }  
      }\right|_{\tup{x}=\tup{x}_0}
      \quad,\quad \slang \simeq p
\\[8pt] \nonumber 
   \inv{\Phi_\tau} 
   &\,\simeq\,
    \left.\fnpmat{ 
    \fnpmat{ 
    \imat_3 \csn{\tau} + \hdge{\htup{\slangup}} \snn{\tau} 
    &  -\tfrac{1}{\slang}\htup{\slangup}\otms\htup{\slangup} \snn{\tau} 
    \\[3pt]
   \slang \htup{\slangup}\otms\htup{\slangup}  \snn{\tau}
    &  \imat_3 \csn{\tau} +   \hdge{\htup{\slangup}} \snn{\tau} 
    } 
    &
    \fnpmat{ \tup{0} & \tup{0} \\[3pt] \tup{0} & \tup{0} }
\\[18pt]
     \fnpmat{ 
         -\big( \tfrac{2\kconst_1}{\slang^2}(1-\csn{\tau}) - \tfrac{1}{\slang}w \snn{\tau} \big) \trn{\tup{q}}
         &\;\;  -\tfrac{1}{\slang}\big( \tfrac{2\kconst_1}{\slang^2}(1-\csn{\tau}) - \tfrac{1}{\slang}w \snn{\tau} \big) \trn{\htup{p}}
     \\[3pt]
           \slang (u + \tfrac{\kconst_1}{\slang^2})\snn{\tau}\,  \trn{\tup{q}}
         &   (u + \tfrac{\kconst_1}{\slang^2})\snn{\tau}\, \trn{\htup{p}}
     }
    &
    \fnpmat{
    \csn{\tau} & -\tfrac{1}{\slang}\snn{\tau}
    \\[3pt]
    \slang \snn{\tau} & \csn{\tau}
    }  
} \right|_{\tup{x}=\tup{x}_\tau}
\end{align}
\end{small}
In all the above, the upper left sub-matrices have used the relation: 
\begin{small}
\begin{align}
    \htup{\slangup}\otms \htup{\slangup} 
    \;=\; \imat_3 - \htup{p}\otms\htup{p} - (\hdge{\htup{\slangup}}\cdot\htup{p})\otms (\hdge{\htup{\slangup}}\cdot\htup{p})
    \;=\; \imat_3 - \htup{q}\otms\htup{q} - (\hdge{\htup{\slangup}}\cdot\htup{q})\otms (\hdge{\htup{\slangup}}\cdot\htup{q})
    \;\simeq\; 
    \imat_3 - \htup{q}\otms\htup{q} - \htup{p}\otms\htup{p}
\end{align}
\end{small}

\begin{notesq}
    \textit{Canonical coordinate STM.} Let us denote by \eq{\Psi_\tau} the STM for the canonical coordinates \eq{\tup{z}=(\tup{q},\tup{p},u,p_\ss{u})}. This may be recovered from the STM \eq{\Phi_\tau} given above for \eq{\tup{x}=(\tup{q},\tup{p},u,w)} using:
    \begin{small}
    \begin{align}
        \Psi_\tau \,=\,
        \pderiv{\tup{z}_\tau}{\tup{z}_0} 
         \,=\, 
         \inv{(\pderiv{\tup{x}_\tau}{\tup{z}_\tau})} \cdot   \pderiv{\tup{x}_\tau}{\tup{x}_0}\cdot   \pderiv{\tup{x}_0}{\tup{z}_0}
         \,=\, 
         \pderiv{\tup{z}_\tau}{\tup{x}_\tau} \cdot   \Phi_{\tau} \cdot   \pderiv{\tup{x}_0}{\tup{z}_0}
    \end{align}
    \end{small}
    where the only transformation is \eq{w=u^2 p_\ss{u} \leftrightarrow p_\ss{u} = w/u^2} such that \eq{\pderiv{\tup{z}}{\tup{x}}} and  \eq{\pderiv{\tup{x}}{\tup{z}}} are given simply by
    \begin{small}
    \begin{align}
          \pderiv{\tup{x}}{\tup{z}} \,=\, 
          \fnpmat{ 
            \imat_6 & 0_{6\times 2} 
            \\
            0_{2 \times 6} 
            & 
            \fnpmat{1 &  0 \\
                   2 u p_\ss{u} & u^2
            } }
     \qquad,\qquad 
             \pderiv{\tup{z}}{\tup{x}}  = \inv{( \pderiv{\tup{x}}{\tup{z}})} 
          \,=\, 
          \fnpmat{ 
            \imat_6 & 0_{6\times 2} 
            \\
            0_{2 \times 6} 
            & 
            \fnpmat{1 &  0 \\
                   -2 p_\ss{u}/u & 1/u^2
            } }
    \end{align}
    \end{small}
\end{notesq}


\paragraph*{Simplified Matrix Relations.}
The above matrix ODEs and solutions are without simplifying anything using the integrals of motion \eq{q=1} or \eq{\lambda=\htup{q}\cdot\tup{p}=0}. 
The simplified expression for the STM in Eq.\eqref{STM_prj_somesimp} is still obtained from the unsimplified Kepler solutions; only the resulting \eq{\Phi_\tau =\pderiv{\tup{x}_\tau}{\tup{x}_0} } is simplified after differentiating.

We will now use \eq{q=1} and \eq{\htup{q}\cdot\tup{p}=0} (and thus \eq{\slang^2 \simeq p^2}) to simply the dynamics and solutions to begin with. 
The Kepler dynamics in Eq.\eqref{qpw_kep_EOM} then simplify to:
\begin{small}
\begin{flalign} \label{uw_eom_kep_simp}
\begin{array}{cc}
     \fnsize{Kepler-type }  \\
      \fnsize{dynamics} 
\end{array} 
&&
\begin{array}{llllll}
     \pdt{\tup{q}}  
     \,\simeq\, \tup{p}
  &,\quad 
  \pdt{\tup{p}}  
  \,\simeq\, - p^2 \tup{q}
  \,\simeq\, - \slang^2 \tup{q}
\\[4pt]
   \pdt{u} \,=\, w
    &,\quad 
     \pdt{w} = -\slang^2 u + \kconst_1  
\end{array}
\qquad\quad \fnsize{or,}  \qquad\quad
 \begin{array}{llllll}
     \rng{\tup{q}} 
     \,\simeq\, \tfrac{1}{\slang}\tup{p} \,\simeq\, \htup{p}
  &,\quad 
  \rng{\tup{p}}  
  \,\simeq\, -\slang \tup{q}
\\[4pt]
   \rng{u} \,=\,  w/\slang 
    &,\quad 
     \rng{w} = -\slang u + {\kconst_1}/{\slang} 
\end{array} 
&&
\end{flalign}
\end{small}
where \eq{\slang=\slang_\zr\simeq p_\zr = p} is preserved. The above is equivalent to the following matrix ODE:
\begin{small}
\begin{flalign} \label{qusol_MAT_2bp}
\quad
\tup{x}=(\tup{q},\tup{p},u,w)
&&
   \boxed{
\begin{array}{rlll}
        &\pdt{\tup{x}} \,\simeq\, \til{\Lambda} \cdot \tup{x} \,+\, \tup{k} 
\\[5pt] 
   \fnsize{or,} & \rng{\tup{x}} \,\simeq\, \tfrac{1}{\slang}\til{\Lambda} \cdot \tup{x} \,+\, \tfrac{1}{\slang}\tup{k} 
\end{array} }
    \quad  \Rightarrow \qquad 
    \tup{x}_\tau \simeq E_\tau \cdot \tup{x}_\zr \,+\, \tup{\sigup}_\tau  
&&
\tau = \slang s
\quad
\end{flalign} 
\end{small}
with \eq{\tup{k}} and \eq{\tup{\sigup}_\tau} the same as in Eq.\eqref{KepMats_nonsimp}, and 
where \eq{\til{\Lambda}=\til{\Lambda}_\zr} is a matrix of integrals of motion given as follows, along with \eq{E_\tau}:
\begin{small}
\begin{align}
        \til{\Lambda} := \smpmat{ 
        \Lambda_6 & 0 \\
        0 & \Lambda_2 }
          \in \spmat{6;2}
\qquad,\qquad 
       E_\tau := \mrm{e}^{\frac{1}{\slang}\til{\Lambda} \tau} =  \mrm{e}^{\til{\Lambda} s}
        = \smpmat{ \mrm{e}^{\frac{1}{\slang}\Lambda_6 \tau}  & 0 \\ 0 & \mrm{e}^{\frac{1}{\slang}\Lambda_2 \tau}  }
          \in \Spmat{6;2}
\qquad,\qquad
\begin{array}{rllll}
     \tup{k} & =\, \fnsize{Eq.\eqref{KepMats_nonsimp}}
     \\[3pt]
     \tup{\sigup}_\tau & =\, \fnsize{Eq.\eqref{KepMats_nonsimp}}
     \\[3pt]
      \tup{\varsigup}_\tau & := -\inv{E_\tau} \cdot \tup{\sigup}_\tau = \tup{\sigup}_{-\tau}
\end{array}
\end{align}
\end{small}
with \eq{\Lambda_{2\en}\in\spmat{2\en}} and \eq{\mrm{e}^{\Lambda_{2\en}\tau/\slang} \in \Spmat{2\en} } defined as in Eq.\eqref{ThisMat_Ngen} for any finite dimension \eq{\en}. 
When written out explicitly,  the solution in Eq.\eqref{qusol_MAT_2bp} agrees with that given previously in Eq.\eqref{qusol_s_2bp_simp}: 
\begin{small}
\begin{flalign}  \label{qusol_MAT_2bp_inv}
\quad
\begin{array}{ccccc}
     \boxed{ \tup{x}_\tau \simeq \theta_\tau(\tup{x}_\zr) = E_\tau \cdot \tup{x}_\zr + \tup{\sigup}_\tau 
    \qquad \leftrightarrow \qquad 
    \tup{x}_\zr \simeq \inv{\theta_\tau}(\tup{x}_\tau) = \inv{E_\tau} \cdot \tup{x}_\tau + \tup{\varsigup}_\tau }  
\\[10pt]
\begin{array}{lllllll}
    \tup{q}_\tau
       \,\simeq\, 
      \tup{q}_\zr \csn{\tau} + \tfrac{1}{\slang_0}\tup{p}_\zr \snn{\tau} 
\\[4pt]
       \tup{p}_\tau
       \,\simeq\, 
      \tup{p}_\zr \csn{\tau} - \slang_\zr \tup{q}_\zr \snn{\tau} 
\\[4pt]
      u_\tau
      \,=\,  u_{\zr} \csn{\tau} + \tfrac{1}{\slang_0}w_{\zr}\snn{\tau} + \tfrac{\kconst_1}{\slang_0^2}(1-\csn{\tau})
\\[4pt]
       w_\tau 
      \,=\,  w_{\zr}\csn{\tau}  -\slang_\zr u_{\zr} \snn{\tau} + \tfrac{\kconst_1}{\slang_0}\snn{\tau}
\end{array}
 \qquad \leftrightarrow \qquad 
\begin{array}{lll}
     \tup{q}_\zr \,\simeq\,  
     \tup{q}_\tau \csn{\tau} - \tfrac{1}{\slang_\tau}\tup{p}_\tau \snn{\tau} 
\\[4pt]
      \tup{p}_\zr \,\simeq\,  \tup{p}_\tau \csn{\tau} + \slang_\tau \tup{q}_\tau \snn{\tau} 
\\[4pt]
     u_\zr \,=\,  u_\tau \csn{\tau} - \tfrac{1}{\slang_\tau}w_\tau\snn{\tau} + \tfrac{\kconst_1}{\slang_\tau^2}(1-\csn{\tau})
\\[4pt]
      w_\zr \,=\,  w_\tau\csn{\tau} + \slang_\tau u_\tau \snn{\tau}  - \tfrac{\kconst_1}{\slang_\tau}\snn{\tau}
\end{array}
\end{array}
&&
\begin{array}{lll}
     \slang_\tau = \slang_\zr  \\[2pt]
     \slang \simeq p
\end{array}
\end{flalign}
\end{small}
where we have defined \eq{\theta_\tau} as the Kepler flow in coordinates \eq{\tup{x}} that has been simplified with \eq{q=1} and \eq{\htup{q}\cdot\tup{p}=0}. 
As before, 
note \eq{\inv{\theta_\tau} =  \theta_{-\tau}}, \eq{\inv{E_\tau} =  E_{-\tau}}, and \eq{\tup{\varsigup}_\tau = \tup{\sigup}_{-\tau}}.
\begin{small}
\begin{itemize}
      \item \textit{Angular momentum dependence.} The comments on angular momentum dependence around Eq.\eqref{KepODE_ang_depend1}-Eq.\eqref{KepODE_ang_depend2} still apply. Yet, for the above simplified relations, this dependence is only on the \textit{magnitude}, \eq{\slang}, which simplifies to \eq{\slang\simeq p}. 
      That is, \eq{\til{\Lambda} = \til{\Lambda}(\slang) }, \eq{ E_\tau = E_\tau(\slang)}, and \eq{\tup{\sigup}_\tau = \tup{\sigup}_\tau(\slang)}, with \eq{\slang\simeq p}. 
      To clarify, we could write the solutions in Eq.\eqref{qusol_MAT_2bp_inv} as:
      \begin{small}
    \begin{flalign}  \label{qusol_2bp_just_p}
    \begin{array}{ccccc}
    \begin{array}{lllllll}
        \tup{q}_\tau
           \simeq 
          \tup{q}_\zr \csn{\tau} + \tfrac{1}{p_0}\tup{p}_\zr \snn{\tau} 
    \\[4pt]
           \tup{p}_\tau
           \simeq 
          \tup{p}_\zr \csn{\tau} - p_\zr \tup{q}_\zr \snn{\tau} 
    \\[4pt]
          u_\tau
          \simeq  u_{\zr} \csn{\tau} + \tfrac{1}{p_0}w_{\zr}\snn{\tau} + \tfrac{\kconst_1}{p_0^2}(1-\csn{\tau})
    \\[4pt]
           w_\tau 
          \simeq  w_{\zr}\csn{\tau}  -p_\zr u_{\zr} \snn{\tau} + \tfrac{\kconst_1}{p_0}\snn{\tau}
    \end{array}
     \quad \leftrightarrow \quad 
    \begin{array}{lll}
         \tup{q}_\zr \simeq  
         \tup{q}_\tau \csn{\tau} - \tfrac{1}{p_\tau}\tup{p}_\tau \snn{\tau} 
    \\[4pt]
          \tup{p}_\zr \simeq  \tup{p}_\tau \csn{\tau} + p_\tau \tup{q}_\tau \snn{\tau} 
    \\[4pt]
         u_\zr \simeq  u_\tau \csn{\tau} - \tfrac{1}{p_\tau}w_\tau\snn{\tau} + \tfrac{\kconst_1}{p_\tau^2}(1-\csn{\tau})
    \\[4pt]
          w_\zr \simeq  w_\tau\csn{\tau} + p_\tau u_\tau \snn{\tau}  - \tfrac{\kconst_1}{p_\tau}\snn{\tau}
    \end{array}
    \end{array}
    &&,&&
    \begin{array}{lll}
        p_\tau = p_\zr  \\[3pt]
         \tup{p}_\tau \neq \tup{p}_\zr 
    \end{array}
    \end{flalign}
    \end{small}
    \item \textit{Back to ``standard'' ordering.} We have been using the modified coordinate ordering \eq{\tup{x} = (\tup{q},\tup{p},u,w)} which is best for illustrating the structure of the \textit{un}simplified dynamics. For the \textit{simplified} dynamics given above, one could just as easily return to the standard coordinate ordering, \eq{\tup{x} = (\tup{q},u,\tup{p},w)}, for which the above simplified matrix ODE takes the form
    \begin{small}
    \begin{flalign} \label{qusol_MAT_2bp_qupw}
    \qquad
    \tup{x}=(\tup{q},u,\tup{p},w)
    &&
    \begin{array}{rlll}
            &\pdt{\tup{x}} \,\simeq\, \Lambda_8 \cdot \tup{x} \,+\, \tup{k} 
    \\[5pt] 
       \fnsize{or,} & \rng{\tup{x}} \,\simeq\, \tfrac{1}{\slang}\Lambda_8 \cdot \tup{x} \,+\, \tfrac{1}{\slang}\tup{k} 
    \end{array}
        \qquad  \Rightarrow \qquad 
        \tup{x}_\tau \simeq E_\tau \cdot \tup{x}_\zr \,+\, \tup{\sigup}_\tau  
    &&
    \end{flalign} 
    \end{small}
    where everything is the ``same'' of what was given in Eq.\eqref{qusol_MAT_2bp}-Eq.\eqref{qusol_MAT_2bp_inv}, up to a re-ordering:
    \begin{small}
    \begin{align}
             \Lambda_8 = \fnsize{ Eq.\eqref{ThisMat_Ngen} }
             \in \spmat{8}
    &&,&& 
           E_\tau = \mrm{e}^{\frac{1}{\slang}\Lambda_8 \tau}
           = \fnsize{ Eq.\eqref{ThisMat_Ngen} }
                \in \Spmat{8}
    &&,&&
        \tup{\sigup}_\tau = \fnpmat{ \tup{0} \\
            \tfrac{\kconst_1}{\slang^2}(1-\csn{\tau}) \\
            \tup{0} \\
              \tfrac{\kconst_1}{\slang}\snn{\tau}  }
    &&,&&
    \tup{\varsigup}_\tau = \tup{\sigup}_{-\tau}
    \end{align}
    \end{small}
    where \eq{\Lambda_8} and \eq{\til{\Lambda}} both have eigenvalues \eq{\pm\mrm{i}\slang} (four of each sign). Note we have re-used some notation (\eq{\tup{x}}, \eq{E_\tau}, \eq{\tup{\sigup}_\tau}, and \eq{\tup{\varsigup}_\tau}) but now with a slightly different meaning (just re-ordering of indices). 
\end{itemize}
\end{small}
\vspace{1ex}

\noindent  Continuing with the modified ordering \eq{\tup{x}=(\tup{q},\tup{p},u,w)}, note that if we use the simplified Kepler flow, \eq{\theta_\tau}, from Eq.\eqref{qusol_MAT_2bp_inv} or Eq.\eqref{qusol_2bp_just_p} — treating \eq{\slang} as \eq{\slang\simeq p} —  to define an STM \eq{\Theta_\tau :=\pderiv{\theta_\tau}{\tup{x}}\big|_{\tup{x}_0} \equiv \pderiv{\tup{x}_\tau}{\tup{x}_0}}, then this would lead to:
\begin{small}
\begin{align}
  \Theta_\tau := \pderiv{\theta_\tau}{\tup{x}}\big|_{\tup{x}_0} \equiv \pderiv{\tup{x}_\tau}{\tup{x}_0}
   &\,\simeq\,
     \left.\fnpmat{ 
    \fnpmat{ 
    \csn{\tau} \, \imat_3
    &  \tfrac{1}{p}(\imat_3-\htup{p}\otms\htup{p})\snn{\tau} 
    \\[3pt]
   -p \snn{\tau} \,\imat_3
    & \csn{\tau} \, \imat_3 - \tup{q}\otms \htup{p} \snn{\tau}
    } 
    &
    \fnpmat{ \tup{0} & \tup{0} \\[3pt] \tup{0} & \tup{0} }
\\[14pt]
   \fnpmat{ 
         0 &
         -\tfrac{1}{p}\big( \tfrac{2\kconst_1}{p^2}(1- \csn{\tau}) + \tfrac{1}{p}w \snn{\tau} \big) \trn{\htup{p}}
         \\[3pt]
          0
         &  -(u + \tfrac{\kconst_1}{p^2})\snn{\tau} \,\trn{\htup{p}}
     }
    &
    \fnpmat{
    \csn{\tau} & \tfrac{1}{p}\snn{\tau}
    \\[3pt]
    -p \snn{\tau} & \csn{\tau}
    }  
} \right|_{\tup{x}=\tup{x}_0}
\\[8pt] \nonumber 
    \inv{\Theta_\tau} =  \pderiv{\tup{x}_0}{\tup{x}_\tau}
     &\,\simeq\, 
      \left.\fnpmat{ 
    \fnpmat{ 
    \csn{\tau} \, \imat_3
    &  -\tfrac{1}{p}(\imat_3-\htup{p}\otms\htup{p})\snn{\tau} 
    \\[3pt]
    p \snn{\tau} \,\imat_3
    & \csn{\tau} \, \imat_3 + \tup{q}\otms \htup{p} \snn{\tau}
    } 
    &
    \fnpmat{ \tup{0} & \tup{0} \\[3pt] \tup{0} & \tup{0} }
\\[14pt]
   \fnpmat{ 
         0 &
         -\tfrac{1}{p}\big( \tfrac{2\kconst_1}{p^2}(1- \csn{\tau}) - \tfrac{1}{p}w \snn{\tau} \big) \trn{\htup{p}}
         \\[3pt]
          0
         &  (u + \tfrac{\kconst_1}{p^2})\snn{\tau} \,\trn{\htup{p}}
     }
    &
    \fnpmat{
    \csn{\tau} & -\tfrac{1}{p}\snn{\tau}
    \\[3pt]
    p \snn{\tau} & \csn{\tau}
    }  
} \right|_{\tup{x}=\tup{x}_\tau}
\end{align}
\end{small}
The above \eq{\Theta_\tau} 
is not the same as the full STM \eq{\Phi_\tau}, from Eq.\eqref{STM_prj_nosimp} — which was obtained from the unsimplified Kepler flow — nor is it the same as the simplified form of \eq{\Phi_\tau} given in Eq.\eqref{STM_prj_somesimp}.

%% file: Mysecs_prj/new_5_prjJ2.tex
\subsection{Perturbed example:~\txi{J2} gravitational perturbation} \label{sec:j2}

\begin{notation}
    Here, the symbol \eq{J_2} is the ``\eq{J_2} constant" (second zonal harmonic coefficient) for the leading-order gravitational correction for oblate bodies — a numeric constant depending on shape and mass distribution of the body. It should not be confused with the canonical symplectic \eq{2\tms 2} matrix (also denoted \eq{J_2} but not used here). 
\end{notation}

\noindent
 As an example of the perturbed Kepler problem, consider the two-body problem for masses \eq{m_a} and \eq{m_b}, where \eq{m_b} is a uniform sphere but \eq{m_a} is some oblate body with equatorial radius \eq{R_a} whose \eq{J_2} term is a significant gravitational perturbation. The inertial cartesian coordinate Hamiltonian, with \eq{V^1} modeling the \eq{J_2} term, is given by
\begin{small}
\begin{flalign} \label{Vj2_r}
&&
     V^{1}(\tup{r})\,=\, \tfrac{1}{3}\tfrac{j_2}{r^3}(3 \hat{r}_3^2 - 1)
\qquad\Rightarrow\qquad
       \mscr{K} = \tfrac{1}{2} \v^2 - \tfrac{\kconst_1}{r} + \tfrac{1}{3}\tfrac{j_2}{r^3}(3 \hat{r}_3^2 - 1)
&&,&&
     j_2 := \tfrac{3}{2} J_2 \kconst_1 R_a^2
\quad 
\end{flalign}
\end{small}
where \eq{\hat{r}_i=r_i/r} and we have absorbed various numeric constants into \eq{j_2} defined above. In the following, we assume no other perturbations are present. That is, \eq{\tup{a}^\nc=0}  such that \eq{\tup{F}=-\pd_{\tup{r}}V^1} and, therefore, \eq{(\tup{\alpha},\alpha_\ss{u})=0} such that \eq{(\tup{f},f_\ss{u})=-(\pd_{\tup{q}} V^1 ,\pd_{u} V^1)}.

As discussed in section \ref{sec:gen_burdet_eom}, there are two approaches we may take when calculating the partials  \eq{\pderiv{V^{1}}{\tup{q}}} and \eq{\pderiv{V^{1}}{u}} appearing in the equations of motion for the projective coordinates. One method is to first obtain the cartesian components \eq{\pderiv{V^1}{\tup{r}}},
which is found from the above as follows:
\begin{small}
\begin{flalign} \label{j2_a}
&&
  \tup{F} \,=\,   -\pderiv{V^{1}}{\tup{r}} 
   \,=\,
    \tfrac{j_2}{r^4}
     \begin{pmatrix}
        (5\hat{r}_3^2 -1)\hat{r}_1 \\[3pt]
         (5\hat{r}_3^2-1)\hat{r}_2 \\[3pt]
          (5\hat{r}_3^2-3)\hat{r}_3 
    \end{pmatrix}
    \,=\,  
    \tfrac{j_2}{r^4} \big( (5 \hat{r}_3^2 -1)\htup{r} - 2 \hat{r}_3 \ibase_3 \big)
    \,=\, 
    j_2  u^4 \big( (5 \hat{q}_3^2 -1)\htup{q} - 2 \hat{q}_3 \ibase_3 \big)
    &&,&&  \hat{q}_i\simeq q_i
    \quad
\end{flalign}
\end{small}
where \eq{\ibase_3 = (0,0,1)} and
where the last expression has been written in terms of \eq{(\tup{q},u)} using  \eq{\hat{r}_i= \hat{q}_i} and \eq{r=1/u} (it could also be simplified using \eq{q=1}).
 We may then find \eq{(\tup{f},f_\ss{u})=-(\pd_{\tup{q}} V^1 ,\pd_{u} V^1) }  
from the above by using Eq.\eqref{sum_dV_2bp}:
\begin{small}
\begin{align}  \label{sum_dV_j2}
    \tup{f} \,=\,  -\pderiv{V^{1}}{\tup{q}} \,=\, -\trn{\pderiv{\tup{r}}{\tup{q}}}\cdot \pderiv{V^{1}}{\tup{r}}  
     \,=\,  \tfrac{1}{u q} (\imat_3 - \htup{q}\otms\htup{q}) \cdot \tup{F}
\quad, && 
   f_\ss{u} \,=\,  -\pderiv{V^{1}}{u} \,=\, -\pderiv{\tup{r}}{u}\cdot \pderiv{V^{1}}{\tup{r}} 
   \,=\, -\tfrac{1}{u^2} \htup{q} \cdot \tup{F}
\end{align}
\end{small}
Alternatively, we may simply use the point transformation, \eq{\tup{r} =\tfrac{1}{u}\htup{q}}, to express \eq{V^{1}} in terms of \eq{(\tup{q},u)} by using  \eq{\hat{r}_i= \hat{q}_i} and \eq{r=1/u}. Substituting this into \eq{V^{1}(\tup{r})} from Eq.\eqref{Vj2_r} leads to \eq{V^{1}(\bartup{q})} for the projective coordinate Hamiltonian:
\begin{small}
\begin{align} \label{dVj2}
V^1(\bartup{q}) = \tfrac{1}{3}j_2 u^3(3\hat{q}_3^2 - 1)
\qquad \Rightarrow \qquad
&\boxed{ \mscr{H} \,=\,
     \tfrac{1}{2}u^2\big( \slang^2  + u^2 p_u^2 \big)  - \kconst_1 u 
     +  \tfrac{1}{3}j_2 u^3(3\hat{q}_3^2 - 1) }
\end{align}
\end{small}
As per Remark \ref{rem:H_simp}, we do \textit{not}
 simplify the potential \eq{V^{1}} using \eq{q=1} as this would lead to an incorrect expression for \eq{\pd_{\tup{q}}{V^1}}. From the above, we find that the correct \eq{J_2} generalized forces are:
\begin{small}
\begin{align} \label{dvdq_j2}
   \tup{f} \,=\, -\pderiv{V^{1}}{\tup{q}} \,=\,
    \tfrac{2}{q}j_2 u^3 ( \hat{q}_3^2 \htup{q} - \hat{q}_3 \ibase_3 ) 
    \,\simeq\,  2j_2  u^3 ( q_3^2\tup{q} - q_3 \ibase_3 )
&&, &&
   f_\ss{u} \,=\,  -\pderiv{V^{1}}{u} \,=\, -j_2  u^2(3\hat{q}_3^2 - 1)
    \,\simeq\, -j_2  u^2(3q_3^2 - 1) 
\end{align}
\end{small}
where only after differentiating do we   simplify the above using \eq{q=1}. It can be verified that substituting Eq.\eqref{j2_a} into Eq.\eqref{sum_dV_j2} leads to the same result as given above.

\paragraph*{Equations of Motion Including $J_2$.}
Assuming only the above \eq{J_2} perturbation,  
the equations of motion for the \eq{J_2}-perturbed Kepler problem in projective coordinates are obtained simply be substituting Eq.\eqref{dvdq_j2} for \eq{(\tup{f},f_\ss{u})} — or Eq.\eqref{j2_a} for \eq{\tup{F}} — wherever these terms appear in the perturbed dynamics given previously in section \ref{sec:2BP}.

For example, 
the extended Hamiltonian and canonical equations of motion for the \eq{s}-parameterized dynamics, given by Eq.\eqref{qp_eom_s_sum}, now take the form:
\begin{small}
\begin{align}  \label{s_J2_2}
\begin{array}{ccccc}
    \wtscr{H} = \tfrac{1}{2}(\slang^2 + u^2 p_\ss{u}^2) - \tfrac{\kconst_1}{u} \,+\, \tfrac{1}{3}j_2 u(3\hat{q}_3^2 - 1) + \tfrac{1}{u^2}p_t
\\[6pt]
\begin{array}{llll}
      \pdt{q}_i  \,=\, -\slang_{ij} q_j
 \\[4pt]
     \pdt{u} \,=\, u^2 p_u 
\\[3pt]
    \vphantom{\tfrac{\kconst_1}{u^2}}
\end{array}
\qquad,\qquad
\begin{array}{llll}
      \pdt{p}_i   &=\, -\slang_{ij} p_j 
     + \tfrac{2}{q}j_2 u(\hat{q}_3^2 \hat{q}_i -  \kd_{i3} \hat{q}_3)
    \\[4pt]
       \pdt{p}_u  
     &=\, -up_u^2 - \tfrac{\kconst_1 }{u^2} - \tfrac{1}{3}j_2 (3 \hat{q}_3^2-1) + \tfrac{2}{u^3}p_t
 \\[3pt]
     &=\, -\tfrac{1}{u}(\slang^2 + 2 u^2 p_\ss{u}^2) + \tfrac{\kconst_1}{u^2} - j_2 (3\hat{q}_3^2 -1)
\end{array} 
\qquad,\qquad
\begin{array}{llll}
      \pdt{t} \,=\, 1/u^2
     \\[4pt]
      \pdt{p}_t \,=\, 0
  \\[3pt]
    \vphantom{\tfrac{\kconst_1}{u^2}}
\end{array} 
\end{array}
\end{align}
\end{small}
where the \eq{\pdt{p}_\ss{u}} equation has been re-written using \eq{p_t=-\mscr{H}}.\footnote{Differentiating \eq{\wtscr{H}} directly gives \eq{\pdt{p}_\ss{u} = -\pderiv{\wtscr{H}}{u} = -up_u^2 - \tfrac{\kconst_1 }{u^2} - \tfrac{1}{3}j_2 (3 \hat{q}_3^2-1) + \tfrac{2}{u^3}p_t }. Substitution of \eq{p_t=-\mscr{H}} leads to the expression for \eq{\pdt{p}_\ss{u}} seen in Eq.\eqref{s_J2_2}.  }
Simplifying the above dynamics using \eq{q=1} and \eq{\lambda=\htup{q}\cdot\tup{p}=0}, and rewriting the \eq{(u,p_\ss{u})} dynamics using \eq{w:=u^2 p_\ss{u}}, we have:
\begin{small}
\begin{flalign}
&&
\begin{array}{llll}
     \pdt{q}_i   \,\simeq\, p_i
     &\quad,\quad 
      \pdt{p}_i  
     \,\simeq\,
    -\slang^2 q_i \,+\, 2j_2 u( q_i q_3 -  \kd_{i3})q_3
 \\[4pt]
     \pdt{u} \,=\, w
      &\quad,\quad 
      \pdt{w} = -\slang^2 u + \kconst_1 \,-\, j_2 u^2(3 \hat{q}_3^2 -1)
\end{array}
&&
\slang^2 \simeq p^2 
\quad
\end{flalign}
\end{small}
The \eq{\tau}-parameterized ODEs are equivalent to using \eq{\rng{\square}=\tfrac{1}{\slang}\pdt{\square}}. 
Similarly, consider the second-order equations for \eq{\bartup{q}=(\tup{q},u)} given in Eq.\eqref{ddqu1}-Eq.\eqref{ddqu2_simp}.
After substitution of the \eq{J_2} generalized forces from Eq.\eqref{dvdq_j2} into said equations (and simplifying the right-hand-side using \eq{q=1} and \eq{\htup{q}\cdot\tup{p}=0}), we obtain:
\begin{small}
\begin{align} \label{ddq_j2}
\begin{array}{llllll}
     \pddt{q}_i + \slang^2q_i 
    \,\simeq\, 2j_2  u( q_i q_3-\kd_{i3})q_3
\\[5pt]
    \pddt{u} + \slang^2 u  - \kconst_1 
    \,\simeq\,  -j_2  u^2 (3q_3^2 -1)
\end{array}
&& \Bigg| &&
\begin{array}{lllll}
        \rrng{q}_{i} + q_i 
        \,\simeq\, 
        2j_2 \tfrac{u}{\slang^2}( q_i q_3  + \rng{q}_i \rng{q}_3 - \kd_{i3}) q_3
\\[5pt]
    \rrng{u} + u - \tfrac{\kconst_1 }{\slang^2} 
    \,\simeq\, -j_2   \tfrac{u}{\slang^2}\big( (3q_3^2-1)u \,-\, 2q_3\rng{q}_3\rng{u} \big)
\end{array}
\end{align}
\end{small}
where the specific angular momentum, \eq{\slang},  is not constant; it is governed by Eq.\eqref{ldot_cen_again}, now leading to
\begin{small}
\begin{align} 
    \dot{\slang}  \simeq \tfrac{1}{\slang}\tup{p}\cdot\tup{f} \simeq -\tfrac{1}{\slang}2j_2 u^3 p_3 q_3
    &&
    \pdt{\slang} = \dot{\slang}/u^2 \simeq -\tfrac{1}{\slang}2j_2 u \pdt{q}_3 q_3
    &&
    \rng{\slang} = \pdt{\slang}/\slang \simeq -\tfrac{1}{\slang}2j_2 u \rng{q}_3 q_3
\end{align}
\end{small}
Though, one does not need to include \eq{\slang} as an additional coordinate; it is given in terms of the projective coordinates by \eq{\slang^2=q^2 p^2-(\tup{q}\cdot\tup{p})^2 \simeq p^2 \simeq \mag{\pdt{\tup{q}}}^2}.

%% file: Mysecs_prj/99_conclusion.tex

\section{Conclusion} 

A method is developed for extending dimension-raising configuration coordinate transformations to canonical/symplectic coordinate transformations. 
Using this method, along with a an evolution parameter transformation, projective point transformations that linearize the Newtonian second-order ODEs for Kepler dynamics
are extended to canonical transformations that linearize the Hamiltonian first-order ODEs for both Kepler and Manev dynamics.
This provides an alternative to the KS and BF coordinates, offering a more intuitive formulation applicable to a broader class of problems.
In addition, the handling of coordinate redundancy is made transparent starting with Hamilton's principle, and the problem of when/how to use the constraints and conserved quantities is largely circumvented — full linearity of Kepler and Manev dynamics is achieved \textit{with or without} observing the extra integrals of motion.   
Closed-form solutions are presented and used to obtain Keplerian state transition matrices.


Within the classical Hamiltonian framework presented here, a typical next step on this path would be the application of Hamilton-Jacobi theory to obtain from the projective coordinates their corresponding action-angle coordinates and investigate their subsequent use in semi-analytic Hamiltonian/symplectic perturbation methods (e.g., Deprit-Hori and Lie series methods). 
The fact that the transformed system is both linear and Hamiltonian makes this formulation an attractive choice for such perturbation methods, as well as for symplectic integration algorithms. 
In a different direction, given the parallels between our projective coordinates and the KS coordinates, quantitative comparisons — e.g., computational efficiency and long-term propagation stability/accuracy — would be valuable.

%% file: Mysecs_prj/new_1apx_prjCoord.tex
\section{Families of BF-type projective transformations} \label{sec:2bp_gen}

\noindent We apply the methods of section \ref{sec:Pxform} to develop a family of canonical transformations for a generalized version of the projective point transformation used by Burdet, Sperling, Vitins, and others. We will do so within the context of arbitrary particle dynamics in 
Euclidean 3-space.\footnote{Though the majority of the developments easily generalize to any finite-dimensional real inner product space. See \cite{peterson2025prjGeomech}.}
The use of such transformations for regularizing and linearizing central-force dynamics is considered in sections \ref{sec:central force} and \ref{sec:2BP}.

    Readers may wish to refer to the present section only for desired derivations and details. 
    In short, this section considers a canonical extension of a family of projective point transformations given by \eq{\tup{r}=u^n q^m \tup{q}} for arbitrary scalar constants \eq{n\neq 0,m\in\mbb{R}}. 
    The results are important for the main developments in sections \ref{sec:central force} and \ref{sec:2BP}. However, the developments of the present section are somewhat cumbersome and not particularly interesting in their own right until we choose desired values of \eq{n} and \eq{m} in sections \ref{sec:central force} and \ref{sec:2BP}. 


\paragraph*{Original Cartesian Coordinate Formulation.} 
We consider an unconstrained particle moving in Euclidean 3-space, subject to arbitrary conservative and nonconservative forces.
Let \eq{\tup{r}\in\mbb{R}^{3}} and \eq{\tup{v}\in\mbb{R}^{3}} denote the cartesian position and velocity coordinates associated to an orthonormal inertial frame.
The set \eq{(\tup{r},\tup{v})\in\mbb{R}^6} will be our starting, minimal, canonical coordinates (what we denoted \eq{(\tup{x},\tup{\piup})} in section \ref{sec:Pxform}). 
The Hamiltonian (per unit mass) and canonical equations of motion — including arbitrary nonconservative forces — for the particle are given in cartesian coordinates simply by:
\begin{small}
\begin{align} \label{Kjoechi}
    &\mscr{K} \,=\, \tfrac{1}{2}\v^2 + V(\tup{r},t)
&&
\begin{array}{ll}
      \dot{\tup{r}} \,=\, \pderiv{\mscr{K}}{\tup{v}} \,=\, \tup{v} 
\\[5pt] 
    \dot{\tup{v}} \,=\, -\pderiv{\mscr{K}}{\tup{r}} + \tup{a}^\nc \,=\,
    -\pderiv{V}{\tup{r}} +  \tup{a}^\nc
\end{array}
\end{align}
\end{small}
with \eq{\v=\mag{\tup{v}}} and where \eq{\tup{a}^\nc\in\mbb{R}^3} are just the inertial cartesian components of the nonconservative forces. 


\subsection{Canonical extension of a family of projective point transformations} \label{sec:gen_burdet}


Starting with the cartesian coordinate description given in  Eq.\eqref{Kjoechi}, we then specify a point transformation, \eq{\mbb{R}^4\ni \bartup{q} \mapsto \tup{r}\in\mbb{R}^3}, 
where the new coordinates, \eq{\bartup{q}\in\mbb{R}^4}, are redundant coordinates 
consisting of three coordinates related to the radial unit vector, \eq{\htup{r}}, and one coordinates related to the radial distance, \eq{r:=\mag{\tup{r}}}. We then use the methods of section \ref{sec:Pxform} to extend this point transformation to the momentum level and obtain a canonical transformation between the minimal phase space with cartesian coordinates \eq{(\tup{r},\tup{v})\in\mbb{R}^6} and Hamiltonian \eq{\mscr{K}(\tup{r},\tup{v},t)}, and a set  of redundant phase space coordinates   \eq{(\bartup{q},\bartup{p})\in\mbb{R}^8} and Hamiltonian, \eq{\mscr{H}(\bartup{q},\bartup{p},t)}.
It will be helpful to partition our redundant coordinates as: 
\begin{small}
\begin{align}
          \bartup{q}\,=\, (\tup{q} , u)   \in\mbb{R}^4  
 \qquad,\qquad
          \bartup{p}\,=\, ( \tup{p} , p_\ss{u}) \in\mbb{R}^4
 \end{align}
 \end{small}

\paragraph*{The Point Transformation.}
The usual generalized projective point transformation would be given as in Eq.\eqref{proj_pt}  by \eq{\tup{r}=u^n\tup{q}} for some \eq{n\neq0}. We instead start with the following further generalized family of projective point transformations: 
\begin{small}
\begin{flalign} \label{PT_Qq}
&&
\begin{array}{lllll}
    \boxed{ \tup{r} \,=\, u^n q^m\tup{q} }
    &,\qquad 
     \fnsize{subject to:} \;\; \varphi(\tup{q}) \,=\, q - \qsc \,=\, 0 
\end{array}
&&
\rnk \pderiv{\tup{r}}{\bartup{q}} = 3
\end{flalign}
\end{small}
where \eq{q=\mag{\tup{q}}} (not \eq{\mag{\bartup{q}}}), where we are free to choose any numbers \eq{n\neq0,m\in\mbb{R}} to construct different transformations, and where \eq{\varphi} is a constraint function decreeing that \eq{\tup{q}} have constant magnitude, \eq{\qsc\in\mbb{R}_\ss{+}}. Though one could consider any number \eq{\qsc >0}, we are only interested in the case \eq{\qsc=1}. Yet, as the number \eq{1} is difficult to track though equations, we will occasionally write \eq{\qsc} to make steps easier to follow. 


 The parameters \eq{n,m\in\mbb{R}} appear throughout the remainder of Appendix \ref{sec:2bp_gen} and will make some equations temporarily cluttered. In sections \ref{sec:central force} and \ref{sec:2BP}, we choose desirable values for these parameters which simplify things considerably.


\begin{notesq}
\textit{The value of $m$.}
The presence of \eq{q^m} in the point transformation \eq{\tup{r} = u^n q^m\tup{q}} may seem strange given that the constraint \eq{\varphi} is specified such that \eq{q=\qsc} is constant. With \eq{\qsc=1} (the case of primary interest), the factor of \eq{q^m} in the transformation might seem altogether pointless.\footnote{Even with \eq{\qsc\neq1}, the term \eq{q^m} would still seem only to contribute a scaling by an arbitrary constant.} 
However, the factor of \eq{q^m} affects the Jacobian matrix \eq{\pderiv{\tup{r}}{\bartup{q}}} which plays a central role in transforming the momenta, Hamiltonian, forces, etc. As such, different values of \eq{m} yield surprisingly differences at the momentum level and in the resulting projective coordinate Hamiltonian system. 
\end{notesq}


\noindent With \eq{\mscr{K}(\tup{r},\tup{v},t)} specified in Eq.\eqref{Kjoechi} and \eq{\tup{r}(\bartup{q})} specified in Eq.\eqref{PT_Qq}, we now wish to find the momenta transformation, \eq{\tup{v}(\bartup{q},\bartup{p})}, and the new Hamiltonian, \eq{\mscr{H}(\bartup{q},\bartup{p},t)}, such that the dynamics in the \eq{(\bartup{q},\bartup{p})} representation satisfy Hamilton's canonical equations, and can be mapped back to the original cartesian coordinate Hamiltonian dynamics for \eq{\mscr{K}(\tup{r},\tup{v},t)}. To accomplish this, we simply apply Eq.\eqref{qpH_new} through Eq.\eqref{Ham_xform} of section \ref{sec:Pxform}. 

\paragraph*{The Momenta Transformation.}
To find the momenta transformation we need the partials of \eq{\tup{r}(\bartup{q})} and \eq{\varphi(\tup{q})} from Eq.\eqref{PT_Qq}. These are collected into a matrix \eq{B\in\Glmat{4}}:
\begin{small}
\begin{align} \label{B_nm}
 B(\bartup{q}) :=\, 
    \begin{pmatrix}
         \pderiv{\tup{r}}{\tup{q}} 
         & \pderiv{\tup{r}}{u} 
         \\[5pt] 
         \pderiv{\varphi}{\tup{q}}
         & \pderiv{\varphi}{u}
    \end{pmatrix} 
  \,=
    \smpmat{
            u^n q^{m}(\imat_3 + m\htup{q}\otms \htup{q}) 
            & \;\;  nu^{n-1}q^{m+1}\htup{q}  
             \\[5pt]
             \trn{\htup{q}} & 0 
            }
        \qquad,\qquad     
        \inv{B}(\bartup{q}) \,=
        \smpmat{
             \tfrac{1}{u^n q^m}(\imat_3 - \htup{q}\otms \htup{q})
             &  \htup{q}
             \\[5pt]
             \tfrac{u^{1-n}}{nq^{m+1}} \trn{\htup{q}}
             & -\tfrac{m+1}{n}\tfrac{u}{q}
          }
\end{align}
\end{small}
The momenta transformation for \eq{\tup{v}(\bartup{q},\bartup{p})} as well as \eq{\lambda(\bartup{q},\bartup{p})} is then obtained from
Eq.\eqref{Pxform_finaly}, or Eq.\eqref{Pxform_div}, leading to:
\begin{small}
\begin{align} \label{P_pq}
    \begin{pmatrix}
     \tup{v} \\ \lambda 
\end{pmatrix} =
     \invtrn{B}\begin{pmatrix}
     \tup{p} \\ p_\ss{u}
\end{pmatrix}
\qquad,\qquad 
    \begin{array}{ll}
         \boxed{ \tup{v}
         \,=\,   \tfrac{1}{u^n q^m}\big[ (\imat_3 - \htup{q}\otms\htup{q})\cdot\tup{p}   \,+\, \tfrac{u}{nq}p_\ss{u}\htup{q}\big] }
    \\[8pt]
     \lambda  \,=\, \tfrac{1}{q}(\tup{q}\cdot\tup{p} - \tfrac{m+1}{n}u p_\ss{u}) 
\end{array}
\end{align}
\end{small}

\begin{notesq}
    \textit{Angular momentum.}
    The inertial cartesian components of the angular momentum (pseudo)vector, \eq{\tup{\slangup}}, the  antisymmetric matrix, \eq{\hdge{\tup{\slangup}}},  and  the magnitude, \eq{\slang}, are given in terms of the cartesian coordinates by the usual relations \eq{\tup{\slangup} =\tup{r}\tms \tup{v} }, \eq{\hdge{\tup{\slangup}} = \tup{r} \wdg \tup{v} }, and \eq{ \slang^2 = r^2 \v^2 - (\tup{r}\cdot\tup{v})^2}.
    Conveniently, substitution of \eq{\tup{r}(\bartup{q})} and \eq{\tup{v}(\bartup{q},\bartup{p})} from Eq.\eqref{PT_Qq} and Eq.\eqref{P_pq} leads to expressions of the precisely the same form:
    \begin{small}
    \begin{align} \label{l_qp_general_again}
    \begin{array}{rllllll}
          \tup{\slangup} &=\, \tup{r}\tms \tup{v} 
         \,=\, \tup{q}\tms \tup{p} 
    \end{array}
     &&,&&
    \begin{array}{rllllll}
           \hdge{\tup{\slangup}} &=\,  \tup{r} \wdg \tup{v}  
        \,=\, \tup{q} \wdg \tup{p} 
    \end{array}
     &&,&&
    \begin{array}{rllllll}
            \slang^2  = \mag{\tup{\slangup}}^2
          &=\, r^2 \v^2  - (\tup{r}\cdot\tup{v})^2
         \,=\, q^2 p^2  - (\tup{q}\cdot\tup{p})^2
    \end{array}
    \end{align}
    \end{small}
    The above does not depend on \eq{u} or \eq{p_\ss{u}}, and holds for any \eq{n\neq 0} and \eq{m} in Eq.\eqref{PT_Qq}. A collection of useful angular momentum relations is given in Appendix \ref{sec:ang_momentum}. 
\end{notesq}

\begin{notesq}
    \textit{The form of \eq{\varphi}.}
    In earlier conference and preprint papers, the authors instead posed the constraint as \eq{\til{\varphi}(\tup{q}) = \tfrac{1}{2}(q^2-\qsc^2)=0} such that \eq{\pderiv{\til{\varphi}}{\tup{q}}=\tup{q}} is globally defined. 
    The matrices \eq{B} and \eq{\inv{B}} of partial derivatives then differ slightly from those given in Eq.\eqref{B_nm}.\footnote{With \eq{\til{\varphi}(\tup{q}) = \tfrac{1}{2}(q^2-\qsc^2)=0}, the matrices \eq{B} and \eq{\inv{B}} in Eq.\eqref{B_nm} change slightly to the following:
    \begin{align} \nonumber 
        \til{\varphi} = \tfrac{1}{2}(q^2-\qsc^2)=0
        && \Rightarrow &&
        \begin{array}{lll}
     B := 
        \begin{pmatrix}
             \pderiv{\tup{r}}{\tup{q}} 
             & \pderiv{\tup{r}}{u} 
             \\[5pt] 
             \pderiv{\til{\varphi}}{\tup{q}}
             & \pderiv{\til{\varphi}}{u}
        \end{pmatrix} 
      \,=
        \begin{pmatrix}
            u^n q^{m}(\imat_3 + m\htup{q}\otms \htup{q}) 
             & \;\; nu^{n-1}q^{m+1}\htup{q} &
         \\[5pt]
             \trn{\tup{q}} & 0 
    \end{pmatrix}
    \qquad,\qquad     
    \inv{B}\,=
        \begin{pmatrix}
             \tfrac{1}{u^n q^m}(\imat_3 - \htup{q}\otms \htup{q})
             &  \tfrac{1}{q^2}\tup{q}
         \\[5pt]
              \tfrac{u^{1-n}}{nq^{m+1}}\trn{\htup{q}}
                  & -\tfrac{m+1}{n}\tfrac{u}{q^2}
        \end{pmatrix}
    \end{array}
    \end{align} }
    This changes little other than the fact that the associated Lagrange multiplier,  \eq{\til{\lambda}}, would then instead be given by \eq{\til{\lambda} = \tfrac{1}{q^2}(\tup{q}\cdot\tup{p} - \tfrac{m+1}{n}u p_\ss{u})}. 
\end{notesq}


\paragraph*{The Hamiltonian.}
 The point transformation and constraint given by Eq.\eqref{PT_Qq}  are not explicit functions of time. Therefore, it is seen from Eq.\eqref{Ham_xform} that the new Hamiltonian is found from direct substitution of \eq{\tup{r}(\bartup{q})} and \eq{\tup{v}(\bartup{q},\bartup{p})} into the original Hamiltonian for the cartesian coordinates, \eq{\mscr{K}}. That is,
\begin{small}
\begin{align} \label{Hsubs}
    \mscr{H}(\bartup{q},\bartup{p},t) \,=\, \mscr{K}\big(\tup{r}(\bartup{q}),\tup{v}(\bartup{q},\bartup{p}),t \big)
\end{align}
\end{small}
where \eq{\mscr{K}} is given by Eq.\eqref{Kjoechi}, \eq{\tup{r}(\bartup{q})}  by Eq.\eqref{PT_Qq}, and  \eq{\tup{v}(\bartup{q},\bartup{p})} by Eq.\eqref{P_pq}. From these last two relations, we see that
\begin{small}
\begin{align} \label{QQPP}
     r^2 = \tup{r}\cdot\tup{r} 
     \,=\, u^{2n}q^{2m+2}
&&,&& 
     \v^2 = \tup{v}\cdot\tup{v} \,=\, \tfrac{1}{u^{2n}q^{2m+2}}\Big( q^2 p^2 -(\tup{q}\cdot\tup{p})^2 \,+\, \tfrac{1}{n^2}u^2 p_\ss{u}^2 \Big)
     \;=\, \tfrac{1}{u^{2n}q^{2m+2}} ( \slang^2 + \tfrac{1}{n^2}u^2 p_\ss{u}^2 )
\end{align}
\end{small}
Substitution of the above \eq{\v^2} into \eq{\mscr{K}} leads to the new Hamiltonian \eq{\mscr{H}}:
\begin{small}
\begin{align} \label{Hqp_general}
\mscr{K} \,=\, \tfrac{1}{2}\v^2 \,+\, V(\tup{r},t)
\qquad \Rightarrow \qquad 
    \boxed{ \mscr{H} \,=\, 
    \tfrac{1}{u^{2n}q^{2m+2}} \tfrac{1}{2}\big( \slang^2 + \tfrac{1}{n^2}u^2 p_\ss{u}^2 \big)  \,+\, V(\bartup{q},t) }
\end{align}
\end{small}
where the angular momentum, \eq{\slang^2 = q^2 p^2 -(\tup{q}\cdot\tup{p})^2}, satisfies the relations in Appx.~\ref{sec:ang_momentum}, and where 
 the abuse of notation \eq{V(\bartup{q},t)} simply indicates the original \eq{V(\tup{r},t)} expressed in terms of \eq{\bartup{q}} by substitution of  \eq{\tup{r} = u^n q^m \tup{q}}.

\begin{remrm} \label{rem:H_simp}
 The Hamiltonian given by Eq.\eqref{Hqp_general} should be constructed and differentiated treating all \eq{(\bartup{q},\bartup{p})} as eight independent coordinates.  That is, the Hamiltonian should \textit{not} be simplified using the constraint \eq{q=\qsc (=1)}. We must treat \eq{q=\mag{\tup{q}}} as a function such that \eq{\pderiv{q}{\tup{q}}=\htup{q}}.
Failing to do this will lead to incorrect partial derivatives of \eq{\mscr{H}} and incorrect equations of motion. However, \textit{after} differentiating \eq{\mscr{H}}, the resulting ODEs themselves can indeed be simplified using \eq{q=1} (as well as \eq{\lambda=0}) for the reasons given later in remark \ref{rem:q=k 1} and \ref{rem:qlam_qp}. 
\end{remrm}

\paragraph*{An Inverse Transformation.}
The transformation \eq{\mbb{R}^8\ni(\bartup{q},\bartup{p})\mapsto(\tup{r},\tup{v})\in\mbb{R}^6}, as given by Eq.\eqref{PT_Qq} and Eq.\eqref{P_pq}, does not have a unique inverse mapping \eq{(\tup{r},\tup{v})\mapsto (\bartup{q},\bartup{p})}. However, 
note the constraint \eq{\varphi=0} in Eq.\eqref{PT_Qq} is equivalent to \eq{q=\qsc}, and that in Eq.\eqref{P_pq} we have a relation for \eq{\lambda(\bartup{q},\bartup{p})}. All together, we then have a map \eq{\mbb{R}^8\ni(\bartup{q},\bartup{p})\mapsto(\tup{r},\qsc,\tup{v},\lambda)\in\mbb{R}^8}:
\begin{small}
\begin{align} \label{prj_nmGen}
\begin{array}{llllll}
     \tup{r} \,=\, u^n q^m \tup{q}
     &,\qquad    \tup{v}
         \,=\,   \tfrac{1}{u^n q^m}\big[ (\imat_3 - \htup{q}\otms\htup{q})\cdot\tup{p}   \,+\, \tfrac{u}{nq}p_\ss{u}\htup{q}\big]
\\[4pt]
     \varphi \,=\, 0 \;\;\; (\fnsize{i.e.,} \; \qsc = q)
     &,\qquad  \lambda  \,=\, \tfrac{1}{q}(\tup{q}\cdot\tup{p} - \tfrac{m+1}{n}u p_\ss{u}) 
\end{array}
\end{align}
\end{small}
which can indeed be inverted to obtain \eq{(\bartup{q},\bartup{p})} in terms of \eq{(\tup{r},\tup{v})}, along with \eq{\qsc} and \eq{\lambda}: 
\begin{small}
\begin{align} \label{prj_inv_nmGen}
\begin{array}{llllll}
     \tup{q} \,=\, \qsc\htup{r}
     &,\qquad  
     \tup{p} \,=\, \tfrac{r}{\qsc}( \imat_3 + m\htup{r}\otms\htup{r} ) \cdot \tup{v} \,+\,  \lambda \htup{r} 
\\[4pt]
    u \,=\, \big(\frac{r}{ \qsc^{m+1} }\big)^{1/n}
     &,\qquad
     p_\ss{u}  \,=\, n \big(\frac{r}{ \qsc^{m+1} }\big)^{-1/n} \tup{r}\cdot\tup{v} 
\end{array}
\end{align}
\end{small}
The above, which is not yet complete, is derived as follows:
\begin{footnotesize}
\begin{itemize}
    \item[] \textit{Derivation of Eq.\eqref{prj_inv_nmGen}.} The point transformation  given by Eq.\eqref{PT_Qq} — which is independent of the momenta and \eq{\lambda} — is easily inverted  for \eq{\tup{q}(\tup{r})} and \eq{u(\tup{r})} 
    by first finding \eq{r(\bartup{q})} and \eq{\htup{r}(\bartup{q})} and then inverting these expressions, using the constraint relation \eq{\qsc=q}, to express \eq{\tup{q}} and \eq{u} in terms of \eq{\tup{r}} and \eq{\qsc}:
    \begin{align} \label{ears}
    \tup{r} = u^n q^m \tup{q} 
    \qquad\Rightarrow\qquad 
    \begin{array}{llll}
        \htup{r} \,=\, \tfrac{1}{r}\tup{r} \,=\, \tfrac{1}{q} \tup{q} = \htup{q}  
     \\[5pt]
          r \,=\, \mag{\tup{r}} \,=\, u^n q^{m+1}   
    \end{array}
    &&\xRightarrow[]{\;\;q\,=\,\qsc\;\;}&&
     \begin{array}{llll}
          \tup{q} \,=\, \tfrac{q}{r} \tup{r}  \,=\, \qsc\htup{r}
      \\[5pt]
        u \,=\, \big(\frac{r}{ q^{m+1} }\big)^{1/n}
         \,=\, \big(\frac{r}{ \qsc^{m+1} }\big)^{1/n}
     \end{array}
    \end{align}
    where \eq{\qsc=1} leads to the the desired \eq{\tup{q}=\htup{r}} and \eq{u=r^{\sfrac{1}{n}}}.
    Now, the relation for the momenta, \eq{\bartup{p}=(\tup{p},p_\ss{u})}, is found as in Eq.\eqref{Pxform_finaly} or Eq.\eqref{Pxform_div} from section \ref{sec:Pxform}. This is simply the inverse of Eq.\eqref{P_pq}:
   \begin{align} \label{zeeky} 
    &\begin{pmatrix}
         \tup{p} \\ p_\ss{u}
    \end{pmatrix} 
    \,=\,\trn{B}\begin{pmatrix}
         \tup{v} \\ \lambda 
    \end{pmatrix} 
    \qquad,\qquad 
    \begin{array}{llll}
        \tup{p} \,=\,    u^n q^{m} ( \imat_3 + m\htup{q}\otms\htup{q} ) \cdot\tup{v} \,+\, \lambda\htup{q}
        & =\,
          \tfrac{r}{\qsc}( \imat_3 + m\htup{r}\otms\htup{r} ) \cdot \tup{v} \,+\, \lambda \htup{r}
    \\[5pt]
        p_\ss{u}  \,=\, nu^{n-1}q^m\tup{q}\cdot\tup{v} 
        & =\,  n\big(\tfrac{\qsc^{m+1}}{r}\big)^\ss{1/n}\tup{r}\cdot\tup{v} 
    \end{array}  
    \end{align}
    where the last equalities follow from substitution of \eq{\bartup{q}(\tup{r})} from Eq.\eqref{ears}.
     The above are precisely what was claimed in Eq.\eqref{prj_inv_nmGen}. 
\end{itemize}
\end{footnotesize}

\noindent 
Now, the ``inverse'' transformation in Eq.\eqref{prj_inv_nmGen} is only unique up to some chosen value \eq{0<\qsc\in\mbb{R}} (we only care about \eq{\qsc=1}). A larger issue is that \eq{\lambda} has no expression in terms of \eq{(\tup{r},\tup{v})} such that Eq.\eqref{prj_inv_nmGen} is still incomplete. 
Yet, this is a nonissue as it turns out that both \eq{q=\mag{\tup{q}}} and \eq{\lambda} are integrals of motion of the new Hamiltonian system whose values we may choose freely. This is clarified in Remarks \ref{rem:q=k 1} and \ref{rem:qlam_qp} below with a chosen, complete, inverse transformation given in Eq.\eqref{qu_k1}.

\begin{remrm}[\textit{Two ``extra'' integrals of motion}] \label{rem:q=k 1}
    \sloppy  It it proven in section \ref{sec:lambda_const} that \eq{q=\mag{\tup{q}}} and \eq{\lambda(\bartup{q},\bartup{p})} are integrals of motion 
    — an abuse of terminology, perhaps\footnote{``Integral of motion'' is perhaps an abuse of terminology in this context; the functions in Eq.\eqref{qlam=const} do not represent a conserved physical quantity for a particle moving in Euclidean 3-space (the actual system we are modeling). They are, rather,  ``kinematic constants'' which arise from artificially increasing the dimension of phase space. However, in a metathetical sense, they are indeed true integrals of motion for the new Hamiltonian system described by \eq{\mscr{H}}.} — 
    of the  projective coordinate Hamiltonian system for \eq{\mscr{H}}.
    It is shown that this holds in the presence of arbitrary conservative and nonconservative forces.
    That is, letting \eq{\bartup{\alphaup}\in\mbb{R}^4} denote generalized nonconservative forces for the projective coordinates (this is detailed later), then: 
    \begin{small}
    \begin{flalign} \label{qlam=const}
         \begin{array}{cc}
         \fnsize{integrals}  \\
         \fnsize{of motion} 
    \end{array} \!\!:
    \left. \;\;\;
    \begin{array}{ll}
        q = \mag{\tup{q}} 
         \\[5pt]
        \lambda =  \tfrac{1}{q}(\tup{q}\cdot\tup{p} - \tfrac{m+1}{n} u p_\ss{u})   
    \end{array} \right.
    ,&&
    \begin{array}{ll}
        \dot{q} = \pbrak{q}{\mscr{H}} + \pderiv{q}{\bartup{p}} \cdot \bartup{\alphaup}  = 0 
       \\[5pt]
        \dot{\lambda}  = \pbrak{\lambda}{\mscr{H}} + \pderiv{\lambda}{\bartup{p}} \cdot \bartup{\alphaup}  = 0
    \end{array} 
    \quad \Rightarrow \quad
    \begin{array}{ll}
         q = q_\zr  
         \\[6pt]
           \lambda = \lambda_\zr  
    \end{array}
    &&
    \fnsize{choose:} \;\;
    \begin{array}{ll}
         q_\zr = 1
         \\[6pt]
           \lambda_\zr  = 0
    \end{array}
    \quad
    \end{flalign}
    \end{small}
    That is, our coordinates \eq{(\bartup{q},\bartup{p})\in\mbb{R}^8}, which are over-parameterized by two, satisfy the two ``extra'' conservation laws, \eq{q=q_\zr} ad \eq{\lambda=\lambda_\zr}, for some values determined by initial conditions. 
    We may choose any values \eq{0<q_\zr} and \eq{\lambda_\zr} that we please since they place no restrictions on the cartesian coordinates.
     We are only interested in the  values \eq{q=q_\zr=1} and \eq{\lambda=\lambda_\zr=0}. 
    Although not necessary, we  may use these relations to simplify the equations of motion \textit{after} they have been found by differentiating the Hamiltonian. 
\end{remrm}

\begin{remrm}[\textit{Inverse transformation}]\label{rem:qlam_qp}
 We do not need to do anything to enforce the above relations that \eq{q=1} and \eq{\lambda=0} other than specify our initial conditions such that these relations hold at the initial time. That is, for some known initial conditions \eq{(\tup{r}_{\zr},\tup{v}_{\zr})}, we specify \eq{(\bartup{q}_{\zr},\bartup{p}_{\zr})} using Eq.\eqref{prj_inv_nmGen} with \eq{q=\qsc=1} and \eq{\lambda=0} (this places no restrictions on the cartesian coordinates):
\begin{small}
\begin{align} \label{qu_k1}
\begin{array}{cc}
\fnsize{choosing} \\
     q=1  \\
    \lambda=0 
\end{array} \;\; \Rightarrow \qquad 
\left\{\quad
\boxed{\begin{array}{llll}
    \tup{q} \,=\, \htup{r} &,
\\[4pt] 
    u \,=\, r^{\sfrac{1}{n}}  &,
\end{array} 
\quad 
\begin{array}{llll}
    \tup{p} \,=\, r ( \imat_3 + m\htup{r}\otms\htup{r} ) \cdot \tup{v} 
\\[4pt]
  p_\ss{u} \,=\,  n r^\ss{-1/n} \tup{r}\cdot\tup{v}
\end{array} } \right.
\begin{array}{lll}
     =\,  r\dot{\tup{r}} + m\dot{r}\tup{r}
 \\[4pt]
     =\, 
nr^\ss{(n-1)/n}\,\dot{r}
\end{array}
\end{align}
\end{small}
The above satisfies \eq{q=1} and \eq{\lambda=0} and is an inverse of Eq.\eqref{prj_nmGen} (for these chosen values). If the above is used for \eq{(\tup{r}_\zr,\tup{v}_\zr)\mapsto(\bartup{q}_\zr,\bartup{p}_\zr)} then \eq{q_\zr=1} and \eq{\lambda_\zr=0} hold automatically and, since \eq{q} and \eq{\lambda} are integrals of motion, it is guaranteed that \eq{q=1} and \eq{\lambda=0} hold for all time along any solution curve for the new Hamiltonian \eq{\mscr{H}} and, therefore, so does the above transformation.  
\end{remrm}

\subsection{Transformed Hamiltonian dynamics and generalized forces}  \label{sec:gen_burdet_eom}

\paragraph*{Canonical EOMs in Projective Coordinates.}
Suppose that, in addition to arbitrary conservative forces modeled by the potential \eq{V},  arbitrary nonconservative forces such as thrust or drag also act on the particle (see Appx.~\ref{app:non_conserv} for an overview of general Hamiltonian analytical dynamics with nonconservative forces). 
Then, as per Appx.~\ref{app:non_conserv}, an additional term is included in the momenta equation such that the canonical equations of motion are given by
\eq{\dot{\bartup{q}} = \pderiv{\mscr{H}}{\bartup{p}} } and
\eq{\dot{\bartup{p}} = -\pderiv{\mscr{H}}{\bartup{q}}+ \bartup{\alphaup}},
where \eq{\bartup{\alphaup}=(\tup{\alphaup},\alpha_\ss{u})\in\mbb{R}^4} are the generalized nonconservative forces for the projective coordinates (discussed soon). 
 With the projective coordinate Hamiltonian given by Eq.\eqref{Hqp_general}, we obtain the following equations of motion:
\begin{small}
\begin{gather}  \nonumber 
  \mscr{H} \,=\, 
    \tfrac{1}{u^{2n}q^{2m+2}} \tfrac{1}{2}\big( \slang^2 + \tfrac{1}{n^2}u^2 p_\ss{u}^2 \big)  \,+\, V(\bartup{q},t)
\\[5pt] \label{sum_qpdot}
\begin{array}{llll}
     \dot{\tup{q}} \,=\, \pderiv{\mscr{H}}{\tup{p}} 
     \,=\,   -\tfrac{1}{u^{2n}q^{2m+2}} \hdge{\tup{\slangup}}\cdot\tup{q} &,
\\[6pt]
    \dot{u} \,=\, \pderiv{\mscr{H}}{p_u} \,=\, \tfrac{1}{u^{2n}q^{2m+2}}\tfrac{1}{n^2}u^2 p_\ss{u} &,
 \end{array}
\qquad 
 \begin{array}{llll}
    \dot{\tup{p}} \,=\,  -\pderiv{\mscr{H}}{\tup{q}} + \tup{\alphaup} 
    \,=\,  - \tfrac{1}{u^{2n}q^{2m+2}} \hdge{\tup{\slangup}}\cdot\tup{p} \,+\,  \tfrac{m+1}{u^{2n}q^{2m+4}}(\slang^2 + \tfrac{1}{n^2}u^2 p_\ss{u}^2)\tup{q} \,-\,   \pderiv{V}{\tup{q}} + \tup{\alphaup} 
\\[6pt]
    \dot{p}_u \,=\, -\pderiv{\mscr{H}}{u} + \alpha_\ss{u} 
    \,=\, \tfrac{n}{u^{2n+1}q^{2m+2}} \big( \slang^2 \,+\, \tfrac{n-1}{n^3}u^2 p_\ss{u}^2 \big)   \,-\, \pderiv{V}{u} + \alpha_\ss{u}
\end{array}
\end{gather}
\end{small}
where the angular momentum, \eq{\slang^2 = q^2 p^2 -(\tup{q}\cdot\tup{p})^2}, 
satisfies Eq.\eqref{angmoment_rels_crd}-Eq.\eqref{angmoment_rels_crd2}.\footnote{In particular, \eq{\pd_{\tup{p}}\tfrac{1}{2}\slang^2 = -\hdge{\tup{\slangup}}\cdot\tup{q}  = q^2\tup{p}-(\tup{q}\cdot\tup{p})\tup{q}= \tup{\slangup}\tms\tup{q}} and  \eq{\pd_{\tup{q}}\tfrac{1}{2}\slang^2 =\hdge{\tup{\slangup}}\cdot\tup{p}  = p^2\tup{q} -(\tup{q}\cdot\tup{p})\tup{p} = -\tup{\slangup}\tms\tup{p}}. }
As per Remark \ref{rem:q=k 1}, the above may be simplified using \eq{q=q_\zr} for any constant \eq{q_\zr>0}. Specifically, for \eq{q=1} this leads to 
\begin{small}
\begin{align} \label{sum_qpdot_k1}
q=1 \;\; \Rightarrow \;\;
\left\{ \qquad 
\begin{array}{llll}
     \dot{\tup{q}}  
     \,=\,   -\tfrac{1}{u^{2n}} \hdge{\tup{\slangup}}\cdot\tup{q} &,
\\[6pt]
    \dot{u}  \,=\, \tfrac{1}{n^2}u^{2-2n} p_\ss{u} &,
 \end{array} \right.
\qquad 
 \begin{array}{llll}
    \dot{\tup{p}} 
    \,=\,  - \tfrac{1}{u^{2n}} \hdge{\tup{\slangup}}\cdot\tup{p} \,+\,  \tfrac{m+1}{u^{2n}}(\slang^2 + \tfrac{1}{n^2}u^2 p_\ss{u}^2)\tup{q} \,-  \pderiv{V}{\tup{q}} + \tup{\alphaup} 
\\[6pt]
    \dot{p}_u 
    \,=\, \tfrac{n}{u^{2n+1}} \big( \slang^2 \,+\, \tfrac{n-1}{n^3}u^2 p_\ss{u}^2 \big)   \,- \pderiv{V}{u} + \alpha_\ss{u}
\end{array}
\end{align}
\end{small}
The above could also be re-expressed using the integral of motion \eq{\lambda=\lambda_\zr=0}, implying \eq{\tup{q}\cdot\tup{p} - \tfrac{m+1}{n}u p_\ss{u}=0}, though this does not yield anything particularly helpful until values of \eq{n} and \eq{m} are chosen.

\begin{notesq}
    Although the above dynamics look rather hideous, with no obvious improvement over the original cartesian coordinate formulation, they are simplified tremendously in sections \ref{sec:central force} and \ref{sec:2BP} where we choose \eq{n=m=-1} and, importantly, where we transform the evolution parameter to something other than time. This yields linear dynamics for certain types of orbital motion.   
\end{notesq}

\paragraph*{The Conservative Force Terms, $\pderiv{V}{\bartup{q}}$} 
For some specified potential, \eq{V}, there are two methods we may use when calculating the conservative force terms, \eq{\pderiv{V}{\bartup{q}}}, appearing in  the equations of motion, Eq.\eqref{sum_qpdot_k1}.

(1.)
We may simply substitute \eq{\tup{r}=u^n q^m\tup{q}} from Eq.\eqref{PT_Qq} into \eq{V(\tup{r},t)} in order to obtain an expression for \eq{V(\bartup{q},t)}, then directly calculate \eq{\pderiv{V}{\bartup{q}}}.  If using this method, note from Remark \ref{rem:H_simp} that \eq{V(\bartup{q},t)} should \textit{not} be simplified using \eq{q=1}.
For example, consider a central-force potential of the form \eq{V(\tup{r})=\kconst  r^b} for some constants \eq{\kconst, b\in\mbb{R} } . 
Substitution of the relation \eq{r = u^n q^{m+1}} from Eq.\eqref{ears} leads to \eq{V(\bartup{q})} given by
\begin{small}
\begin{align} \label{V_examples}
    \begin{array}{rl}
    \fnsize{example}:
  &\quad V(r)\,=\, \kconst  r^b \;\,=\,\; V(\bartup{q})\,=\, \kconst  (u^{n}q^{m+1})^b  
\end{array}
\end{align}
\end{small}
We may then directly calculate \eq{\pderiv{V}{\tup{q}}} and \eq{\pderiv{V}{u}} from the above. As per Remark \ref{rem:q=k 1}, we may simplify the derivatives using \eq{q=1} \textit{after} differentiating \eq{V}, but we cannot simplify \eq{V} itself.\footnote{For the given example \eq{V = \kconst  r^b = \kconst  (u^{n}q^{m+1})^b}, this leads to: \\
\eq{\qquad\qquad  \pderiv{V}{\tup{q}} \,=\, \tfrac{b(m+1)}{q^2}\kconst  (u^{n}q^{m+1})^b \tup{q} \,=\, b(m+1)\kconst  u^{nb}\tup{q}  
\qquad,\qquad 
\pderiv{V}{u} \,=\, \tfrac{b n}{u}\kconst  (u^{n}q^{m+1})^b \,=\, nb\kconst  u^{nb-1} }. \\
Where \eq{q=1} is only used \textit{after} differentiating. If we had simplified the potential itself in Eq.\eqref{V_examples} using \eq{q=1} as \eq{V(\bartup{q})=\kconst  u^{nb}}, then we would have incorrectly obtained \eq{\pderiv{V}{\tup{q}}=\tup{0}}. }

(2.)
Alternatively, if we already know the gradient  \eq{\pderiv{V}{\tup{r}}} with respect to the original cartesian coordinates, \eq{\tup{r}}, then we may express 
\eq{\pderiv{V}{\tup{q}}} and \eq{\pderiv{V}{u}} as follows:
\begin{small}
\begin{align} \label{dvdq_dvdr}
\begin{array}{lllll}
      \pderiv{V}{\tup{q}} \,=\, \trn{\pderiv{\tup{r}}{\tup{q}}}  \cdot \pderiv{V}{\tup{r}}
      &=\,  u^n q^m  ( \imat_3 +m \htup{q}\otms\htup{q})  \cdot \pderiv{V}{\tup{r}}
      & \simeq\, u^n  ( \imat_3 + m \tup{q}\otms\tup{q})  \cdot \pderiv{V}{\tup{r}}
\\[6pt]
      \pderiv{V}{u} \,=\, \pderiv{\tup{r}}{u} \cdot \pderiv{V}{\tup{r}}
      &=\, n u^{n-1}q^{m+1} \htup{q} \cdot \pderiv{V}{\tup{r}}
      &\simeq\, n u^{n-1} \tup{q} \cdot \pderiv{V}{\tup{r}}
\end{array}
\end{align}
\end{small}
Where, as per Remark \ref{rem:q=k 1}, the last expressions have been simplified using \eq{q=1}. 
 For example, for the central-force potential \eq{V(\tup{r})=\kconst  r^b} from Eq.\eqref{V_examples}, 
 we obtain:
\begin{small}
\begin{align} \label{dVdr_examples}
    \begin{array}{rl}
    \fnsize{example }:
  &\qquad \pderiv{V}{\tup{r}} \,=\, b\kconst  r^{b-1}\htup{r}  \,=\, b\kconst  (u^n q^{m+1})^{b-1}\htup{q}
  \,\simeq\, b\kconst  u^{n(b-1)}\tup{q}
\end{array}
\end{align}
\end{small}
where we have used \eq{r=u^n q^{m+1}} and \eq{\htup{r}=\htup{q}}. 
The expressions on the far right-hand-side above have been simplified using \eq{q=1}. 
It can be verified that substituting the  above into Eq.\eqref{dvdq_dvdr} leads to the same result we would get if we differentiated \eq{V(\bartup{q})} in Eq.\eqref{V_examples}.

\paragraph*{The Nonconservative Force Terms, $\bartup{\alphaup}$.} 
Let \eq{\tup{a}^\nc\in\mbb{R}^3} be the inertial cartesian components of the total nonconservative forces (per unit mass).
Then, as seen in Appendix \ref{app:non_conserv}, the generalized forces (per unit mass), \eq{\bartup{\alphaup}=(\tup{\alphaup},\alpha_\ss{u})\in\mbb{R}^4}, appearing in the equation for \eq{\dot{\bartup{p}}} are given in terms of \eq{\tup{a}^\nc} by
\begin{small}
\begin{align} \label{g_nc} 
    \bartup{\alphaup}  \,=\, \trn{\pderiv{\tup{r}}{\bartup{q}}} \cdot \tup{a}^\nc 
&&\Rightarrow &&
\begin{array}{llll}
    \tup{\alphaup}  \,=\, \trn{\pderiv{\tup{r}}{\tup{q}}} \cdot \tup{a}^\nc
    \,=\,  u^n q^{m}(\imat_3 + m\htup{q}\otms \htup{q}) \cdot \tup{a}^\nc
    & \simeq\,  u^n (\imat_3 + m\tup{q}\otms \tup{q}) \cdot \tup{a}^\nc
\\[5pt]
   \alpha_\ss{u} \,=\, \pderiv{\tup{r}}{u}\cdot\tup{a}^\nc \,=\, nu^{n-1}q^{m+1}\htup{q} \cdot \tup{a}^\nc
    & \simeq\, nu^{n-1}\tup{q} \cdot \tup{a}^\nc
\end{array}
\end{align}
\end{small}
where the the partials were given in Eq.\eqref{B_nm} and, as before, the far right-hand-side has been simplified using \eq{q=1}. Note the expressions for \eq{\bartup{\alphaup}} in terms of the \textit{nonconservative} forces, \eq{\tup{a}^\nc}, are the same as
 the expressions for \eq{-\pderiv{V}{\bartup{q}}} in terms of the \textit{conservative} forces, \eq{-\pderiv{V}{\tup{r}}}, given in Eq.\eqref{dvdq_dvdr}. Readers may recognize it as the standard transformation rule for any 1-form/covector.

\subsection{Comparing \txi{m = 0} (the BF transformation) and \txi{m = –1} } \label{sec:BF_vs_mine}

 For the  parameter \eq{m} in the point transformation \eq{\tup{r}= u^n q^m \tup{q}} of Eq.\eqref{PT_Qq}, there are two values of particular interest: \eq{m=0} and \eq{m=-1}. The former corresponds to the classic general form of the projective point transformation,  \eq{\tup{r}=u^n\tup{q}}, which was used by Ferrándiz \cite{ferrandiz1987general} to develop a canonical extension of Burdet's projective point transformation (specifically, \eq{m=0} and \eq{n=\pm 1}). 
  The latter case (\eq{m=-1}) corresponds to \eq{\tup{r}=u^n\htup{q}}, which may seem redundant given that we impose \eq{q=1} and thus \eq{m} does not affect the physical meaning of \eq{\tup{q}=\htup{r}} or \eq{u=r^{\sfrac{1}{n}}}. However, the resulting Hamiltonian systems for \eq{m=0} and \eq{m=-1} turn out to have notably different properties and, although either can be used to linearize Kepler-type dynamics, the process is greatly simplified for the \eq{m=-1} case.  
  
The canonical projective transformation with \eq{m=0} (i.e., the BF transformation with \eq{n} arbitrary) is: 
\begin{small}
\begin{align} \label{BF_xform_narb}
m=0 \;\; \left\{ \quad 
\begin{array}{lllll}
      \tup{r} \,=\,  u^n \tup{q} 
  \\[5pt]
     \tup{v} \,=\,  \tfrac{1}{u^n} \big[ (\imat_3 - \htup{q}\otms\htup{q})\cdot\tup{p}
         \,+\, \tfrac{1}{n q} u p_\ss{u} \htup{q} \big]
\end{array}  \right.
&&\leftrightarrow &&
\begin{array}{llll}
     \tup{q} \,=\,  \htup{r} 
      &,\quad \tup{p}  \,=\,  r \tup{v}
\\[5pt]
  u \,=\, r^\ss{1/n}  
   &, \quad 
   p_\ss{u} \,=\, 
      n r^\ss{-1/n } \tup{r}\cdot\tup{v} 
\end{array}
\end{align}
\end{small}
where the relations on the right follow from restricting to \eq{q=1} and \eq{\lambda = \tfrac{1}{q}( \tup{q}\cdot\tup{p} - \tfrac{1}{n} u p_\ss{u}) = 0}. 
Compare the above (especially \eq{\tup{p}}) to
the canonical projective transformation with \eq{m=-1} given as follows: 
\begin{small}
\begin{align} \label{jtap_xform_narb}
m=-1 \;\; \left\{ \quad
\begin{array}{lllll}
      \tup{r} \,=\,  u^n \htup{q} 
  \\[5pt]
     \tup{v} \,=\,  \tfrac{q}{u^n} \big[ (\imat_3 - \htup{q}\otms\htup{q})\cdot\tup{p}
         \,+\, \tfrac{1}{n q} u p_\ss{u} \htup{q} \big]
\end{array}  \right.
&&\leftrightarrow &&
\begin{array}{llll}
     \tup{q} \,=\,  \htup{r} 
      &,\quad \tup{p}  \,=\,  r(\imat_3-\htup{r}\otms\htup{r})\cdot\tup{v}
     \,=\, -\hdge{\tup{\slangup}} \cdot \htup{r} 
\\[5pt]
  u \,=\, r^\ss{1/n}  
   &, \quad 
   p_\ss{u} \,=\, 
      n r^\ss{-1/n } \tup{r}\cdot\tup{v} 
\end{array}
\end{align}
\end{small}
where the relation on the right follow from restricting to \eq{q=1} and \eq{\lambda = \htup{q}\cdot\tup{p} = 0}.

    As seen above, and in Eq.\eqref{ears}-Eq.\eqref{qu_k1}, \eq{m}  \textit{does} effect the  meaning of \eq{\tup{p}} and \eq{\lambda}, as well as the dependence of \eq{r(\bartup{q})} and \eq{\mscr{H}(\bartup{q},\bartup{p},t)} on \eq{q=\mag{\tup{q}}}.
    For example, some differences following from \eq{m=0} versus \eq{m=-1}:
    \begin{small}
    \begin{align} \label{m_compare}
    \fnsize{if \eq{m=0}} \;\;
    \left\{\begin{array}{llll}
        \tup{r} = u^n\tup{q}
        &,\quad  \tup{p} = r\tup{v} 
    \\[5pt]
          r = u^n q
          &,\quad  p^2 = \slang^2 + r^2\dot{r}^2
    \\[5pt]
        \pd_{\tup{q}} f(r)   \neq 0
        &,\quad
        \lambda = \tfrac{1}{q}(\tup{q}\cdot\tup{p} - \tfrac{1}{n} u p_\ss{u}) 
    \end{array} \right.
    &&,&&
    \fnsize{if \eq{m=-1}} \;\;
    \left\{\begin{array}{llll}
        \tup{r} = u^n\htup{q}
        &,\quad
         \tup{p}   = \tup{\slangup}\tms \htup{r}  
    \\[5pt]
         r = u^n
         &,\quad
         p^2 = \slang^2
    \\[5pt]
         \pd_{\tup{q}} f(r)   = 0
        &,\quad
        \lambda = \htup{q}\cdot\tup{p}  
    \end{array} \right.
    \end{align}
    \end{small}
    where, in both cases, it holds that \eq{\tup{\slangup}=\tup{q}\tms\tup{p}}. 
    The properties for \eq{m=-1} on the right-hand-side of the above — that \eq{r} is only a function of \eq{u}, that \eq{p=\slang}, and that \eq{(\tup{q},\tup{p},\tup{\slangup})} are all mutually orthogonal — will lead to several desirable features compared to the \eq{m=0} case. 
For instance, with \eq{m=-1}, then \eq{\tup{q}} and \eq{\tup{p}} alone directly define the local vertical, local horizontal (LVLH) basis associated with the particle's instantaneous position: 
\begin{small}
\begin{flalign}
\quad
    \fnsize{if  $m=-1:$} 
    &&
 \begin{array}{llll}
    \begin{array}{cc}
         \fnsize{inertial cartesian coordinate}  \\
         \fnsize{representation of LVLH basis}
    \end{array}  
        =\, \{\htup{t}_r,\htup{t}_\tau, \htup{t}_\slang \}
        =  \{\tup{q},-\hdge{\htup{\slangup}}\cdot\htup{q},\htup{\slangup}\} \simeq \{\htup{q},\htup{p},\htup{\slangup}\}
\end{array}
 &&
\end{flalign}
\end{small}
with \eq{ \htup{\slangup} = \htup{q}\tms \htup{p}}, \eq{q^2 \simeq 1}, and \eq{\slang^2 \simeq p^2}.  

\paragraph*{Dynamics with $m=0$ vs.~$m=-1$.} 
Let us quickly compare the Hamiltonian dynamics in Eq.\eqref{sum_qpdot} for the cases \eq{m=0} and \eq{m=-1}. We will do so for the case that the potential function is of the form \eq{V=V^0(r) + V^1(\tup{r},t)} where \eq{V^0} accounts for central-forces and \eq{V^1} for all other arbitrary perturbations. That is, the cartesian coordinate Hamiltonian is given by
\begin{small}
\begin{align}
    \mscr{K} = \tfrac{1}{2}\v^2 + V^{0}(r) + V^{1}(\tup{r},t)
    &&,&&
    \dot{\tup{r}} = \tup{v} 
    \quad,\quad 
    \dot{\tup{v}} \,=\, -\pderiv{V^{0}}{r}\htup{r} \,+\, \tup{F}
     &&,&&
     \tup{F} \,:=\, -\pderiv{V^{1}}{\tup{r}} + \tup{a}^\nc
\end{align}
\end{small}
The above is transformed by the projective transformation with \eq{m=0} or \eq{m=-1} as given below:

\begin{notesq}
\textit{Dynamics with $m=0$.} 
For the cases \eq{m=0}, Eq.\eqref{sum_qpdot} leads to
\begin{small}
\begin{align}
m=0 \;\; \left\{ \;\;
\begin{array}{llllll}
     \mscr{H}  \,=\, 
    \tfrac{1}{u^{2n}q^2} \tfrac{1}{2} \big(\slang^2 + \tfrac{1}{n^2} u^2 p_\ss{u}^2 \big)  \,+\, V^{0}(\bartup{q}) \,+\, V^{1}(\bartup{q},t) 
\\[6pt] 
  \begin{array}{lll}
      \dot{q}_i = -\tfrac{1}{u^{2n}} \slang_{ij}q_j &,
\\[6pt]
    \dot{u}  = \tfrac{1}{n^2 u^{2n-2}}p_\ss{u} &,
 \end{array}
\quad
 \begin{array}{ll}
     \dot{p}_i = -\tfrac{1}{u^{2n}} \slang_{ij} p_j
     + \tfrac{1}{u^{2n}} ( \slang^2 + \tfrac{1}{n^2} u^2 p_\ss{u}^2) q_i
      - \pderiv{V^{0}}{q_i} + f_i
\\[6pt]
     \dot{p}_\ss{u}   =  \tfrac{n}{u^{2n+1}}\big(\slang^2 +  \tfrac{n-1}{n^3} u^2 p_\ss{u}^2 \big)   -  \pderiv{V^{0}}{u}  +  f_{u}
\end{array}  
&&\qquad
 \begin{array}{ll}
     \tup{f} 
     \,=\,  u^n  \tup{F}
\\[6pt]
    f_{u} 
    \,=\,  n u^{n-1}\tup{q}\cdot \tup{F}
\end{array}  
\end{array} \right.
\end{align}
\end{small}
where \eq{\tup{f}:=-\pderiv{V^1}{\tup{q}}+\tup{\alphaup}} and \eq{f_\ss{u}:=-\pderiv{V^1}{u}+\alpha_\ss{u}}, and
 where the above ODEs have been simplified with \eq{q=1}  \textit{after} differentiating \eq{\mscr{H}}. One could also rewrite the above using \eq{\lambda = \tfrac{1}{q}(\tup{q}\cdot\tup{p} - \tfrac{u}{n}p_\ss{u}) = 0}.  
 For \eq{m=0}, it does \textit{not} hold that \eq{\slang^2=p^2} nor that \eq{\tup{q}\cdot\tup{p}=0}. 
\end{notesq}

\begin{notesq}
\textit{Dynamics with $m=-1$.} 
For the cases \eq{m=-1}, Eq.\eqref{sum_qpdot} leads to
\begin{small}
\begin{align}
m=-1\;\; \left\{ \;\;
 \begin{array}{lllllll}
     \mscr{H}  \,=\, 
    \tfrac{1}{u^{2n}} \tfrac{1}{2} \big(\slang^2 + \tfrac{1}{n^2} u^2 p_\ss{u}^2 \big) + V^{0}(u) + V^{1}(\bartup{q},t)   
\\[6pt] 
   \begin{array}{lll}
     \dot{q}_i = -\tfrac{1}{u^{2n}} \slang_{ij}q_j&,
\\[6pt]
    \dot{u}  = \tfrac{1}{n^2 u^{2n-2}}p_\ss{u} &,
 \end{array}
\quad
 \begin{array}{ll}
    \dot{p}_i =
    -\tfrac{1}{u^{2n}} \slang_{ij} p_j
     \,+\, f_i 
\\[6pt]
     \dot{p}_\ss{u}   =
    \tfrac{n}{u^{2n+1}}\big(\slang^2 +  \tfrac{n-1}{n^3} u^2 p_\ss{u}^2 \big)   -  \pderiv{V^{0}}{u}  +  f_{u} 
    \end{array}  
&&\qquad
 \begin{array}{ll}
     \tup{f} 
     \,=\,  \tfrac{u^n}{q}(\imat_3 - \htup{q}\otms\htup{q})\cdot  \tup{F}
\\[6pt]
    f_{u} 
    \,=\,  n u^{n-1}\htup{q}\cdot \tup{F}
\end{array}  
\end{array} \right.
\end{align}
\end{small}
where \eq{\tup{f}:=-\pderiv{V^1}{\tup{q}}+\tup{\alphaup}} and \eq{f_\ss{u}:=-\pderiv{V^1}{u}+\alpha_\ss{u}}.
The above has \textit{not} been simplified with \eq{q=1} or \eq{\lambda=\htup{q}\cdot\tup{p}=0}. Even so, the above dynamics, especially for \eq{\tup{p}}, are much simpler than in the \eq{m=0} case.
For  instance: 
\begin{small}
\begin{itemize}
    \item For the case \eq{m=-1} then \eq{r = u^n}  such that \eq{\pderiv{V^{0}}{\tup{q}}=0}; any central-forces drop out of the ODE for \eq{\tup{p}}, now appearing only on the ODE for \eq{p_\ss{u}} (not true for the \eq{m=0} case). 
    \item Even without simplifying using \eq{q=1} and \eq{\lambda=0}, the above dynamics are, after a transformation of the evolution parameter, are already \textit{linear} for unperturbed Kepler or Manev dynamics. The ``\eq{(\tup{q},\tup{p})}-part'' will be fully linear for \textit{any} central-force dynamics. This is detailed in sections \ref{sec:central force} and \ref{sec:2BP}. 
    \item With \eq{m=-1}, we conveniently have that the integral \eq{\lambda} is now \eq{\lambda=\htup{q}\cdot\tup{p}} such that \eq{\lambda=0} is equivalent to \eq{\tup{q}\cdot\tup{p} = 0}. This can be used to simplify many expressions.  For instance, simplifying the above ODEs with \eq{q=1} and \eq{\lambda=\htup{q}\cdot\tup{p}=0} (which also means \eq{\slang^2 \simeq p^2}) leads to:
    \begin{small}
    \begin{align}
       \begin{array}{lll}
         \dot{\tup{q}} \simeq \tfrac{1}{u^{2n}}\tup{p} &,
    \\[6pt]
        \dot{u}  = \tfrac{1}{n^2 u^{2n-2}}p_\ss{u} &,
     \end{array}
    \quad
     \begin{array}{ll}
        \dot{\tup{p}} \simeq
        -\tfrac{1}{u^{2n}}  p^2\tup{q} 
         \,+\, \tup{f} 
         \;\simeq\;  -\tfrac{1}{u^{2n}}  \slang^2\tup{q} 
         \,+\, \tup{f} 
    \\[6pt]
         \dot{p}_\ss{u}   \simeq
        \tfrac{n}{u^{2n+1}}\big(p^2 +  \tfrac{n-1}{n^3} u^2 p_\ss{u}^2 \big)   -  \pderiv{V^{0}}{u}  + f_{u} 
        \end{array}  
    \end{align}
    \end{small}
    \item  In fact, the above are equivalent to what we would get if we used \eq{\lambda=\htup{q}\cdot\tup{p}=0} to simplify \eq{\mscr{H}} itself \textit{before} differentiating (and then used \eq{q=1} after obtaining the ODEs). That is, with \eq{m=-1} we could, if we wished, use the simplified Hamiltonian:
\begin{small}
\begin{align}
      \mscr{H}_\ss{\mrm{simplified}} =
    \tfrac{1}{u^{2n}} \tfrac{1}{2}\big( q^2 p^2  + \tfrac{1}{n^2} u^2 p_\ss{u}^2 \big) + V^{0}(u) + V^{1}(\bartup{q},t)
\end{align}
\end{small}
Yet, in contrast to the BF transformation, we find little use for the above simplified Hamiltonian in this work. 
\end{itemize}
\end{small}
\end{notesq}


\subsection{Derivation of two ``extra'' integrals of motion} \label{sec:lambda_const} 
It was claimed in Remark \ref{rem:q=k 1}, Eq.\eqref{qlam=const}, that \eq{q=\mag{\tup{q}}} and \eq{\lambda=q^\ss{-2}  (\tup{q}\cdot\tup{p} - \tfrac{m+1}{n}u p_\ss{u})} are integrals of motion of the general projective coordinate Hamiltonian given in Eq.\eqref{Hqp_general} as \eq{\mscr{H} =  \tfrac{1}{u^{2n}q^{2m+2}} \tfrac{1}{2}\big( \slang^2 + \tfrac{1}{n^2}u^2 p_\ss{u}^2 \big)  + V}, where \eq{\slang^2 = q^2 p^2  - (\tup{q}\cdot\tup{p})^2}. 
We will now prove that \eq{q} and \eq{\lambda} are integrals of motion by showing that  \eq{\dot{q}=0} and \eq{\dot{\lambda}=0} using the Poisson bracket expression for any function \eq{f(\bartup{q},\bartup{p},t)}:
\begin{small}
\begin{align} \label{eom_pbrack}
    \dot{f} \,=\, \pbrak{f}{\mscr{H}} \,+\, \pderiv{f}{\bartup{p}}\cdot\bartup{\alphaup} \,+\, \pderiv{f}{t}
\end{align}
\end{small}
where the above follows from Eq.\eqref{eom_pbrack1} in Appendix \ref{app:non_conserv} on Hamiltonian dynamics with nonconservative forces.

\paragraph*{Proof that \eq{q=\mag{\tup{q}}} is constant.}
Noting that \eq{\pderiv{q}{\bartup{p}}=\bartup{0}} and  \eq{\pderiv{q}{t}=0}, it is strait forward to see that Eq.\eqref{eom_pbrack} leads to \eq{\dot{q}=0}:
\begin{small}
\begin{align} \label{q=k}
\begin{array}{rl}
       \dot{q} \,=\, \pbrak{q}{\mscr{H}} 
    \,=\, \pderiv{q}{\bartup{q}}\cdot\pderiv{\mscr{H}}{\bartup{p}}  \,-\, \cancel{\pderiv{q}{\bartup{p}}}\cdot\pderiv{\mscr{H}}{\bartup{q}}
    \,=\,
    \htup{q}\cdot\pderiv{\mscr{H}}{\tup{p}}
    \,=\,
    \htup{q}\cdot\big(
    \tfrac{1}{u^{2n}q^{2m+2}} \big[ q^2\tup{p}-(\tup{q}\cdot\tup{p})\tup{q} \big]
    \big) 
    \,=\,  -\tfrac{1}{u^{2n}q^{2m+2}} \htup{q}\cdot\hdge{\tup{\slangup}}\cdot\tup{q}
    \,=\, 0
\end{array}
\end{align}
\end{small}
where \eq{\mscr{H}} is given  by Eq.\eqref{Hqp_general}. 
Thus, the constraint \eq{\varphi=\tfrac{1}{2}(q^2-\qsc^2)=0} 
— i.e., \eq{q=\qsc} for constant \eq{\qsc>0} —   
 is, in fact,  an integral of motion (regardless of the forces present). The above leads immediately to \eq{\pbrak{\varphi}{\mscr{H}}=0}.

\paragraph*{Proof that \eq{\lambda} is constant.}
Consider \eq{\dot{\lambda}} as given by the Poisson bracket expression of Eq.\eqref{eom_pbrack}:
\begin{small}
\begin{align} \label{d_lambda_1}
    \dot{\lambda} \,=\, \pbrak{\lambda}{\mscr{H}} \,+\, \pderiv{\lambda}{\bartup{p}}\cdot\bartup{\alphaup}
    \,=\, 
    \pderiv{\lambda}{\bartup{q}}\cdot\pderiv{\mscr{H}}{\bartup{p}}  \,-\, \pderiv{\lambda}{\bartup{p}}\pderiv{\mscr{H}}{\bartup{q}} \,+\, \pderiv{\lambda}{\bartup{p}}\cdot\bartup{\alphaup}
\end{align}
\end{small}
 Where \eq{\lambda} is given by Eq.\eqref{P_pq} as \eq{\lambda=q^\ss{-2}  (\tup{q}\cdot\tup{p} - \tfrac{m+1}{n}u p_\ss{u})}. The above then leads to \eq{\dot{\lambda}} given by
\begin{small}
\begin{align} \label{lam_dot_long}
    \dot{\lambda} &\,=\, \pbrak{\tfrac{1}{q^2}(\tup{q}\cdot\tup{p} - \tfrac{m+1}{n}u p_\ss{u})}{\mscr{H}}
    \,+\, \bartup{\alphaup}\cdot\pderiv{}{\bartup{p}}(\tup{q}\cdot\tup{p} - \tfrac{m+1}{n}u p_\ss{u})\tfrac{1}{q^2}
\\ \nonumber
    &\,=\,
    \tfrac{1}{q^2} \Big( (\tup{p}- 2\lambda\tup{q})\cdot\pderiv{\mscr{H}}{\tup{p}} \,-\, \tfrac{m+1}{n}p_\ss{u}\pderiv{\mscr{H}}{p_u} \,-\, \tup{q} \cdot \pderiv{\mscr{H}}{\tup{q}} \,+\, \tfrac{m+1}{n}u\pderiv{\mscr{H}}{u}
    \,+\, \tup{q}\cdot\tup{\alphaup} \,-\,  \tfrac{m+1}{n}u\alpha_\ss{u} \Big)
\\ \nonumber
    &\,=\,  \tfrac{1}{q^2}\Bigg\{
    \tfrac{1}{u^{2n}q^{2m+2}} \Big( \big[ q^2 p^2 -(\tup{q}\cdot\tup{p})^2   -  \tfrac{m+1}{n}\tfrac{1}{n^2}u^2 p_\ss{u}^2 \big]
    -2\lambda\big[ \cancel{ q^2\tup{q}\cdot\tup{p}-\tup{q}\cdot\tup{p}q^2 } \big] \Big)
\\ \nonumber
    &\quad\;\;\dots \, + \,
     \tfrac{1}{u^{2n}q^{2m+2}}\big[mq^2 p^2  - m(\tup{q}\cdot\tup{p})^2  +  \tfrac{m+1}{n^2}u^2 p_\ss{u}^2 \big] 
        -  \pderiv{V}{\tup{q}}\cdot\tup{q}
\\ \nonumber
    &\quad\;\;\dots\, -\,
    \tfrac{m+1}{u^{2n}q^{2m+2}}\big[q^2 p^2 -(\tup{q}\cdot\tup{p})^2  +  \tfrac{n-1}{n^3}u^2 p_\ss{u}^2 \big]    +  \tfrac{m+1}{n}u\pderiv{V}{u}
      \,+\, \tup{q}\cdot\tup{\alphaup} -  \tfrac{m+1}{n}u\alpha_\ss{u} \Bigg\}
\end{align}
\end{small}
Simplifying and collecting like terms, we obtain
\begin{small}
\begin{align} \nonumber
    \dot{\lambda} &\,=\,  
   \tfrac{1}{u^{2n}q^{2m+4}}\Big\{ 
    \cancel{ [1+m - (m+1)] }\big[q^2 p^2 -(\tup{q}\cdot\tup{p})^2\big]
    \,+\, \tfrac{(m+1)}{n^2}u^2 p_\ss{u}^2 \cancel{ (1-\tfrac{1}{n}  - \tfrac{n-1}{n}) }   \Big\} 
\\[4pt] \nonumber
    &\qquad\;\;\dots\; 
    +\;  \tfrac{1}{q^2}(\tfrac{m+1}{n}u\pderiv{V}{u} \,-\, \pderiv{V}{\tup{q}}\cdot\tup{q} ) \,-\,   \tfrac{1}{q^2}(\tfrac{m+1}{n}u\alpha_\ss{u} \,-\, \tup{q}\cdot\tup{\alphaup} )
\end{align}
\end{small}
It is seen that the entire first line in brackets — the terms arising purely from the kinematics — reduces to zero, leaving only the terms involving \eq{\pderiv{V}{\bartup{q}}} and \eq{\bartup{\alphaup}} — the terms arising from conservative and nonconservative forces, respectively.  However, from Eq.\eqref{dvdq_dvdr} and Eq.\eqref{g_nc}, we see that these terms also cancel:
\begin{small}
\begin{align}
    \tfrac{m+1}{n}u\pderiv{V}{u}  - \pderiv{V}{\tup{q}}\cdot\tup{q} 
    &\;\,=\,\;
    \tfrac{m+1}{n}u u^{n-1}q^m\pderiv{V}{\tup{r}}\cdot\tup{q}  \,-\, u^n q^m(m+1)\pderiv{V}{\tup{r}}\cdot\tup{q} \;\,=\,\; 0
\\[5pt]
    \tfrac{m+1}{n}u\alpha_\ss{u} \,-\, \tup{q}\cdot\tup{\alphaup} &\;\,=\,\;  \tfrac{m+1}{n}u u^{n-1}q^m\tup{a}^\nc\cdot\tup{q}  \,-\, u^n q^m(m+1)\tup{a}^\nc\cdot\tup{q} \;\,=\,\; 0
\end{align}
\end{small}
The entire right-hand-side of Eq.\eqref{d_lambda_1} then reduces to zero. I.e., \eq{\lambda} is an integral of motion:
\begin{small}
\begin{align} \label{lambdadot_simp}
        \dot{\lambda} \,=\, \pbrak{\lambda}{\mscr{H}} \,+\, \pderiv{\lambda}{\bartup{p}}\cdot\bartup{\alphaup} \,=\, 0
\end{align}
\end{small}
Note that this holds regardless of whether the forces involved are conservative or nonconservative and that the derivation did not require us to assume that \eq{q} is constant at any point. 

%% file: Mysecs_prj/apx_angMoment.tex
\subsection{Useful angular momentum relations} \label{sec:ang_momentum}

We collect some useful relations involving projective coordinates and the angular momentum, which plays a central role throughout this  work. 
Recall that the angular momentum coordinate vector, \eq{\tup{\slangup}\in\mbb{R}^3}, the  antisymmetric matrix, \eq{\hdge{\tup{\slangup}}\in\somat{3}},  and  the magnitude, \eq{\slang}, are given in terms of inertial cartesian coordinates \eq{(\tup{r},\tup{v})\in\mbb{R}^6} by the usual relations:
\begin{small}
\begin{align} \label{l_rv_general}
      \tup{\slangup} \,=\, \tup{r}\tms \tup{v} 
 \qquad,\qquad
       \hdge{\tup{\slangup}} \,=\,  \tup{r} \wdg \tup{v}  
 \qquad,\qquad
        \slang^2  \,=\, \mag{\tup{\slangup}} 
      \,=\, r^2 \v^2  - (\tup{r}\cdot\tup{v})^2
\end{align}
\end{small}
where \eq{\tup{\slangup}} and \eq{\hdge{\tup{\slangup}}} are Hodge duals of one another: \eq{\slang_{ij}  = \lc_{ijk}\slang_k \leftrightarrow \slang_i = \tfrac{1}{2}\lc_{ijk}\slang_{jk}}.\footnote{We write the components of \eq{\hdge{\tup{\slangup}}} simply as \eq{\slang_{ij}} rather than \eq{\hdge{\slang}_{ij}}.}
For the following, it may be helpful to recall the relations for the Hodge dual on \eq{\mbb{R}^3} given in Eq.\eqref{hodge_cord}.

\paragraph*{General Relations.}
The functions in Eq.\eqref{l_rv_general} are invariant under a family of ``canonically extended'' projective transformations, \eq{\mbb{R}^8\ni(\tup{q},u,\tup{p},p_\ss{u})\mapsto (\tup{r},\tup{v})\in \mbb{R}^6}, of the general form developed in section \ref{sec:gen_burdet}:
\begin{small}
\begin{align} \label{prj_family}
\; \left. \;\;
\begin{array}{lllll}
      \tup{r} \,=\, u^n q^m \tup{q} 
\\[4pt]
     \tup{v} \,=\,  \tfrac{1}{u^n q^m} \big(  (\imat_3 - \htup{q}\otms\htup{q})\cdot\tup{p} \,+\, \tfrac{u}{n q} p_\ss{u}\htup{q}\big)
\end{array} \right.
&&,&&
  \left(\begin{array}{lllll}
         \varphi = q-1  = 0
     \\[3pt]
         \lambda = \tfrac{1}{q}(\tup{q}\cdot\tup{p} - \tfrac{m+1}{n} u p_\ss{u} )
    \end{array}\right)
\end{align}
\end{small}
for any real numbers \eq{n\neq0,m\in\mbb{R}}. 
It can be verified by direct substitution of the above \eq{\tup{r}(\tup{q},u)} and \eq{\tup{v}(\tup{q},u,\tup{p},p_\ss{u})} that the angular momentum functions in Eq.\eqref{l_rv_general} are invariant under any such transformation:
\begin{small}
\begin{align} \label{l_qp_general}
\begin{array}{rllllll}
      \tup{\slangup} &\,=\, \tup{r}\tms \tup{v} \,=\, \hdge{\tup{v}}\cdot\tup{r}
  \\[3pt]
     &\,=\, \tup{q}\tms \tup{p}  \,=\, \hdge{\tup{p}}\cdot\tup{q}
\end{array}
 &&,&&
\begin{array}{rllllll}
       \hdge{\tup{\slangup}} &\,=\,  \tup{r} \wdg \tup{v}  \,=\, \tup{r}\otms \tup{v}-\tup{v}\otms \tup{r}
  \\[3pt]
     &\,=\,  \tup{q} \wdg \tup{p}  \,=\, \tup{q}\otms \tup{p}-\tup{p}\otms \tup{q}
\end{array}
 &&,&&
\begin{array}{rllllll}
        \slang^2  = \mag{\tup{\slangup}}^2
      &\,=\, r^2 \v^2  - (\tup{r}\cdot\tup{v})^2
  \\[3pt]
     &\,=\, q^2 p^2  - (\tup{q}\cdot\tup{p})^2
\end{array}
\end{align}
\end{small}
We further note the following useful relations (where \eq{\slang_{ij}=q_i p_j-p_i q_j=-\slang_{ji}}):
\begin{small}
\begin{flalign} \label{angmoment_rels_crd}
\;\;
\begin{array}{lllllll}
    \slang^2 =  \slang_i \slang_i = \tfrac{1}{2} \slang_{ij} \slang_{ij} 
    = \tup{q}\cdot \hdge{\tup{\slangup}} \cdot\tup{p}  = q^2 p^2 - (\tup{q}\cdot\tup{p})^2
\\[4pt]
     0 =  \slang_{ij}q_i q_j = \slang_{ij}p_i p_j  = \tup{q}\cdot \hdge{\tup{\slangup}} \cdot\tup{q} = \tup{p}\cdot \hdge{\tup{\slangup}} \cdot\tup{p}
\end{array} 
&&,&&
\begin{array}{lll}
     \hdge{\tup{\slangup}} \cdot\hdge{\tup{\slangup}} \cdot\tup{q} = -\slang^2 \tup{q}
\\[4pt]
     \hdge{\tup{\slangup}} \cdot\hdge{\tup{\slangup}} \cdot\tup{p} = -\slang^2 \tup{p}
\end{array}
&&,&&
\begin{array}{lll}
     \hdge{\tup{\slangup}} \cdot\tup{q} = -(q^2 \imat_3 - \tup{q}\otms\tup{q}) \cdot\tup{p} 
     =  \hdge{\tup{q}} \cdot \hdge{\tup{q}} \cdot\tup{p}
\\[4pt]
     \hdge{\tup{\slangup}} \cdot\tup{p} = (p^2 \imat_3 - \tup{p}\otms\tup{p}) \cdot\tup{q} 
     = -\hdge{\tup{p}}\cdot\hdge{\tup{p}}\cdot\tup{q}
\end{array}
\end{flalign}
\end{small}         
as well as the following relations for various partial derivatives (where \eq{\pd^2_{\tup{b}\tup{a}} := \pd_{\tup{b}}\pd_{\tup{a}}}):
\begin{small}
\begin{align} \label{angmoment_rels_crd2}
 &\begin{array}{llllll}
      \pd_{\tup{q}}
     \tfrac{1}{2}\slang^2 \,=\,  \hdge{\tup{\slangup}}\cdot\tup{p} 
\\[4pt]
      \pd_{\tup{p}}
     \tfrac{1}{2}\slang^2 \,=\,   -\hdge{\tup{\slangup}}\cdot\tup{q}
\end{array}
&&,&&
\begin{array}{rllllll}
    \pd^2_{\tup{q}\tup{q}} \tfrac{1}{2}\slang^2  =&
  \pd_{\tup{q}} (\hdge{\tup{\slangup}}\cdot\tup{p})
  \,=\,  p^2 \imat_3 - \tup{p}\otms\tup{p} 
  &=\,  - \hdge{\tup{p}}\cdot\hdge{\tup{p}}
\\[4pt] 
    -\pd^2_{\tup{p}\tup{p}} \tfrac{1}{2}\slang^2 =&
    \pd_{\tup{p}} (\hdge{\tup{\slangup}}\cdot\tup{q}) 
    \,=\,   \tup{q}\otms\tup{q} - q^2 \imat_3
    &=\, \hdge{\tup{q}}\cdot\hdge{\tup{q}}  
\\[4pt]
     -\pd^2_{\tup{q}\tup{p}} \tfrac{1}{2}\slang^2 =&
   \pd_{\tup{q}}
      (\hdge{\tup{\slangup}}\cdot\tup{q}) \,=\, \hdge{\tup{\slangup}} - \hdge{\tup{q}}\cdot\hdge{\tup{p}} 
  \\[4pt]
    \pd^2_{\tup{p}\tup{q}} \tfrac{1}{2}\slang^2 =&
     \pd_{\tup{p}} (\hdge{\tup{\slangup}}\cdot\tup{p}) \,=\, \hdge{\tup{\slangup}} +  \hdge{\tup{p}}\cdot\hdge{\tup{q}}
\end{array}
\\[6pt] \nonumber 
&\begin{array}{llllll}
  \pd_{\tup{q}}
     \slang \,=\,   \hdge{\htup{\slangup}}\cdot\tup{p}
\\[4pt]
      \pd_{\tup{p}} \slang \,=\,   -\hdge{\htup{\slangup}}\cdot\tup{q}
\end{array}
&&,&& 
\begin{array}{rllllll}
      \pd^2_{\tup{q}\tup{q}} \slang =&
     \pd_{\tup{q}} (\hdge{\htup{\slangup}}\cdot\tup{p}) 
     \,=\,
      \tfrac{p^2}{\slang} \big( \imat_3 - \htup{p}\otms\htup{p} - (\hdge{\htup{\slangup}}\cdot\htup{p})\otms (\hdge{\htup{\slangup}}\cdot\htup{p}) \big)
     &= \tfrac{p^2}{\slang} \htup{\slangup} \otms \htup{\slangup}
\\[4pt] 
    -\pd^2_{\tup{p}\tup{p}} \slang =&
   \pd_{\tup{p}} (\hdge{\htup{\slangup}}\cdot\tup{q})
   \,=\, 
   - \tfrac{q^2}{\slang} \big( \imat_3 - \htup{q}\otms\htup{q} - (\hdge{\htup{\slangup}}\cdot\htup{q})\otms (\hdge{\htup{\slangup}}\cdot\htup{q}) \big)
    &=- \tfrac{q^2}{\slang} \htup{\slangup} \otms \htup{\slangup}
\\[4pt]
     -\pd^2_{\tup{q}\tup{p}} \slang =&
   \pd_{\tup{q}}
      (\hdge{\htup{\slangup}}\cdot\tup{q}) \,=\, \tfrac{1}{\slang} \big( \hdge{\tup{\slangup}} - \hdge{\tup{q}}\cdot\hdge{\tup{p}} -   (\hdge{\htup{\slangup}}\cdot\tup{q})\otms (\hdge{\htup{\slangup}}\cdot\tup{p}) \big)
\\[4pt]
    \pd^2_{\tup{p}\tup{q}} \slang =&
     \pd_{\tup{p}} (\hdge{\htup{\slangup}}\cdot\tup{p}) \,=\, \tfrac{1}{\slang} \big( \hdge{\tup{\slangup}} +  \hdge{\tup{p}}\cdot\hdge{\tup{q}} +   (\hdge{\htup{\slangup}}\cdot\tup{p})\otms (\hdge{\htup{\slangup}}\cdot\tup{q}) \big)
 \end{array}
\end{align}
\end{small}
where we have 
used\footnote{In more detail:
\begin{align} \nonumber
\begin{array}{rllll}
     \pd^2_{\tup{q}\tup{q}} \slang =&
     \pd_{\tup{q}} (\hdge{\htup{\slangup}}\cdot\tup{p}) \,=\,
     -\tfrac{1}{\slang} \big(  \hdge{\tup{p}}\cdot\hdge{\tup{p}} +  (\hdge{\htup{\slangup}}\cdot\tup{p})\otms (\hdge{\htup{\slangup}}\cdot\tup{p}) \big)
     &=\, 
     \tfrac{p^2}{\slang} \big( \imat_3 - \htup{p}\otms\htup{p} - (\hdge{\htup{\slangup}}\cdot\htup{p})\otms (\hdge{\htup{\slangup}}\cdot\htup{p}) \big)
     \,=\,  \tfrac{p^2}{\slang} \htup{\slangup} \otms \htup{\slangup}
\\[4pt]
     -\pd^2_{\tup{p}\tup{p}} \slang =&
   \pd_{\tup{p}} (\hdge{\htup{\slangup}}\cdot\tup{q}) \,=\,
   \tfrac{1}{\slang} \big(  \hdge{\tup{q}}\cdot\hdge{\tup{q}} + (\hdge{\htup{\slangup}}\cdot\tup{q})\otms (\hdge{\htup{\slangup}}\cdot\tup{q}) \big)
   &=\, 
   - \tfrac{q^2}{\slang} \big( \imat_3 - \htup{q}\otms\htup{q} - (\hdge{\htup{\slangup}}\cdot\htup{q})\otms (\hdge{\htup{\slangup}}\cdot\htup{q}) \big)
   \,=\,  - \tfrac{q^2}{\slang} \htup{\slangup} \otms \htup{\slangup}
\end{array}
\end{align}
}
the below relation, following from the fact that \eq{\{\htup{q}, -\hdge{\htup{\slangup}}\cdot\htup{q}, \htup{\slangup}\} } and \eq{\{\htup{p}, -\hdge{\htup{\slangup}}\cdot\htup{p}, \htup{\slangup}\} } are orthonormal bases: 
\begin{small}
\begin{align} \label{angmoment_rels_iden}
\imat_3 \,=\,  
     \htup{q}\otms\htup{q} + (\hdge{\htup{\slangup}}\cdot\htup{q})\otms (\hdge{\htup{\slangup}}\cdot\htup{q}) +  \htup{\slangup} \otms \htup{\slangup}
    \,=\,  \htup{p}\otms\htup{p} + (\hdge{\htup{\slangup}}\cdot\htup{p})\otms (\hdge{\htup{\slangup}}\cdot\htup{p}) +  \htup{\slangup} \otms \htup{\slangup}
\end{align}
\end{small}

\noindent For the above angular momentum relations in Eq.\eqref{l_qp_general}-Eq.\eqref{angmoment_rels_iden}, we note:
\begin{small}
\begin{enumerate}
    \item Everything holds exactly the same with \eq{(\tup{q},\tup{p})} replaced by inertial cartesian coordinates \eq{(\tup{r},\tup{v})}. 
    \item More generally, if, for any \eq{\tup{x},\tup{y}\in\mbb{R}^3} we define \eq{\hdge{\tup{z}}:=\tup{x}\wdg\tup{y}} (i.e, \eq{\tup{z}:=\tup{x}\tms\tup{y}=\hdge{\tup{y}}\cdot\tup{x}}), then everything in Eq.\eqref{angmoment_rels_crd}-Eq.\eqref{angmoment_rels_iden} holds with \eq{(\tup{q},\tup{p},\tup{\slangup})} replaced by \eq{(\tup{x},\tup{y},\tup{z})}.  
    \item Nothing depends on \eq{(u,p_\ss{u})} nor the choices of \eq{n\neq 0,m\in\mbb{R}} in the family of transformations in Eq.\eqref{prj_family}.
    \item Nothing involves the functions \eq{\varphi} and \eq{\lambda} in Eq.\eqref{prj_family}. Nothing is simplified using \eq{q=1} or \eq{\lambda=0}. 
\end{enumerate}
\end{small}

\paragraph*{Simplified Relations (in the case \eq{m=-1}).}
 As mentioned, the above angular momentum relations do not involve the functions \eq{\varphi} or \eq{\lambda} in Eq.\eqref{prj_family}. 
 Yet, one \textit{could}, if they wished, rewrite some of the above using \eq{\varphi=0} (i.e., \eq{q=1}) and \eq{\lambda=0}. The reasons for this is the same as
 previously discussed\footnote{It was shown in section \ref{sec:lambda_const} that \eq{q=\mag{\tup{q}}} and \eq{\lambda} are integrals of motion of the Hamiltonian system in any projective coordinates defined by Eq.\eqref{prj_family} for any \eq{n,m\in\mbb{R}} (note \eq{q} being an integral of motion is equivalent to \eq{\varphi} being an integral of motion). This holds for any arbitrary forces, conservative or nonconservative. Furthermore, we may choose to limit consideration to the values \eq{q=1} and \eq{\lambda=0} as this places no restrictions on the inertial cartesian coordinates. }. 
This does not lead to anything particularly interesting until specific values of \eq{n,m\in\mbb{R}} are chosen. In particular, only for the choice \eq{m=-1} does rewriting Eq.\eqref{l_qp_general}-Eq.\eqref{angmoment_rels_iden} using \eq{q=1} and \eq{\lambda=0} lead to notable simplifications, detailed below.

This work ultimately prefers a projective transformation which corresponds to \eq{m=n=-1} in Eq.\eqref{prj_family}.\footnote{That is, much of this work uses a transformation given as follows: 
\begin{align} \label{prj__family}
\begin{array}{lllll}
      \tup{r} \,=\, \tfrac{1}{u} \htup{q} 
      &,\quad 
      \varphi = q-1 = 0 
\\[4pt]
     \tup{v} \,=\,  u q (\imat_3 - \htup{q}\otms\htup{q})\cdot\tup{p} - u^2 p_\ss{u}\htup{q} 
     &,\quad  \lambda = \htup{q}\cdot\tup{p} 
\end{array} 
&& \longleftrightarrow &&
\begin{array}{llll}
     \tup{q} \,=\,  \htup{r} 
     &,\quad 
      u \,=\, 1/r  
 \\[4pt]
     \tup{p}  
     \,=\, -(\tup{r}\wdg\tup{v})\cdot\htup{r}
     \,=\, -\hdge{\tup{\slangup}}\cdot\htup{r}
     &,\quad
      p_\ss{u} \,=\,  - r^2 \htup{r}\cdot\tup{v}
\end{array}
\end{align}
}
Although the value of \eq{n} is inconsequential for the present discussion, the value of \eq{m} is not: with \eq{m=-1}, the multiplier \eq{\lambda} is then given by \eq{\lambda=\htup{q}\cdot\tup{p}} such that limiting consideration to \eq{\lambda=0} means that \eq{(\tup{q},\tup{p},\tup{\slangup})} are all mutually orthogonal (this holds only for \eq{m=-1}). Together with \eq{q=1}, this means \eq{\tup{q}} and \eq{\tup{p}} satisfy the below relations (from Eq.\eqref{qpl_rels}): 
\begin{small}
\begin{flalign} \label{qpl_rels_apx}
    \qquad
    \begin{array}{lll}
          \tup{\slangup}=\tup{q}\tms \tup{p}
        \\[4pt]
          \tup{q}\cdot\tup{\slangup} = \tup{p}\cdot\tup{\slangup}  = 0
       \\[4pt]
        ` \slang^2 = q^2 p^2-(\tup{q}\cdot\tup{p})^2 
    \end{array}
    &&,&&
    \begin{array}{cc}
        \fnsize{using:} \\
         q=1  \\[1pt]
         \lambda=\htup{q}\cdot\tup{p}=0 
    \end{array}
     \Rightarrow \;\; 
    \left\{ \quad
    \begin{array}{ll}
        \tup{q} \simeq \htup{q} \simeq \htup{p}\tms \htup{\slangup} 
      \\[4pt]
        \tup{p} \simeq \tup{\slangup}\tms \htup{q} 
    \\[4pt]
          \slang^2 \simeq p^2  
    \end{array} \right.
    \;\;,\quad 
     \begin{array}{lll}
        \hdge{\tup{\slangup}}\cdot\tup{q} \simeq -q^2\tup{p} \simeq -\tup{p}
        \\[4pt]
        \hdge{\tup{\slangup}}\cdot\tup{p} \simeq p^2\tup{q} \simeq \slang^2 \tup{q} 
        \\[4pt]
          \imat_3 \simeq \htup{q} \otms\htup{q} + \htup{p} \otms\htup{p} + \htup{\slangup} \otms\htup{\slangup}
    \end{array}
    \quad
\end{flalign}
\end{small}
(where ``\eq{\simeq}'' indicates relations that follow from \eq{q=1} and \eq{\lambda=0}). 
In particular, \eq{ \{\htup{q},\htup{p},\htup{\slangup}\} \simeq \{\htup{t}_r,\htup{t}_\tau, \htup{t}_\slang \}} are the inertial cartesian components of the LVLH basis. 
Using the above relations, several derivatives in Eq.\eqref{angmoment_rels_crd2} simplify to:
\begin{small}
\begin{align} \label{angmoment_rels_simplified}
 &\begin{array}{llllll}
      \pd_{\tup{q}}
     \tfrac{1}{2}\slang^2 \,=\, \hdge{\tup{\slangup}}\cdot\tup{p} 
     &\simeq\, \slang^2 \tup{q}
\\[4pt]
      \pd_{\tup{p}}
     \tfrac{1}{2}\slang^2 \,=\, -\hdge{\tup{\slangup}}\cdot\tup{q}
     &\simeq\,  \tup{p}
\end{array}
&&,&&
\begin{array}{llllll}
      \pd_{\tup{q}} (\hdge{\tup{\slangup}}\cdot\tup{p}) \,=\,  - \hdge{\tup{p}}\cdot\hdge{\tup{p}} 
\\[4pt] 
   \pd_{\tup{p}} (\hdge{\tup{\slangup}}\cdot\tup{q}) \,=\, \hdge{\tup{q}}\cdot\hdge{\tup{q}} 
\\[4pt]
   \pd_{\tup{q}}
      (\hdge{\tup{\slangup}}\cdot\tup{q}) \,=\, \hdge{\tup{\slangup}} - \hdge{\tup{q}}\cdot\hdge{\tup{p}} 
      &\simeq\, \hdge{\tup{\slangup}} - \tup{p}\otms\tup{q}
\\[4pt]
     \pd_{\tup{p}} (\hdge{\tup{\slangup}}\cdot\tup{p}) \,=\, \hdge{\tup{\slangup}} +  \hdge{\tup{p}}\cdot\hdge{\tup{q}}
     &\simeq\, \hdge{\tup{\slangup}} + \tup{q}\otms\tup{p}
\end{array}
\\[6pt] \nonumber 
&\begin{array}{llllll}
 \pd_{\tup{q}}
     \slang \,=\, \hdge{\htup{\slangup}}\cdot\tup{p}
     &\simeq\, \slang \tup{q}
\\[4pt]
      \pd_{\tup{p}} \slang \,=\, -\hdge{\htup{\slangup}}\cdot\tup{q} 
      &\simeq\, \htup{p}
\end{array}
&&,&& 
\begin{array}{llllll}
      \pd_{\tup{q}} (\hdge{\htup{\slangup}}\cdot\tup{p}) 
      \,=\, \tfrac{p^2}{\slang}\htup{\slangup}\otms\htup{\slangup}
       &\simeq\,
       \slang  \htup{\slangup}\otms\htup{\slangup}
        &\simeq\, 
      \slang ( \imat_3  - \htup{q}\otms \htup{q}  - \htup{p}\otms\htup{p} )
\\[4pt] 
         \pd_{\tup{p}} (\hdge{\htup{\slangup}}\cdot\tup{q}) 
         \,=\, - \tfrac{q^2}{\slang}\htup{\slangup}\otms\htup{\slangup}
    &\simeq\,  - \tfrac{1}{\slang} \htup{\slangup} \otms \htup{\slangup}
     &\simeq\,   
       - \tfrac{1}{\slang} ( \imat_3 - \htup{q}\otms\htup{q} - \htup{p}\otms \htup{p} )
\\[4pt]
  \pd_{\tup{q}}
      (\hdge{\htup{\slangup}}\cdot\tup{q}) 
      \,=\, \fnsize{Eq.\eqref{angmoment_rels_crd2}}
       &\simeq\,  \hdge{\htup{\slangup}}
\\[4pt]
     \pd_{\tup{p}} (\hdge{\htup{\slangup}}\cdot\tup{p}) 
     \,=\, \fnsize{Eq.\eqref{angmoment_rels_crd2}}
     &\simeq\, \hdge{\htup{\slangup}}
 \end{array}
\end{align}
\end{small}
(the last two relations are shown in the footnote\footnote{Using the relations in Eq.\eqref{qpl_rels_apx}: 
\begin{align} \nonumber 
\begin{array}{llll}
     \pd_{\tup{q}}
      (\hdge{\htup{\slangup}}\cdot\tup{q}) \,=\, \tfrac{1}{\slang} \big( \hdge{\tup{\slangup}} - \hdge{\tup{q}}\cdot\hdge{\tup{p}} -   (\hdge{\htup{\slangup}}\cdot\tup{q})\otms (\hdge{\htup{\slangup}}\cdot\tup{p}) \big)
      &\simeq\, 
      \hdge{\htup{\slangup}} - \hdge{\tup{q}}\cdot\hdge{\htup{p}} + \htup{p}\otms \tup{q}
       &=\,  \hdge{\htup{\slangup}} + (\tup{q}\cdot\htup{p})\imat_3 
       &\simeq\,  \hdge{\htup{\slangup}}
\\[4pt]
     \pd_{\tup{p}} (\hdge{\htup{\slangup}}\cdot\tup{p}) \,=\, \tfrac{1}{\slang} \big( \hdge{\tup{\slangup}} +  \hdge{\tup{p}}\cdot\hdge{\tup{q}} +   (\hdge{\htup{\slangup}}\cdot\tup{p})\otms (\hdge{\htup{\slangup}}\cdot\tup{q}) \big)
     &\simeq\, 
     \hdge{\htup{\slangup}} +  \hdge{\htup{p}}\cdot\hdge{\tup{q}}
     - \tup{q}\otms\htup{p}
     &=\,  \hdge{\htup{\slangup}}  - (\tup{q}\cdot\htup{p})\imat_3
     &\simeq\, \hdge{\htup{\slangup}} 
\end{array}
\end{align} }).
The above relations denoted with ``\eq{\simeq}'' follow from using \eq{\lambda=\htup{q}\cdot\tup{p}=0} and \eq{q=1} (and thus \eq{\slang^2\simeq p^2}). \textit{These relations are specific to projective coordinates defined using \eq{m=-1} in Eq.\eqref{prj_family}.} 
We further note that these relations denoted with ``\eq{\simeq}'' are valid specifically along solution curves of the projective coordinate Hamiltonian system that start with initial conditions satisfying \eq{q_\zr=1} and \eq{\lambda_\zr =0} (we are always free to limit consideration to such initial conditions).

%% file: Mysecs_prj/apx_nonCon_coord.tex

\section{Some aspects of analytical Hamiltonian mechanics} \label{app:Ham_mech_crd}

\begin{notation}
  In the following, \eq{\tup{r}=(r_1,\dots,r_\ii{N})\in\mbb{R}^{\en}} denotes some inertial cartesian coordinates, and \eq{(\tup{x},\tup{\piup})\in\mbb{R}^{2n}} and \eq{(\tup{q},\tup{p})\in\mbb{R}^{2n}} denote some generalized configuration and momenta coordinates, where \eq{n\leq \en} is the degrees of freedom (i.e., dimension of configuration space).
  For everything in this Appendix, \eq{(\tup{x},\tup{\piup})} and \eq{(\tup{q},\tup{p})} are arbitrary with no particular relation to the projective coordinates.  
\end{notation}

\subsection{Hamilton's equations with nonconservative forces} \label{app:non_conserv}

We will outline the process — within the classic analytical dynamics framework — for including nonconservative forces in Hamilton's canonical equations of motion, and for transforming such forces between different canonical/symplectic coordinate sets.

\paragraph*{Generalized Forces for Arbitrary Canonical Coordinates, $(\tup{x},\tup{\piup})$.} 
 Suppose some nonconservative forces, \eq{\vecbs{a}^\nc}, acts on the particle, or system of particles, of interest. In this work, \eq{\vecbs{a}^\nc} is usually the force per unit reduced mass.
For some arbitrary canonical/symplectic coordinates, \eq{(\tup{x},\tup{\piup})}, Hamilton's canonical equations of motion including \eq{\vecbs{a}^\nc}
 are given by
\begin{small}
\begin{align} \label{dQP_nc}
\begin{array}{lll}
     & \dot{\tup{x}} \,=\, \pderiv{\mscr{K}}{\tup{\piup}} \,-\, \tup{A}_{x}  &\quad,
     \\[5pt]
     & \dot{\tup{\piup}} \,=\, -\pderiv{\mscr{K}}{\tup{x}} \,+\, \tup{A}_{\pi}
     &\quad, 
\end{array} 
\qquad\qquad
\begin{array}{ll}
     &  \tup{A}_{x} \,:=\, \trn{\pderiv{\tup{r}}{\tup{\piup}}} \tup{a}^\nc 
     \\[5pt]
     & \tup{A}_{\pi} \,:=\,  \trn{\pderiv{\tup{r}}{\tup{x}}} \tup{a}^\nc
\end{array}
\end{align}
\end{small}
where \eq{\tup{r},\tup{a}^\nc \in \mbb{R}^3} are the cartesian component tuples of the vectors \eq{\vecbs{r}} and \eq{\vecbs{a}^\nc} in an inertial basis. 
The governing equation for any function (possibly time-dependent) on phase space, \eq{f(\tup{x},\tup{\piup},t)}, may be obtained from the Hamiltonian \eq{\mscr{K}}  using the Poisson bracket as
\begin{small}
\begin{align}\label{eom_pbrack1}
\begin{array}{rlll}
     \dot{f} 
     \,=\, \pbrak{f}{\mscr{K}} 
     \,+\, \pderiv{f}{\tup{\piup}}\cdot\tup{A}_{\pi} \,-\, \pderiv{f}{\tup{x}}\cdot\tup{A}_{x} \,+\, \pderiv{f}{t} 
\end{array}
\end{align}
\end{small}
Note that it is often the case that \eq{\tup{x}} are
true configuration-level coordinates (e.g., cartesian position coordinates, spherical coordinates, Euler angles, etc.) and \eq{\tup{\piup}} are velocity/momentum-level coordinates such that \eq{\tup{r}=\tup{r}(\tup{x},t)} is not a function of \eq{\tup{\piup}}. When this is the case, then \eq{\tup{A}_{x}=0} and the dynamics simplify as follows 
(where we now write \eq{\tup{A}=\tup{A}_{\pi}}):\footnote{Note in the special case that \eq{(\tup{x},\tup{\piup})} are simply inertial cartesian position and velocity/momenta coordinates (if  \eq{\tup{x}=\tup{r}}), then \eq{\tup{A}=\tup{a}^\nc} are simply the inertial cartesian components the nonconservative forces in that same basis.}
\begin{small}
\begin{align} \label{dQP_nc_2}
\begin{array}{lll}
       \dot{\tup{x}} \,=\, \pderiv{\mscr{K}}{\tup{\piup}} 
    \qquad,\qquad 
      \dot{\tup{\piup}} \,=\, -\pderiv{\mscr{K}}{\tup{x}} \,+\, \tup{A}
\end{array} 
\end{align}
\end{small}
Although various derivations of Eq.\eqref{dQP_nc} and Eq.\eqref{dQP_nc_2} can be found in many texts on analytical mechanics \cite{lanczos2012variational,schaub2003analytical}, we briefly summarize the details below.

\begin{small}
\begin{itemize}
    \item[] \textit{Derivation of Eq.\eqref{dQP_nc}.}  
We may account for the effects of some nonconservative force on the equations of motion by adding the work due to this force, denoted \eq{W^\nc},  to the action such that Hamilton's principle leads to
\begin{small}
\begin{align} \label{dIQP}
\begin{array}{rlllll}
   0\,=\,\delta I & \,=\, \delta\int_{t_\iio}^{t_f} (\tup{\piup}\cdot\dot{\tup{x}} - \mscr{K} + W^\nc)\,\mrm{d}t  
   \,=\, \int_{t_\iio}^{t_f} ( \dot{\tup{x}}\cdot\delta\tup{\piup} \,+\, \tup{\piup}\cdot\delta\dot{\tup{x}}
    \,-\, \delta\mscr{K} 
    \,+\, \delta W^\nc)\,\mrm{d} t 
\\[5pt] 
  & \,=\,
    \int_{t_\iio}^{t_f} ( \dot{\tup{x}}\cdot\delta\tup{\piup} \,-\, \dot{\tup{\piup}}\cdot\delta\tup{x}
    \,-\, \delta\mscr{K} 
    \,+\, \delta W^\nc)\,\mrm{d} t 
\\[5pt] 
    &\,=\,
    \int_{t_\iio}^{t_f} \big( (\dot{\tup{x}}-\pderiv{\mscr{K}}{\tup{\piup}})\cdot\delta\tup{\piup} \,+\,  (-\dot{\tup{\piup}}-\pderiv{\mscr{K}}{\tup{x}})\cdot\delta\tup{x} \,+\, \tup{a}^\nc\cdot\delta\tup{r} \big)\mrm{d} t \,=\, 0 
 \end{array}
\end{align}
\end{small}
where we have used \eq{\delta W^\nc=\vecbs{a}^\nc\cdot\delta \vecbs{r}} as well as fact that the variations vanish at the boundaries:  \eq{\eval{\delta\tup{x}}{t_\iio} \!\!=\eval{\delta\tup{\piup}}{t_\iio}\!\!=\tup{0}}, \eq{\delta t_\iio=0}, and likewise at \eq{t_f}. 
Integration by parts was used to replace \eq{\tup{\piup}\cdot\delta\dot{\tup{x}}} with \eq{-\dot{\tup{\piup}}\cdot\delta\tup{x}} inside the integral.
Next, it is often assumed that \eq{\tup{x}} are true configuration/position-level coordinates and that the momenta, \eq{\tup{\piup}}, are velocity-level coordinates such that \eq{\tup{r}=\tup{r}(\tup{x},t)} is a function only of the configuration coordinates and, possibly, time.
Yet, this is not the most general case as a canonical transformation can render this assumption inaccurate.\footnote{Take the Delaunay variables for example. These are canonical coordinates for orbital motion for which \eq{\vecbs{r}} is a function of both the configuration and momentum coordinates.}  
Thus, taking the most general case that \eq{\tup{r}=\tup{r}(\tup{x},\tup{\piup},t)}, we find that \eq{\delta W^\nc} is given by 
\begin{small}
\begin{align} \label{dW_QP}
    \delta W^\nc \,=\, \tup{a}^\nc\cdot\delta\tup{r} \,=\,
    \tup{a}^\nc\cdot(\pderiv{\tup{r}}{\tup{x}}\delta \tup{x} \,+\, \pderiv{\tup{r}}{\tup{\piup}}\delta\tup{\piup})
    \,=\,
    (\trn{\pderiv{\tup{r}}{\tup{x}}} \tup{a}^\nc)\cdot\delta\tup{x}
    \,+\, 
    ( \trn{\pderiv{\tup{r}}{\tup{\piup}}} \tup{a}^\nc)\cdot\delta\tup{\piup}
\end{align}
\end{small}
Such that Eq.\eqref{dIQP} becomes
\begin{small}
\begin{align} \label{dIQP2}
    \delta I \,=\, 
    \int_{t_\iio}^{t_f} \Big( (\dot{\tup{x}} - \pderiv{\mathscr{K}}{\tup{\piup}}+\trn{\pderiv{\tup{r}}{\tup{\piup}}} \tup{a}^\nc )\cdot\delta\tup{\piup} \,-\,
    (\dot{\tup{\piup}} + \pderiv{\mathscr{K}}{\tup{x}} - \trn{\pderiv{\tup{r}}{\tup{x}}} \tup{a}^\nc ) \cdot\delta\tup{x} 
    \Big)\mrm{d} t = 0
\end{align}
\end{small}
for the above to hold for all arbitrary \eq{\delta\tup{x}} and \eq{\delta\tup{\piup}}  requires that their coefficients be zero. This leads to the canonical equations of motion given by Eq.\eqref{dQP_nc}. 
\end{itemize}
\end{small}


\paragraph*{Transformation of Generalized Forces.} 
Now suppose that we perform some general canonical transformation from \eq{(\tup{x},\tup{\piup})} to new canonical coordinates, \eq{(\tup{q},\tup{\pf})}.\footnote{The developments of this section apply the same weather \eq{\tup{q}} and \eq{\tup{\pf}} are  minimal or non-minimal canonical coordinates.} 
Since \eq{\tup{q}} and \eq{\tup{\pf}} are also canonical coordinates, the derivation of the equations of motion including nonconservative forces exactly parallels Eq.\eqref{dQP_nc}-Eq.\eqref{dIQP2}. In general, we assume \eq{\tup{r}=\tup{r}(\tup{q},\tup{\pf},t)} such that \eq{\delta W^\nc} is now given by Eq.\eqref{dW_QP} with \eq{\tup{q}} and \eq{\tup{\pf}} taking the place of \eq{\tup{x}} and \eq{\tup{\piup}}, respectively. 
 The requirement that \eq{\delta I =0} then leads to canonical equations of motion which are precisely the same \textit{form} as Eq.\eqref{dQP_nc}:
\begin{small}
\begin{align} \label{dqp_nc}
    \begin{array}{lll}
     & \dot{\tup{q}} \,=\, \pderiv{\mscr{H}}{\tup{q}} \,-\, \tup{\alphaup}_{q}  &\quad,
     \\[5pt]
     & \dot{\tup{\pf}} \,=\, -\pderiv{\mscr{H}}{\tup{q}} \,+\, \tup{\alphaup}_{\pf}
     &\quad,
\end{array} 
\qquad\qquad
\begin{array}{ll}
     &  \tup{\alphaup}_{q} \,:=\, \trn{\pderiv{\tup{r}}{\tup{\pf}}} \tup{a}^\nc 
     \\[5pt]
     & \tup{\alphaup}_{\pf} \,:=\,  \trn{\pderiv{\tup{r}}{\tup{q}}} \tup{a}^\nc
\end{array}
\end{align}
\end{small}
where \eq{\mscr{H}} is the Hamiltonian for the \eq{(\tup{q},\tup{\pf})} phase space. 
What we would like to know is, given some canonical transformation between \eq{(\tup{x},\tup{\piup})} and  \eq{(\tup{q},\tup{\pf})}, what then is the relation between the above  \eq{\tup{\alphaup}} and the \eq{\tup{A}}  defined in Eq.\eqref{dQP_nc}? We will present the most general case as well a more specific, yet common, case. 
\begin{small}
\begin{itemize}
\item \textit{Case 1.} Consider the most general case:
    \begin{small}
    \begin{itemize}[nosep]
        \item \eq{(\tup{x},\tup{\piup})} are arbitrary canonical/symplectic coordinates such that, in general, \eq{\tup{r}=\tup{r}(\tup{x},\tup{\piup},t)}. 
        \item  \eq{(\tup{x},\tup{\piup})} and  \eq{(\tup{q},\tup{\pf})} are related by a general canonical transformation such that, in general,  \eq{\tup{x}=\tup{x}(\tup{q},\tup{\pf},t)} and \eq{\tup{\piup}=\tup{\piup}(\tup{q},\tup{\pf},t)}.
    \end{itemize}
    \end{small}
Then the relation between \eq{\tup{A}} and \eq{\tup{\alphaup}} can be found from Eq.\eqref{dqp_nc} as
\begin{small}
\begin{align} \label{genForce_CT}
\begin{array}{lll}
       \tup{\alphaup}_{q}\,:=\, \trn{\pderiv{\tup{r}}{\tup{\pf}}} \tup{a}^\nc
     \,=\,  \trn{\big( \pderiv{\tup{r}}{\tup{x}} \pderiv{\tup{x}}{\tup{\pf}} + \pderiv{\tup{r}}{\tup{\piup}} \pderiv{\tup{\piup}}{\tup{\pf}} \big)} \tup{a}^\nc
      &\,=\,
      \trn{\pderiv{\tup{x}}{\tup{\pf}}} \trn{\pderiv{\tup{r}}{\tup{x}}} \tup{a}^\nc  + \trn{\pderiv{\tup{\piup}}{\tup{\pf}}}  \trn{\pderiv{\tup{r}}{\tup{\piup}}} \tup{a}^\nc  
    &\,=\,
     \trn{\pderiv{\tup{x}}{\tup{\pf}}} \tup{A}_{\pi}  + \trn{\pderiv{\tup{\piup}}{\tup{\pf}}}  \tup{A}_{x}
  \\[10pt]
    \tup{\alphaup}_{\pf}\,:=\, \trn{\pderiv{\tup{r}}{\tup{q}}} \tup{a}^\nc
      \,=\, \trn{\big( \pderiv{\tup{r}}{\tup{x}} \pderiv{\tup{x}}{\tup{q}} + \pderiv{\tup{r}}{\tup{\piup}} \pderiv{\tup{\piup}}{\tup{q}} \big)} \tup{a}^\nc
      &\,=\,
        \trn{\pderiv{\tup{x}}{\tup{q}}} \trn{\pderiv{\tup{r}}{\tup{x}}} \tup{a}^\nc  + \trn{\pderiv{\tup{\piup}}{\tup{q}}}  \trn{\pderiv{\tup{r}}{\tup{\piup}}} \tup{a}^\nc  
       &\,=\,
       \trn{\pderiv{\tup{x}}{\tup{q}}} \tup{A}_{\pi}  + \trn{\pderiv{\tup{\piup}}{\tup{q}}}  \tup{A}_{x} 
\end{array}
\end{align}
\end{small}
where we have used the definitions of \eq{\tup{A}_{x}} and  \eq{\tup{A}_{\pi}} from Eq.\eqref{dQP_nc}. 
Thus, once one knows the relation between the phase space coordinates, one may then use the above to calculate the relation between the generalized forces. 
\item \textit{Case 2.}  Now consider a more specific, common, case:
    \begin{small}
    \begin{itemize}[nosep]
        \item \eq{(\tup{x},\tup{\piup})} are configuration-level and velocity-level coordinates, respectively, such that \eq{\tup{r}=\tup{r}(\tup{x},t)} is \textit{not} a function of \eq{\tup{\piup}}
        \item  \eq{(\tup{x},\tup{\piup})} and  \eq{(\tup{q},\tup{\pf})} are related by a time-dependent \textit{point} transformation such that \eq{\tup{x}=\tup{x}(\tup{q},t)} and \eq{\tup{\piup}=\tup{\piup}(\tup{q},\tup{\pf},t)}.
    \end{itemize}
    \end{small}
The first of the above leads to \eq{\tup{A}_Q=\tup{0}} such that the equations of motion for \eq{\tup{x}} and \eq{\tup{\piup}} are given by Eq.\eqref{dQP_nc_2}. 
The second condition — that the two phase spaces are related by a canonical \textit{point} transformation —  means that \eq{\pderiv{\tup{x}}{\tup{\pf}}=0}. 
Eq.\eqref{genForce_CT} then simplifies to
\begin{small}
\begin{align} \label{dqp_nc_pt}
 \begin{array}{lll}
     & \dot{\tup{q}} \,=\, \pderiv{\mscr{H}}{\tup{q}}   
     \\[5pt]
     & \dot{\tup{\pf}} \,=\, -\pderiv{\mscr{H}}{\tup{q}} \,+\, \tup{\alphaup}
\end{array} 
\quad, \qquad\qquad
\begin{array}{ll}
     & \tup{\alphaup}_{q} \,=\,  \tup{0}
     \\[5pt]
     & \tup{\alphaup} \,\equiv\, \tup{\alphaup}_{\pf}  \,=\,  \trn{\pderiv{\tup{r}}{\tup{q}}} \tup{a}^\nc \,=\, \trn{\pderiv{\tup{x}}{\tup{q}}} \tup{A}
\end{array} 
\end{align}
\end{small}
\end{itemize}
\end{small}

%% file: Mysecs_prj/apx_Pext_coord.tex
\subsection{Extended phase space coordinates} \label{sec:ext_coord}

\begin{small}
\begin{notation}
    In the following, 
    a ``hat'' \eq{\,\wh{\sblt}\,} does \textit{not} denote a normalized unit vector or unit tuple but will, instead, denote objects associated with extended space.
\end{notation}
\end{small}

\noindent Suppose we wish to use some evolution parameter, \eq{\varep}, that is related to the time through some specified differential relation:
\begin{small}
\begin{align}
    \diff{t}{\varep} \,=\, \rng{t} 
    \qquad, \qquad \mrm{d} t \,=\, \rng{t} \mrm{d} \varep
\end{align}
\end{small}
where \eq{\rng{(\,)}} will denote differentiation with respect to \eq{\varep} and where, in general, \eq{\rng{t}} may be some function of the generalized coordinates, \eq{\tup{q}}, the generalized velocities, \eq{\dot{\tup{q}}}, and time, \eq{t}. Now, the classic action integral written in terms of the Lagrangian is given by \eq{I=\int_{t_\iio}^{t_f}\mscr{L}\mrm{d} t}. Using the above differential relation, this same integral with \eq{\varep} as the evolution parameter is given by \eq{I=\int_{\varep_\iio}^{\varep_f}\mscr{L}\rng{t}\mrm{d} \varep}.  We then define extended generalize coordinate and velocity vectors, and the \textit{extended Lagrangian}, \eq{\whscr{L}}, as
\begin{small}
\begin{align}
  \whtup{q}:=(\tup{q},t) \qquad,\qquad \rng{\whtup{q}}:=(\rng{\tup{q}},\rng{t})
  \qquad,\qquad
  \whscr{L}(\whtup{q},\rng{\whtup{q}}) 
  \,:=\, \rng{t}\mscr{L}(\tup{q},\dot{\tup{q}},t)
\end{align}
\end{small}
Where the configuration space has increased by one dimension to include the time as another generalized coordinate with generalized velocity \eq{\rng{t}}.

Now, the conjugate momenta coordinate are classically defined as \eq{\tup{p}=\pderiv{\mscr{L}}{\dot{\tup{q}}}}. In an analogous manner, we now define the momenta conjugate to the time as
\begin{small}
\begin{align} \label{pt_deff}
     p_t \,:=\, \pderiv{\whscr{L}}{\rng{t}} \,=\, \mscr{L} \,-\, \tup{p} \cdot \dot{\tup{q}} \,=\, -\mscr{H}
\end{align}
\end{small}
where the equality \eq{\pderiv{\whscr{L}}{\rng{t}} = \mscr{L} - \tup{p} \cdot \dot{\tup{q}}} is derived in
the footnote\footnote{From \eq{\whscr{L}:=\rng{t}\mscr{L}}, then  \eq{\pderiv{\whscr{L}}{\rng{t}}  = \mscr{L} + \rng{t}\pderiv{\mscr{L}}{\rng{t}}}. Noting that \eq{\mscr{L}(\tup{q},\dot{\tup{q}},t)} only depends on \eq{\rng{t}} through the relation \eq{\dot{\tup{q}}=\tfrac{1}{\rng{t}}\rng{\tup{q}}}, it follows that \eq{\pderiv{\mscr{L}}{\rng{t}} =\pderiv{\mscr{L}}{\dot{\tup{q}}}\cdot \pderiv{\dot{\tup{q}}}{\rng{t}} = - \tup{p} \cdot \tfrac{\rng{\tup{q}}}{\rng{t}^2} = -\tfrac{1}{\rng{t}} \tup{p} \cdot \dot{\tup{q}} }. Thus, \eq{\pderiv{\whscr{L}}{\rng{t}} = \mscr{L} + \rng{t}\pderiv{\mscr{L}}{\rng{t}} = \mscr{L} - \tup{p} \cdot \dot{\tup{q}}}, which is precisely the classical definition of \eq{-\mscr{H}}. }.
From the above, we see that \eq{p_t} is equal to the \textit{negative} of the classic Hamiltonian, \eq{\mscr{H}}, where \eq{\mscr{H}} itself is often (but not always) equal to the total energy.  
We then define the \textit{extended Hamiltonian}, \eq{\whscr{H}(\whtup{q},\whtup{p})} — where \eq{(\whtup{q},\whtup{p})} are extended coordinate and momenta vectors for the extended phase space — using the a Legendre transformation in the usual manner:
\begin{small}
\begin{align} \label{Hext_app}
    \whtup{q} :=\, (\tup{q},t) 
\qquad,\qquad
    \whtup{p} :=\, (\tup{p},p_t) 
\qquad,\qquad
  \whscr{H} :=\, \whtup{p}\cdot \rng{\whtup{q}}  \,-\, \whscr{L} \,=\, \rng{t}(\mscr{H}\,+\, p_t)
    \,=\, 0
\end{align}
\end{small}
where \eq{\tup{q}} and \eq{\tup{p}} are the usual coordinates and momenta for the non-extended phase space. 
 Note from the above that the value of \eq{\whscr{H}} is always equal to zero on account of  Eq.\eqref{pt_deff}. Then — as we will soon show —  the canonical equations of motion with  \eq{\varep} as the evolution parameter are obtained from the extended Hamiltonian as
 \begin{small}
 \begin{align} \label{q'p'_app}
 \begin{array}{lllll}
      \rng{\tup{q}} \,=\, \pderiv{\whscr{H}}{\tup{p}} \,-\, \rng{t}\tup{\alphaup}_{q}  
  \\[5pt]
   \rng{t} \,=\, \pderiv{\whscr{H}}{p_t} 
\end{array}
\qquad,\qquad 
\begin{array}{llll}
   \rng{\tup{p}} \,=\, -\pderiv{\whscr{H}}{\tup{q}} \,+\, \rng{t}\tup{\alphaup}_{p}
 \\[5pt]
      \rng{p}_t \,=\, -\pderiv{\whscr{H}}{t} \,+\, \rng{t}\alpha_t
\end{array} 
\end{align}
\end{small}
where \eq{\tup{\alphaup}} are the generalized forces defined in section \ref{app:non_conserv} and where we have also defined an additional generalized ``force'', \eq{\alpha_t} (which really has units of power). If \eq{\tup{a}^\nc\in\mbb{R}^\ii{N}} denotes the components of the nonconservative force vector in the inertial \eq{\tup{r}}-basis, then these generalized forces for arbitrary \eq{\tup{q}} and \eq{\tup{p}} are given by:
\begin{small}\begin{align} \label{gt_qp_def}
     \tup{\alphaup}_{q} = \trn{\pderiv{\tup{r}}{\tup{p}}}  \cdot  \tup{a}^\nc 
&&,&&
    \tup{\alphaup}_{p} =  \trn{\pderiv{\tup{r}}{\tup{q}}} \cdot  \tup{a}^\nc
&&,&&
    \alpha_t \,:=\,   -( \tup{\alphaup}_p \cdot  \dot{\tup{q}}   + \tup{\alphaup}_q \cdot  \dot{\tup{p}})
  \,=\, -\tup{a}^\nc\cdot ( \pderiv{\tup{r}}{\tup{q}}\cdot  \pderiv{\mscr{H}}{\tup{p}} - \pderiv{\tup{r}}{\tup{p}}\cdot  \pderiv{\mscr{H}}{\tup{q}} ) \,=\, -\tup{a}^\nc\cdot  (\dot{\tup{r}} - \pd_t \tup{r} )
\end{align}
\end{small}
For the case that \eq{\tup{q}} and \eq{\tup{p}} are true configuration-level and momentum-level variables, respectively — such that \eq{\tup{r}=\tup{r}(\tup{q},t)} is \textit{not} a function of the momenta — then the canonical equations of motion in the extended phase space are given by
\begin{small}
\begin{align} \label{qp'_coord}
\begin{array}{lllll}
      \rng{\tup{q}} \,=\, \pderiv{\whscr{H}}{\tup{p}} 
  \\[5pt]
   \rng{t} \,=\, \pderiv{\whscr{H}}{p_t} 
\end{array}
\qquad,\qquad 
\begin{array}{llll}
   \rng{\tup{p}} \,=\, -\pderiv{\whscr{H}}{\tup{q}} \,+\, \rng{t}\tup{\alphaup}
 \\[5pt]
      \rng{p}_t \,=\, -\pderiv{\whscr{H}}{t} \,+\, \rng{t}\alpha_t
\end{array} 
&&
 \fnsize{where:}
 \quad 
\begin{array}{llll}
    \tup{\alphaup} \equiv \tup{\alphaup}_p \,=\, \tup{a}^\nc \cdot \pderiv{\tup{r}}{\tup{q}}
 \\[5pt]
     \alpha_t \,=\, - \tup{\alphaup} \cdot  \dot{\tup{q}} \,=\, - \tup{\alphaup} \cdot  \pderiv{\mscr{H}}{\tup{p}}  
\end{array} 
\end{align}
\end{small}
The aforementioned relation \eq{p_t=-\mscr{H}}, such that \eq{\whscr{H}=0}, may give the incorrect impression that the above partials of \eq{\whscr{H}} vanish. We will explain why this is not the case soon. First, we derive Eq.\eqref{q'p'_app}:
 \begin{small}
\begin{itemize}
    \item[] \textit{Derivation of Eq.\eqref{q'p'_app}.} 
Hamilton's principle, including the work from nonconservative forces, \eq{W^\nc}, is given in the extended phase space by
\begin{small}
\begin{align} \label{dI_ext}
\begin{array}{rl}
   0 \,=\,   \delta I  \,=\, \delta\int_{\varep_\iio}^{\varep_f} (\tup{p}\cdot \dot{\tup{q}} - \mscr{H} + W^\nc )\rng{t} \,\mrm{d} \varep
       &\,=\, \delta\int_{\varep_\iio}^{\varep_f} (\whtup{p}\cdot \rng{\whtup{q}} - \whscr{H} + \rng{t}W^\nc )\,\mrm{d}\varep  
\\[8pt] 
      & \,=\, \int_{\varep_\iio}^{\varep_f} \big(\delta\whtup{p}\cdot \rng{\whtup{q}} - \rng{\whtup{p}} \cdot  \delta\whtup{q} - \delta\whscr{H} + \delta(\rng{t} W^\nc ) \big) \, \mrm{d}\varep 
\end{array}
\end{align}
\end{small}
where we have used integration by parts along with the fact that the variations vanish at the boundaries.  Using the same process, the term involving \eq{W^\nc } leads to
\begin{small}
\begin{align} \label{intyint}
\begin{array}{lll}
     \int_{\varep_\iio}^{\varep_f} \delta(\rng{t}W^\nc) \mrm{d} \varep \,=\, \int_{\varep_\iio}^{\varep_f} (\rng{t}\delta W^\nc - \rng{W}^\nc\delta t) \mrm{d} \varep
      \,=\, \int_{\varep_\iio}^{\varep_f} \rng{t}(\delta W^\nc - \dot{W}^\nc\delta t) \mrm{d} \varep
\end{array}
\end{align}
\end{small}
where we have used \eq{\rng{(\;)}=\rng{t}\dot{(\;)}}. 
Let us take \eq{\tup{q}} and \eq{\tup{p}} to be arbitrary canonical coordinates such that, in general, \eq{\tup{r}=\tup{r}(\tup{q},\tup{p},t)}.  
With  \eq{t} now treated as any other generalized coordinate, \eq{\dot{W}^\nc} and \eq{\delta W^\nc} are now given
by\footnote{Note in  Eq.\eqref{wdot_ext} that \eq{\pd_t \tup{r}} is now included in \eq{\delta W^\nc} but that \eq{\pderiv{\tup{r}}{p_t}=0} .}
\begin{small}
\begin{align}  \label{wdot_ext}
\begin{array}{rlll}
       \dot{W}^\nc   
       \,=\, \tup{a}^\nc \cdot \dot{\tup{r}}  \,=\,  \tup{a}^\nc \cdot \big( \pderiv{\tup{r}}{\whtup{q}}\dot{\whtup{q}}  +  \pderiv{\tup{r}}{\whtup{p}}\dot{\whtup{p}} \big)
       & \,=\,  \tup{a}^\nc \cdot  \big( \pderiv{\tup{r}}{\tup{q}}\dot{\tup{q}} + \pd_t \tup{r}  +  \pderiv{\tup{r}}{\tup{p}}\dot{\tup{p}}  \big)
\\[6pt]
     \delta W^\nc   
     \,=\, \tup{a}^\nc \cdot \delta \tup{r}  \,=\,  \tup{a}^\nc\cdot \big( \pderiv{\tup{r}}{\whtup{q}}\delta\whtup{q}  +  \pderiv{\tup{r}}{\whtup{p}}\delta\whtup{p} \big)
     & \,=\,  \tup{a}^\nc \cdot  \big( \pderiv{\tup{r}}{\tup{q}}\delta\tup{q} + \pd_t \tup{r} \delta t  +  \pderiv{\tup{r}}{\tup{p}}\delta\tup{p} \big)
\end{array}
\end{align}
\end{small}
The term in parenthesis in  Eq.\eqref{intyint} is then given by
\begin{small}
\begin{align} \label{saucy}
\begin{array}{rl}
       \delta W^\nc - \dot{W}^\nc\delta t  \;\,=\,\; \tup{a}^\nc\cdot \delta \tup{r} \,-\, \tup{a}^\nc\cdot \dt{\tup{r}}\delta t 
       &\,=\,
       \tup{a}^\nc\cdot  \big( \pderiv{\tup{r}}{\tup{q}}\delta\tup{q} + \pderiv{\tup{r}}{\tup{p}}\delta\tup{p}
       \,-\, (\pderiv{\tup{r}}{\tup{q}}\dot{\tup{q}}   +  \pderiv{\tup{r}}{\tup{p}}\dot{\tup{p}})\delta t \big)
\\[5pt]
     & \,=\,  \tup{\alphaup}_p\cdot  \delta\tup{q}  \,+\,  \tup{\alphaup}_q\cdot \delta\tup{p}  \,-\, (\tup{\alphaup}_p \cdot  \dot{\tup{q}}   \,+\,  \tup{\alphaup}_q \cdot  \dot{\tup{p}})\delta t
\end{array}
\end{align}
\end{small}
where we have used \eq{\tup{\alphaup}_p:=\tup{a}^\nc\cdot  \pderiv{\tup{r}}{\tup{q}}} and \eq{\tup{\alphaup}_q:=\tup{a}^\nc\cdot  \pderiv{\tup{r}}{\tup{p}}}. 
Substitution of Eq.\eqref{intyint} and Eq.\eqref{saucy} into  Eq.\eqref{dI_ext} then leads to
\begin{small}
\begin{align} \label{dI_ext=0}
    \delta I \,=  \int_{\varep_\iio}^{\varep_f} 
    \Bigg \{ (\dot{\tup{q}}-\pderiv{ \whscr{H}}{\tup{p}} + \rng{t}\tup{\alphaup}_q)\cdot \delta\tup{p}
    \,+\,  (-\dot{\tup{p}}-\pderiv{ \whscr{H}}{\tup{q}} + \rng{t}\tup{\alphaup}_p)\cdot \delta\tup{q} 
    \,+\,  (\rng{t} - \pderiv{ \whscr{H}}{p_t})\delta p_t  \,+\,  \big[-\rng{p}_t-\pderiv{ \whscr{H}}{t} - \rng{t}(\tup{\alphaup}_p \cdot  \dot{\tup{q}}   \,+\,  \tup{\alphaup}_q \cdot  \dot{\tup{p}}) \big]\delta t
    \Bigg\} \,\mrm{d}\varep   \,=\, 0
\end{align}
\end{small}
In order for the above to hold for all arbitrary \eq{\delta \whtup{q}} and \eq{\delta \whtup{p}}, their coefficients must vanish. This leads to the canonical equations of motion with respect to the evolution parameter \eq{\varep} as seen in  Eq.\eqref{q'p'_app} with the generalized forces defined in  Eq.\eqref{gt_qp_def}. 
\end{itemize}
\end{small}

\vspace{1ex}
\noindent Now, we indicated in Eq.\eqref{pt_deff} that \eq{p_t=-\mscr{H}} such that \eq{\whscr{H}=0}. This may give the incorrect impression that the partials of \eq{\whscr{H}} vanish in the extended phase space canonical equations of motion seen in Eq.\eqref{q'p'_app}-Eq.\eqref{qp'_coord}. They do not. The relation \eq{p_t=-\mscr{H}} must be interpreted more carefully as follows: along any solution curve \eq{(\whtup{q}_\varep,\whtup{p}_\varep)=(\tup{q}_\varep,\tup{p}_\varep,t_\varep,p_{t_\varep})} of the extended phase space ODEs, it holds that the \textit{value} of \eq{p_t} along the curve is related to the \textit{value} of \eq{\mscr{H}} along the curve by: 
\begin{small}
\begin{align} \label{pt_H_specific}
    p_{t_\varep} = -\mscr{H}(\tup{q}_\varep,\tup{p}_\varep,t_\varep)
\end{align}
\end{small}
When we write \eq{p_t=-\mscr{H}}, it has the above interpretation. It does \textit{not} mean that \eq{p_t} should be treated as a function of the other other coordinates. It is important to keep this in mind when differentiating the extended Hamiltonian \eq{\whscr{H}} to obtain the ODEs;  \eq{p_t} must be treated as an independent coordinate. 
\textit{However}, after obtaining the equations of motion, we may use the fact that Eq.\eqref{pt_H_specific} holds along any solution curve. That is, the relation \eq{ p_t  = -\mscr{H} } (interpreted as above) may be substituted into the resulting ODEs themselves. This is not necessary, but it is permissible and it leads to the dynamics for \eq{(\tup{q},\tup{p})} in the form \eq{(\rng{\tup{q}},\rng{\tup{p}}) = \rng{t}(\dot{\tup{q}},\dot{\tup{p}})}. 
That is, the extended phase space ODE for \eq{(\tup{q},\tup{p})} in Eq.\eqref{q'p'_app} are found to be equivalent to:
\begin{small}
\begin{align} \label{qpdot_ext_alt_apx}
\begin{array}{llll}
      \rng{\tup{q}} \,=\, \pderiv{\whscr{H}}{\tup{p}} - \rng{t}\bs{\alphaup}_q
      & =\; \rng{t}(\pderiv{\mscr{H}}{\tup{p}} -\bs{\alphaup}_q )
      \;=\; \rng{t}\dot{\tup{q}}
\\[5pt]
     \rng{\tup{p}} \,=\, -\pderiv{\whscr{H}}{\tup{q}} +\rng{t}\tup{\alphaup}_p
     & =\;  \rng{t}(-\pderiv{\mscr{H}}{\tup{q}} + \tup{\alphaup}_p )
     \;=\; \rng{t}\dot{\tup{p}}
\end{array}
\end{align}
\end{small}
where the second equalities follow from \eq{ p_t = -\mscr{H} }, which leads to \eq{\pderiv{\whscr{H}}{\tup{p}} = \rng{t}\pderiv{\mscr{H}}{\tup{p}}} and \eq{\pderiv{\whscr{H}}{\tup{q}} = \rng{t}\pderiv{\mscr{H}}{\tup{q}}}.
To verify these relations, let \eq{z_i\in (q_1,\dots,q_n,p_1,\dots,p_n)} be any phase space coordinate. We then find: 
\begin{small}
\begin{align}
\begin{array}{llll}
    \pderiv{\whscr{H}}{z_i} = \pderiv{}{z_i} \big( \rng{t}(\mscr{H}+p_t) \big) 
    \,=\, \rng{t} \pderiv{\mscr{H}}{z_i} \,+\,  (\mscr{H}+p_t) \pderiv{}{z_i} \rng{t} 
\end{array}
\;, &&  p_t = -\mscr{H}
\quad \Rightarrow \quad 
 \pderiv{\whscr{H}}{z_i} \,=\, \rng{t} \pderiv{\mscr{H}}{z_i}
\end{align}
\end{small}

%% file: Mysecs_prj/old_prjPlots.tex
\section{Numerical verification:~Kepler dynamics with \txi{J2} perturbation} \label{sec:J2_plots_prj}


The following plots serve as numerical verification of the dynamics given in section \ref{sec:2BP} (specifically, section \ref{sec:j2}).
The following are simulations run for the $J_2$-perturbed Kepler problem using the projective coordinates of section \ref{sec:2BP}, including the Earth's $J_2$ gravitational perturbation term, with initial conditions given in terms of classic orbit elements \eq{(a,e,i,\omega,\Omega,f= \tau)}, with \eq{\tau} the true anomaly, by \eq{a_\zr = 8.59767038 \cdot 10^{3} (\mrm{km})}, \eq{e_\zr=0.2}, \eq{i_\zr=20 (\mrm{deg.})}, \eq{\omega_\zr=70 (\mrm{deg.})},  \eq{\Omega_\zr = 135 (\mrm{deg.})}, and true anomaly \eq{f_\zr = \tau_\zr = 0}. The numerical integration is carried out using ode45 in \sc{matlab} with units scaled such that the Earth's radius is \eq{R_\mrm{e}=1} and the Earth's gravitational parameter is \eq{\kconst = 1}.

\begin{figure}[!h]
	\centering
	\includegraphics[scale=0.45]{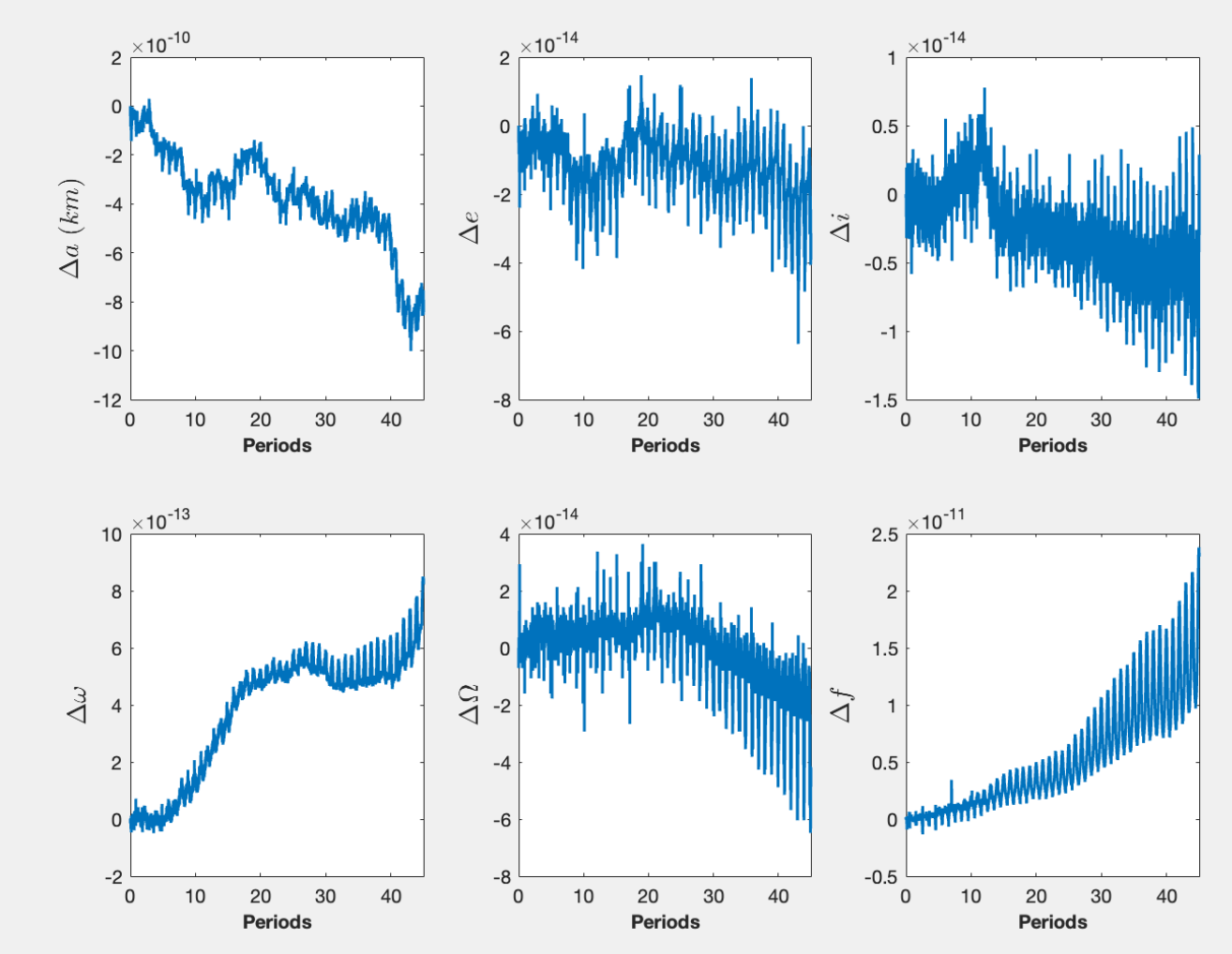}
	\caption{Verifying accuracy of the canonical equations of motion for $(\tup{q},u,\tup{p},p_u)$ given in Eq.\eqref{s_J2_2}, showing errors in the classic orbit elements (COEs), as compared to propagating the modified equinoctial elements (MEEs). ($J_2$ gravitational term included.)}

\vspace{1cm}

	\centering
	\includegraphics[scale=0.45]{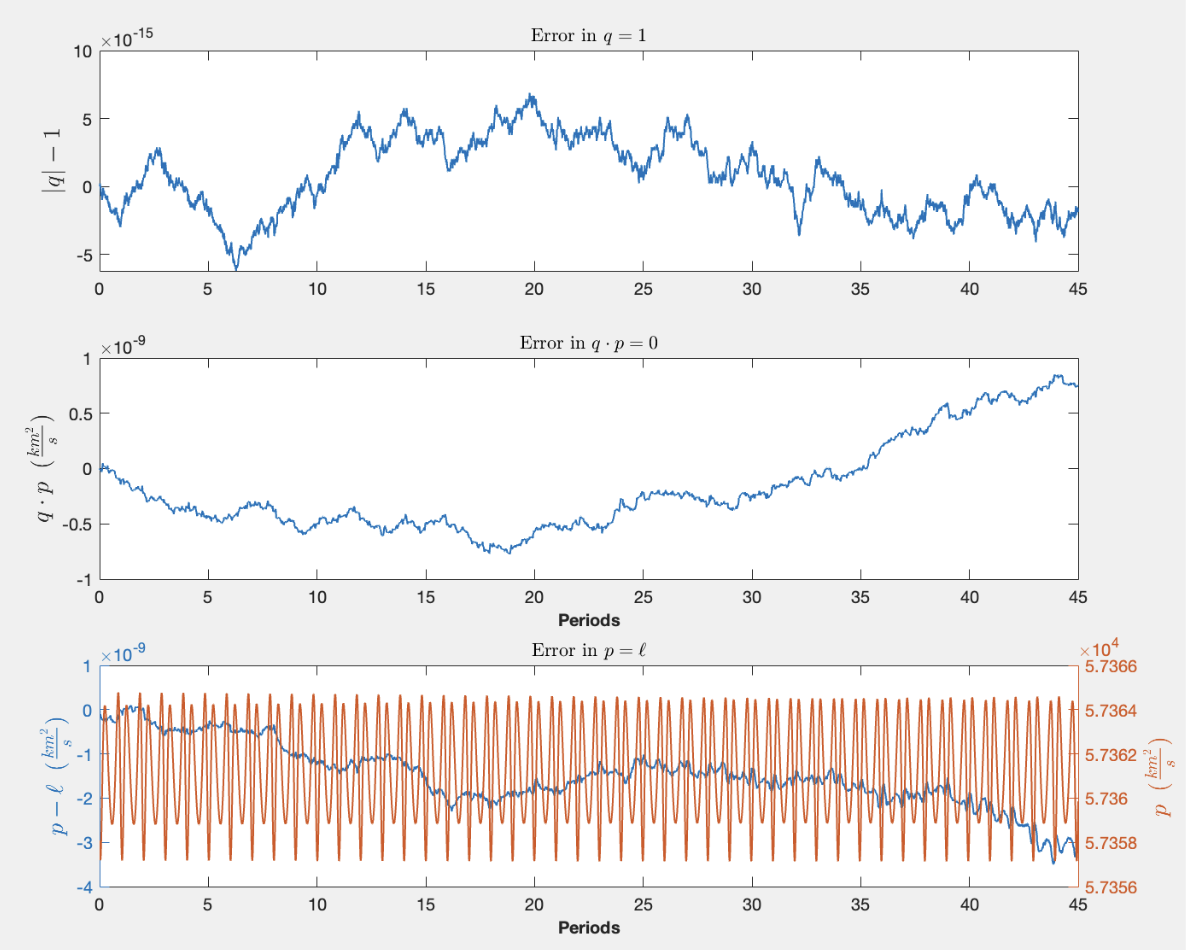}
	\caption{Verification that $q=\mag{\tup{q}}=1$ and \eq{\tup{q}\cdot\tup{p}=0} are integrals of motion and that \eq{p=\mag{\tup{p}}=\slang} is equal to the (specific) angular momentum magnitude. ($J_2$ gravitational term included.)}
\end{figure}


\begin{figure}[!h]
	\centering
	\includegraphics[scale=0.5]{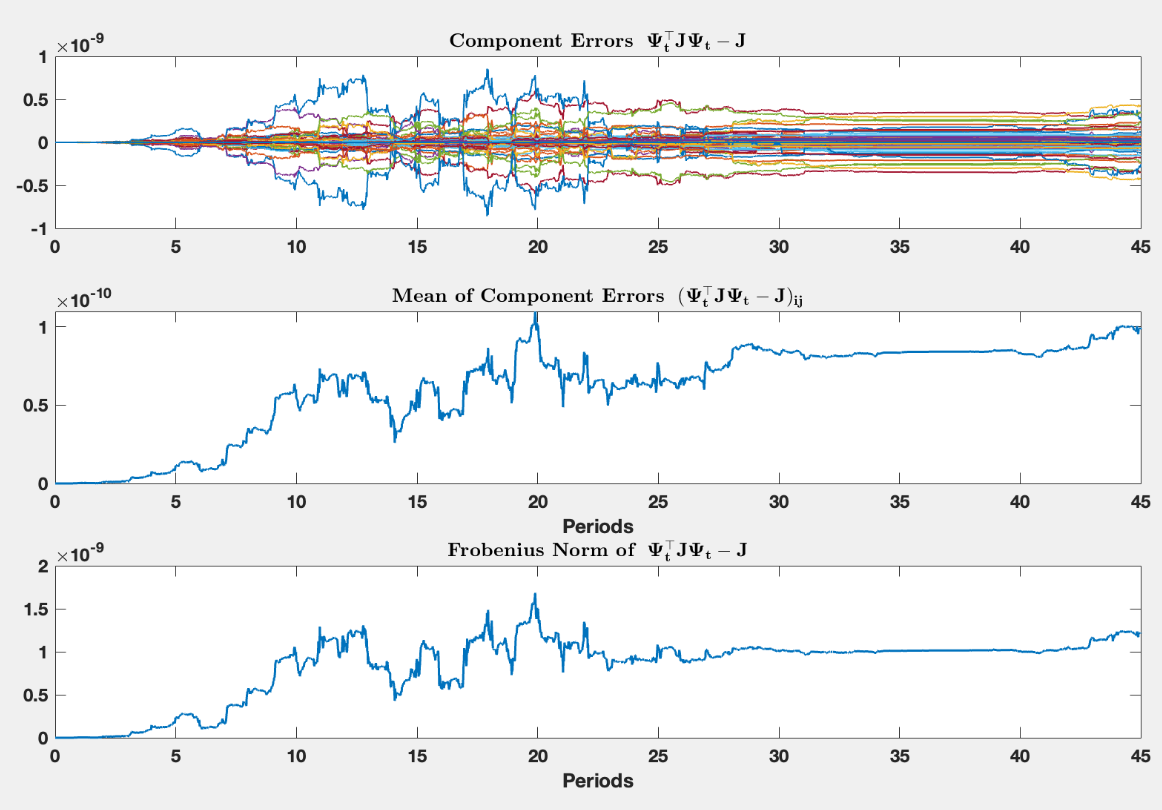}
	\caption{Verification that the STM for $(\tup{q},u,\tup{p},p_u)$, with \eq{t} as the evolution parameter, is a symplectic matrix. ($J_2$ gravitational term included.)  
    Results obtained by propagating the general ODE for the STM from Eq.\eqref{STM_ODE}.  }

\vspace{1cm}

\centering    
\begin{subfigure}[b]{0.48\textwidth}
\centering    
    \includegraphics[width = \textwidth]{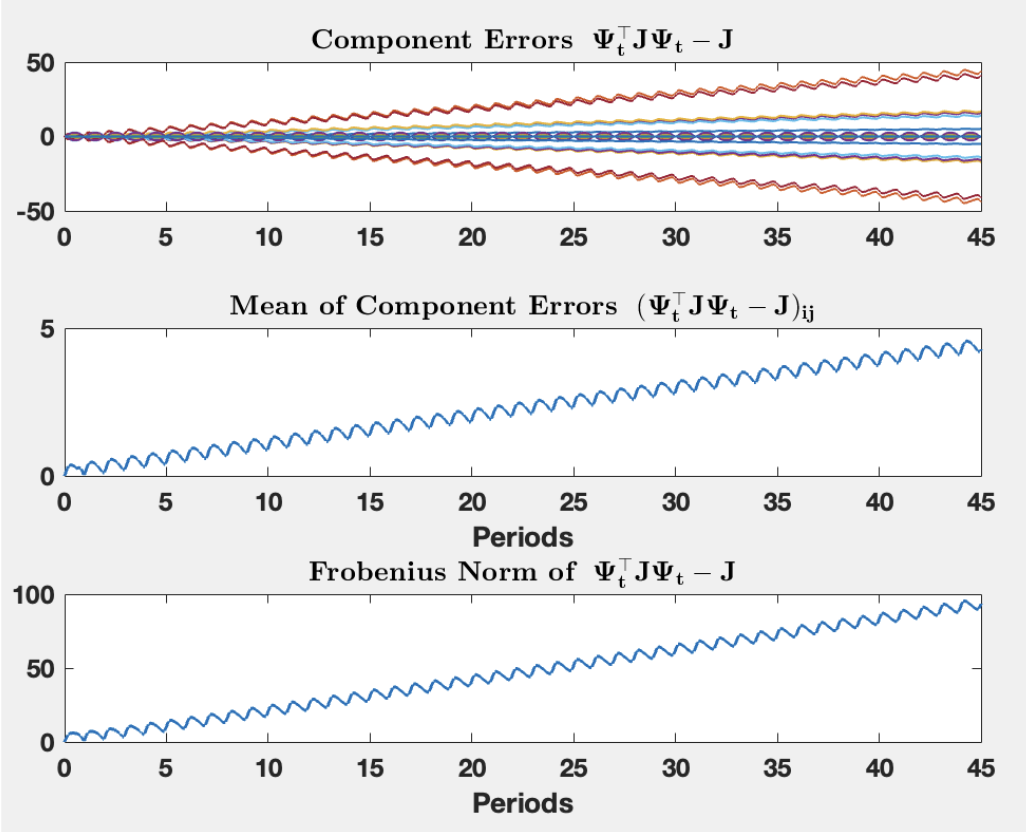}
   \caption{Verification that the STM for $(\tup{q},u,\tup{p},p_u)$, with \eq{s} as the evolution parameter, is \textit{not} a symplectic matrix ($J_2$ gravitational term included.)
   Results obtained by propagating the general ODE for the STM from Eq.\eqref{STM_ODE}. }
    \end{subfigure}
    \hfill
    \begin{subfigure}[b]{0.48\textwidth}
    \centering
    \includegraphics[width =\textwidth]{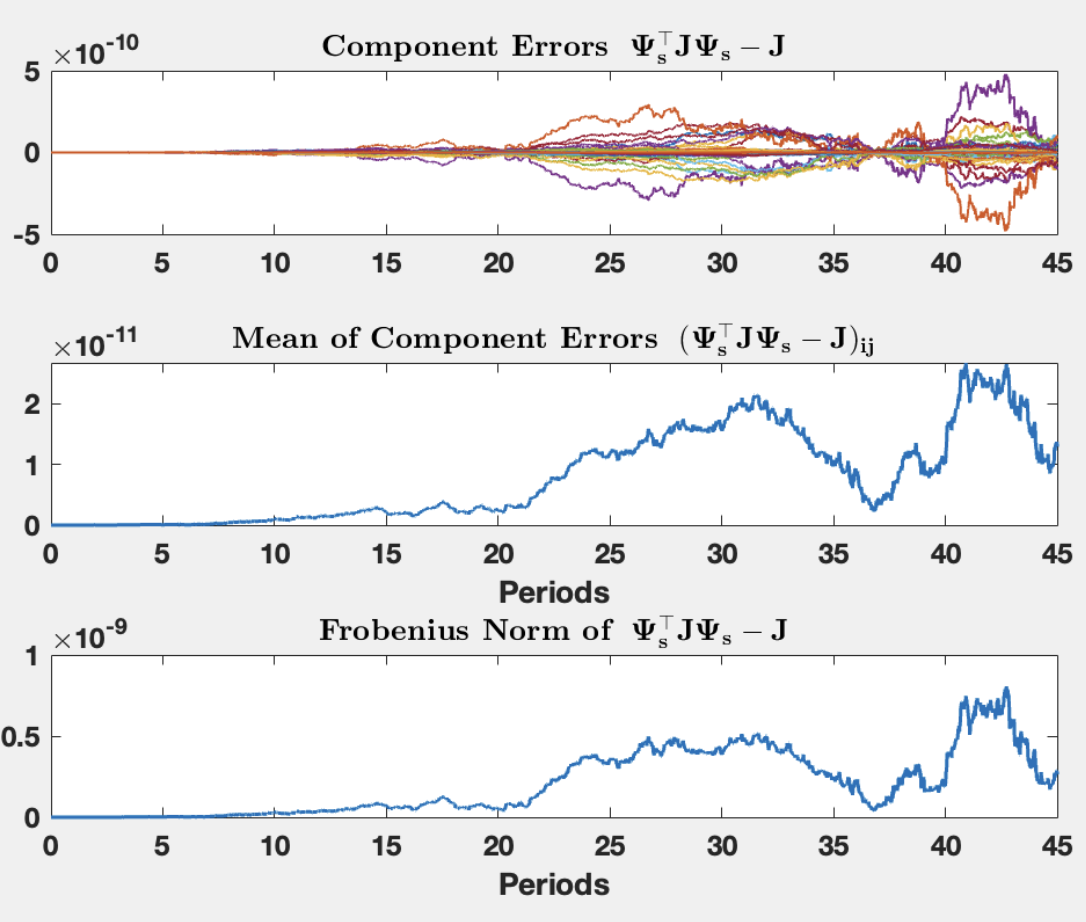}
    \caption{Verification that the STM for $(\tup{q},u,t,\tup{p},p_u,p_t)$, with \eq{s} as the evolution parameter, is again a symplectic matrix ($J_2$ gravitational term included.)
    Results obtained by propagating the general ODE for the STM from Eq.\eqref{STM_ODE}. } 
    \end{subfigure}
\end{figure}




